%% file: article_HAL.tex
\documentclass[final,3p]{elsarticle}

\usepackage{graphicx,color}
\usepackage{amsmath,amsfonts,amssymb}
\usepackage{mathrsfs,mathtools,stmaryrd}
\usepackage{fancybox}
\usepackage[english,francais]{babel}
\usepackage{subcaption}
\usepackage{float}
\usepackage{xspace}
\usepackage{enumitem}
\setlist[itemize,1]{label={--}}
\usepackage{pgf,tikz,pgfplots}
\pgfplotsset{compat=newest}
\pgfplotsset{
xlabel near ticks,
ylabel near ticks,
tick label style={font=\footnotesize},
label style={font=\small},
yticklabel style={/pgf/number format/fixed,,/pgf/number format/precision=4},
legend style={font=\small},
try min ticks=7,
legend pos=outer north east
}
\usepgfplotslibrary{external} 
\newlength\figureheight
\newlength\figurewidth
\newtheorem{e-proposition}[theorem]{Proposition}

\newtheorem{e-definition}[theorem]{Definition\rm}
\newdefinition{remark}{Remark}

\newcommand{\Om}{\Omega}
\newcommand{\om}{\omega}

\newcommand{\dO}{\textrm{d} \Omega}
\newcommand{\dom}{\textrm{d} \omega}
\newcommand{\dS}{\textrm{d} S}
\newcommand{\dt}{\textrm{d} t}
\newcommand{\dpp}{\textrm{d} \pb}

\newcommand{\intO}{\int_{\Omega}}
\newcommand{\into}{\int_{\omega}}
\newcommand{\intI}{\int_I}

\newcommand{\intP}{\int_P}
\newcommand{\CRE}{\mathrm{CRE}}
\newcommand{\PGD}{\mathrm{PGD}}
\newcommand{\dis}{\mathrm{dis}}
\newcommand{\corr}{\mathrm{corr}}

\newcommand{\bb}{\boldsymbol{b}}
\newcommand{\nb}{\boldsymbol{n}}
\newcommand{\pb}{\boldsymbol{p}}
\newcommand{\qb}{\boldsymbol{q}}
\newcommand{\rb}{\boldsymbol{r}}
\newcommand{\ub}{\boldsymbol{u}}
\newcommand{\xb}{\boldsymbol{x}}
\newcommand{\psib}{\boldsymbol{\psi}}
\newcommand{\taub}{\boldsymbol{\tau}}
\newcommand{\Bb}{\mathbf{B}}
\newcommand{\Fb}{\boldsymbol{F}}
\newcommand{\Gb}{\mathbf{G}}
\newcommand{\Nb}{\mathbf{N}}

\newcommand{\Kc}{\mathcal{K}}
\newcommand{\Lc}{\mathcal{L}}
\newcommand{\Mc}{\mathcal{M}}

\newcommand{\Pc}{\mathcal{P}}

\newcommand{\Sc}{\mathcal{S}}
\newcommand{\Tc}{\mathcal{T}}

\newcommand{\Vc}{\mathcal{V}}
\newcommand{\Wc}{\mathcal{W}}

\def\Abb{\mathbb{A}}

\def\Rbb{\mathbb{R}}
\def\Nbb{\mathbb{N}}

\newcommand{\forallin}[2]{\forall #1 \in #2}

\newcommand{\abs}[1]{\lvert#1\rvert}
\newcommand{\labs}[1]{\left\lvert#1\right\rvert}
\newcommand{\norm}[1]{\lVert#1\rVert}

\newcommand{\trinorm}[1]{\interleave#1\interleave} 

\newcommand{\grandO}[1]{O\mathopen{}\left(#1\right)}
\newcommand{\set}[1]{\{#1\}}

\newcommand{\setst}[2]{\{#1\mathrel{;}#2\}}

\newcommand{\restrictto}{\mathclose{}|\mathopen{}}

\DeclareMathOperator{\gradd}{grad\kern-.5em{grad}}

\DeclareMathOperator*{\argmax}{arg\,max}

\DeclareMathOperator*{\spann}{span}
\newcommand{\spanset}[1]{\ensuremath\spann\set{#1}}

\newcommand{\interval}[4]{\mathopen{#1}#2 \mathclose{}\mathpunct{},#3 \mathclose{#4}}
\newcommand{\intervalcc}[2]{\interval{[}{#1}{#2}{]}}
\newcommand{\intervaloc}[2]{\interval{]}{#1}{#2}{]}}

\renewcommand{\(}{\left(}
\renewcommand{\)}{\right)}
\renewcommand{\[}{\left[}
\renewcommand{\]}{\right]}

\renewcommand{\geq}{\geqslant}
\renewcommand{\leq}{\leqslant}
\let\oldtimes\times
\renewcommand{\times}{\!\oldtimes\!}

\newcommand{\ie}{i.e.\xspace}
\newcommand{\eg}{e.g.\xspace}

\journal{Computer Methods in Applied Mechanics and Engineering}

\begin{document}

\begin{frontmatter}



\title{\textit{A posteriori} error estimation and adaptive strategy for PGD model reduction applied to parametrized linear parabolic problems}


\author[authorlabel1]{Ludovic Chamoin}
\ead{chamoin@lmt.ens-cachan.fr}
\author[authorlabel2]{Florent Pled}
\ead{florent.pled@univ-paris-est.fr}
\author[authorlabel1]{Pierre-Eric Allier}
\ead{allier@lmt.ens-cachan.fr}
\author[authorlabel1]{Pierre Ladev\`eze\corref{cor1}}
\ead{ladeveze@lmt.ens-cachan.fr}

\cortext[cor1]{Corresponding author}

\address[authorlabel1]{
LMT-Cachan (ENS Cachan/CNRS/Universit\'e Paris-Saclay) \\
61 Avenue du Pr\'esident Wilson, 94235 Cachan Cedex, France}
\address[authorlabel2]{
Universit\'e Paris-Est, Laboratoire Mod\'elisation et Simulation Multi Echelle, MSME UMR 8208 CNRS, 5 bd Descartes, 77454 Marne-la-Vall\'ee, France \\
}

\begin{abstract}
We define an \textit{a posteriori} verification procedure that enables to control and certify PGD-based model reduction techniques applied to parametrized linear elliptic or parabolic problems. Using the concept of constitutive relation error, it provides guaranteed and fully computable global/goal-oriented error estimates taking both discretization and PGD truncation errors into account. Splitting the error sources, it also leads to a natural greedy adaptive strategy which can be driven in order to optimize the accuracy of PGD approximations. The focus of the paper is on two technical points: (i) construction of equilibrated fields required to compute guaranteed error bounds; (ii) error splitting and adaptive process when performing PGD-based model reduction. Performances of the proposed verification and adaptation tools are shown on several multi-parameter mechanical problems. 
\end{abstract}


\end{frontmatter}


\section{Introduction}
\label{section:introduction}

With continuous advances in modeling and computing methods, numerical simulation has progressively become a common tool for analysis and design in engineering activities. Nowadays, it enables to deal with complex (multiscale, multiphysics, multi-parameter, \dots) problems that include modeling with finer and finer features of the real world. Nevertheless, numerical simulation tools remain limited or even powerless for some categories of problems, in particular when considering complex multidimensional models with many fluctuating parameters. Such high-dimensional problems are encountered in several branches of computational science and engineering, such as parametric modeling (control, optimization, inverse analysis, \dots) or stochastic modeling (uncertainty quantification and propagation, risk assessment, sensitivity analysis, \dots). Classical numerical methods, known as brute force (mesh-based) discretization methods, then require huge and often unreasonable computational costs and storage requirements, as the number of degrees of freedom (dofs) grows exponentially with respect to the number of dimensions of the resulting approximation spaces; this is the so-called curse of dimensionality \cite{Bel61} related to computational intractability for high-dimensional problems. Consequently, new robust approximation methods need to be introduced to address multi-parameter models and efficiently compute numerical approximations for high-dimensional problems.

In this context, model reduction is an attractive alternative approach which has been widely developed during the last decade. It leans on the fact that the (full-order) solution of complex numerical models can often be accurately approximated by the (reduced-order) solution of surrogate models; this latter is obtained through the projection of the initial model onto a low-dimensional (reduced) subspace spanned by global basis functions, so that dimensionality can be drastically reduced. The various model (or complexity) reduction methods, such as Reduced Basis (RB) approaches \cite{Roz08} or Proper Orthogonal Decomposition (POD) \cite{Kun01,Lia02a}, distinguish themselves by the way of constructing and selecting the basis functions.

A promising model reduction method, denoted Proper Generalized Decomposition (PGD), has recently gained much attention and is currently the topic of numerous research works (see \cite{Chi10} for an overview) following pioneering ideas developed in \cite{Lad89,Lad99d}. It is a low-rank tensor approximation method which consists of a representation of the solution as a linear combination of separated variables functions (called modes), after defining all model parameters as extra-coordinates of the problem. The specificity of this spectral approach comes from the fact that no \textit{a priori} partial knowledge of the solution is required, contrary to POD which requires a preliminary offline stage, called learning phase and based on appropriate snapshots allowing to efficiently explore the parameter domain by means of \textit{a priori} error estimation procedures. 
In the PGD framework, modes are first computed on the fly, once at all, in an offline and iterative phase that provides an approximate solution of the model for any value of the parameters; this solution can then be used in an online phase, with cheap and fast computations on light computing platforms, in order to perform real-time parametric or stochastic analysis, sensitivity analysis, design or shape optimization, inverse identification and optimal control. The PGD method allows to circumvent the curse of dimensionally, as the number of dofs grows linearly with respect to the number of dimensions which enables considerable savings in terms of computing time and memory storage and leads to affordable simulations of complex engineering problems. Performances of the PGD method have been illustrated in several applications with linear or nonlinear problems, and including model variabilities of many kinds (material, loading, initial or boundary conditions, geometry, \dots) \cite{Amm06,
Nou07,
Nou10b,Lad10ter,Nou11b,Che13a,Gon12,Bou13,Chi14,Ner15,Lad16}.
 
However, an intensive use of PGD capabilities for numerical analysis as well as the transfer of PGD solvers in industrial activities still face several major difficulties. One of them is the control of PGD reduced-order models, in order to certify the accuracy of PGD-based numerical solutions for applications of increasing complexity, that represents a fundamental and well-identified concern for robust design and decision making. This scientific issue requires to construct dedicated \textit{a posteriori} error estimation tools as well as adaptive strategies in order to define a suitable PGD approximation in terms of required number of terms in the modal representation of the solution, but also in terms of the discretization meshes used to compute modes. Model verification is a pillar of simulation-based engineering for which J.T.~Oden has been a pioneer and a main contributor \cite{Ain00}. Verification of PGD reduced-order models has been addressed in very few works until now, in opposition to the vast literature dedicated to the control of RB techniques (see \cite{Mac01,Gre05,Rov06,Boy09,Qua11} for instance, where \textit{a priori} or \textit{a posteriori} explicit residual methods are used). Preliminary works were shown in \cite{Amm10} using residual-based techniques for goal-oriented error estimation; in this work, mainly devoted to adaptivity, only the error coming from the truncation of the modal representation was controlled and the error bounds were not guaranteed (the adjoint solution was approximated with a finer PGD decomposition). Another recent work \cite{Alf15} extended this concept to the nonlinear context by using a linearized version of the problem to define the adjoint problem, before using a weighted residuals method with a higher number of PGD modes to represent the adjoint solution and catch PGD truncation error. Even though this approach is cheap, it cannot deliver guaranteed error bounds in which all error sources are taken into account. A first robust verification approach, using the concept of Constitutive Relation Error (CRE) \cite{Lad83,Lad04}, was introduced in \cite{Lad11,Cha12bis,Lad12c} to control and assess the numerical quality of PGD reduced-order models. Based on the derivation of equilibrated fields from a post-processing of the approximate PGD solution, it provides guaranteed error bounds involving both discretization error and truncation error in the PGD modal representation. A similar approach was recently proposed in \cite{Moi13}, even though equilibrated fields were here obtained using a dual PGD approach. Other sources of errors resulting \eg from numerical integration (quadrature) and round-off (machine precision) are assumed to be negligible.

In this paper, we go one step forward by investigating three technical points. First, we present an automatic manner to construct equilibrated fields from a post-processing of the PGD approximation. Such a construction, required to apply CRE concepts, is actually the only way to get guaranteed error bounds. Second, as engineering design and optimization require the prediction of selected outputs as a function of specified inputs, we extend the verification procedure to the goal-oriented framework using similar tools as in the classical Finite Element Method (FEM) context \cite{Lad06,Lad08}. We thus define strict error bounds on given outputs of interest from the approximate PGD solution of an adjoint problem, and show how this latter can be computed effectively using PGD tools. Third, we introduce an adaptive strategy, based on a greedy algorithm involving the relative contributions of various error sources, in order to drive the PGD computations effectively and reach a prescribed error tolerance. Performances of the proposed verification method, valid for linear problems with potentially numerous fluctuating parameters, are illustrated through several numerical experiments carried out on elliptic (elasticity) or parabolic (transient thermal conduction) problems in one, two or three space dimensions. Advection problems are not considered here, even though the philosophy we present in the paper could be extended to such non-symmetric problems with minor changes \cite{Par09bis,Ern10bis}.

The paper is organized as follows. After this introduction, Section~\ref{section:refpb} presents the reference problem and notations. An overview of the technique used to get the PGD approximation, as well as its post-processing to obtain an equilibrated field (in a weak sense), is given in Section~\ref{section:PGDsolution}. Section~\ref{section:globalerror} introduces the global error estimation method, using the CRE concept and applied to PGD computations, as well as the associated adaptive strategy. An extension of the method to goal-oriented error estimation and model adaptation is proposed in Section~\ref{section:GOerror}. Numerical results are reported in Section~\ref{section:results}. Eventually, conclusions and prospects are drawn in Section~\ref{section:conclusions}.

\section{Reference problem}
\label{section:refpb}

We consider a transient linear diffusion-reaction problem defined on an open bounded domain $\Omega \subset \Rbb^d$ ($d=1, 2$ or $3$ being the space dimension), with boundary $\partial \Omega$, over a time interval $I = \intervalcc{0}{T}$. We assume that prescribed homogeneous Dirichlet boundary conditions are imposed on part $\partial_u \Omega \neq \emptyset$ of $\partial \Omega$, whereas time-dependent Neumann boundary conditions (given flux $g_d(\xb,t)$) are applied on the complementary part $\partial_q \Omega$ of $\partial \Omega$, with $\partial_u \Omega \cap \partial_q \Omega = \emptyset$ and $\overline{\partial_u \Omega \cup \partial_q \Omega} = \partial \Omega$. A given source term $f_d(\xb,t)$ may also be active in domain $\Omega$. For simplicity reasons, initial conditions are set to zero.

For the sake of simplicity of the technical aspects addressed in the paper, we consider variabilities in material properties only; however, the PGD method as well as the verification approach we present can be extended to other variabilities occurring in the loading, geometry, and initial or boundary conditions. The material behavior is assumed to be fluctuating due to possible heterogeneities and/or uncertainties characterized by a set of $n_p$ input design parameters gathered in a vector $\pb = [p_1,\dots,p_{n_p}] \in P = P_1 \times \cdots \times P_{n_p}$, where the input parameter domains $P_j$ are bounded subsets of $\Rbb$ defining the range of variations of parameters $p_j$. \\

Denoting by $c(\xb,\pb)$ (resp. $k(\xb,\pb)$ and $r(\xb,\pb)$) the evolution (resp. diffusion and reaction) coefficient, the reference mathematical problem consists in finding $u(\xb,t,\pb)$ (and associated flux $\qb=k\boldsymbol{\nabla}u$), with $(\xb,t,\pb) \in \Omega \times I \times P$, solution of the following partial differential equation (PDE):
\begin{subequations}\label{eq:initpb}
\begin{equation}\label{eq:initPDE}
c \frac{\partial u}{\partial t} - \boldsymbol{\nabla} \cdot (k \boldsymbol{\nabla} u) + r u = f_d \quad \text{on } \Omega \times I \times P,
\end{equation}
with given initial and boundary conditions:
\begin{equation}
u = 0 \quad \text{on } \Omega \times \set{0} \times P, \quad
u = 0 \quad \text{on } \partial_u \Omega \times I \times P, \quad
k \boldsymbol{\nabla} u \cdot \nb = g_d \quad \text{on } \partial_q \Omega \times I \times P,
\end{equation}
\end{subequations}
$\nb$ denoting the unit outgoing normal vector to $\Omega$. \\

In the following, a function $v(\xb,t)$ defined on $\Omega\times I$ will be equivalently considered as a function $v(t)$ defined on $I$ with values in the Hilbert space $\Vc = H_0^1(\Omega) = \setst{v \in H^1(\Omega)}{v_{\restrictto{\partial_u \Omega}} = 0}$, corresponding to the Sobolev space of functions vanishing on $\partial_u \Omega$. Similarly, a function $v(\xb,t,\pb)$ defined on $\Omega\times I\times P$ will be equivalently considered as a function $v(t,\pb)$ defined on $I\times P$ with values in $\Vc$. For the sake of readability, the dependence of functions on variables in space $\xb$, time $t$ and parameters $\pb$ will be often omitted.

Introducing $\Tc = L^2(I)$ the Lebesgue space associated to $I$, identifying the Bochner space $L^2(I;\Vc)$ (corresponding to the Hilbert space of square integrable functions defined from $I$ into $\Vc$) with the tensor product space $\Vc \otimes \Tc$: $L^2(I;\Vc) \simeq \Vc \otimes \Tc$, and denoting $\Vc^{\ast} = H^{-1}(\Omega)$ the topological dual space of $\Vc$, the space-time weak formulation of \eqref{eq:initpb} reads for all $\pb \in P$:
\begin{equation}\label{eq:spacetimeweakform}
\text{find $u \in L^2(I;\Vc)$, with $\dfrac{\partial u}{\partial t} \in L^2(I;\Vc^{\ast})$, such that }
b(u,v) = l(v) \quad \forallin{v}{L^2(I;\Vc)},
\end{equation}
with $u_{\restrictto{t=0}} = 0$, and where the bilinear form $b(\cdot,\cdot)$ and linear form $l(\cdot)$ are defined on $L^2(I;\Vc)$ by:
\begin{equation}
\begin{aligned}
b(u,v) &\coloneqq \intI \intO \(c \frac{\partial u}{\partial t} v + k \boldsymbol{\nabla} u \cdot \boldsymbol{\nabla} v + r u v\) \dO \dt
, 
\\ 
l(v) &\coloneqq \intI \(\intO f_d v \dO + \int_{\partial_q \Omega} g_d v \dS \)\dt, 
\end{aligned}
\end{equation}
where $\dfrac{\partial u}{\partial t}(t)\in \Vc^{\ast}$ is assimilated to its Riesz representation in $\Vc$ \cite{Nou10b} and the duality pairing between $\Vc$ and $\Vc^{\ast}$ is considered as a continuous extension of the natural inner product on Lebesgue space $L^2(\Omega)$ \cite{Rou05}.

Then, introducing $\Pc_j = L^2(P_j)$ the Lebesgue space associated to $P_j$ and identifying the Bochner space $L^2(I,P;\Vc)$ (corresponding to the Hilbert space of square integrable functions defined from $I\times P$ into $\Vc$) with the tensor product space $\Vc \otimes \Tc \otimes \Pc$: $L^2(I,P;\Vc) \simeq \Vc \otimes \Tc \otimes \Pc$, with $\Pc = \otimes_{j=1}^{n_p} \Pc_j$, the full weak formulation of \eqref{eq:initpb} reads:
\begin{equation}\label{eq:weakform}
\text{find $u \in L^2(I,P;\Vc)$, with $\dfrac{\partial u}{\partial t} \in L^2(I,P;\Vc^{\ast})$, such that }
B(u,v) = L(v) \quad \forallin{v}{L^2(I,P;\Vc)},
\end{equation}
where $B(\cdot,\cdot)$ and $L(\cdot)$ are bilinear and linear forms defined on $L^2(I,P;\Vc)$ by:
\begin{equation}
B(u,v) \coloneqq \intP b(u,v) \dpp, \quad L(v) \coloneqq \intP l(v) \dpp.
\end{equation}
The solution $u$ of \eqref{eq:weakform}, which depends on $D=(d+1+n_p)$ dimensions, may be approximated with a classical brute force (grid-based) approach only for a moderate number of parameters $n_p$. Indeed, considering a regular Cartesian domain $\Omega$ and $n$ discretization points in each dimension, a full tensor product approximation of the solution $u \in L^2(I,P;\Vc)$ reads:
\begin{equation}
u \approx \sum_{i_1=1}^{n} \dots \sum_{i_D=1}^n a_{i_1\dots i_D} \otimes_{\ell=1}^{D} \varphi^{\ell}_{i_{\ell}},
\end{equation}
where $\varphi^{\ell}_{i_{\ell}}$ are canonical basis functions in dimension $\ell$ and $a \in \Rbb^{n^D}$ gathers the components $a_{i_1\dots i_D}$ of the full tensor product approximation on the canonical basis. 
It leads to an exponential blow up of the number of dofs (or complexity) $N$ with respect to the number of dimensions $D$ ($N = n^D$), often referred to as the so-called curse of dimensionality. Furthermore, assuming that solution $u \in C^{s}(\Omega \times I \times P)$ with order $s \in \Nbb$, the uniform accuracy $\varepsilon$ measured in the $L^{\infty}(I,P;\Vc)$-norm is such that $\varepsilon(N,D) = \grandO{n^{-s}} =\grandO{N^{-s/D}}$ as $N \to +\infty$ \cite{Che66}. Consequently, the associated numerical complexity $N$ is such that $N(\varepsilon,D) = \grandO{\varepsilon^{-D/s}}$ as $\varepsilon \to 0$, which corresponds to a very poor convergence rate $-s/D$ for high-dimensional problems, \ie involving a large number of dimensions $D$. In the remainder of the paper, we do not perform any separation of space variables in order to avoid restrictions in the shape of the domain $\Omega$ (\eg geometrical symmetry or periodicity conditions). We then consider general physical domains $\Omega$ so that the $d$ space dimensions are kept together. This way, we get $D=2+n_p$ and the number $n$ of discretization points in the physical domain $\Omega$ corresponds to the number of degrees of freedom over the physical mesh.

A classical model reduction technique, known as the Reduced Basis (RB) method, consists of properly selecting a discrete subset $S_J = \set{\pb_j}_{j=1}^J$ of parameters values and computing the solution $u_j(\xb,t)=u(\xb,t,\pb_j)$ ($1 \leq j \leq J$) of the parametrized variational problem \eqref{eq:spacetimeweakform} for each element $\pb_j$ of $S_J$ (involving $N = J \times n^{d+1}$ dofs). The $J$-dimensional Lagrangian reduced-basis subspace $\Wc_J = \spanset{u_j}_{j=1}^J \subset \Pc$ is then defined, and an approximation of the solution $u(\xb,t,\pb)$, for any input parameter vector $\pb \in P$, can be obtained by a standard Galerkin projection onto the reduced space $\Wc_J$. The performances of RB-like approaches, which depend on the Kolmogorov $J$-width 
 of the set $\Kc = \set{u(\cdot,\cdot,\pb)}_{\pb \in P} \subset L^2(I;\Vc)$, as well as \textit{a posteriori} verification tools to define optimal elements in $S_J$, have been investigated in several works (see \cite{Rov06,Qua11} for instance). The effectivity of RB methods strongly depends on the way to construct the low-dimensional subspace $\Wc_J$ and especially on the choice of snapshots $\set{u_j}_{j=1}^J$. In the next section, an alternative model reduction approach based on PGD is detailed.

\section{PGD reduced-order approximation: computation and post-processing}
\label{section:PGDsolution}

The basic idea in model reduction using PGD is to \textit{a priori} construct an approximation of the solution $u(\xb,t,\pb)$ as a separated variables representation defined in tensor product spaces, \ie a finite sum of products of separable functions. We consider here such a representation associated with a canonical format and (low-rank) separated structure for variables in space $\xb$, time $t$ and parameters $p_j \in \pb$; it reads:
\begin{equation}\label{eq:lowrank}
u \approx \sum_{i=1}^m \otimes_{\ell=1}^D u^{\ell}_i, \quad \text{with } u^{\ell}_i = \sum_{i_\ell=1}^n a^{\ell}_{i,i_\ell} \varphi^{\ell}_{i_\ell},
\end{equation}
where $m \in \Nbb^{\ast}$ is the rank (or order) of the PGD approximation and $u^{\ell}_i$ are global reduced basis functions in dimension $\ell$.
Consequently, the computational cost is drastically reduced with this specific low-rank structure when $D$ increases, as it now leads to a linear scaling (or linear growth) of $N$ with respect to $D$ ($N=mDn$). Besides, assuming that the solution $u$ can be approximated by a rank-$m$ separated representation \eqref{eq:lowrank} (also called order~$m$ canonical tensor decomposition) with $u^{\ell}_i \in C^{s}(P_{\ell})$, where $P_{\ell}$ denotes here the domain of the generic coordinate in dimension $\ell$, the numerical complexity $N$ to achieve accuracy $\varepsilon$ is such that $N(\varepsilon,D) = \grandO{m^{1/s} D^{1+1/s} \varepsilon^{-1/s}}$ as $\varepsilon \to 0$. 
In the continuous setting, the rank-$m$ PGD approximation \eqref{eq:lowrank} defines a separated modal representation with space, time and parameters functions.

Consequently, the continuous order~$m$ PGD approximation $u_m$ is defined as: 
\begin{equation}\label{eq:PGDrep}
u_m(\xb,t,\pb) = \sum_{i=1}^m \psi_i(\xb) \lambda_i(t) \prod_{j=1}^{n_p} \gamma_{j,i}(p_j),
\end{equation}
where the $\psi_i \in \Vc$ (resp. $\lambda_i \in \Tc$ and $\gamma_{j,i} \in \Pc_j$) form a low-dimensional reduced basis composed of space (resp. time and parameters) functions. The discretized version $u^{h,\Delta t}_m$, denoting by $h$ (resp. $\Delta t$) the space mesh size (resp. time step) which is used, reads:
\begin{equation}\label{eq:discretizedPGDrep}
u^{h,\Delta t}_m(\xb,t,\pb) = \sum_{i=1}^m \psi^h_i(\xb) \lambda^{\Delta t}_i(t) \prod_{j=1}^{n_p} \gamma_{j,i}(p_j),
\end{equation}
where the $\psi^h_i$ (resp. $\lambda^{\Delta t}_i$) functions are discretized counterparts of $\psi_i$ (resp. $\lambda_i$) lying in a finite-dimensional approximation subspace $\Vc_h \subset \Vc$ (resp. $\Tc_{\Delta t} \subset \Tc$).

Even though discretizations in the parameters $p_j$ dimensions are required for computational purposes, we do not explicitly consider them in the error analysis as they are not associated with a given numerical approximation method that we wish to adapt. In the numerical results presented in Section~\ref{section:results}, we will consider a very fine discretization grid in each parameter domain $P_j$ in order to safely neglect the error due to the numerical approximation over the parametric space $\Pc$.

\subsection{Computation of the PGD approximation}\label{section:PGDcomputation}

In the case where the solution $u$ is known (at least partially by means of snapshots), an optimal low-rank separated representation may be searched by minimizing the distance to the exact solution with respect to a given metric on the tensor product $\Vc \otimes \Tc \otimes \Pc$; the classical POD approach (known as SVD technique in matrix computations) corresponds to a particular case where a $L^2$-norm is used for two variables functions ($D=2$), leading to eigenvalue problems. Here, we wish to compute modes on the fly with no \textit{a priori} knowledge on the solution $u$.

Among the various methods which have been introduced so far to compute PGD modes (separable functions) and build an approximate separated representation of the solution $u$ (see \cite{Chi10,Nou10b} for instance), we use the classical one referred to as progressive Galerkin-based PGD method; it operates in an iterative strategy based on the progressive construction of successive order~$1$ corrections, and defines Galerkin approximations in tensor product spaces from the full weak formulation \eqref{eq:weakform}. In particular, it involves a space-time weak formulation \cite{Lio72} of parabolic problem \eqref{eq:initpb} in order to separate the space and time modes. Assuming that an order~$(m-1)$ PGD decomposition $u_{m-1}$ is known, the order~$m$ decomposition $u_m$ is searched as:
\begin{equation}
u_m(\xb,t,\pb) = u_{m-1}(\xb,t,\pb) + \psi(\xb) \lambda(t) \prod_{j=1}^{n_p} \gamma_j(p_j),
\end{equation}
where the new space function $\psi$, time function $\lambda$ and parameter functions $\gamma_j$ are the unknown functions to be determined and respectively belong to $\Vc$, $\Tc$ and $\Pc_j$.

The Galerkin approach then implies that these functions verify the following Galerkin orthogonality conditions:
\begin{equation}\label{eq:Galerkinorthcond}
\begin{aligned}
B(u_{m-1}+\psi \lambda \prod_{j=1}^{n_p} \gamma_j,v^{\ast}) = L(v^{\ast}), \quad &\text{with } v^{\ast}=\psi^{\ast} \lambda \prod_{j=1}^{n_p} \gamma_j + \psi \lambda^{\ast} \prod_{j=1}^{n_p} \gamma_j + \sum_{j_0=1}^{n_p} \psi \lambda \gamma_{j_0}^{\ast} \prod_{\substack{j=1\\ j \neq j_0}}^{n_p} \gamma_j, \\
&\forallin{\psi^{\ast}}{\Vc}, \; \forallin{\lambda^{\ast}}{\Tc}, \; \forallin{\gamma_{j_0}^{\ast}}{\Pc_{j_0}}, \; j_0=1,\dots,n_p, \\
\end{aligned}
\end{equation}
or, equivalently, the following stationarity conditions:
\begin{subequations}\label{eq:Galerkinorth}
\begin{alignat}{2}
&B(u_{m-1}+\psi \lambda \prod_{j} \gamma_j,\psi^{\ast} \lambda \prod_{j} \gamma_j) &&= L(\psi^{\ast} \lambda \prod_{j} \gamma_j) \quad \forallin{\psi^{\ast}}{\Vc}, \label{eq:spaceGalerkincond}\\
&B(u_{m-1}+\psi \lambda \prod_{j} \gamma_j,\psi \lambda^{\ast} \prod_{j} \gamma_j) &&= L(\psi \lambda^{\ast} \prod_{j} \gamma_j) \quad \forallin{\lambda^{\ast}}{\Tc}, \label{eq:timeGalerkincond}\\
&B(u_{m-1}+\psi \lambda \prod_{j} \gamma_j,\psi \lambda \gamma_{j_0}^{\ast} \prod_{j \neq j_0} \gamma_j) &&= L(\psi \lambda \gamma_{j_0}^{\ast} \prod_{j \neq j_0} \gamma_j) \quad \forallin{\gamma_{j_0}^{\ast}}{\Pc_{j_0}}, \; j_0=1,\dots,n_p. \label{eq:paramGalerkincond}
\end{alignat}
\end{subequations}
Problem \eqref{eq:Galerkinorth} is a complex nonlinear multidimensional problem which can be interpreted as a pseudo-eigenvalue problem; it may thus be solved using specific iterative algorithms inspired from classical power iterations algorithms dedicated to eigenvalue problems or dominant subspace methods \cite{Nou10b}. In the present work, we use the fixed-point iteration method (also called alternating direction algorithm in an optimization context). For each mode $m \in \Nbb^{\ast}$, starting from an \textit{ad hoc} initialization $(\psi^{(0)},\lambda^{(0)},\gamma^{(0)}_1,\dots,\gamma^{(0)}_{n_p})$ at iteration 0, one builds a sequence $\set{(\psi^{(k)},\lambda^{(k)},\gamma^{(k)}_1,\dots,\gamma^{(k)}_{n_p})}_{k \in \Nbb^{\ast}}$ with the following power sub-iterations algorithm requiring the solution of a sequence of simple low-dimensional problems at each sub-iteration $k \in \Nbb^{\ast}$:
\begin{itemize}
\item Compute $\lambda^{(k)} \in \Tc$ such that:
\begin{equation}\label{eq:timepb}
B(u_{m-1}+\psi^{(k-1)} \lambda^{(k)} \prod_{j} \gamma^{(k-1)}_j, \psi^{(k-1)} \lambda^{\ast} \prod_{j} \gamma^{(k-1)}_j) = L(\psi^{(k-1)} \lambda^{\ast} \prod_{j} \gamma^{(k-1)}_j) \quad \forallin{\lambda^{\ast}}{\Tc}.
\end{equation}
Problem \eqref{eq:timepb} is the weak formulation of a scalar ordinary differential equation (ODE) in time which may be solved in practice using the FEM with associated finite-dimensional approximation subspace $\Tc_{\Delta t} \subset \Tc$. An alternative method would consist in using an incremental time integration scheme dedicated to first-order differential equations;
\item For $j_0=1,\dots,n_p$, compute $\gamma_{j_0}^{(k)} \in \Pc_{j_0}$ such that:
\begin{equation}\label{eq:parampb}
\begin{aligned}
&B(u_{m-1}+\psi^{(k-1)} \lambda^{(k)} \gamma_{j_0}^{(k)} \prod_{j<j_0} \gamma_j^{(k)} \prod_{j>j_0} \gamma_j^{(k-1)}, \psi^{(k-1)} \lambda^{(k)} \gamma^{\ast} \prod_{j<j_0} \gamma_j^{(k)} \prod_{j>j_0} \gamma_j^{(k-1)}) \\
&= L(\psi^{(k-1)} \lambda^{(k)} \gamma^{\ast} \prod_{j<j_0} \gamma_j^{(k)} \prod_{j>j_0} \gamma_j^{(k-1)}) \quad \forallin{\gamma^{\ast}}{\Pc_{j_0}}.
\end{aligned}
\end{equation}
Problem \eqref{eq:parampb} is a scalar algebraic equation and provides an explicit definition of parameter function $\gamma_{j_0}^{(k)} \in \Pc_{j_0}$;
\item Compute $\psi^{(k)} \in \Vc$ such that:
\begin{equation}\label{eq:spacepb}
B(u_{m-1}+\psi^{(k)} \lambda^{(k)} \prod_{j} \gamma^{(k)}_j, \psi^{\ast} \lambda^{(k)} \prod_{j} \gamma^{(k)}_j) = L(\psi^{\ast} \lambda^{(k)} \prod_{j} \gamma^{(k)}_j) \quad \forallin{\psi^{\ast}}{\Vc}.
\end{equation}
Problem \eqref{eq:spacepb} is the weak formulation of a (time-independent) partial differential equation (PDE) in space which may be solved in practice using the FEM with associated finite-dimensional approximation subspace $\Vc_h \subset \Vc$.
\end{itemize}

\vspace{1em}

Sub-iterations may be performed until the convergence is reached at a given tolerance. In practice, the power sub-iterations algorithm converges quite fast and generally does not require more than a few iterations to obtain a good approximation of $(\psi^{(k)},\lambda^{(k)},\gamma^{(k)}_1,\dots,\gamma^{(k)}_{n_p})$. Here, we choose to stop sub-iterations after a given number $k_{\mathrm{max}}$ (we take $k_{\mathrm{max}}=4$ iterations in the numerical experiments shown in Section~\ref{section:results}). Furthermore, time function $\lambda^{(k)}$ and parameter functions $\gamma^{(k)}_j$ are normalized at each sub-iteration $k$ so that the magnitude of PGD mode $m$ is supported by space function $\psi^{(k)}$ alone.

Several possible variants, which will not be considered in the numerical experiments, can be introduced in the progressive Galerkin-based PGD approach in order to capture a good approximation of the optimal decomposition, which would be obtained by directly computing all modes simultaneously (and not progressively). In particular, at any mode $m \in \Nbb^{\ast}$: 
\begin{itemize}
\item space function $\psi_m$ may be orthogonalized with respect to the existing space basis $\set{\psi_i}_{i=1}^{m-1}$ in order to decouple the possible dependencies between functions and thus improve the condition number of the space problem \eqref{eq:spacepb}; 
\item time functions $\lambda_i$ and parameter functions $\gamma_{j,i}$ associated with previously computed modes $i \leq m-1$ may be updated before starting the power sub-iterations algorithm for computing mode $m$, in order to satisfy a stronger Galerkin orthogonality condition leading to an order~$m$ decomposition $u_m$ of better quality; conversely, space functions $\psi_i$ are usually conserved as they generally require most of the computational cost. This preliminary stage actually corresponds to a POD step. We thus compute time functions $\lambda_i^{\mathrm{up}}$ and parameter functions $\gamma_{j,i}^{\mathrm{up}}$, with $1 \leq i \leq m-1$, such that:
\begin{equation}\label{eq:Galerkinorthcondup}
\begin{aligned}
B(\sum_{i=1}^{m-1} \psi_i \lambda_i^{\mathrm{up}} \prod_j \gamma_{j,i}^{\mathrm{up}},\sum_{i=1}^{m-1} \psi_i v_i^{\ast}) = L(\sum_{i=1}^{m-1} \psi_i v_i^{\ast}), \quad &\text{with } v_i^{\ast}=\lambda_i^{\ast} \prod_j \gamma_{j,i}^{\mathrm{up}} + \sum_{j_0=1}^{n_p} \lambda_i^{\mathrm{up}} \gamma_{j_0,i}^{\ast} \prod_{j \neq j_0} \gamma_{j,i}^{\mathrm{up}}, \\
&\forallin{\lambda_i^{\ast}}{\Tc}, \; \forallin{\gamma_{j_0,i}^{\ast}}{\Pc_{j_0}}, \; j_0=1,\dots,n_p, \\
\end{aligned}
\end{equation}
or, equivalently, for all $i_0 \in \{1,\dots,m-1\}$:
\begin{subequations}\label{eq:Galerkinorthup}
\begin{alignat}{2}
&B(\sum_{i=1}^{m-1} \psi_i \lambda_i^{\mathrm{up}} \prod_j \gamma_{j,i}^{\mathrm{up}},\psi_{i_0} \lambda_{i_0}^{\ast} \prod_j \gamma_{j,i_0}^{\mathrm{up}}) &&= L(\psi_{i_0} \lambda_{i_0}^{\ast} \prod_j \gamma_{j,i_0}^{\mathrm{up}}) \quad \forallin{\lambda_{i_0}^{\ast}}{\Tc}, \\
&B(\sum_{i=1}^{m-1} \psi_i \lambda_i^{\mathrm{up}} \prod_j \gamma_{j,i}^{\mathrm{up}},\psi_{i_0} \lambda_{i_0}^{\mathrm{up}} \gamma_{j_0,i_0}^{\ast} \prod_{j \neq j_0} \gamma_{j,i_0}^{\mathrm{up}}) &&= L(\psi_{i_0} \lambda_{i_0}^{\mathrm{up}} \gamma_{j_0,i_0}^{\ast} \prod_{j \neq j_0} \gamma_{j,i_0}^{\mathrm{up}}) \quad \forallin{\gamma_{j_0,i_0}^{\ast}}{\Pc_{j_0}}, \; j_0=1,\dots,n_p.
\end{alignat}
\end{subequations}
The nonlinear problem \eqref{eq:Galerkinorthup} may also be solved using a fixed-point iteration method (with initialization $(\lambda_i,\gamma_{1,i},\dots,\gamma_{n_p,i})$ and normalization of parameter functions $\gamma_{j,i}^{\mathrm{up}}$); it requires low-cost numerical approximation methods as only $(m-1)$ ODEs in time need to be solved, leading to inexpensive computations. Such improvements based on updating steps allow recovering good convergence properties of approximate separated representations in many applications \cite{Nou08b,
Nou10b}.
\end{itemize}

\subsection{Post-processing of the PGD approximation}

Due to the fact that it uses the full weak formulation \eqref{eq:weakform} of the problem, the progressive Galerkin-based PGD technique presented in Section~\ref{section:PGDcomputation} (or any other variant of this PGD technique) provides an approximate solution $u_m^{h,\Delta t}$ (and associated flux $\qb(u_m^{h,\Delta t}) = k \boldsymbol{\nabla} u_m^{h,\Delta t}$) that satisfies the kinematic constraints and initial conditions but fails to verify the equilibrium equations \eqref{eq:initPDE} in any weak sense in space. In order to overcome this drawback and be able to use error estimation tools defined in the next sections (inspired from those used in the FEM context), we propose here a method to recover a FE-equilibrated PGD approximation from $u_m^{h,\Delta t}$ and the prescribed data alone. In other words, we construct a flux $\hat{\qb}_m^h$ that satisfies equilibrium in the FE sense with $u_m^{h,\Delta t}$ for all $(t,\pb)\in I \times P$:
\begin{equation}\label{eq:FEequilibrium}
\intO \(c \frac{\partial u_m^{h,\Delta t}}{\partial t} v^h + \hat{\qb}_m^h \cdot \boldsymbol{\nabla}v^h + r u_m^{h,\Delta t} v^h\) \dO = \intO f_d v^h \dO + \int_{\partial_q \Omega} g_d v^h \dS \quad \forallin{v^h}{\Vc_h}.
\end{equation}
The idea is to exploit the following property satisfied by any computed order~$m_0$ PGD approximation $u^{h,\Delta t}_{m_0}$ ($1 \leq m_0 \leq m$) when the fixed-point iteration method used to solve the nonlinear problem \eqref{eq:Galerkinorth} is stopped after solving the space problem \eqref{eq:spacepb} in the approximation space $\Vc_h$:
\begin{equation}\label{prop:spacepb}
B(u^{h,\Delta t}_{m_0},\psi^{\ast}\lambda^{\Delta t}_{m_0}\prod_j \gamma_{j,m_0})=L(\psi^{\ast}\lambda^{\Delta t}_{m_0}\prod_j \gamma_{j,m_0}) \quad \forall \psi^{\ast} \in \Vc_h.
\end{equation}
We assume that the functional $L$ in the right-hand side of \eqref{prop:spacepb} can be put under the following separated form in terms of space and time variables (known as radial approximation):
\begin{equation}\label{eq:radialapproxloading}
L(v) = \intP \intI \(\sum_{s=1}^S \alpha_s(t) L_s(v)\) \dt \dpp,
\end{equation}
where $\alpha_s(t)$ are time functions, and $L_s(v)$ are time-independent linear forms involving space functions alone. Note that this assumption is generally met in practical applications, since the loading $(f_d,g_d)$ is usually given in terms of some products of space and time functions. Then, solving, using a standard FEM, the $S$ steady-state (time-independent) space problems of the form:
\begin{equation}\label{eq:spacestaticpb}
\text{find $\hat{\rb}^h_s \in \Sc_h$} such that 
\intO \hat{\rb}^h_s \cdot \boldsymbol{\nabla}v^h \dO = L_s(v^h) \quad \forallin{v^h}{\Vc_h}, \quad s=1,\dots,S
\end{equation}
where $\Sc_h$ is a subspace of $\Lc^2(\Omega)$, and defining the flux $\hat{\rb}^h \in \Sc_h \otimes \Tc$ as $\hat{\rb}^h(\xb,t) = \sum_{s=1}^S \alpha_s(t) \hat{\rb}^h_s(\xb)$, the property \eqref{prop:spacepb} can be recast as:
\begin{equation}\label{eq:spaceproperty}
\intO \(\sum_{i=1}^{m_0} a_{m_0i} \psi_i^h \psi^{\ast} + \hat{\bb}^h_{m_0} \cdot \boldsymbol{\nabla}\psi^{\ast}\) \dO = 0 \quad \forallin{\psi^{\ast}}{\Vc_h},
\end{equation}
with
\begin{equation*}
\begin{aligned}
a_{m_0i} &= \intP \intI \lambda_{m_0}^{\Delta t} \prod_j \gamma_{j,m_0} \(c \dot{\lambda}_i^{\Delta t} + r \lambda_i^{\Delta t}\) \prod_j \gamma_{j,i} \dt \dpp, \\
\hat{\bb}^h_{m_0} &= \intP \intI \lambda^{\Delta t}_{m_0} \prod_j \gamma_{j,m_0} \(k \boldsymbol{\nabla}u_{m_0}^{h,\Delta t} - \hat{\rb}^h\) \dt \dpp.
\end{aligned}
\end{equation*}

\begin{remark}
A solution $\hat{\rb}^h_s$ to \eqref{eq:spacestaticpb} can be computed using a primal approach, \ie defining a primal field $w^h$ over the FE approximation space $\Vc_h$ and taking the corresponding dual field $\hat{\rb}^h_s = \boldsymbol{\nabla}w^h$. This way, solving \eqref{eq:spacestaticpb} comes back to the solution of a classical static problem (with a constant elliptic operator) using standard FEM.
\end{remark}

Then, writing the property \eqref{eq:spaceproperty} for all $m_0 \in \{1,\dots,m\}$ yields the following triangular system of equations:
\begin{equation}\label{eq:globalspaceproperty}
\intO \(\Abb_m \boldsymbol{\Psi}_m^h \psi^{\ast} + \hat{\Bb}_m^h \cdot \boldsymbol{\nabla}\psi^{\ast}\) \dO = 0 \quad \forallin{\psi^{\ast}}{\Vc_h}
\end{equation}
with $(\Abb_m)_{ij}=a_{ij}$, $\boldsymbol{\Psi}_m^h=[\psi_1^h,\dots,\psi_m^h]^T$ and $\hat{\Bb}_m^h = [\hat{\bb}^h_1,\dots,\hat{\bb}^h_m]^T$. \\
Assuming that material parameters $c$ and $r$ are homogeneous (uniform in space) over the domain $\Omega$, the matrix $\Abb_m$ becomes constant and the system \eqref{eq:globalspaceproperty} can be easily inverted. It follows, introducing $\hat{\boldsymbol{\Phi}}_m^h = [\hat{\boldsymbol{\phi}}^h_1,\dots,\hat{\boldsymbol{\phi}}^h_m]^T = \Abb_m^{-1}\hat{\Bb}_m^h$ which is fully computable, that for all $(t,\pb) \in I \times P$:
\begin{equation}
\intO \(\boldsymbol{\Psi}_m^h \psi^{\ast} + \hat{\boldsymbol{\Phi}}_m^h \cdot \boldsymbol{\nabla}\psi^{\ast}\) \dO = 0 \quad \forallin{\psi^{\ast}}{\Vc_h}
\end{equation}
or, equivalently, for all $i \in \{1,\dots,m\}$,
\begin{equation}
\intO \(\psi_i^h \psi^{\ast} + \hat{\boldsymbol{\phi}}^h_i \cdot \boldsymbol{\nabla}\psi^{\ast}\) \dO = 0 \quad \forallin{\psi^{\ast}}{\Vc_h}.
\end{equation}
Note that, incorporating the flux $\hat{\rb}^h$ into the FE equilibrium equations \eqref{eq:FEequilibrium}, the researched FE-equilibrated flux $\hat{\qb}_m^h$ should satisfy for all $(t,\pb) \in I \times P$:
\begin{equation}\label{eq:FEequilibriumbis}
\intO \(\hat{\qb}_m^h - \hat{\rb}^h\) \cdot \boldsymbol{\nabla}v^h \dO = -\intO \(c \frac{\partial u_m^{h,\Delta t}}{\partial t} v^h + r u_m^{h,\Delta t} v^h\) \dO = -\sum_{i=1}^m \(c \dot{\lambda}_i^{\Delta t} + r \lambda_i^{\Delta t}\) \prod_j \gamma_{j,i} \intO \psi_i^h v^h \dO \quad \forallin{v^h}{\Vc_h}.
\end{equation}
Consequently, a flux $\hat{\qb}_m^h$ which verifies \eqref{eq:FEequilibrium} (or \eqref{eq:FEequilibriumbis}) can be defined as:
\begin{equation}
\hat{\qb}_m^h(\xb,t,\pb) = \hat{\rb}^h(\xb,t) + \sum_{i=1}^m \hat{\boldsymbol{\phi}}^h_i(\xb) \(c \dot{\lambda}^{\Delta t}_i(t) + r \lambda^{\Delta t}_i(t)\) \prod_j \gamma_{j,i}(p_j).
\end{equation}

The FE-equilibration procedure used to recover the flux $\hat{\qb}_m^h$ can be easily extended to the case where the external loading $(f_d,g_d)$ depends on parameters $\pb \in P$. Indeed, in that case, assuming that the functional $L$ can be written as $L(v) = \intP \intI \(\sum_{s=1}^S \alpha_s(t,\pb) L_s(v)\) \dt \dpp$, we simply need to define the flux $\hat{\rb}^h\in \Vc_h \otimes \Tc \otimes \Pc$ as $\hat{\rb}^h(\xb,t,\pb) = \sum_{s=1}^S \alpha_s(t,\pb) \hat{\rb}^h_s(\xb)$.

Besides, in the case where the geometry of the space domain $\Omega$ depends on parameters $\pb \in P$, different strategies have been proposed in order to reformulate the weak problem \eqref{eq:weakform} on a reference fixed (parameter-independent) domain $\Omega_{\mathrm{ref}}$, by introducing a suitable (parameter-dependent) mapping to a fixed domain, or by using a fictitious domain method. In \cite{Amm14}, a specific geometrical transformation $\Mc(\pb) \colon \Omega_{\mathrm{ref}} \to \Omega(\pb)$ maps a fixed domain $\Omega_{\mathrm{ref}}$ into the geometrically parametrized domain $\Omega(\pb)$, and then allows defining the weak problem \eqref{eq:weakform} in a tensor product space and applying the PGD method developed in Section~\ref{section:PGDcomputation}. In \cite{Nou11b,Che13a}, fictitious domain approaches, which consist in embedding the parametrized domain $\Omega(\pb)$ into a fixed domain $\Omega_{\mathrm{ref}}$, are combined with tensor-based methods such as the PGD method described in Section~\ref{section:PGDcomputation}, in order to reformulate the weak problem \eqref{eq:weakform} on a fixed fictitious domain and then construct optimal separated representations of the solution.

The proposed procedure can also be readily extended to the case where material parameters $c(\xb,\pb)$ and/or $r(\xb,\pb)$ are heterogeneous (non-uniform in space) inside the domain $\Omega$ and can be expressed as separated representations of the form $\sum_j \chi_j(\xb) \delta_j(\pb)$. In that case, space functions $\chi_j$ should be associated with space modes $\psi_i^h$ in \eqref{eq:spaceproperty} before inverting the system.

\begin{remark}
For steady-state diffusion problems, the system \eqref{eq:globalspaceproperty} simply reads:
\begin{equation}
\intO \hat{\Bb}_m^h \cdot \boldsymbol{\nabla}\psi^{\ast} \dO = 0 \quad \forallin{\psi^{\ast}}{\Vc_h},
\end{equation}
so that fluxes $\hat{\bb}^h_i$ are self-equilibrated (in a FE sense). Consequently, a flux $\hat{\qb}_m^h$ verifying the FE equilibrium \eqref{eq:FEequilibrium} can be merely defined as:
\begin{equation}
\hat{\qb}_m^h(\xb,\pb) = \hat{\rb}^h(\xb) + \sum_{i=1}^m \hat{\bb}^h_i(\xb) \prod_j \eta_{j,i}(p_j),
\end{equation}
where parameter functions $\eta_{j,i}$ may \textit{a priori} be defined arbitrarily; nevertheless, they could be adequately chosen and optimized by minimizing an appropriate quadratic functional related to the CRE measure introduced in the next section.
\end{remark}

An alternative approach allowing to construct a statically admissible flux $\hat{\qb}_m^h$ in the FE sense was proposed in \cite{Lad11}. It relies on an additional PGD approach based on a nonclassical static (dual) formulation of the problem \eqref{eq:initpb}. Consequently, it requires to perform another power iterations algorithm resulting in more costly computations compared to the proposed approach.

\section{Global error estimation and adaptive strategy}
\label{section:globalerror}

In this section, we wish to define tools in order to assess the error between the exact solution $u$ and the approximate solution $u_m^{h,\Delta t}$ computed by means of the PGD reduction method. On the one hand, \textit{a priori} error estimates defined in \cite{Lad99b} can be used to assess convergence properties of the PGD approximation but do not give quantitative and useful information for design purposes. On the other hand, classical and relatively simple \textit{a posteriori} error estimates using norms of the discretized residual may not be enough accurate and enable to assess the PGD truncation error alone (without taking discretization error into account). In the following, we define a robust error estimate based on the Constitutive Relation Error (CRE). 

\subsection{Basics on Constitutive Relation Error}

The CRE concept, explained in full details in \cite{Lad04} and shared with various methods in the literature (\eg equilibrated residual method \cite{Ban85,Ain93ter}, flux-free approach \cite{Par06,Par09,Moi09,Cot09}), is currently the only way to get guaranteed and computable error bounds on the solution $u$. For the diffusion-reaction problem we consider, it applies to an admissible triplet solution $(\hat{u},\hat{\qb},\hat{z})$ of the space weak formulation of \eqref{eq:initpb}, \ie a solution verifying the following initial conditions, boundary conditions, as well as equilibrium equations:
\begin{subequations}\label{eq:spaceweakform}
\begin{align}
\hat{u} &= 0 \quad \text{on } \Omega \times \set{0} \times P, \\
\hat{u} &= 0 \quad \text{on } \partial_u \Omega \times I \times P, \\
\intO \(c \frac{\partial \hat{u}}{\partial t} v + \hat{\qb} \cdot \boldsymbol{\nabla}v + \hat{z} v\) \dO &= \intO f_d v \dO + \int_{\partial_q \Omega} g_d v \dS \quad \forallin{v}{\Vc}, \quad \text{on } I \times P. \label{eq:spacefullequilibrium}
\end{align}
\end{subequations}
Only constitutive relations $\qb=k\boldsymbol{\nabla}u$ and $z=ru$ associated to diffusion and reaction mechanisms, respectively, are relaxed for an admissible triplet solution $(\hat{u},\hat{\qb},\hat{z})$ of \eqref{eq:spaceweakform}.
The associated CRE measure $E_{\CRE}$, depending on $\pb$ and computed from any admissible solution $(\hat{u},\hat{\qb},\hat{z})$, is then defined as:
\begin{equation}
\begin{aligned}
E_{\CRE}^2 = \intI e_{\CRE}^2 \dt, \quad \text{with }
e_{\CRE}^2 &= \intO \(\frac{1}{k}(\hat{\qb}-k\boldsymbol{\nabla}\hat{u})^2 + \frac{1}{r}(\hat{z}-r\hat{u})^2\) \dO \\
&= \norm{\hat{\qb}-k\boldsymbol{\nabla}\hat{u}}^2_{k^{-1}} + \norm{\hat{z}-r\hat{u}}^2_{r^{-1}},
\end{aligned}
\end{equation}
where $\norm{\cdot}_{k^{-1}}$ and $\norm{\cdot}_{r^{-1}}$ are energy norms (or equivalently weighted $L^2$-norms) in the space domain $\Omega$ defined by:
\begin{equation}
\norm{\qb}^2_{k^{-1}} = \intO \frac{1}{k} \qb^2 \dO \quad \text{and} \quad \norm{z}^2_{r^{-1}} = \intO \frac{1}{r} z^2 \dO.
\end{equation}
Similarly, we define the following energy norms $\norm{\cdot}_k$ and $\norm{\cdot}_r$ as:
\begin{equation}
\norm{u}^2_k = \intO k (\boldsymbol{\nabla}u)^2 \dO \quad \text{and} \quad \norm{u}^2_r = \intO r u^2 \dO.
\end{equation}
Noticing that
\begin{equation}
\begin{aligned}
\norm{\hat{\qb}-k\boldsymbol{\nabla}\hat{u}}^2_{k^{-1}} &= \norm{\qb-\hat{\qb}}^2_{k^{-1}} + \norm{u-\hat{u}}^2_k - 2 \intO (\qb-\hat{\qb}) \cdot \boldsymbol{\nabla}(u-\hat{u}) \dO, \\
\norm{\hat{z}-r\hat{u}}^2_{r^{-1}} &= \norm{z-\hat{z}}^2_{r^{-1}} + \norm{u-\hat{u}}^2_r - 2 \intO (z-\hat{z}) (u-\hat{u}) \dO \\
\text{and} \quad \intO c \frac{\partial(u-\hat{u})^2}{\partial t} \dO &= 2 \intO c \frac{\partial(u-\hat{u})}{\partial t} (u-\hat{u}) \dO \\
&= - 2 \intO (\qb-\hat{\qb}) \cdot \boldsymbol{\nabla}(u-\hat{u}) \dO - 2 \intO (z-\hat{z}) (u-\hat{u}) \dO,
\end{aligned}
\end{equation}
we get an extension of the well-known Prager-Synge theorem \cite{Pra47} to parametrized evolution (time-dependent) problems, linking the CRE measure $E_{\CRE}$ to a global measure of the discretization error $u-\hat{u}$ over the space-time domain $\Omega \times I$:
\begin{equation}\label{eq:Pragerinequality}
\begin{aligned}
E_{\CRE}^2 &= \trinorm{\qb-\hat{\qb}}^2_{k^{-1}} + \trinorm{z-\hat{z}}^2_{r^{-1}} + \trinorm{u-\hat{u}}^2_k + \trinorm{u-\hat{u}}^2_r + \norm{u-\hat{u}}^2_c \\
&\geq \trinorm{u-\hat{u}}^2_k + \trinorm{u-\hat{u}}^2_r + \norm{u-\hat{u}}^2_c \eqqcolon \trinorm{u-\hat{u}}^2, \\
\end{aligned}
\end{equation}
where $\trinorm{\cdot}^2_{\square} = \intI \norm{\cdot}^2_{\square} \dt$ are energy norms in the space-time domain $\Omega \times I$, and $\norm{u-\hat{u}}^2_c = \intO c(u - \hat{u})^2_{\restrictto{t=T}} \dO$. Similarly, an extension of the Prager-Synge equality (see \cite{Lad04}) can be derived from \eqref{eq:Pragerinequality} by introducing fluxes $\hat{\qb}^{\ast}=\dfrac{1}{2}\[\hat{\qb}+k\boldsymbol{\nabla}\hat{u}\]$ and $\hat{z}^{\ast}=\dfrac{1}{2}\[\hat{z}+r\hat{u}\]$:
\begin{equation}\label{eq:Pragerequality}
\frac{1}{2}E_{\CRE}^2 = 2 \trinorm{\qb-\hat{\qb}^{\ast}}^2_{k^{-1}} + 2 \trinorm{z-\hat{z}^{\ast}}^2_{r^{-1}} + \norm{u-\hat{u}}^2_c.
\end{equation}
\begin{remark}
For steady-state diffusion problems, the previous equalities simply read:
\begin{equation}
E_{\CRE}^2 = \norm{\qb-\hat{\qb}}^2_{k^{-1}} + \norm{u-\hat{u}}^2_k = 4 \norm{\qb-\hat{\qb}^{\ast}}^2_{k^{-1}}.
\end{equation}
\end{remark}

\vspace{1em}

Consequently, the CRE measure $E_{\CRE}$ defines an upper bound (guaranteed estimate) of the global error (measured in the energy norm) between the exact solution $u$ and its approximation $\hat{u}$.

\subsection{Global error estimator on the PGD approximation}

The technical point in the CRE framework is the computation of an admissible triplet denoted by $(\hat{u}_m,\hat{\qb}_m,\hat{z}_m)$. Such an admissible solution can be defined from the data and the approximate solution fields $(u_m^{h,\Delta t},\hat{\qb}_m^{h})$ provided by the PGD method. On the one hand, $u_m^{h,\Delta t}$ satisfies initial and Dirichlet boundary conditions, so that we usually choose $\hat{u}_m=u_m^{h,\Delta t}$. On the other hand, considering $\hat{\qb}_m^{h}$ which verifies equilibrium with $u_m^{h,\Delta t}$ in the FE sense, and applying classical equilibration techniques used in the FEM context (see \cite{Lad04,Lad10bis,Ple11} and the references therein for an overview of the topic), it is possible to construct an admissible flux $\hat{\qb}_m$ that satisfies full equilibrium with $u_m^{h,\Delta t}$ for all $(t,\pb) \in I \times P$:
\begin{equation}\label{eq:fullequilibrium}
\intO \(c \frac{\partial u_m^{h,\Delta t}}{\partial t} v + \hat{\qb}_m\cdot \boldsymbol{\nabla}v + r u_m^{h,\Delta t} v\) \dO = \intO f_d v \dO + \int_{\partial_q \Omega} g_d v \dS \quad \forallin{v}{\Vc}.
\end{equation}
The admissible flux $\hat{\qb}_m$ then takes the following separated form:
\begin{equation}
\hat{\qb}_m(\xb,t,\pb) = \hat{\rb}(\xb,t) + \sum_{i=1}^m \hat{\boldsymbol{\phi}}_i(\xb) \(c \dot{\lambda}^{\Delta t}_i(t) + r \lambda^{\Delta t}_i(t)\) \prod_j \gamma_{j,i}(p_j), \quad \text{with } \hat{\rb}(\xb,t) = \sum_{s=1}^S \alpha_s(t) \hat{\rb}_{s}(\xb),
\end{equation}
where $\hat{\rb}_{s}$ (resp. $\hat{\boldsymbol{\phi}}_i$) is obtained from a post-processing of $\hat{\rb}^h_{s}$ (resp. $\hat{\boldsymbol{\phi}}^h_i$), and $\alpha_s$ are the time functions involved in the separated form~\eqref{eq:radialapproxloading} of functional $L$. Here, we use hybrid-flux (also referred to as Element Equilibration Technique (EET) or Element Equilibration + Star-Patch Technique EESPT in \cite{Ple11}) strategies which lean on an energy relation, called \textit{extension condition} (or \textit{prolongation condition}), that takes the following form \cite{Lad04}:
\begin{equation}
\int_E \(\hat{\taub}-\hat{\taub}^h\) \cdot \boldsymbol{\nabla}\varphi^h_i \dO = 0 \quad \forallin{E}{\Omega_h}, \quad \forallin{i}{E},
\end{equation}
where $E$ is an element of the space mesh $\Omega_h$, $\varphi^h_i$ is the FE shape function associated to any node $i$ connected to element $E$, and $\taub$ is a flux of interest ($\rb_{s}$ or $\boldsymbol{\phi}_j$). $\hat{\taub}^h$ (resp. $\hat{\taub}$) refers to a given (resp. researched) approximation of $\taub$ satisfying FE (resp. full) equilibrium with $u_m^{h,\Delta t}$ in the sense of \eqref{eq:FEequilibrium} (resp. \eqref{eq:fullequilibrium}). This condition, in addition to FE properties of $\hat{\taub}^h$, enables to determine equilibrated fluxes on element edges by solving well-posed local algebraic systems associated with each node of the space mesh $\Omega_h$. Eventually, the equilibrated flux $\hat{\taub}$ is computed at the element level by solving Neumann problems that involve the pre-computed fluxes over element boundaries. Further details on the practical construction of admissible fluxes can be found in \cite{Ple11}.\\

Consequently, according to \eqref{eq:Pragerinequality}, the CRE measure $E_{\CRE}$ computed from the admissible triplet solution $(u_m^{h,\Delta t}, \hat{\qb}_m, \hat{z}_m=ru_m^{h,\Delta t})$ verifies for all $\pb \in P$:
\begin{equation}\label{eq:Prager}
\begin{aligned}
E_{\CRE}^2 &= \trinorm{\hat{\qb}_m-k\boldsymbol{\nabla}u_m^{h,\Delta t}}^2_{k^{-1}} \\
&\geq \trinorm{u-u_m^{h,\Delta t}}^2_k + \trinorm{u-u_m^{h,\Delta t}}^2_r + \norm{u-u_m^{h,\Delta t}}^2_c = \trinorm{u-u_m^{h,\Delta t}}^2
\end {aligned}
\end{equation}
and therefore defines, for any $\pb \in P,$ a guaranteed estimate (upper bound) of the global error (measured in the energy norm) between the exact solution $u$ and the approximate PGD solution $u_m^{h,\Delta t}$ of the problem.

\subsection{Adaptive algorithm}\label{section:adaptiveglobal}
 
The error estimate $E_{\CRE}^2$ previously defined in \eqref{eq:Prager} takes into account the various error sources inherent to the PGD approach, \ie including:
\begin{itemize}
\item the PGD truncation error (indicated by subscript $m$) due to the restriction of the PGD modal representation to a limited (finite) number of modes $m$;
\item the discretization error (indicated by superscripts $h$ and $\Delta t$) due to the use of numerical methods (FEM here) to compute PGD modes, this error source being itself split into space discretization and time discretization errors.
\end{itemize}

\vspace{1em}

The other possible error sources (algebraic error due to the use of iterative solvers, numerical integration error, roundoff error due to machine precision, \dots) are assumed to be negligible compared to both discretization and PGD truncation errors.

In this section, we wish to assess the relative contribution of each error source to the error estimate $E_{\CRE}^2$ in order to efficiently drive an adaptive algorithm based on a greedy approach. For that purpose, we follow a natural procedure which has been already employed in previous works dealing with error estimation and adaptivity within the CRE framework \cite{Gal96,Lad98c,Pel00}. The idea consists of introducing specific error indicators based on the CRE concept applied to admissible fields in the sense of intermediate reference problems (weaker sense compared to the initial reference problem \eqref{eq:initpb}). In practice, we consider a new (intermediate) reference problem defined as the discrete space-time weak formulation of the initial reference problem \eqref{eq:initpb} using the same discretization space $\Vc_h \otimes \Tc_{\Delta t}$ as for the PGD approximation $u_m^{h,\Delta t}$; it reads for all $\pb \in P$:
\begin{equation}\label{eq:discreteweakform}
\text{find $u^{h,\Delta t} \in \Vc_h \otimes \Tc_{\Delta t}$ such that }
b(u^{h,\Delta t},v) = l(v) \quad \forallin{v}{\Vc_h \otimes \Tc_{\Delta t}},
\end{equation}
with $u^{h,\Delta t}_{\restrictto{t=0}} = 0$. Using the following Galerkin orthogonality condition defined for all $\pb \in P$:
\begin{equation}
b(u-u^{h,\Delta t},v) = \intI \intO \(c \frac{\partial(u-u^{h,\Delta t})}{\partial t} v + k \boldsymbol{\nabla}(u-u^{h,\Delta t}) \cdot \boldsymbol{\nabla}v + r (u-u^{h,\Delta t}) v\) \dO \dt = 0 \quad \forallin{v}{\Vc_h \otimes \Tc_{\Delta t}},
\end{equation}
with $v=u^{h,\Delta t}-u_m^{h,\Delta t}$, we obtain:
\begin{equation}\label{eq:splierror}
\trinorm{u-u_m^{h,\Delta t}}^2 = \trinorm{u-u^{h,\Delta t}}^2 + \trinorm{u^{h,\Delta t}-u_m^{h,\Delta t}}^2 + \intI \intO c (u-u^{h,\Delta t}) \frac{\partial(u^{h,\Delta t}-u_m^{h,\Delta t})}{\partial t} \dO \dt,
\end{equation}
where the norm $\trinorm{\cdot}$ is defined in \eqref{eq:Pragerinequality}. The global error $\Delta^2 = \trinorm{u-u_m^{h,\Delta t}}^2$ can then be split into:
\begin{equation}
\Delta^2 = \Delta_{\PGD}^2 + \Delta_{\dis}^2,
\end{equation}
where $\Delta_{\PGD}^2 = \trinorm{u^{h,\Delta t}-u_m^{h,\Delta t}}^2$ quantifies the error coming from the PGD truncation alone ($\Delta_{\PGD} \to 0$ when $m \to +\infty$) whereas $\Delta_{\dis}^2 = \Delta^2 - \Delta_{\PGD}^2 = \trinorm{u-u_m^{h,\Delta t}}^2 - \trinorm{u^{h,\Delta t}-u_m^{h,\Delta t}}^2$ quantifies the error coming from the space-time discretization alone ($\Delta_{\dis} \to 0$ when $h \to 0$ and $\Delta t \to 0$). \\

Both contributions $\Delta_{\PGD}$ and $\Delta_{\dis}$ can be assessed from the CRE property \eqref{eq:Pragerinequality} and a direct post-processing of available approximate fields $u_m^{h,\Delta t}$, $\hat{\qb}_m^h$ and $\hat{\qb}_m$:
\begin{itemize}
\item we first compute the CRE measure $E^{h,\Delta t}_{\CRE}$ applied to a pair $(u_m^{h,\Delta t},\hat{\qb}^{h,\Delta t}_m)$ which is admissible in the FE sense of the intermediate reference problem \eqref{eq:discreteweakform}, \ie satisfying the following equilibrium equations for all $\pb \in P$:
\begin{equation}\label{eq:discreteFEequilibrium}
\intI \(\intO \(c \frac{\partial u_m^{h,\Delta t}}{\partial t} v + \hat{\qb}^{h,\Delta t}_m \cdot \boldsymbol{\nabla}v + r u_m^{h,\Delta t} v\) \dO - \intO f_d v \dO - \int_{\partial_q \Omega} g_d v \dS\) \dt = 0 \quad \forallin{v}{\Vc_h \otimes \Tc_{\Delta t}}.
\end{equation}
Noticing that the flux $\hat{\qb}_m^h$ has been constructed in such a way that it verifies the FE equilibrium \eqref{eq:FEequilibrium} for all $t \in I$, and denoting by $\Nb = [N_0(t),N_{1}(t),\dots,N_N(t)]^T$ ($N$ being the number of time steps: $T=N\Delta t$) the vector of shape functions used for the FEM in the time domain $I$, the flux $\hat{\qb}^{h,\Delta t}_m$ may be recovered as a simple post-processing (linear interpolation) of $\hat{\qb}_m^h$ over $I = \intervalcc{0}{T}$:
\begin{equation}\label{eq:pondN}
\hat{\qb}^{h,\Delta t}_m = \(\intI \Nb^T \Nb\dt \)^{-1} \Nb^T \Gb, \quad \text{with } \Gb = [\Gb_0,\Gb_{1},\dots,\Gb_N]^T \text{ and } \Gb_{j} = \intI \hat{\qb}_m^h N_{j} \dt.
\end{equation}
We then define an indicator $\eta_{\PGD}$ of the PGD truncation error $\Delta_{\PGD}$ as:
\begin{equation}
\Delta_{\PGD}^2 \approx {E^{h,\Delta t}_{\CRE}}^2 = \trinorm{\hat{\qb}^{h,\Delta t}_m-k\boldsymbol{\nabla}u_m^{h,\Delta t}}^2_{k^{-1}} \eqqcolon \eta^2_{\PGD};
\end{equation}
\item we eventually deduce an indicator $\eta_{\dis}$ of the discretization error $\Delta_{\dis}$ as:
\begin{equation}
\Delta_{\dis}^2 \approx E_{\CRE}^2-\eta^2_{\PGD} = \trinorm{\hat{\qb}_m-k\boldsymbol{\nabla}u_m^{h,\Delta t}}^2_{k^{-1}} - \trinorm{\hat{\qb}^{h,\Delta t}_m-k\boldsymbol{\nabla}u_m^{h,\Delta t}}^2_{k^{-1}} \eqqcolon \eta_{\dis}^2.
\end{equation}
Furthermore, we also propose consistent error indicators $\eta_{h}$ and $\eta_{\Delta t}$ in order to quantify the error contributions coming from space and time discretizations, respectively; they read:
\begin{equation}
\eta_{h}^2 = \trinorm{\hat{\qb}_m-\hat{\qb}_m^h}^2_{k^{-1}} \quad \text{and} \quad \eta_{\Delta t}^2 = \eta_{\dis}^2 - \eta_{h}^2.
\end{equation}
It is then possible to identify whether the source of the discretization error is due to the spatial or the time discretization.
\end{itemize}

\vspace{1em}

\begin{remark}
In the case where an incremental numerical method (time integration scheme) is used to solve evolution problems \eqref{eq:timepb} and compute time functions $\lambda_i$ ($1 \leq i \leq m$), it is possible to define an equivalent variational formulation by introducing some weight functions \cite{Zie05}. Then, one can use the simple recovering technique described in \eqref{eq:pondN} to build an admissible field $\hat{\qb}^{h,\Delta t}_m$ in the sense of the incremental numerical method.
\end{remark}

\vspace{1em}

Previously defined error indicators can then be used as stopping criteria or adaptation indicators in a greedy algorithm in order to adaptively construct a suitable PGD approximation $u_m^{h,\Delta t}$ and thus control the computation process, searching for the highest error contributions (between PGD truncation and modes discretizations) and minimizing them. In practice, after computing each mode $i \in \{1,\dots,m\}$, the adaptive strategy is conducted as follows:
\begin{enumerate}
\item we identify the parameter set $\pb_{\mathrm{max}} = \argmax_{\pb\in P} E_{\CRE}(\pb)$ (worst case scenario) as well as associated error indicators;
\item we compare the relative contributions of $\eta_{\PGD}(\pb_{\mathrm{max}})$ and $\eta_{\dis}(\pb_{\mathrm{max}})$ to the error estimate $E_{\CRE}(\pb_{\mathrm{max}})$:
\begin{itemize}
\item if $\eta_{\PGD}(\pb_{\mathrm{max}}) \geq \eta_{\dis}(\pb_{\mathrm{max}})$, the mode $(i+1)$ is computed keeping the same space-time discretization as for the mode $i$;
\item otherwise, \ie if $\eta_{\PGD}(\pb_{\mathrm{max}}) < \eta_{\dis}(\pb_{\mathrm{max}})$, discretization parameters $(h,\Delta t)$ are modified in order to compute the mode $i$ and next modes with better accuracy. For that, optimal mesh adaptation techniques based on local contributions of $\eta_{h}$ and $\eta_{\Delta t}$, as well as on asymptotic convergence rates predicted by \textit{a priori} error estimates, are used to reach a given discretization error threshold \cite{Lad04}. In practice, this error threshold is chosen as $\alpha \, \eta_{\PGD}(\pb_{\mathrm{max}})$ with $\alpha\in \intervaloc{0}{1}$ a scalar parameter to set. We emphasize that the $(i-1)$ modes previously computed with coarser meshes are kept unchanged.
\end{itemize}
\end{enumerate}
The adaptive procedure is performed until $\max_{\pb\in P} E_{\CRE}(\pb)\leq \gamma_{\mathrm{tol}}$ for a given mode $m_{\mathrm{tol}}$, where $\gamma_{\mathrm{tol}}$ is a predefined error tolerance.

The identification of the parameter set $\pb_{\mathrm{max}}$ (corresponding to the parameter values which maximize the error estimate $E_{\CRE}(\pb)$ over $\pb\in P$) may be a computational issue for high-dimensional parametric approximation space (\ie for large dimensions $n_p$). For low-dimensional parametric approximation space, the screening of the parameter domain $P$ can be easily performed over a multidimensional numerical grid as the number of dimensions $n_p$ is limited. Conversely, for high-dimensional parametric approximation space, it would require dedicated strategies and specific algorithms such as those used for solving optimization problems in large dimension or in the empirical interpolation method (EIM) \cite{Bar04}.

The adaptive process could be optimized comparing $\eta_{\PGD}(\pb)$ and $\eta_{\dis}(\pb)$ for each value $\pb\in P$, leading to parametrized adaptivity and refined meshes. However, this strategy looks complex to implement and use in practice.

\begin{remark}
After performing mesh adaptation, the intermediate reference problem (with discrete space-time weak formulation) used to compute error indicators is changed. Consequently, for modes computed before mesh adaptation, part of the discretization error (the one which can be captured with the new finer mesh) is transferred into the indicator $\eta_{\PGD}$ of the PGD truncation error. This procedure is consistent with the definition of discretization error (which should tend to zero when the mesh size used for the current PGD mode goes to zero).
\end{remark}

\section{Extension to goal-oriented error estimation}
\label{section:GOerror}

Error measured in a global (energy) norm is clearly not the best criterion for control and adaptivity when one is interested in specific outputs of interest of the problem. In this section, we define a goal-oriented error estimation method, based on the classical extraction (adjoint-based) technique \cite{Bec96,Par97,Per98,Pru99,Bec01} and the CRE concept, which is dedicated to the accurate and robust computation of outputs of interest from a PGD approximation. 

\subsection{Adjoint problem and associated PGD approximation} 

We consider a functional output of interest $Q$, possibly depending explicitly on $\pb$, linear and continuous with respect to $u$, which is defined globally (over the space-time domain $\Omega \times I$) by means of an extraction pair $(\qb_{\Sigma},f_{\Sigma})$:
\begin{equation}
Q(u)= \intI \intO \(\qb_{\Sigma} \cdot \boldsymbol{\nabla}u + f_{\Sigma} u\) \dO \dt.
\end{equation}
Space-time functions $\qb_{\Sigma}$ and $f_{\Sigma}$, referred to as extraction operators or extractors, may be defined explicitly or implicitly (depending on the quantity of interest $Q$) possibly using Dirac distributions in the case of pointwise outputs. \\

\begin{remark}
For nonlinear functionals $Q$ with respect to $u$, a classical approach consists in using linearization techniques and yields unguaranteed local error bounds \cite{Joh92,Mad99}. An alternative approach, valid only for non-linear pointwise quantities in space, relies on projection procedures and allows to recover strict local error bounds \cite{Lad10}.
\end{remark}

\vspace{1em}

Then, following the optimal control approach \cite{Bec01} based on duality arguments, we now introduce the space-time weak formulation of the adjoint problem, defined on the space-time domain $\Omega \times I$, associated with functional $Q$; it reads for all $\pb \in P$:
\begin{equation}\label{eq:spacetimeweakformadj}
\text{find $\tilde{u} \in L^2(I;\Vc)$, with $\dfrac{\partial \tilde{u}}{\partial t} \in L^2(I;\Vc^{\ast})$, such that }
b(v,\tilde{u}) = b^{\ast}(\tilde{u},v) = Q(v) \quad \forallin{v}{L^2(I;\Vc)},
\end{equation}
with $\tilde{u}_{\restrictto{t=T}} = 0$, where bilinear form $b^{\ast}(\cdot,\cdot)$ is the adjoint operator defined on $L^2(I;\Vc)$ by \cite{Gil02}:
\begin{equation}
b^{\ast}(\tilde{u},v) \coloneqq \intI \intO \(-c \frac{\partial \tilde{u}}{\partial t} v + k \boldsymbol{\nabla} \tilde{u} \cdot \boldsymbol{\nabla} v + r \tilde{u} v\) \dO \dt
.
\end{equation}

The adjoint problem \eqref{eq:spacetimeweakformadj} is linear and reverse in time. Furthermore, performing the change of variable $t \to T-t$, this dual problem becomes a diffusion-reaction problem similar to the primal problem \eqref{eq:spacetimeweakform} but with another loading composed of the extractors $(\qb_{\Sigma},f_{\Sigma})$ defining the quantity of interest $Q$. In practice, primal and dual (adjoint) problems may be solved in parallel for computational efficiency.

An approximate solution $\tilde{u}_{m'}^{h',\Delta t'}$, with potentially order $m' \neq m$ and discretization parameters $(h',\Delta t') \neq (h,\Delta t)$, is first computed using the PGD method. An admissible triplet solution $(\hat{\tilde{u}},\hat{\tilde{\qb}},\hat{\tilde{z}}) = (\tilde{u}_{m'}^{h',\Delta t'},\hat{\tilde{\qb}}_{m'},\tilde{z}_{m'}^{h',\Delta t'} = r\tilde{u}_{m'}^{h',\Delta t'})$ is then derived using the equilibration technique presented in Section~\ref{section:globalerror}. In particular, this solution satisfies the following full equilibrium for all $(t,\pb) \in I \times P$:
\begin{equation}\label{eq:fullequilibriumadj}
\intO \(-c \frac{\partial \tilde{u}_{m'}^{h',\Delta t'}}{\partial t} v + \hat{\tilde{\qb}}_{m'} \cdot \boldsymbol{\nabla}v + \tilde{z}_{m'}^{h',\Delta t'} v\) \dO = \intO \(\qb_{\Sigma} \cdot \boldsymbol{\nabla}v + f_{\Sigma} v \) \dO \quad \forallin{v}{\Vc}.
\end{equation}
Consequently, the CRE measure $\tilde{E}_{\CRE}$ expressed as $\tilde{E}_{\CRE} = \trinorm{\hat{\tilde{\qb}}_{m'}-k\boldsymbol{\nabla}\tilde{u}_{m'}^{h',\Delta t'}}_{k^{-1}}$ defines, for any $\pb \in P,$ an upper bound (guaranteed estimate) of the global error (measured in the energy norm) between the exact solution $\tilde{u}$ and the approximate PGD solution $\tilde{u}_{m'}^{h',\Delta t'}$ of the adjoint problem \eqref{eq:spacetimeweakformadj}. \\

\begin{remark}
A non-intrusive approach, henceforth known as handbook techniques, can be introduced for the approximate solution of the adjoint problem \eqref{eq:spacetimeweakformadj} \cite{Cha08,Cha09,Lad10,Wae12}. Noticing that the adjoint loading $(\qb_{\Sigma},f_{\Sigma})$ usually applies on a local subdomain of $\Omega \times I$ (provided the quantity of interest $Q$ refers to local features in space and time of the solution $u$), and therefore leads to an adjoint solution with localized high gradients in the space-time domain, the idea is to introduce local enrichment functions in the vicinity of the space-time region of interest where the quantity $Q$ is defined. This enrichment is particularly well-suited to handle pointwise quantities of interest and yields accurate local error bounds without requiring any regularization (\eg mollification \cite{Pru99}) of the functional being considered or any specific local remeshing technique.
\end{remark}

\subsection{Local error estimator on a functional output of interest}

Using the linearity assumption for the functional $Q$, a representation of the local error $\Delta Q=Q(u)-Q(u_m^{h,\Delta t})$ between the exact value $Q(u)$ and the approximate PGD value $Q(u_m^{h,\Delta t})$ of the output of interest $Q$ can be defined as a weighted residual from the adjoint solution $\tilde{u}$ \cite{Pru99}; it reads for all $\pb \in P$:
\begin{equation}\label{eq:residual}
\Delta Q = Q(u-u_m^{h,\Delta t}) = b^{\ast}(\tilde{u},u-u_m^{h,\Delta t}) = b(u-u_m^{h,\Delta t},\tilde{u}) = l(\tilde{u})-b(u_m^{h,\Delta t},\tilde{u}) \eqqcolon R(\tilde{u}),
\end{equation}
where $R(\cdot)$ is the weak residual functional associated with the primal space-time weak formulation \eqref{eq:spacetimeweakform} of the reference problem. \\

\begin{remark}
A simple PGD truncation error indicator based on \eqref{eq:residual} is defined in \cite{Amm10}; it merely consists in replacing the exact adjoint solution $\tilde{u}$ with an accurate PGD approximation $\tilde{u}_{m'}^{h',\Delta t'}$ of order $m' \geq m$, and potentially with finer space-time discretization $(h',\Delta t') \leq (h,\Delta t)$. Despite its computational efficiency, this method does not provide robust and guaranteed bounds of the local error $\Delta Q$ on the functional $Q$ as the term $R(\tilde{u}-\tilde{u}_{m'}^{h',\Delta t'})$ is neglected. In the following, we use a technique which takes this term into account in the error estimation procedure and thus leads to strict local error bounds.
\end{remark}

\vspace{1em}

Introducing the approximate PGD adjoint solution $\tilde{u}_{m'}^{h',\Delta t'}$ in \eqref{eq:residual}, we get:
\begin{equation}\label{eq:residualbis}
\Delta Q - R(\tilde{u}_{m'}^{h',\Delta t'}) = R(\tilde{u}-\tilde{u}_{m'}^{h',\Delta t'}) = b(u-u_m^{h,\Delta t},\tilde{u}-\tilde{u}_{m'}^{h',\Delta t'}) = b^{\ast}(\tilde{u}-\tilde{u}_{m'}^{h',\Delta t'},u-u_m^{h,\Delta t}),
\end{equation}
where $R(\tilde{u}_{m'}^{h',\Delta t'})$ is a computable correction term involving both approximate PGD solutions $u_m^{h,\Delta t}$ and $\tilde{u}_{m'}^{h',\Delta t'}$ of primal and dual (adjoint) problems, respectively. Following an approach similar to the one introduced in the FEM context in \cite{Lad06,Cha07,Lad08,Cha08,Lad09}, we define an upper bound of $\abs{\Delta Q - R(\tilde{u}_{m'}^{h',\Delta t'})}$ using properties of admissible solutions $(u_{m}^{h,\Delta t},\hat{\qb}_{m},z_{m}^{h,\Delta t})$ and $(\tilde{u}_{m'}^{h',\Delta t'},\hat{\tilde{\qb}}_{m'},\tilde{z}_{m'}^{h',\Delta t'})$, as well as associated CRE measures $E_{\CRE}$ and $\tilde{E}_{\CRE}$. \\
Indeed, noticing that both admissible solution $(\tilde{u}_{m'}^{h',\Delta t'},\hat{\tilde{\qb}}_{m'},\hat{\tilde{z}}_{m'}^{h',\Delta t'})$ and exact solution $(\tilde{u},\tilde{\qb}=k\boldsymbol{\nabla}\tilde{u},\tilde{z}=r\tilde{u})$ verify the full equilibrium \eqref{eq:fullequilibriumadj}, and that $u-u_m^{h,\Delta t} \in \Vc$, for all $(t,\pb) \in I \times P$, we obtain:
\begin{equation}\label{eq:equalfond}
\begin{aligned}
&b^{\ast}(\tilde{u}-\tilde{u}_{m'}^{h',\Delta t'},u-u_m^{h,\Delta t}) \\
&= \intI \intO \(-c \frac{\partial(\tilde{u}-\tilde{u}_{m'}^{h',\Delta t'})}{\partial t} (u-u_m^{h,\Delta t}) + k \boldsymbol{\nabla}(\tilde{u}-\tilde{u}_{m'}^{h',\Delta t'}) \cdot \boldsymbol{\nabla}(u-u_m^{h,\Delta t}) + r (\tilde{u}-\tilde{u}_{m'}^{h',\Delta t'}) (u-u_m^{h,\Delta t})\)\dO \dt \\
&= \intI \intO \((\hat{\tilde{\qb}}_{m'} - k \boldsymbol{\nabla}\tilde{u}_{m'}^{h',\Delta t'}) \cdot \boldsymbol{\nabla}(u-u_m^{h,\Delta t}) + (\hat{\tilde{z}}_{m'}^{h',\Delta t'} - r \tilde{u}_{m'}^{h',\Delta t'}) (u-u_m^{h,\Delta t})\) \dO \dt \\
&= \intI \intO \(\frac{1}{k} (\hat{\tilde{\qb}}_{m'} - k \boldsymbol{\nabla}\tilde{u}_{m'}^{h',\Delta t'}) \cdot (\qb - k \boldsymbol{\nabla}u_m^{h,\Delta t}) + \frac{1}{r} (\hat{\tilde{z}}_{m'}^{h',\Delta t'} - r \tilde{u}_{m'}^{h',\Delta t'}) (z-ru_m^{h,\Delta t})\) \dO \dt.
\end{aligned}
\end{equation}
Consequently, a first bounding can be easily derived by applying the classical space-time Cauchy-Schwarz inequality:
\begin{equation}\label{eq:bounds1}
\begin{aligned}
\abs{\Delta Q - R(\tilde{u}_{m'}^{h',\Delta t'})} &\leq \sqrt{\trinorm{\hat{\tilde{\qb}}_{m'} - k \boldsymbol{\nabla}\tilde{u}_{m'}^{h',\Delta t'}}^2_{k^{-1}} + \trinorm{\hat{\tilde{z}}_{m'}^{h',\Delta t'} - r \tilde{u}_{m'}^{h',\Delta t'}}^2_{r^{-1}}} \times \sqrt{\trinorm{u-u_m^{h,\Delta t}}^2_{k} + \trinorm{u-u_m^{h,\Delta t}}^2_{r}} \\
&\leq \tilde{E}_{\CRE} E_{\CRE}.
\end{aligned}
\end{equation}
Observing that both admissible solution $(u_{m}^{h,\Delta t},\hat{\qb}_{m},\hat{z}_{m}^{h,\Delta t})$ and exact solution $(u,\qb=k\boldsymbol{\nabla}u,z=ru)$ satisfy the full equilibrium \eqref{eq:fullequilibrium}, and that $\tilde{u}_{m'}^{h',\Delta t'} \in \Vc$, for all $(t,\pb) \in I \times P$, the correction term $R(\tilde{u}_{m'}^{h',\Delta t'})$ can be rewritten in terms of both approximate PGD solutions $u_m^{h,\Delta t}$ and $\tilde{u}_{m'}^{h',\Delta t'}$ as well as the admissible flux pair $(\hat{\qb}_{m},\hat{z}_m^{h,\Delta t})$ only:
\begin{equation}\label{eq:approxresidual}
\begin{aligned}
R(\tilde{u}_{m'}^{h',\Delta t'}) &= b(u-u_m^{h,\Delta t},\tilde{u}_{m'}^{h',\Delta t'}) \\
&= \intI \intO \(c \frac{\partial(u-u_m^{h,\Delta t})}{\partial t} \tilde{u}_{m'}^{h',\Delta t'} + k \boldsymbol{\nabla}(u-u_m^{h,\Delta t}) \cdot \boldsymbol{\nabla}\tilde{u}_{m'}^{h',\Delta t'} + r (u-u_m^{h,\Delta t}) \tilde{u}_{m'}^{h',\Delta t'}\) \dO \dt \\
&= \intI \intO \((\hat{\qb}_m - k \boldsymbol{\nabla}u_m^{h,\Delta t}) \cdot \boldsymbol{\nabla}\tilde{u}_{m'}^{h',\Delta t'} + (\hat{z}_m^{h,\Delta t} - r u_m^{h,\Delta t}) \tilde{u}_{m'}^{h',\Delta t'}\) \dO \dt.
\end{aligned}
\end{equation}

A more accurate bounding can also be deduced from \eqref{eq:equalfond} by introducing the pair $(\hat{\qb}^{\ast}_m,\hat{z}_m^{h,\Delta t\ast}) = \dfrac{1}{2}\[(\hat{\qb}_m,\hat{z}_m^{h,\Delta t})+(k\boldsymbol{\nabla}u_m^{h,\Delta t},ru_m^{h,\Delta t})\]$ and using the extended Prager-Synge equality \eqref{eq:Pragerequality}; it reads:
\begin{equation}\label{eq:bounds2}
\begin{aligned}
\abs{\Delta Q - Q_{\corr}} &= \labs{\intI \intO \(\frac{1}{k} (\hat{\tilde{\qb}}_{m'}-k\boldsymbol{\nabla}\tilde{u}_{m'}^{h',\Delta t'}) \cdot (\qb-\hat{\qb}^{\ast}_m)+\frac{1}{r} (\hat{\tilde{z}}_{m'}^{h',\Delta t'}-r\tilde{u}_{m'}^{h',\Delta t'})(z-\hat{z}_m^{h,\Delta t\ast}) \)\dO\dt} \\
&\leq \sqrt{\trinorm{\hat{\tilde{\qb}}_{m'}-k\boldsymbol{\nabla}\tilde{u}_{m'}^{h',\Delta t'}}^2_{k^{-1}}+\trinorm{\hat{\tilde{z}}_{m'}^{h',\Delta t'}-r\tilde{u}_{m'}^{h',\Delta t'}}^2_{r^{-1}}} \times 
\sqrt{\trinorm{\qb-\hat{\qb}_m^{\ast}}^2_{k^{-1}}+\trinorm{z-\hat{z}_m^{h,\Delta t\ast}}^2_{r^{-1}}} \\
&\leq \frac{1}{2} \tilde{E}_{\CRE} E_{\CRE},
\end{aligned} 
\end{equation}
where $Q_{\corr}$ is a computable correction term defined as:
\begin{equation}
\begin{aligned}
Q_{\corr} &= R(\tilde{u}_{m'}^{h',\Delta t'}) + 
\intI \intO \(\frac{1}{k} (\hat{\tilde{\qb}}_{m'}-k\boldsymbol{\nabla}\tilde{u}_{m'}^{h',\Delta t'}) \cdot (\hat{\qb}^{\ast}_m-k\boldsymbol{\nabla}u_m^{h,\Delta t}) + \frac{1}{r} (\hat{\tilde{z}}_{m'}^{h',\Delta t'}-r\tilde{u}_{m'}^{h',\Delta t'}).(\hat{z}_m^{h,\Delta t\ast}-ru_m^{h,\Delta t})\) \dO \dt \\
&= \intI \intO \(\frac{1}{k} (\hat{\qb}_m-k\boldsymbol{\nabla}u_m^{h,\Delta t}) \cdot \hat{\tilde{\qb}}^{\ast}_{m'} + \frac{1}{r} (\hat{z}_m^{h,\Delta t}-ru_m^{h,\Delta t})\hat{\tilde{z}}_{m'}^{h',\Delta t'\ast}\) \dO \dt.
\end{aligned}
\end{equation}
Consequently, \eqref{eq:bounds2} provides guaranteed bounds $\rho$ on the local error $\Delta Q$ (or directly on the exact output $Q(u)$) for any value of parameters $\pb\in P$:
\begin{equation}
\abs{\Delta Q} \leq \abs{Q_{\corr} \pm \frac{1}{2} \tilde{E}_{\CRE} E_{\CRE}} \eqqcolon \rho.
\end{equation}
It is worthwhile to point out the effects of the correcting term $Q_{\corr}$: having at hand an accurate approximate solution $(\tilde{u}_{m'}^{h',\Delta t'},\hat{\tilde{\qb}}_{m'},\hat{\tilde{z}}_{m'}^{h',\Delta t'})$ of the adjoint problem enables to compute very sharp bounds on $\Delta Q$ (or on $Q(u)$) for any $\pb \in P$ as $\tilde{E}_{\CRE}$ tends to zero when $(\tilde{u}_{m'}^{h',\Delta t'},\hat{\tilde{\qb}}_{m'},\hat{\tilde{z}}_{m'}^{h',\Delta t'})$ converges toward the exact adjoint solution $(\tilde{u},\tilde{\qb},\tilde{z})$ and $Q_{\corr}$ tends to an asymptotic value $R(\tilde{u})$ equal to the exact local error $\Delta Q$. Indeed, we have $Q_{\corr}=R(\tilde{u})=\Delta Q$ when the dual (adjoint) problem is solved exactly. We also notice that $Q_{\corr}=0=\Delta Q$ when the primal problem is solved exactly. \\

\begin{remark}
Sharper bounds, based on alternative strategies to the classical space-time Cauchy-Schwarz inequality \eqref{eq:bounds1} or \eqref{eq:bounds2}, can still be obtained using tools introduced in previous works. In \cite{Lad13}, a bounding technique based on Saint-Venant's principle and properties on homothetic domains was proposed to improve bounding in the space domain $\Omega$. In \cite{Cha07}, weighted CRE functionals and Legendre-Fenchel inequality were used to improve bounding in the time domain $I$. Nevertheless, these bounds will not be considered in the following.
\end{remark}

\subsection{Adaptive algorithm}

Here again, it is still possible to assess the relative contributions of error sources that contribute to the total local error $\Delta Q$ on a given quantity of interest $Q$ (\ie PGD truncation error and discretization error). Indeed, considering a linear functional $Q$ and introducing the approximate FE value $Q(u^{h,\Delta t})$ of the output of interest $Q$ obtained from the approximate FE solution $u^{h,\Delta t}$ of the discrete space-time weak formulation \eqref{eq:discreteweakform} of the reference problem, the local error $\Delta Q = Q(u)-Q(u_m^{h,\Delta t})$ can be split into:
\begin{equation}
\Delta Q = \Delta Q_{\PGD} + \Delta Q_{\dis},
\end{equation}
where $\Delta Q_{\PGD} = Q(u^{h,\Delta t})-Q(u_m^{h,\Delta t})$ is part of the local error due to truncation in the PGD decomposition alone ($\Delta Q_{\PGD} \to 0$ when $m \to +\infty$) whereas $\Delta Q_{\dis} = \Delta Q - \Delta Q_{\PGD} = Q(u)-Q(u^{h,\Delta t})$ is the one due to space-time discretization alone ($\Delta Q_{\dis} \to 0$ when $h \to 0$ and $\Delta t \to 0$). \\

Both contributions $\Delta Q_{\PGD}$ and $\Delta Q_{\dis}$ can be easily assessed from previously defined techniques:
\begin{itemize}
\item using the method introduced in Section~\ref{section:adaptiveglobal}, we first construct an admissible pair $(u_m^{h,\Delta t},\hat{\qb}^{h,\Delta t}_m)$ (resp. $(\tilde{u}_{m'}^{h',\Delta t'},\hat{\tilde{\qb}}^{h',\Delta t'}_{m'})$) in the discretized (FE) sense \eqref{eq:discreteweakform} and compute the associated CRE measure $E^{h,\Delta t}_{\CRE}$ (resp. $\tilde{E}^{h,\Delta t}_{\CRE}$) for the primal/reference (resp. dual/adjoint) problem. We then define an indicator $\rho_{\PGD}$ of the PGD truncation error $\Delta Q_{\PGD}$ as:
\begin{equation}
\abs{\Delta Q_{\PGD}} \approx \abs{Q_{\corr}^{h,\Delta t} \pm \frac{1}{2} \tilde{E}^{h,\Delta t}_{\CRE} E^{h,\Delta t}_{\CRE}} \eqqcolon \rho_{\PGD},
\end{equation}
where the correcting term $Q_{\corr}^{h,\Delta t}$ as well as the CRE measures $E^{h,\Delta t}_{\CRE}$ and $\tilde{E}^{h,\Delta t}_{\CRE}$ are computed from $(u_m^{h,\Delta t},\hat{\qb}^{h,\Delta t}_m)$ and $(\tilde{u}_{m'}^{h',\Delta t'},\hat{\tilde{\qb}}^{h',\Delta t'}_{m'})$.
\item we finally deduce an indicator $\rho_{\dis}$ of the discretization error $\Delta Q_{\dis}$ as:
\begin{equation}
\abs{\Delta Q_{\dis}} \approx \rho - \rho_{\PGD} = \abs{Q_{\corr} \pm \frac{1}{2} \tilde{E}_{\CRE} E_{\CRE}} - \abs{Q_{\corr}^{h,\Delta t} \pm \frac{1}{2} \tilde{E}^{h,\Delta t}_{\CRE} E^{h,\Delta t}_{\CRE}} \eqqcolon \rho_{\dis}.
\end{equation}
\end{itemize}
These error indicators are then used to build up a greedy algorithm that aims at controlling the error on the PGD approximation in a goal-oriented manner. The adaptive strategy is similar to the one developed for the control of global error in Section~\ref{section:adaptiveglobal}, with a comparison between $\rho_{\PGD}$ and $\rho_{\dis}$ after computing each mode $i \in \{1,\dots,m\}$, and the use of the local contributions to the discretization error indicators $\eta_{\dis}$ and $\tilde{\eta}_{\dis}$ to drive space and time mesh adaptations.

\section{Numerical results}
\label{section:results}

We present in this section one-, two-, and three-dimensional numerical experiments which illustrate the proposed error estimation method and adaptive strategy. 

%
%
%
%
%

\subsection{One-dimensional beam problem with time-dependent traction loading}

\begin{figure}[h!]
\centering
\includegraphics[height=20mm]{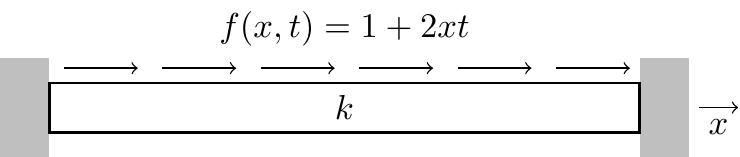}
\caption{1D beam problem.}\label{fig:1D_problem}
\end{figure}

We first consider a 1D transient diffusion problem on a beam structure (see Figure~\ref{fig:1D_problem}). The structure, of length $L=1$~m, is clamped at both ends $x=0$ and $x=L$, and is subjected to an evolutive affine source term $f(x,t) = 1+2xt$ over the time period $I=\intervalcc{0}{T}$ with $T=1$~s. The diffusion coefficient $k \in P_k=\intervaloc{0}{100}$ is considered as an extra-coordinate in the PGD representation $u_m(x,t,k)$ of the solution.

\begin{figure}[h!]
\centering
\includegraphics[height=100mm]{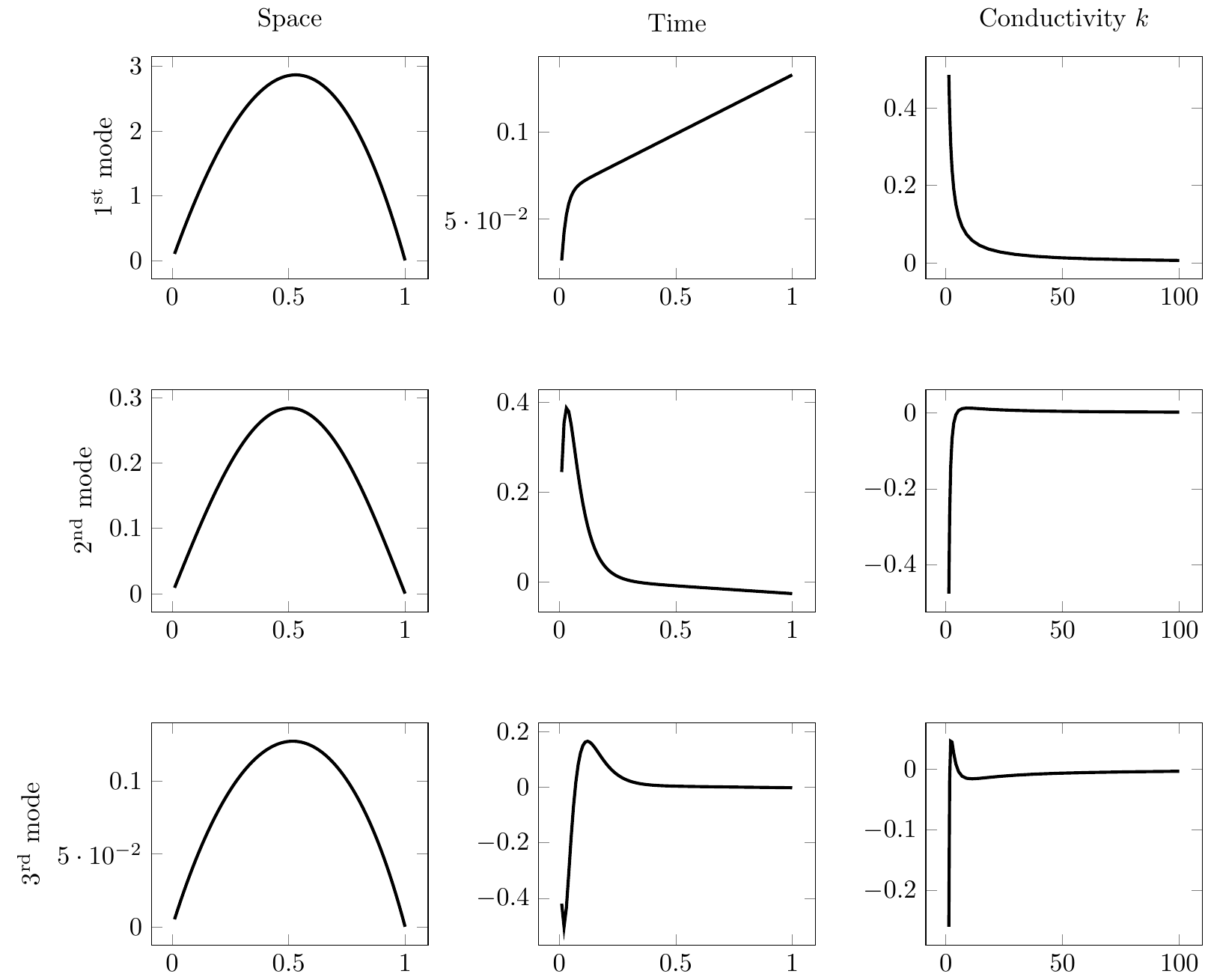}
\caption{Space functions $\psi_m(x)$, time functions $\lambda_m(t)$ and parameter functions $\gamma_m(k)$ (from left to right) obtained for order $m=1,2,3$ (from top to bottom).}\label{fig:1D_modes}
\end{figure}

\begin{figure}[h!]
\centering
\includegraphics[height=40mm]{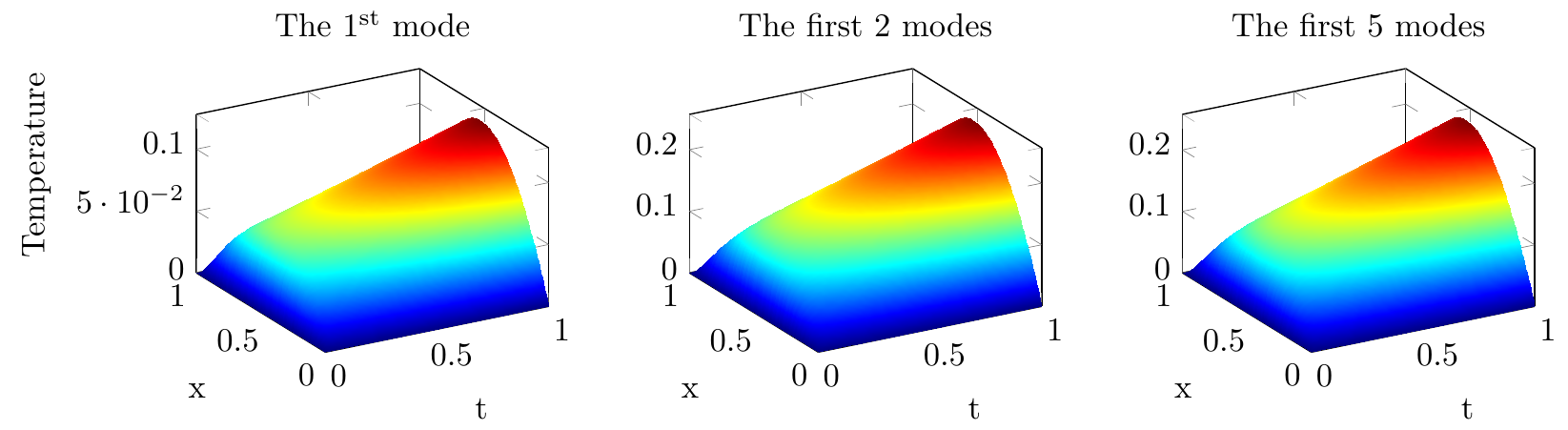}
\caption{Space-time mapping of the approximate PGD solution $u_m^{h,\Delta t}$ given for order $m=1,2,5$ (from left to right).}\label{fig:1D_PGDsol}
\end{figure}

The initial discretization used to compute the PGD approximate solution $u_m^{h,\Delta t}$ is made of $20$ $2$-nodes bar elements ($21$ dofs) in space and $10$ $2$-nodes bar elements in time (FEM in space and time). The first three PGD modes are given in Figure~\ref{fig:1D_modes}, whereas a representation of the space-time PGD solution $u_m^{h,\Delta t}$ for various numbers $m$ of PGD modes and for $k=2.07$ is given in Figure~\ref{fig:1D_PGDsol}.

\begin{figure}[h!]
\centering
\includegraphics[width=0.5\linewidth]{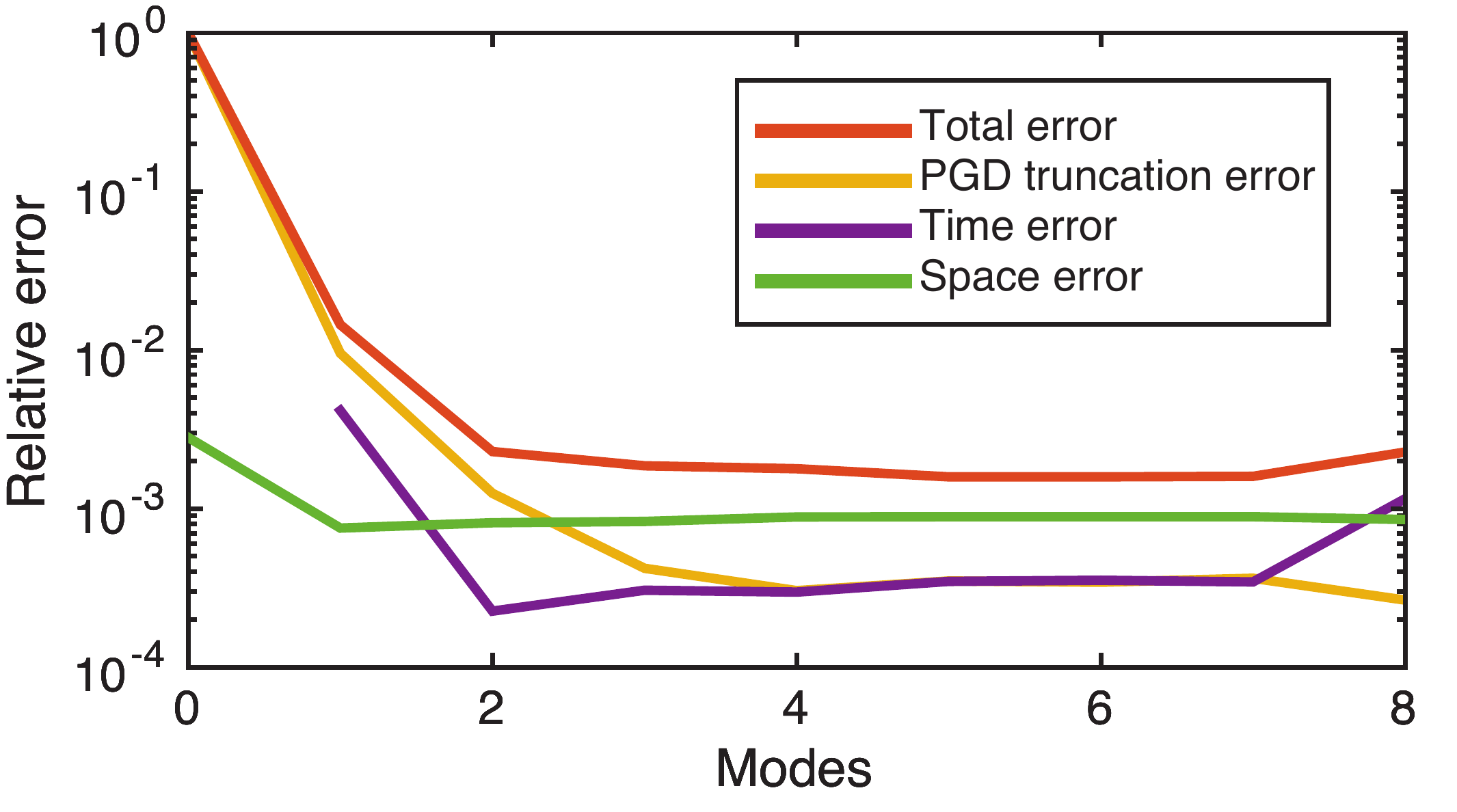}
\caption{Evolutions of the error estimate and error indicators with respect to the number $m$ of PGD modes without any adaptive strategy.}\label{fig:1D_estimates}
\end{figure}

\begin{figure}[h!]
\centering
\includegraphics[width=0.6\linewidth]{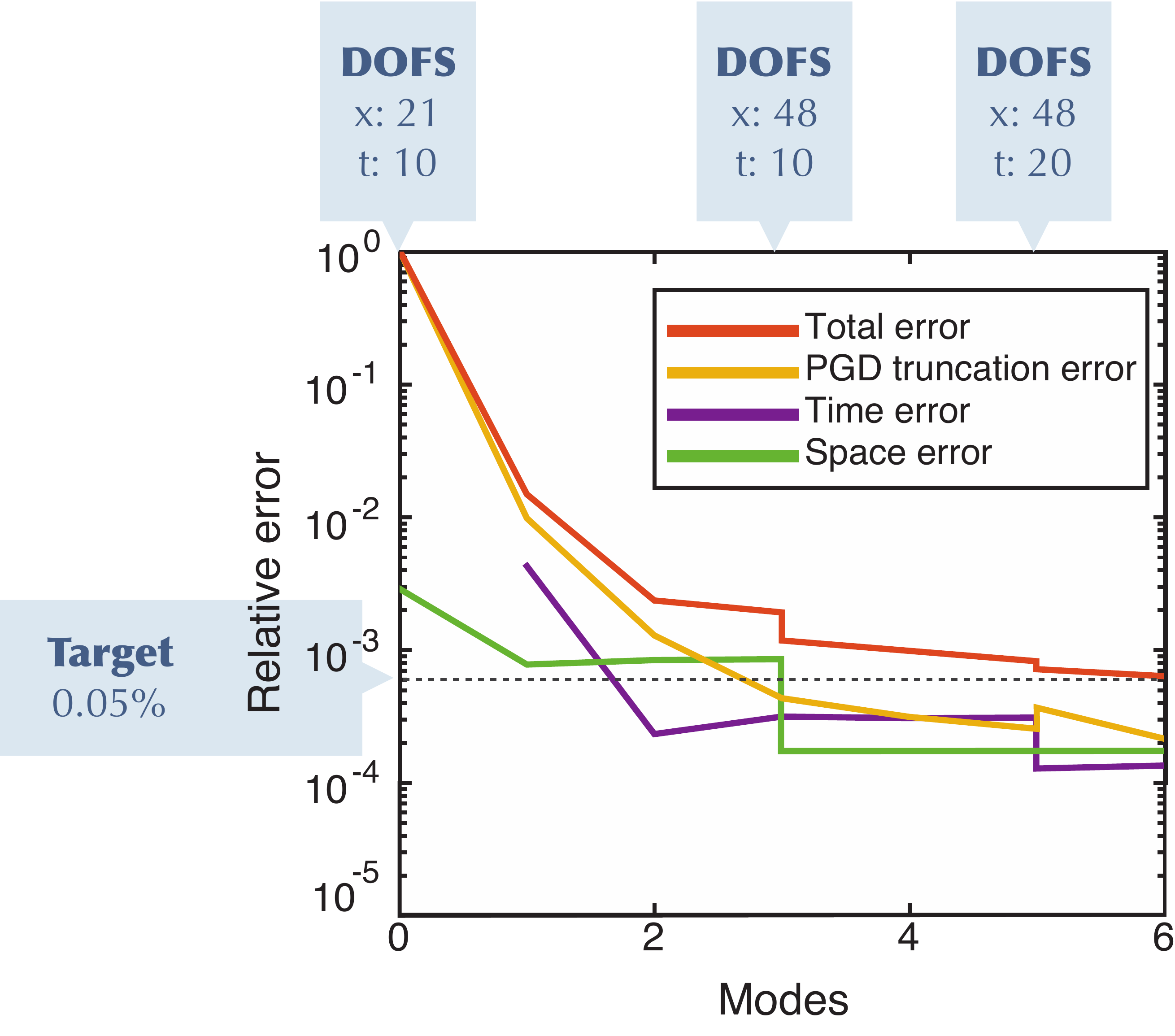}
\caption{Evolutions of the error estimate and error indicators with respect to the number $m$ of PGD modes with adaptive strategy.}\label{fig:1D_adapt_estimates}
\end{figure}

For fixed discretization parameters, Figure~\ref{fig:1D_estimates} shows the evolutions of the error estimate $E_{\CRE}$ and associated error indicators $\eta_{\PGD}$, $\eta_{h}$ and $\eta_{\Delta t}$ with respect to the number $m$ of PGD modes and for the maximal value obtained with $k\in P_k$. We observe that the PGD decomposition is converged at order $m = 4$, so that the error estimate and associated error indicators (measured in energy norm) become very small and sensitive to numerical noise. Such numerical effects could explain the slight increase of some error indicators after mode $m=4$. Conversely, Figure~\ref{fig:1D_adapt_estimates} shows the convergence of the error estimate and associated indicators when performing the adaptive strategy. Discontinuities in the curves correspond to adaptations of the space and time discretizations.

\subsection{Two-dimensional transient thermal problem}

\begin{figure}[h!]
\centering
\includegraphics[height=40mm]{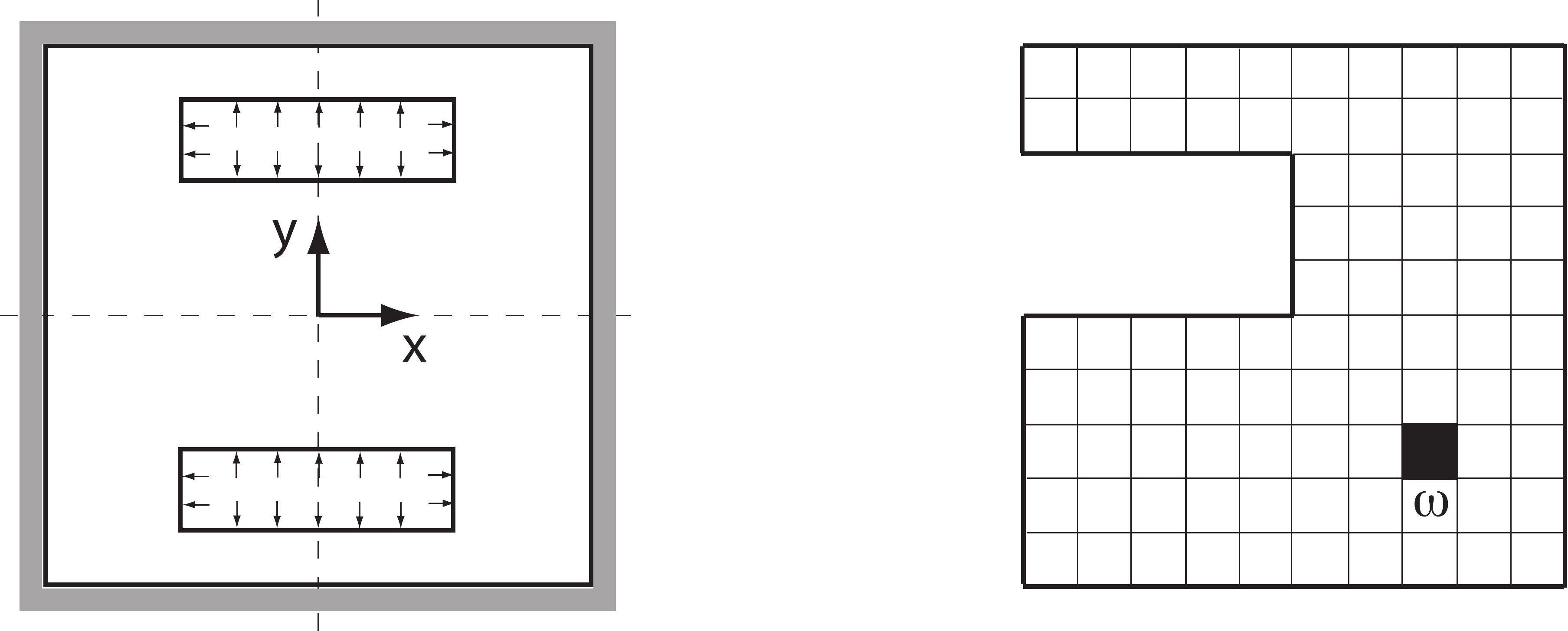}
\caption{2D thermal problem: space domain (left) and associated FE mesh (right).}\label{fig:2D_problem}
\end{figure}

We consider a transient thermal problem on the structure represented in Figure~\ref{fig:2D_problem}, which contains two symmetric rectangular holes in which a fluid circulates; exploiting symmetries, we study only one quarter, denoted $\Omega$, of the whole 2D domain. It is clamped on the external boundary, and is subjected over the time interval $I=\intervalcc{0}{T}$ (with $T=10$~s) to a given unit flux $g_d = -1$ applied on the hole boundary, homogeneous Neumann boundary conditions $g_d = 0$ on the symmetry planes, and a time-independent source term $f_d(x,y) = 200 x y$ in domain $\Omega$.

Extra parameters in the PGD decomposition are related to diffusion coefficient $k \in P_k=\intervalcc{1}{10}$ and thermal capacity $c \in P_c=\intervalcc{1}{10}$ (assumed to be uniform in space domain $\Omega$), so that the PGD representation of the solution reads $u_m(\xb,t,k,c)$.

\begin{figure}[h!]
	\centering
	\begin{subfigure}[c]{0.30\textwidth}
		\centering
		\includegraphics[height=35mm]{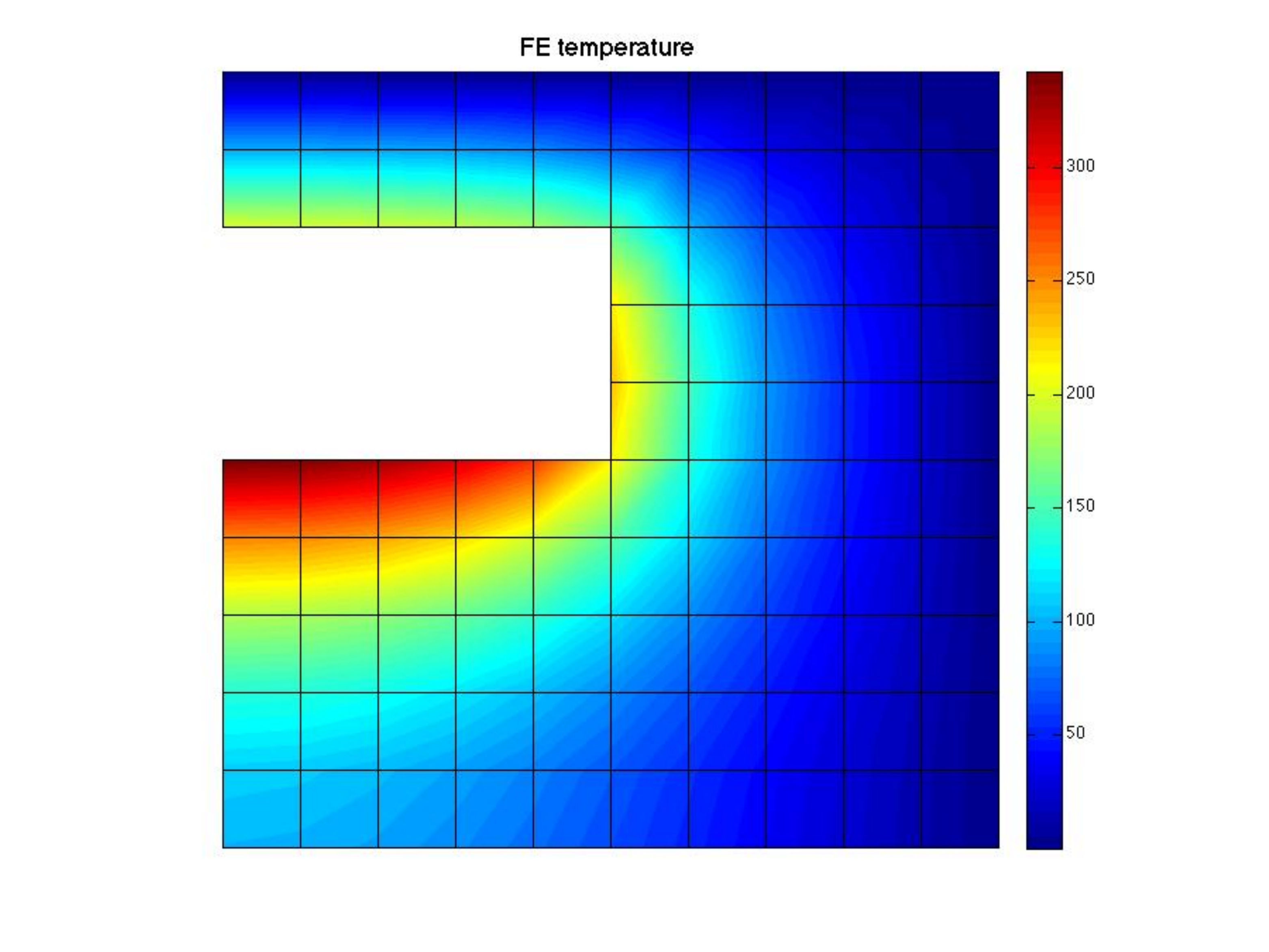}
	\end{subfigure}\hfill
	\begin{subfigure}[c]{0.30\textwidth}
		\centering
		\includegraphics[height=35mm]{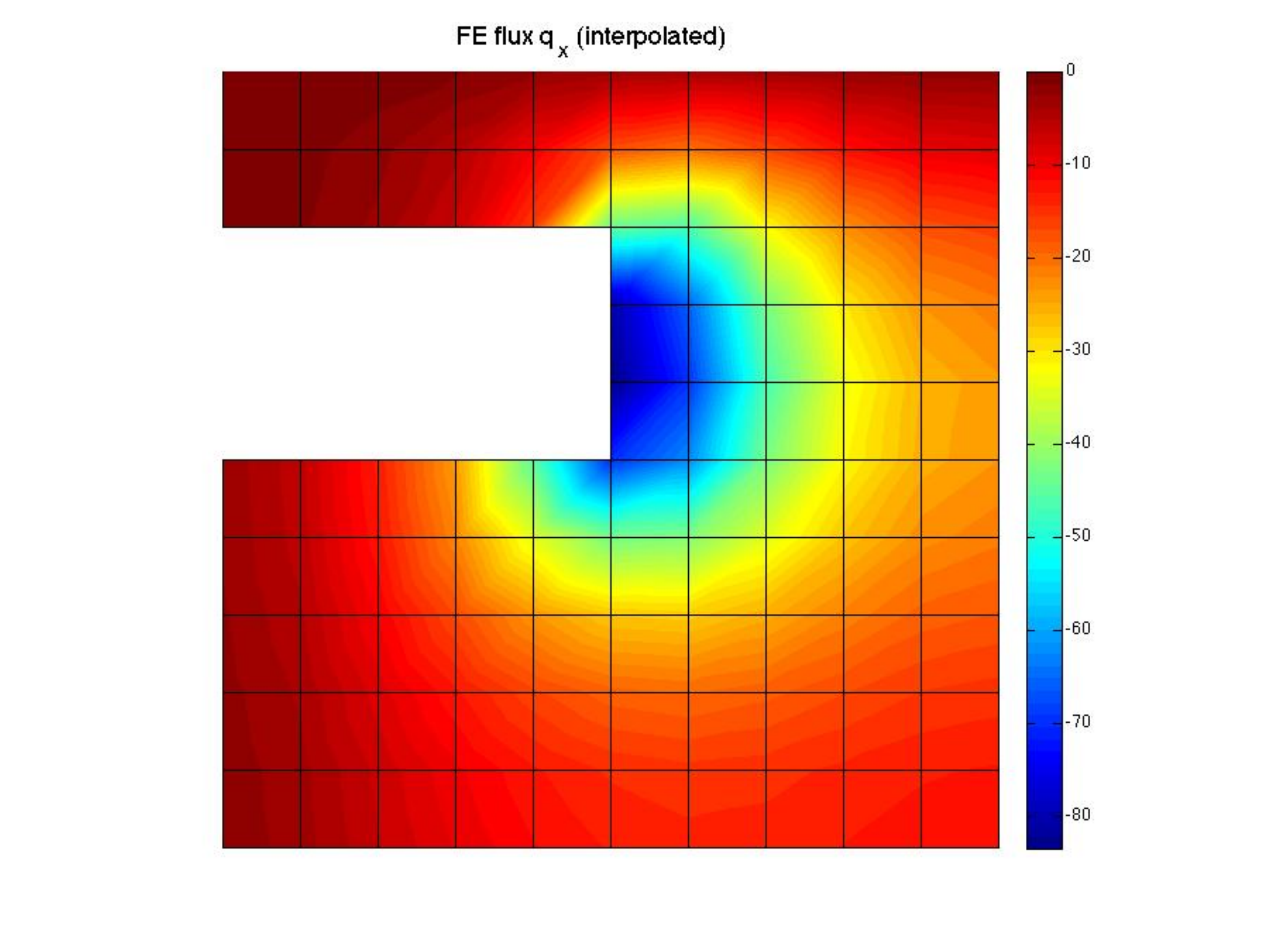}
	\end{subfigure}\hfill
	\begin{subfigure}[c]{0.30\textwidth}
		\centering
		\includegraphics[height=35mm]{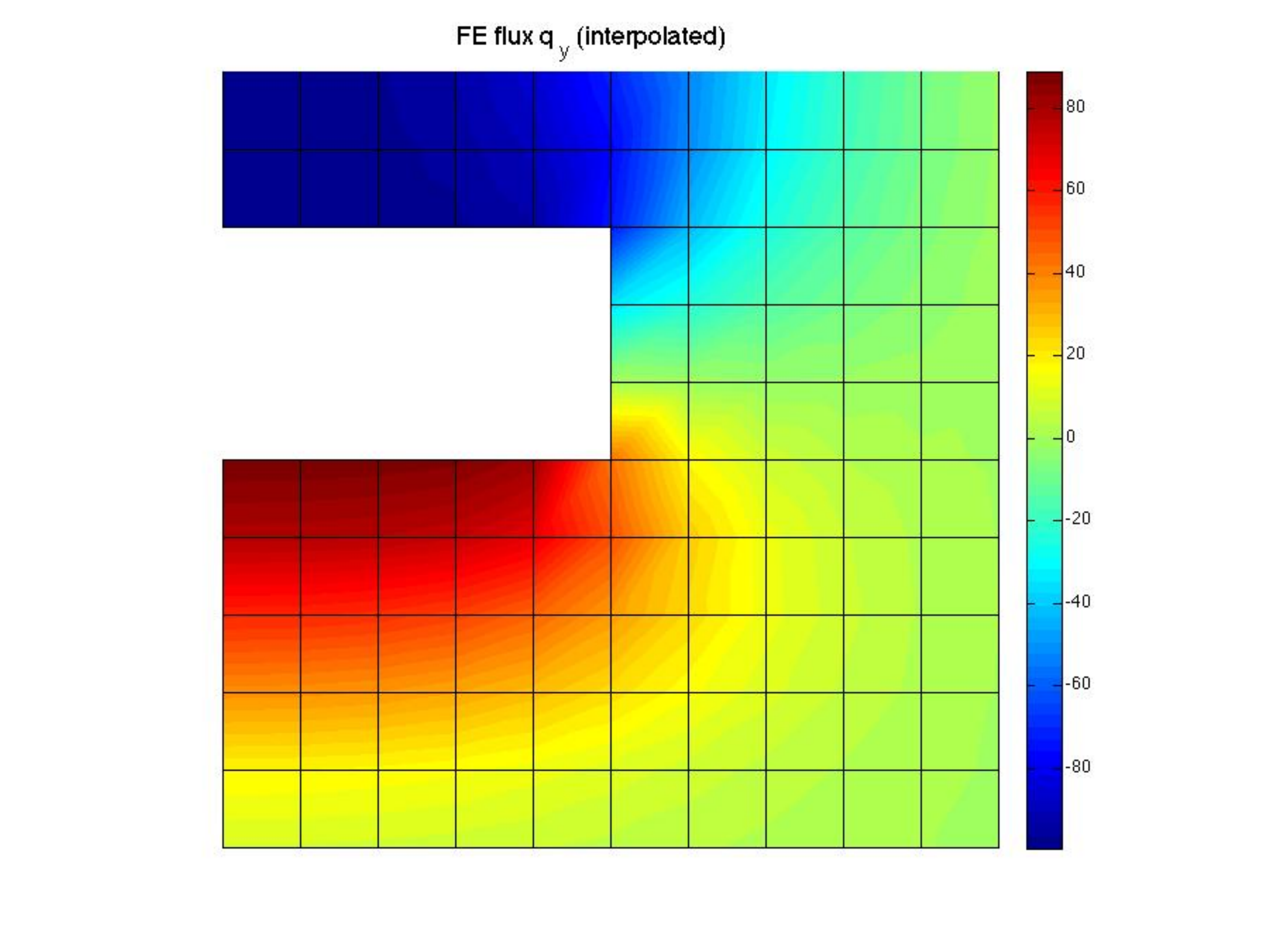}
	\end{subfigure}
\caption{FE solution $u^{h,\Delta t}$ (left) and associated flux components $q_x^{h,\Delta t}$ (center) and $q_y^{h,\Delta t}$ (right).}\label{fig:2D_FEsol}
\end{figure}

The initial mesh of the space domain $\Omega$ is made of $N_h=85$ regular $4$-nodes quadrangular elements of uniform size $h=0.1$, and a forward Euler time scheme is used with $N_{\Delta T}=1\,000$ time steps of uniform size $\Delta t=T/N_{\Delta T}=0.01$~s. We give in Figure~\ref{fig:2D_FEsol} the FE solution $u^{h,\Delta t}$ and associated flux $\qb^{h,\Delta t}$ for $(k,c)=(1,1)$ and at final time $t=T$.

\begin{figure}[h!]
	\centering
	\rotatebox[origin=c]{90}{$1\textsuperscript{st}$ mode}\hfill
	\begin{subfigure}[c]{0.22\textwidth}
		\centering
		\includegraphics[height=30mm]{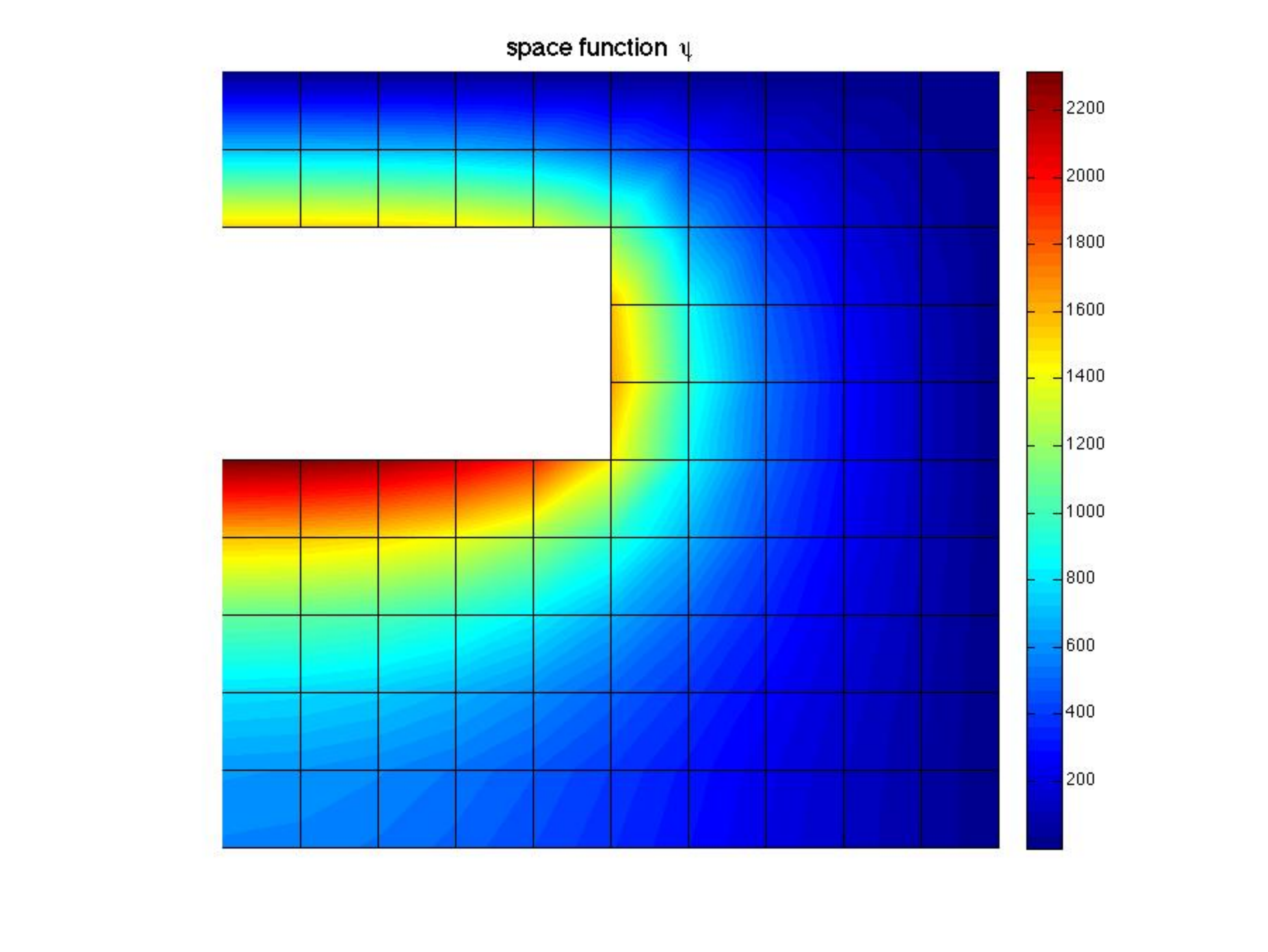}
	\end{subfigure}\hfill
	\begin{subfigure}[c]{0.22\textwidth}
		\centering
		\includegraphics[height=30mm]{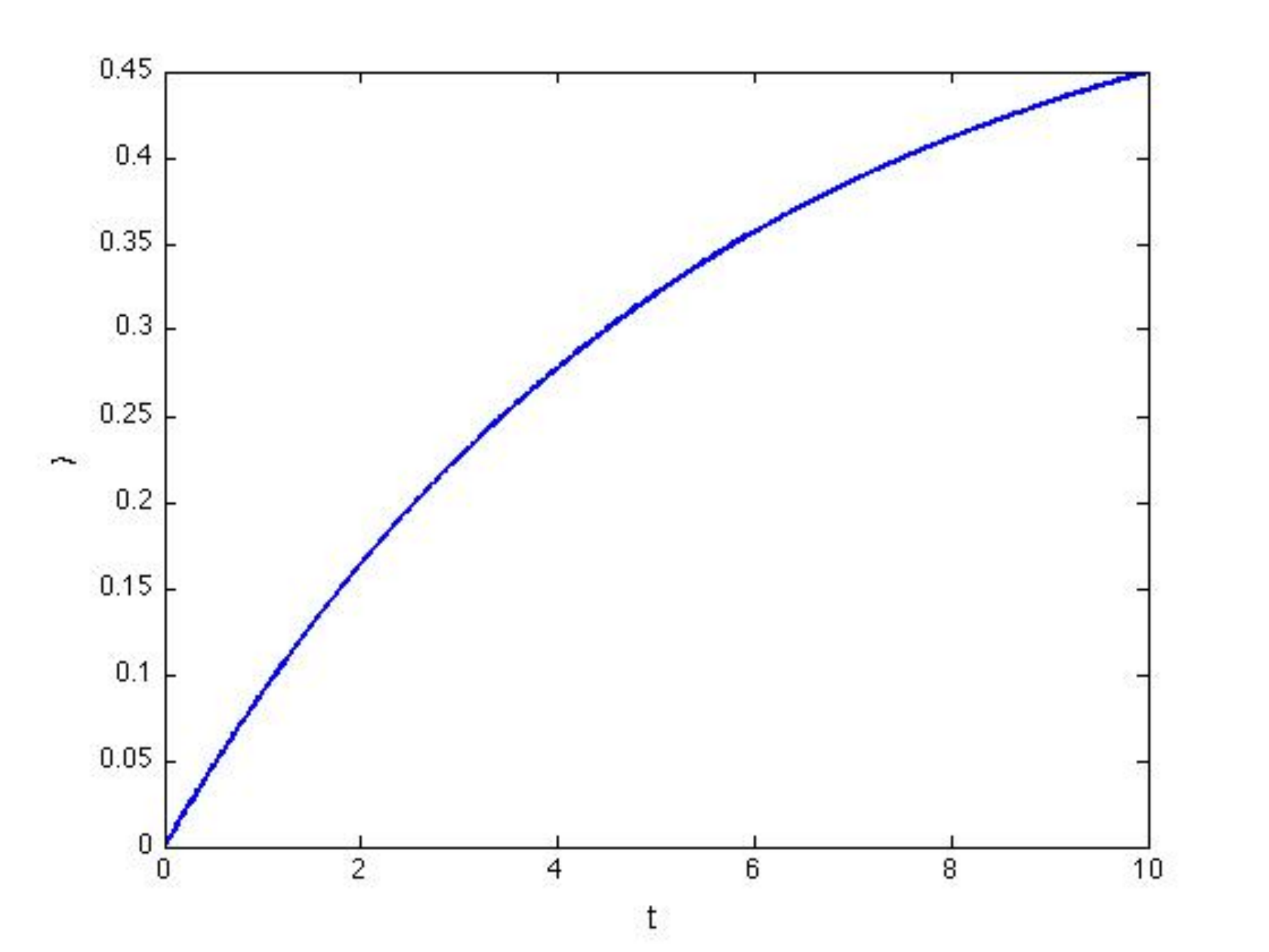}
	\end{subfigure}\hfill
	\begin{subfigure}[c]{0.22\textwidth}
		\centering
		\includegraphics[height=30mm]{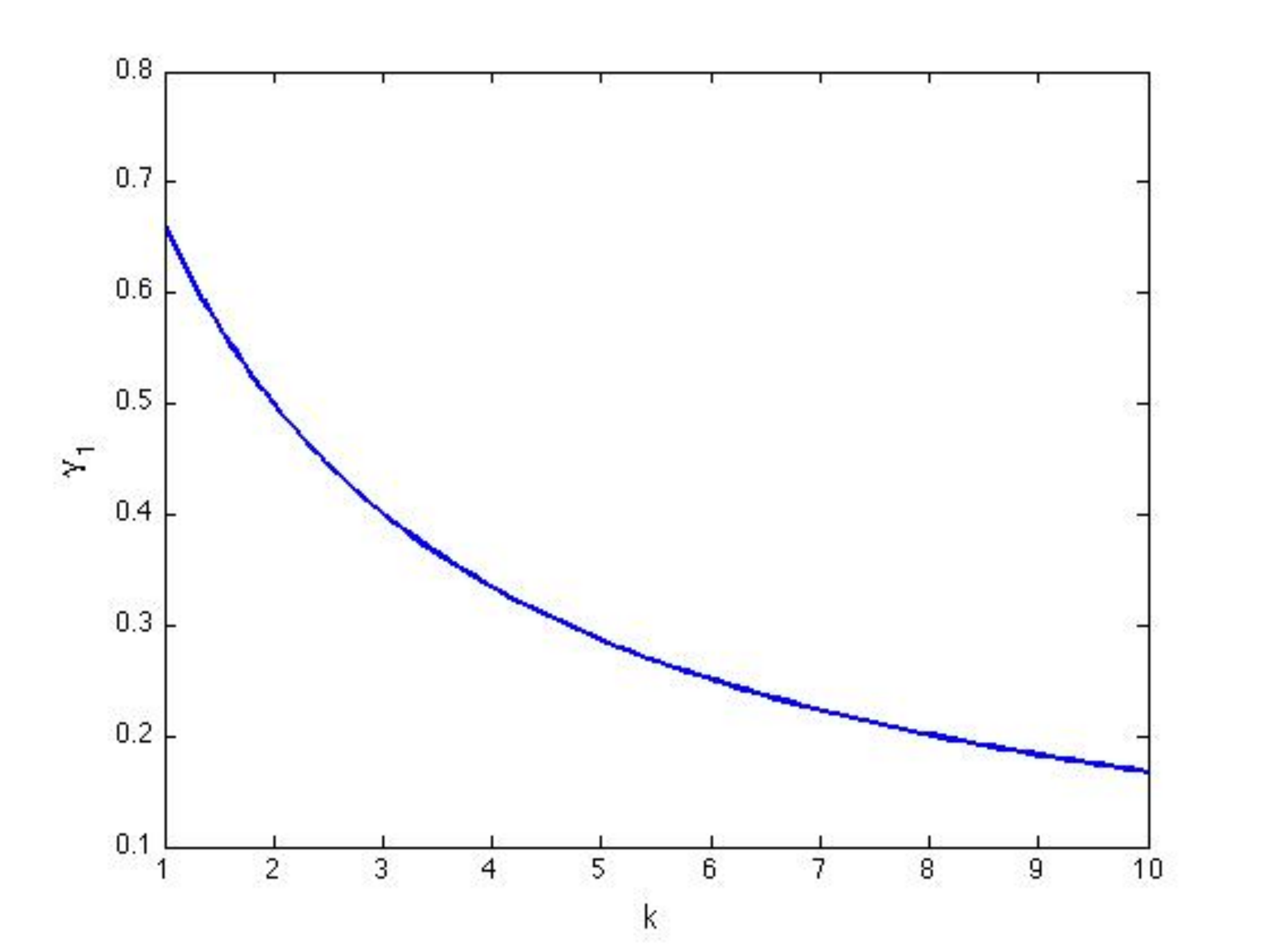}
	\end{subfigure}\hfill
	\begin{subfigure}[c]{0.22\textwidth}
		\centering
		\includegraphics[height=30mm]{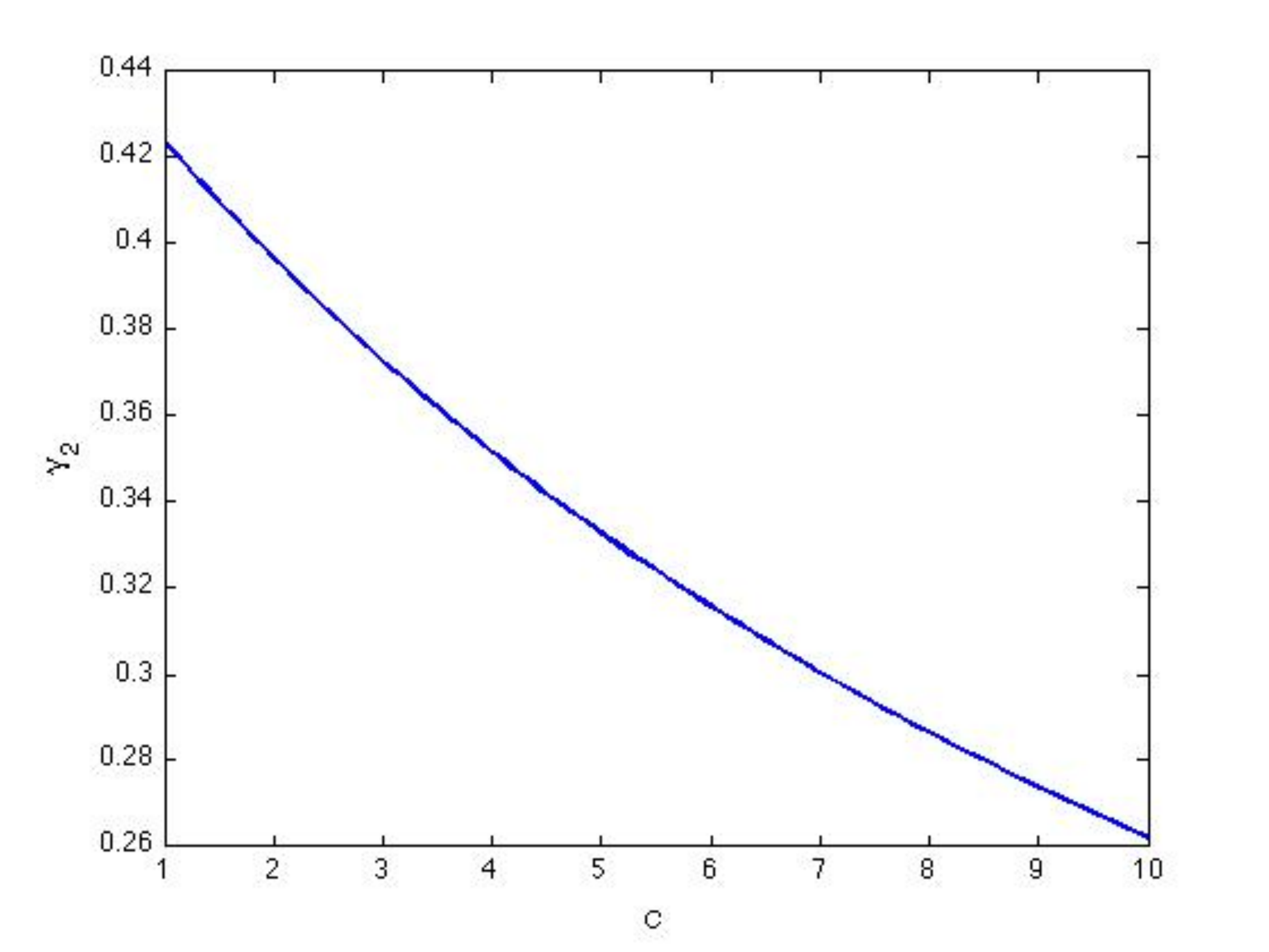}
	\end{subfigure}
	\\
	\rotatebox[origin=c]{90}{$2\textsuperscript{nd}$ mode}\hfill
	\begin{subfigure}[c]{0.22\textwidth}
		\centering
		\includegraphics[height=30mm]{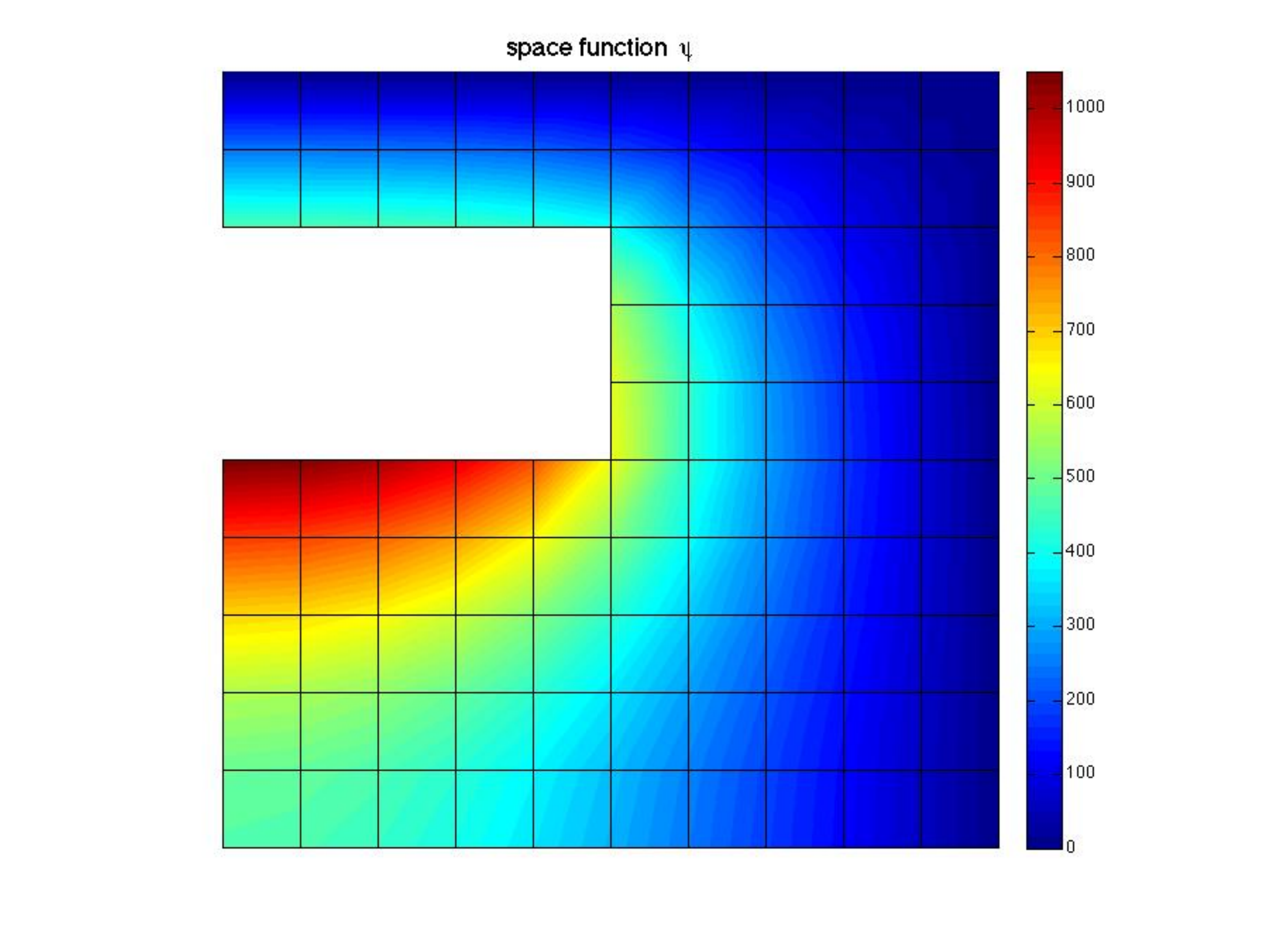}
	\end{subfigure}\hfill
	\begin{subfigure}[c]{0.22\textwidth}
		\centering
		\includegraphics[height=30mm]{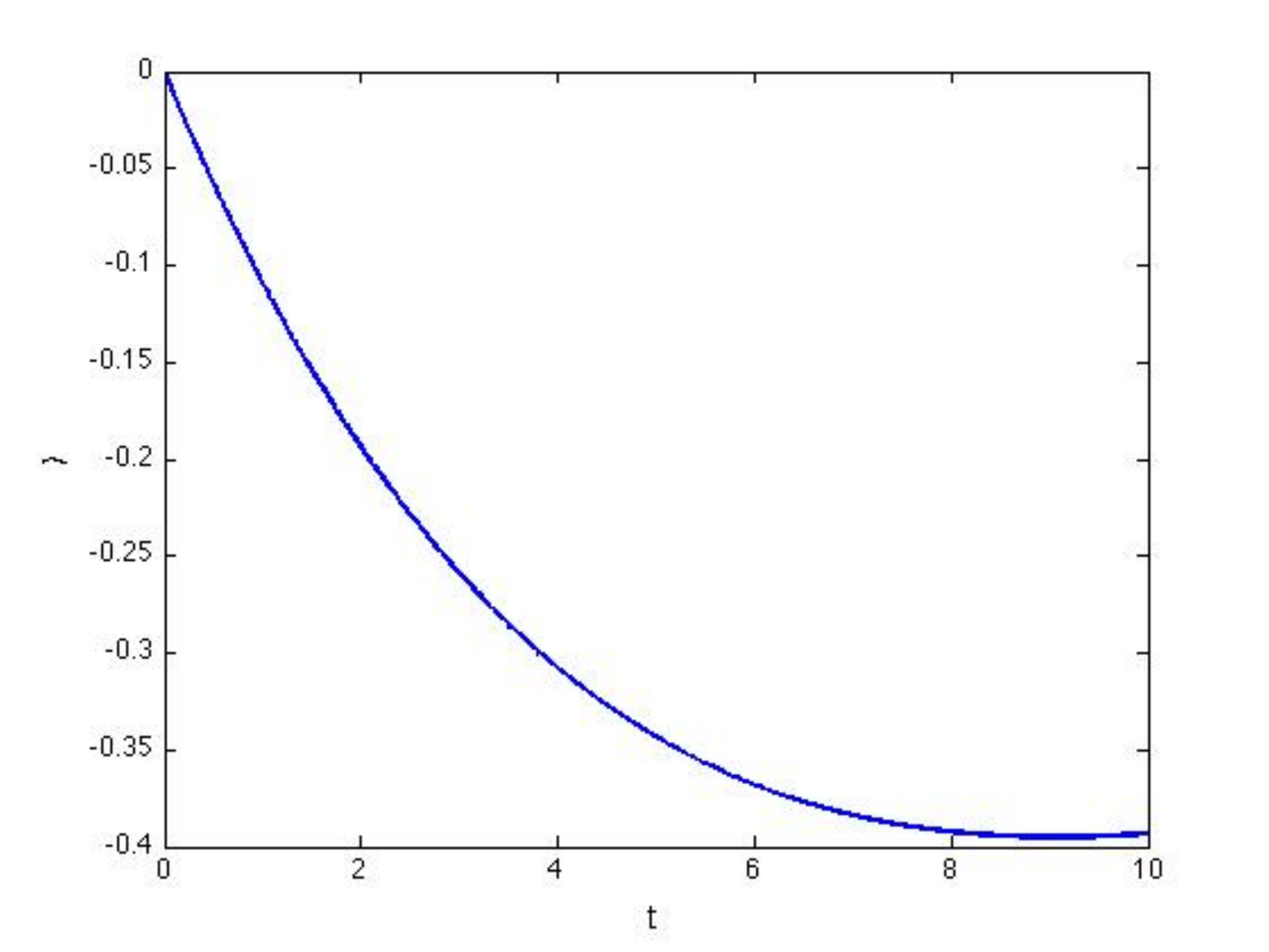}
	\end{subfigure}\hfill
	\begin{subfigure}[c]{0.22\textwidth}
		\centering
		\includegraphics[height=30mm]{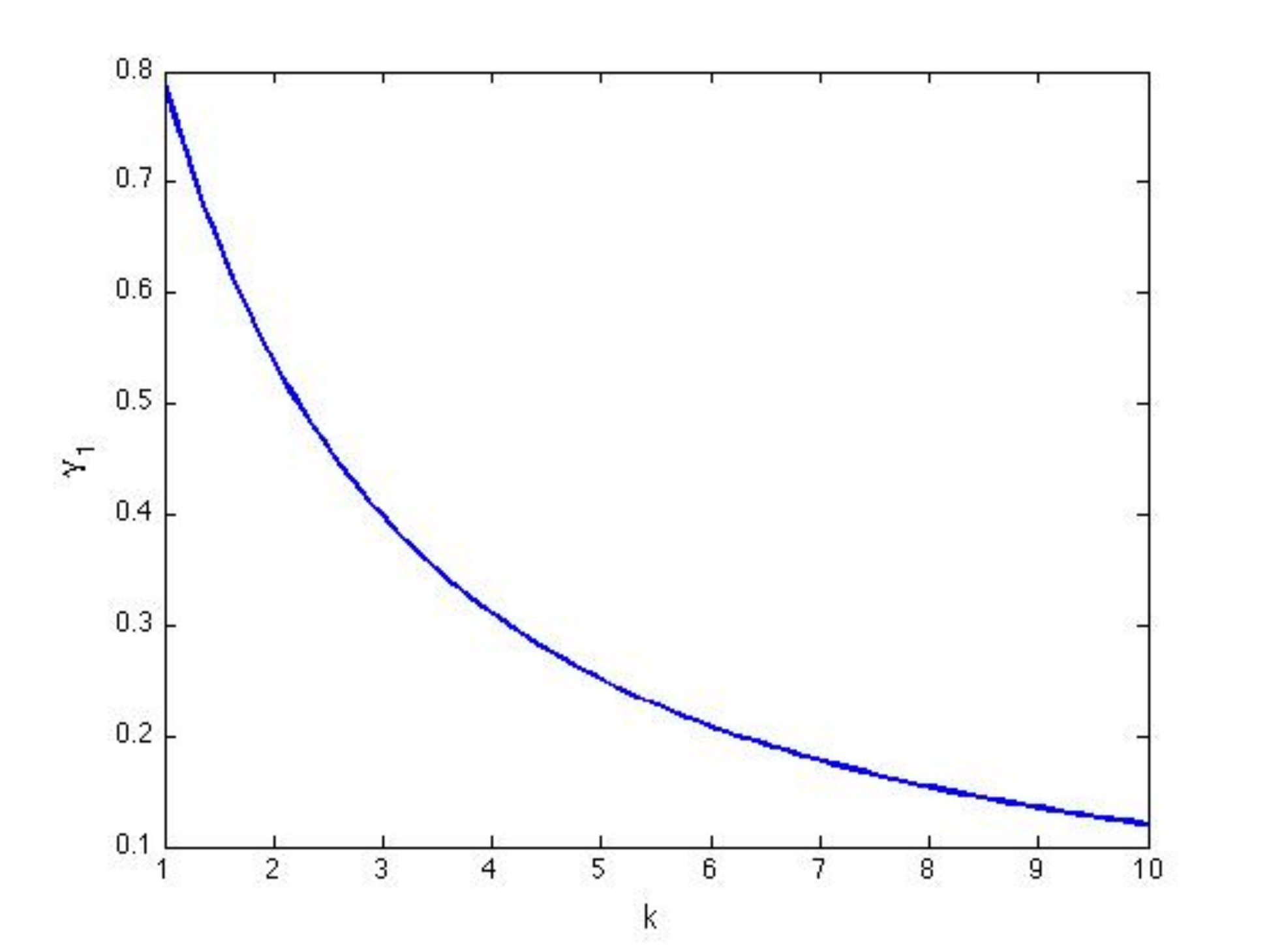}
	\end{subfigure}\hfill
	\begin{subfigure}[c]{0.22\textwidth}
		\centering
		\includegraphics[height=30mm]{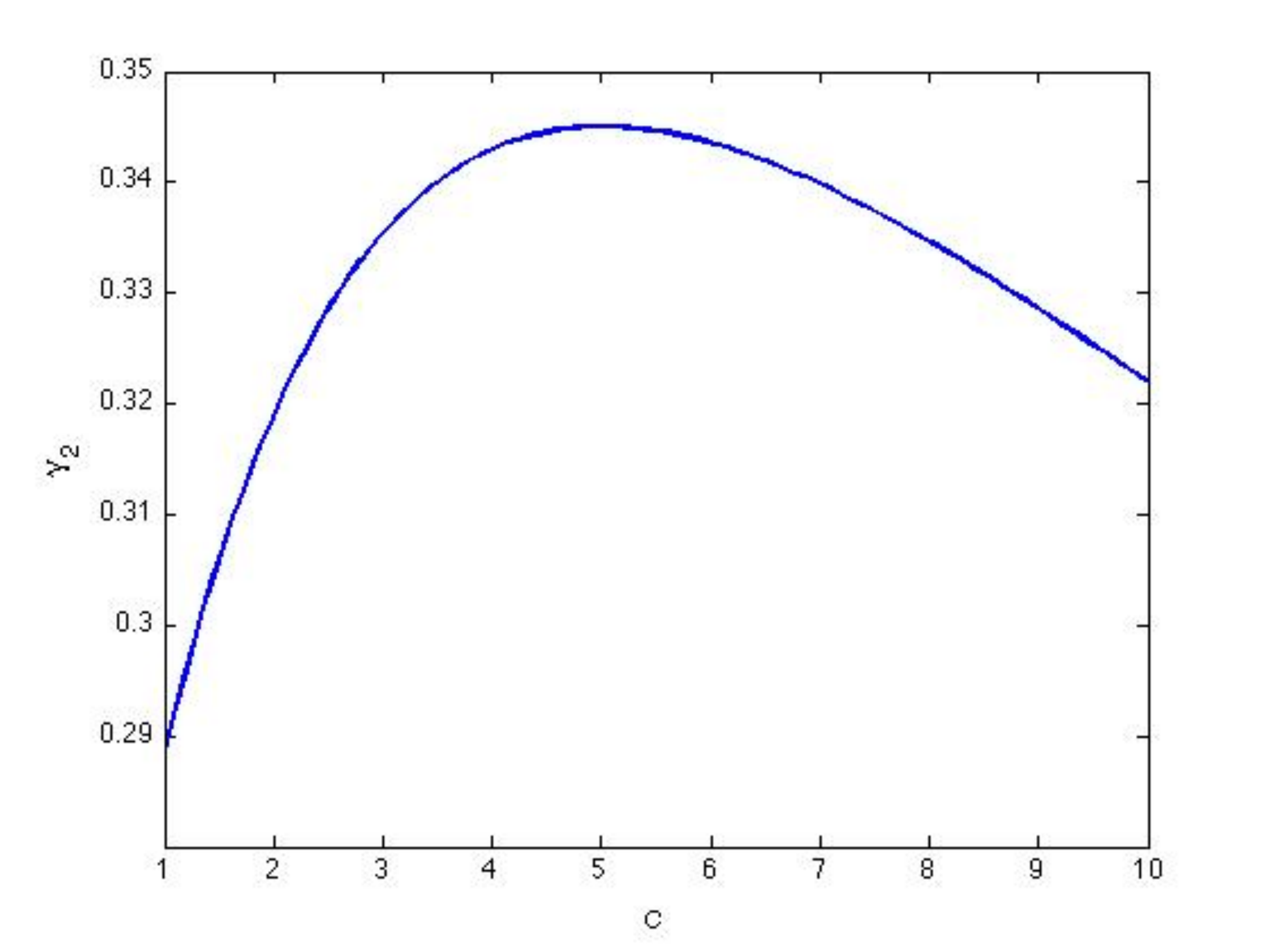}
	\end{subfigure}\\
	\rotatebox[origin=c]{90}{$3\textsuperscript{rd}$ mode}\hfill
	\begin{subfigure}[c]{0.22\textwidth}
		\centering
		\includegraphics[height=30mm]{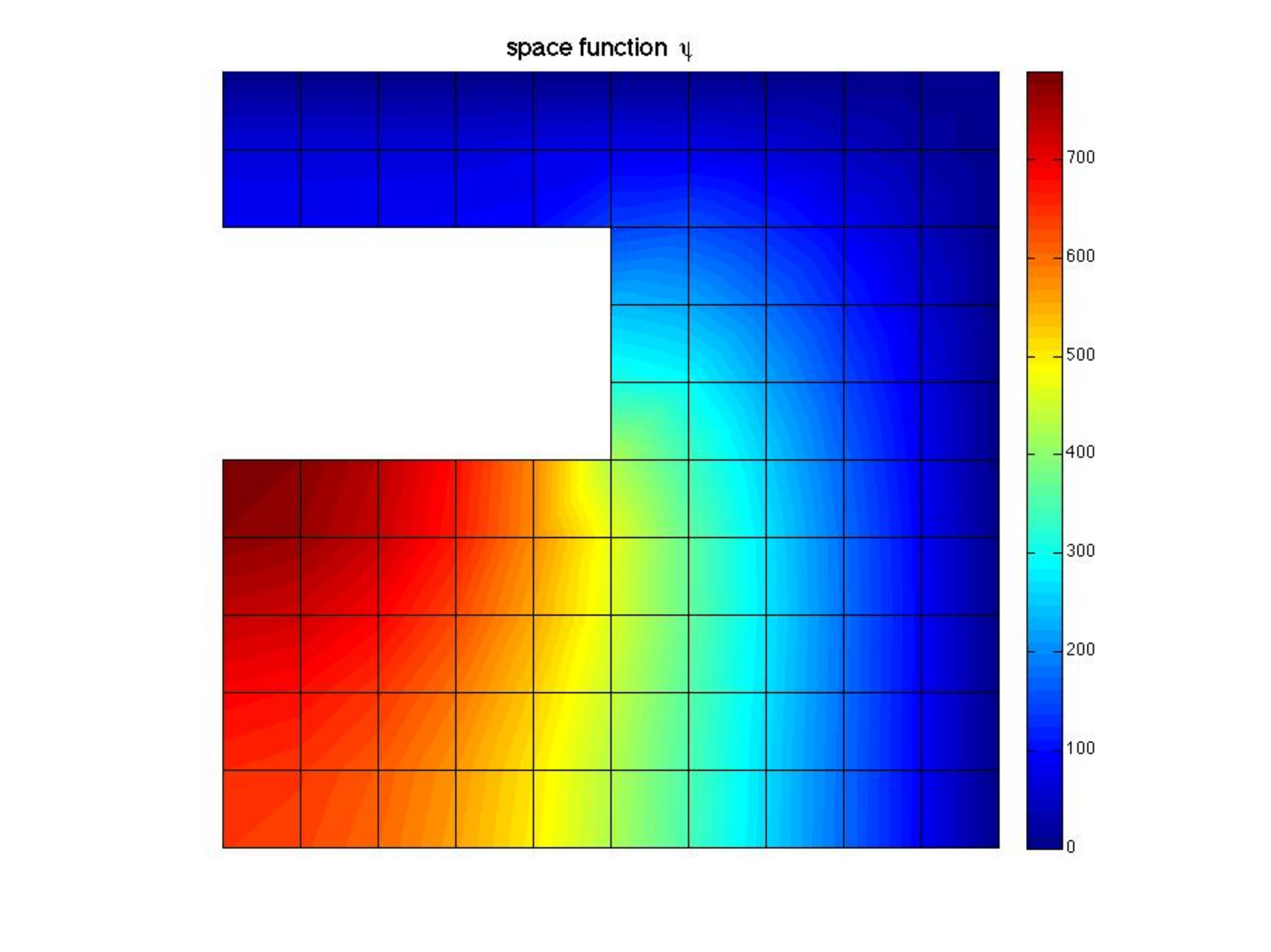}
	\end{subfigure}\hfill
	\begin{subfigure}[c]{0.22\textwidth}
		\centering
		\includegraphics[height=30mm]{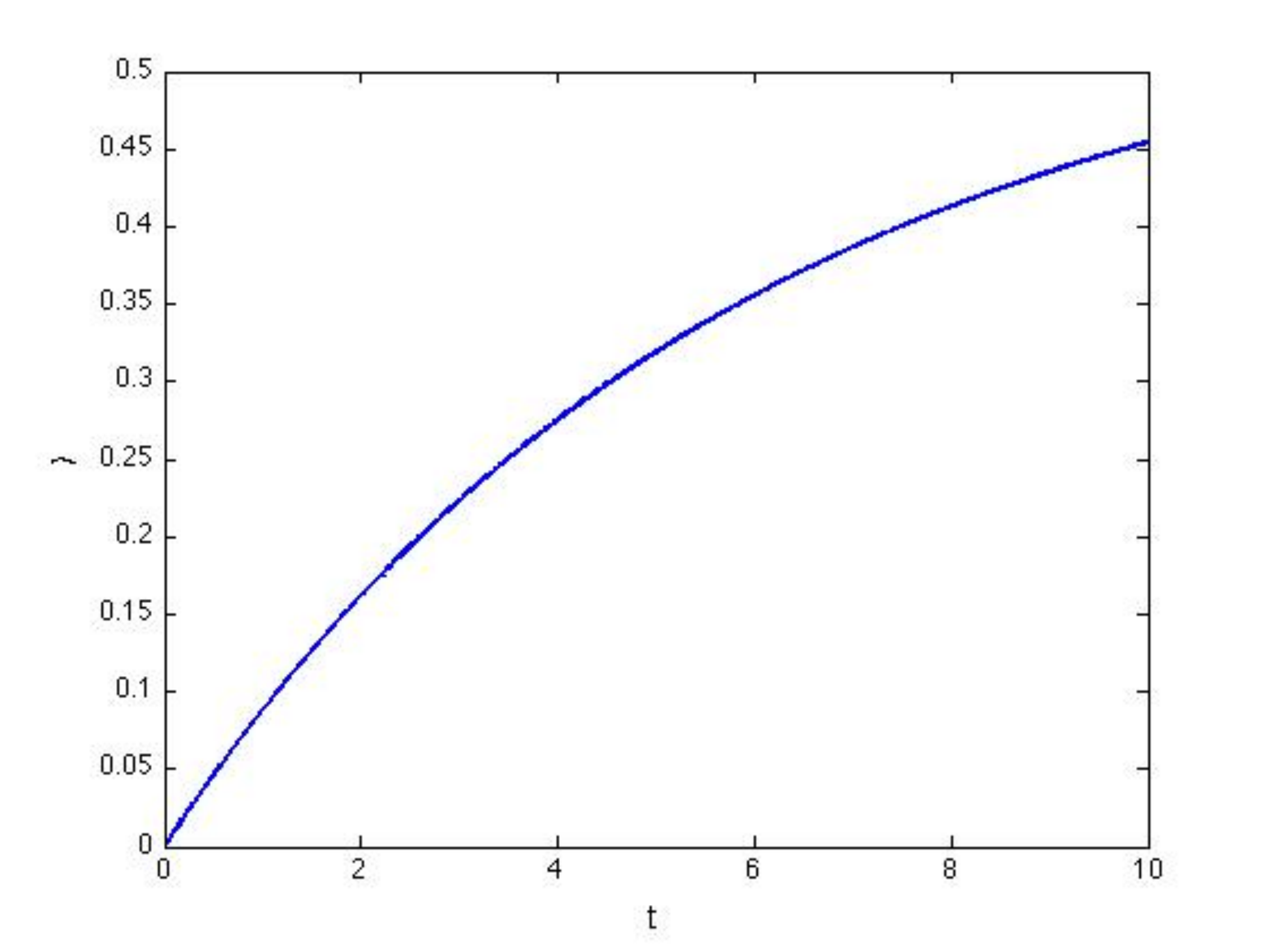}
	\end{subfigure}\hfill
	\begin{subfigure}[c]{0.22\textwidth}
		\centering
		\includegraphics[height=30mm]{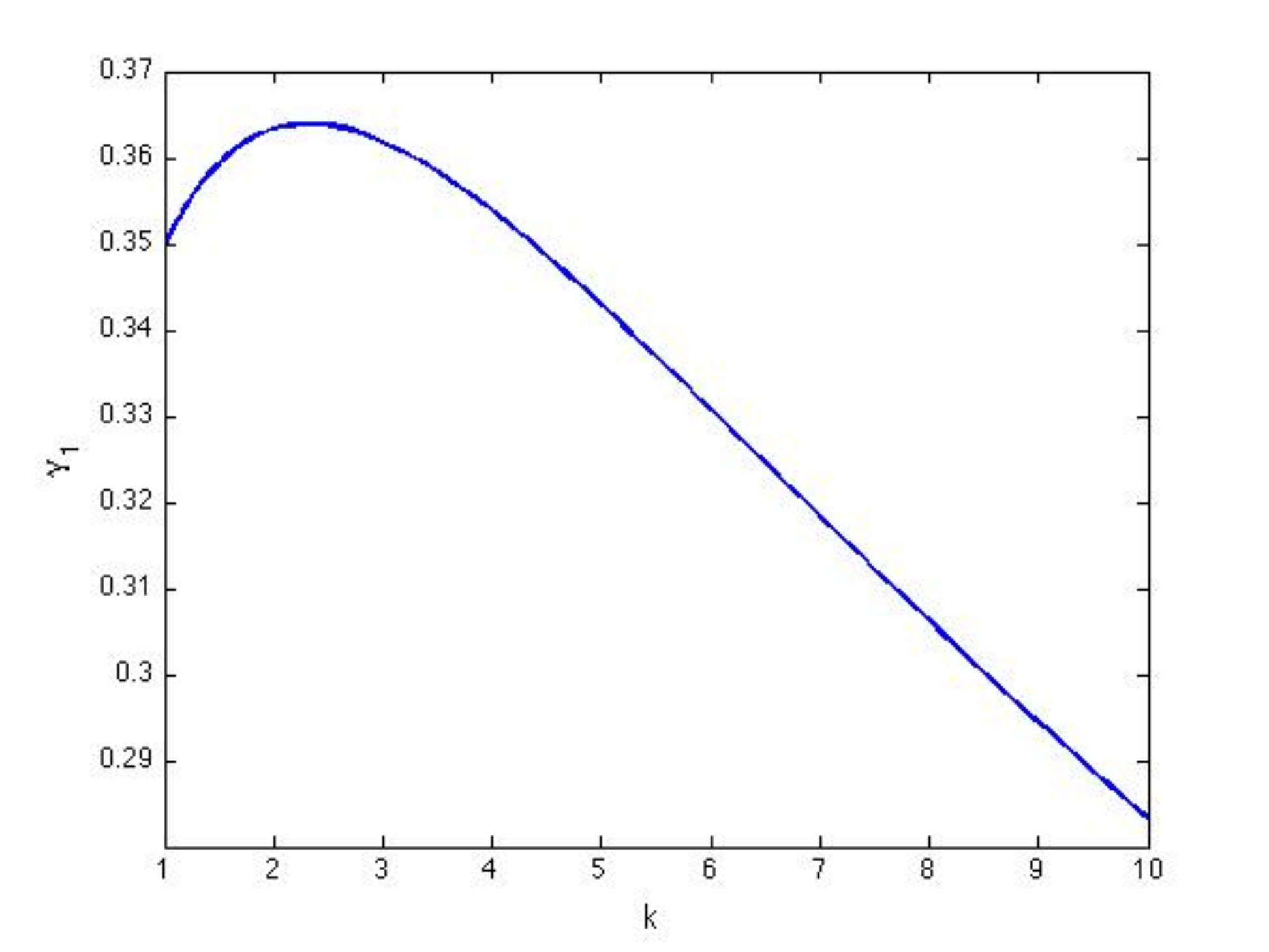}
	\end{subfigure}\hfill
	\begin{subfigure}[c]{0.22\textwidth}
		\centering
		\includegraphics[height=30mm]{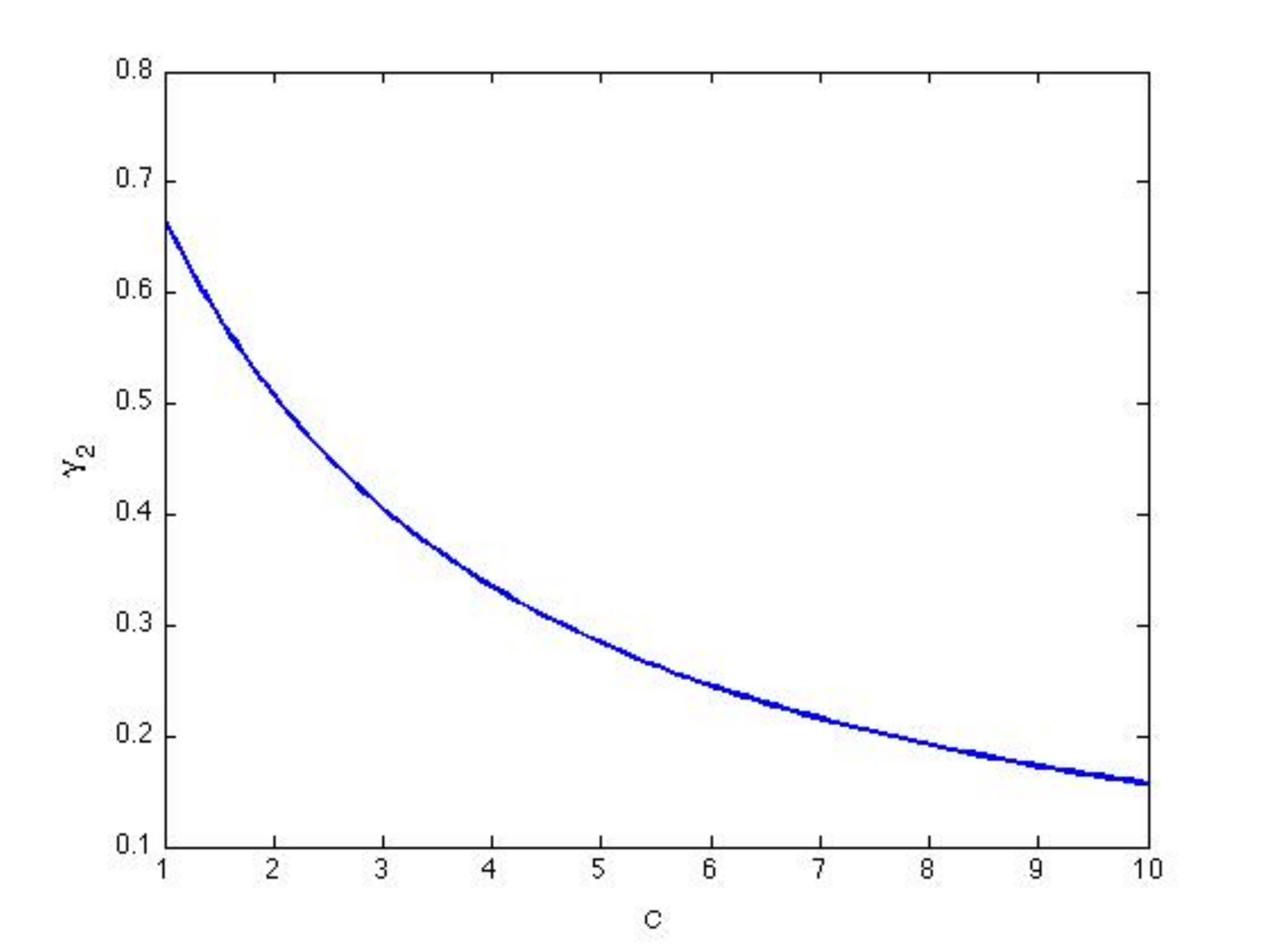}
	\end{subfigure}
	\\
	\rotatebox[origin=c]{90}{$4\textsuperscript{th}$ mode}\hfill
	\begin{subfigure}[c]{0.22\textwidth}
		\centering
		\includegraphics[height=30mm]{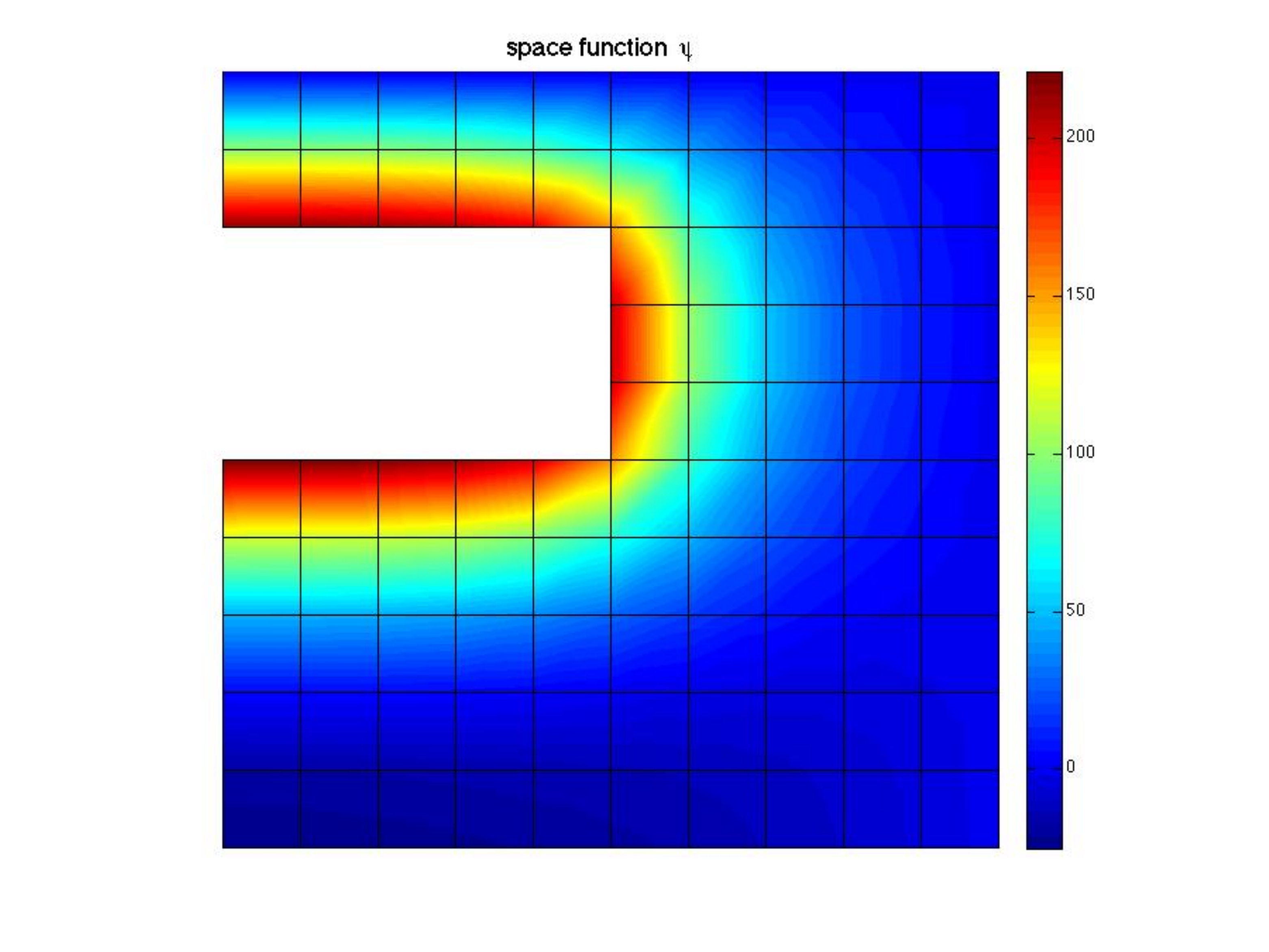}
	\end{subfigure}\hfill
	\begin{subfigure}[c]{0.22\textwidth}
		\centering
		\includegraphics[height=30mm]{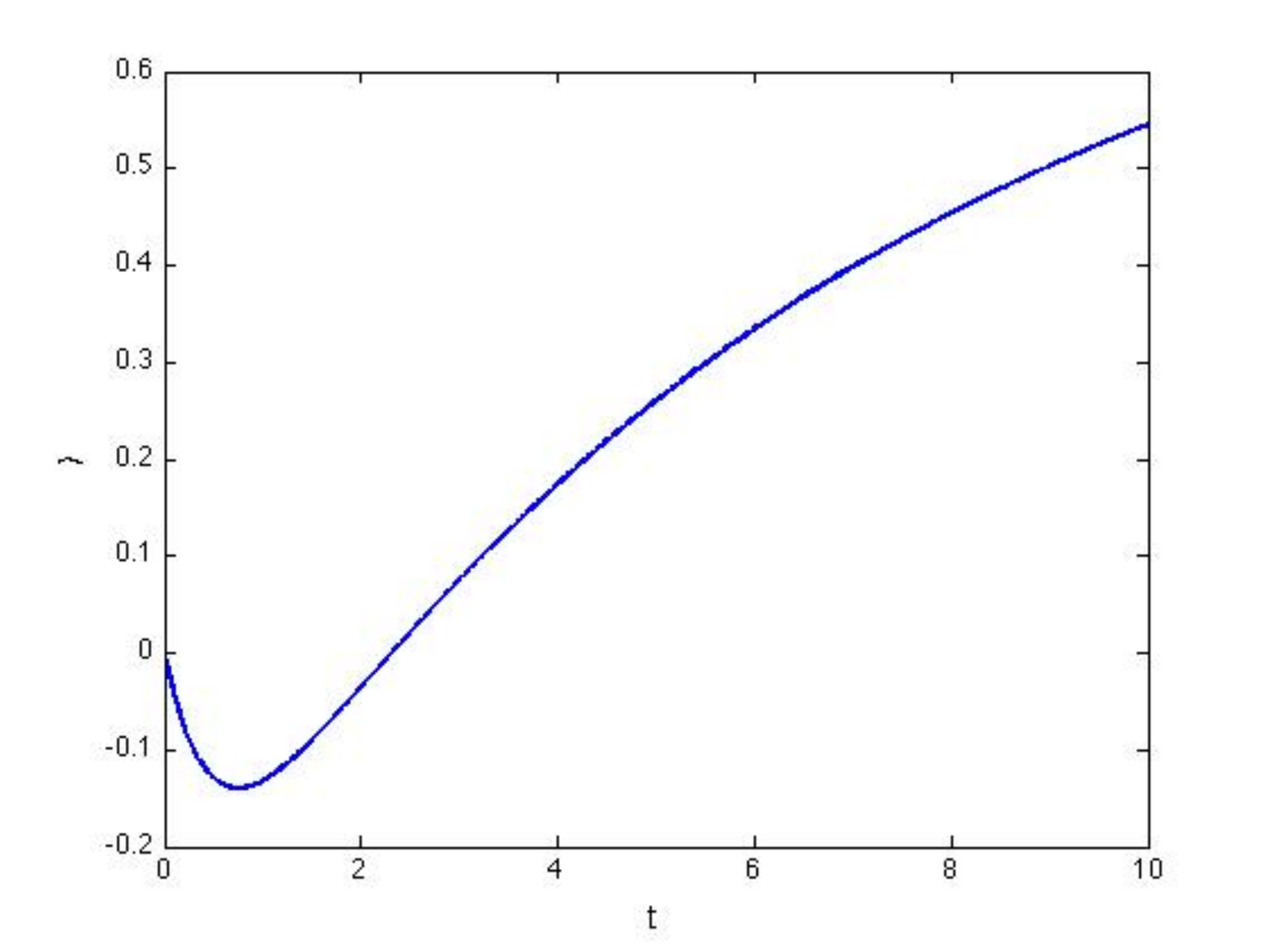}
	\end{subfigure}\hfill
	\begin{subfigure}[c]{0.22\textwidth}
		\centering
		\includegraphics[height=30mm]{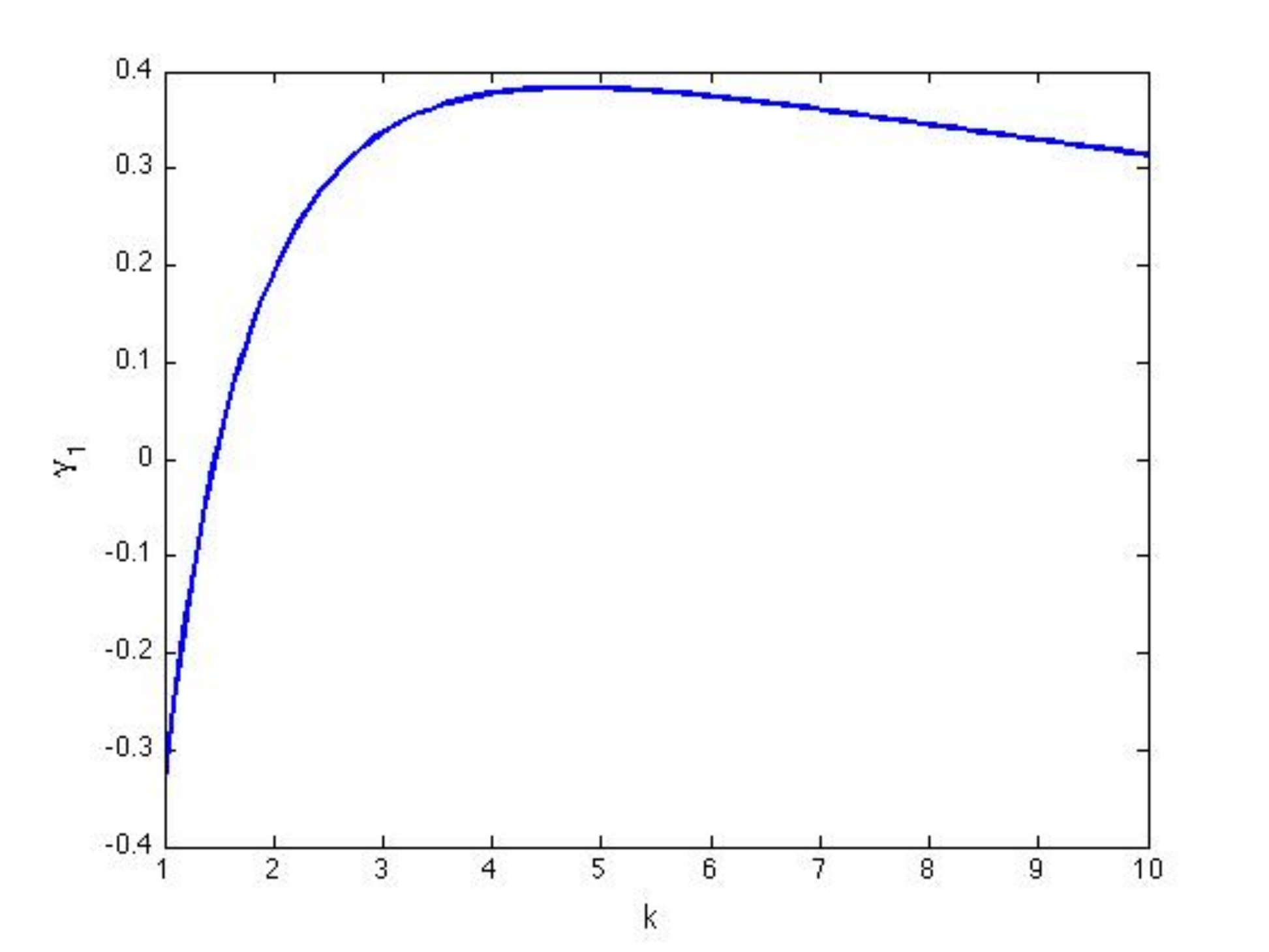}
	\end{subfigure}\hfill
	\begin{subfigure}[c]{0.22\textwidth}
		\centering
		\includegraphics[height=30mm]{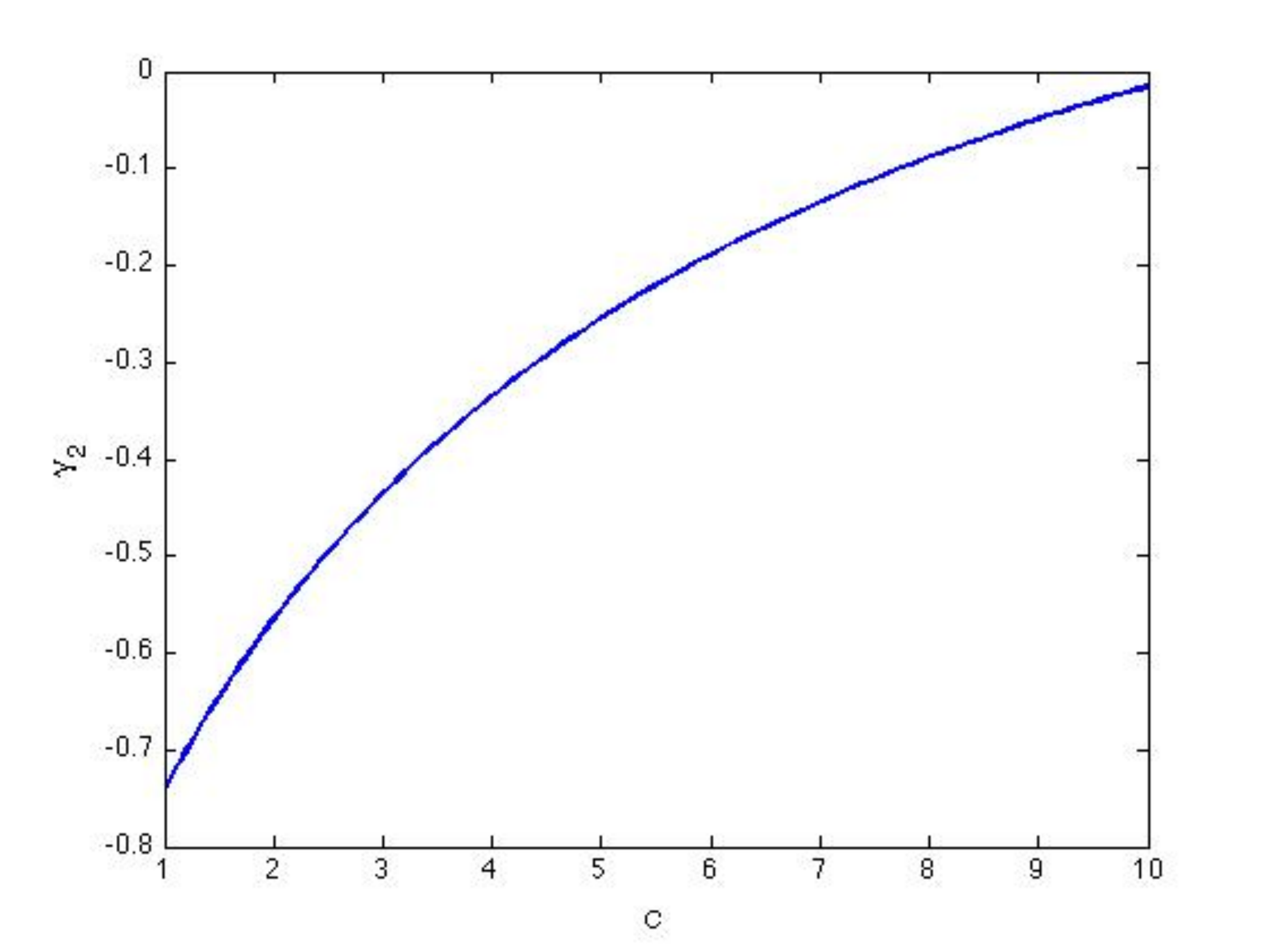}
	\end{subfigure}
	\\
	\rotatebox[origin=c]{90}{$5\textsuperscript{th}$ mode}\hfill
	\begin{subfigure}[c]{0.22\textwidth}
		\centering
		\includegraphics[height=30mm]{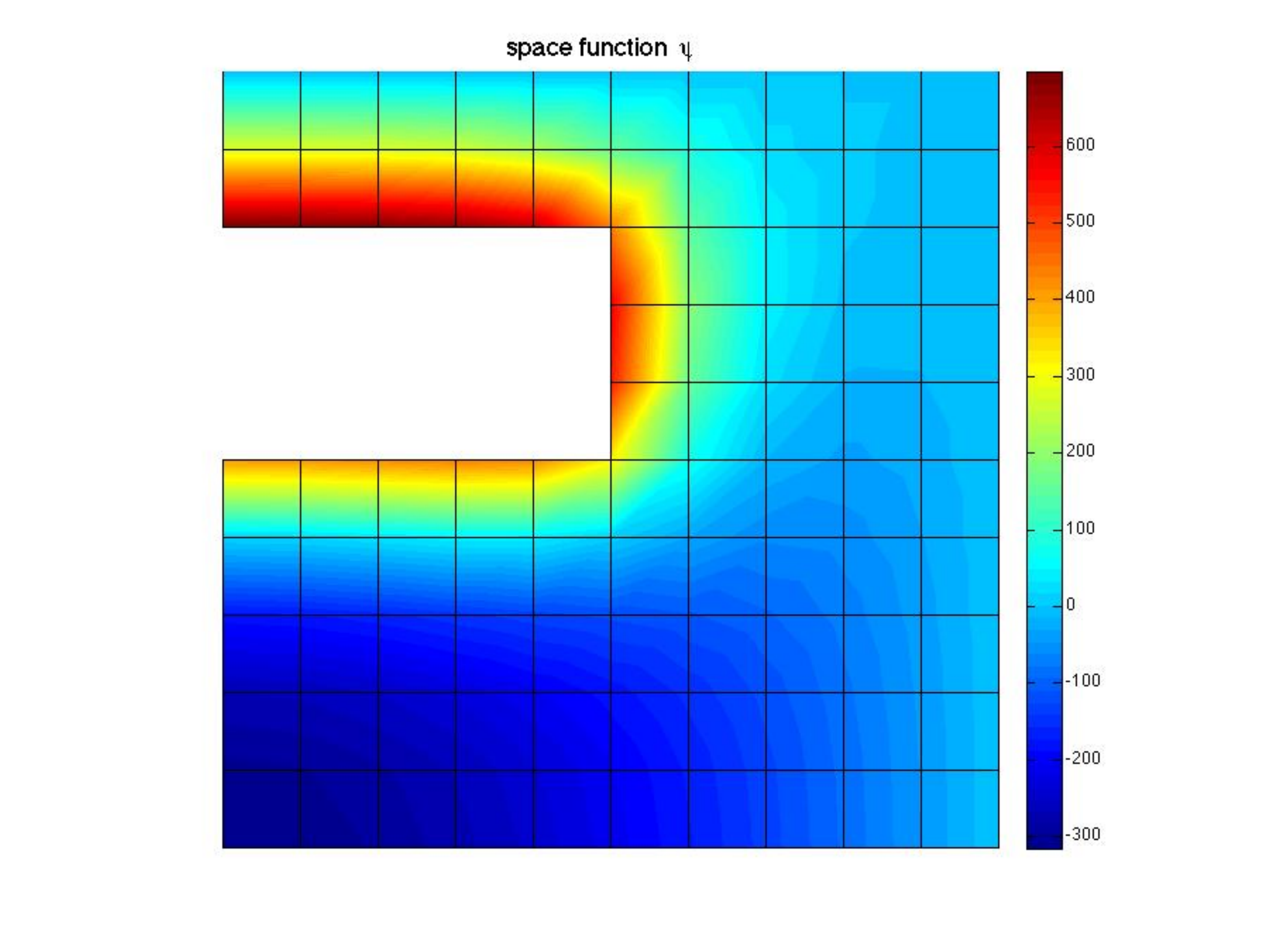}
	\end{subfigure}\hfill
	\begin{subfigure}[c]{0.22\textwidth}
		\centering
		\includegraphics[height=30mm]{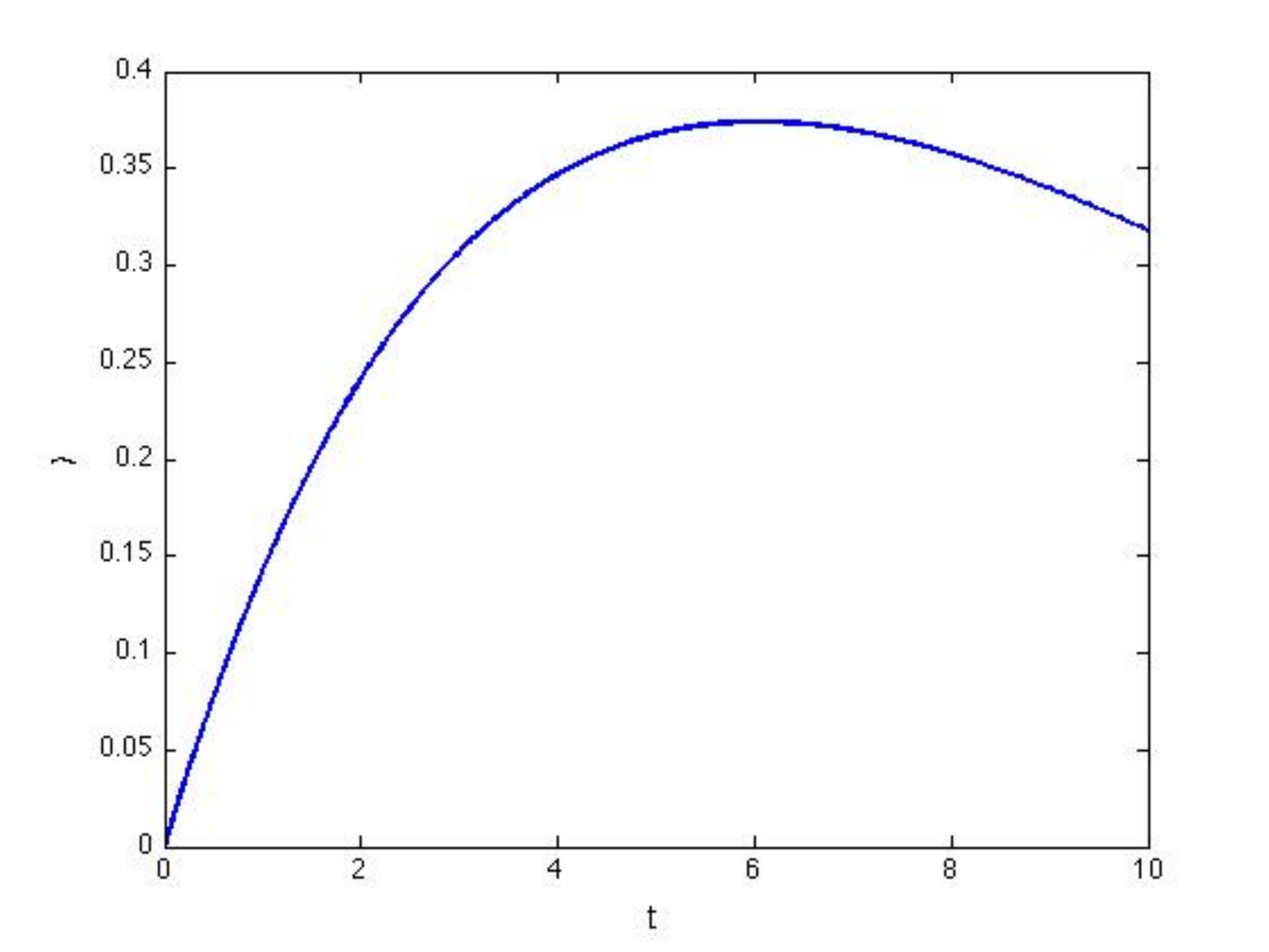}
	\end{subfigure}\hfill
	\begin{subfigure}[c]{0.22\textwidth}
		\centering
		\includegraphics[height=30mm]{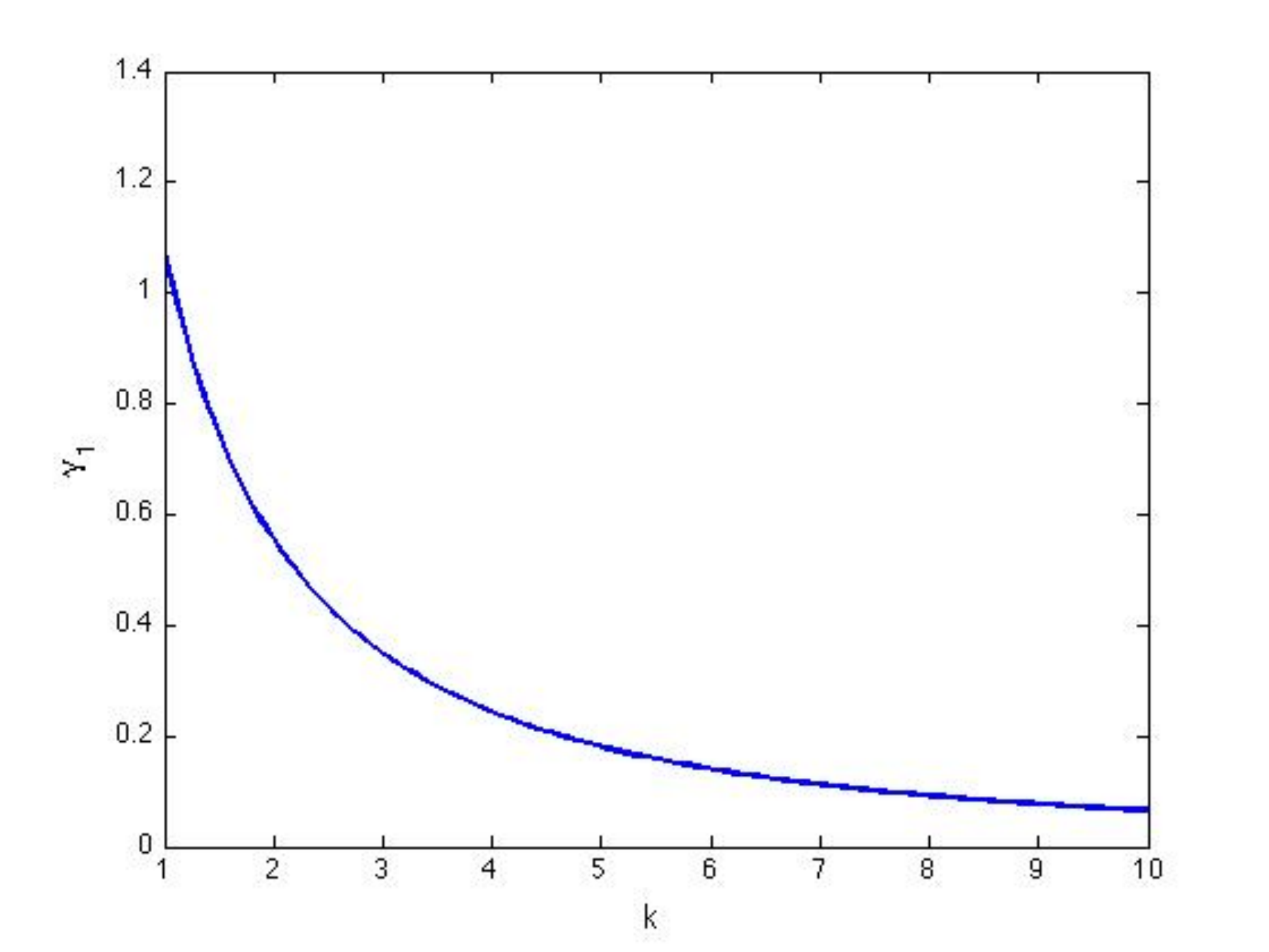}
	\end{subfigure}\hfill
	\begin{subfigure}[c]{0.22\textwidth}
		\centering
		\includegraphics[height=30mm]{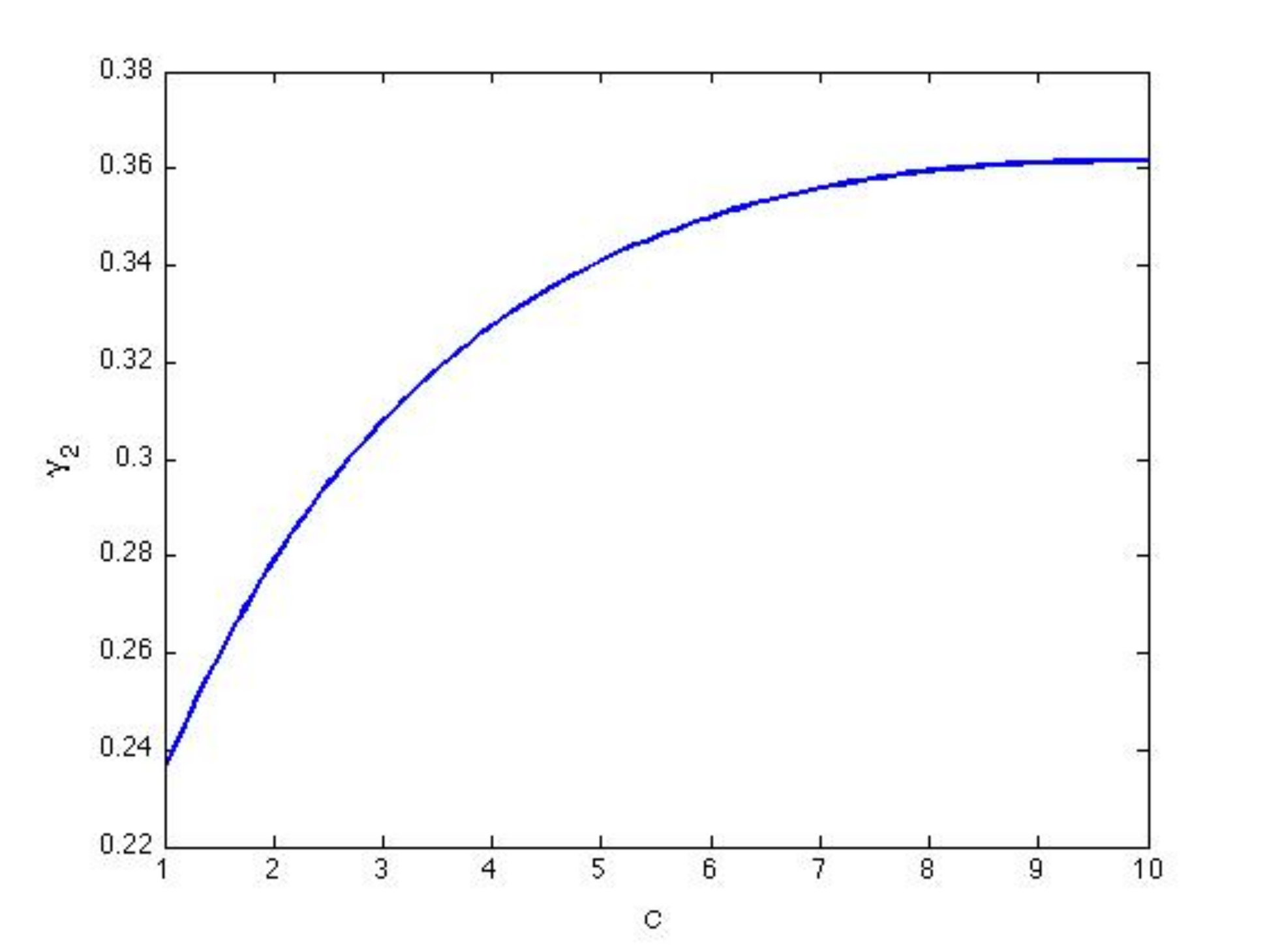}
	\end{subfigure}
\caption{Space functions $\psi_m(\xb)$, time functions $\lambda_m(t)$ and parameter functions $\gamma_{1,m}(k)$ and $\gamma_{2,m}(c)$ (from left to right) obtained for order $m=1,\dots,5$ (from top to bottom).}\label{fig:2D_modes}
\end{figure}

We represent in Figure~\ref{fig:2D_modes} the first five PGD modes obtained without any adaptive strategy.

\begin{figure}[h!]
\centering
\includegraphics[width=8.cm]{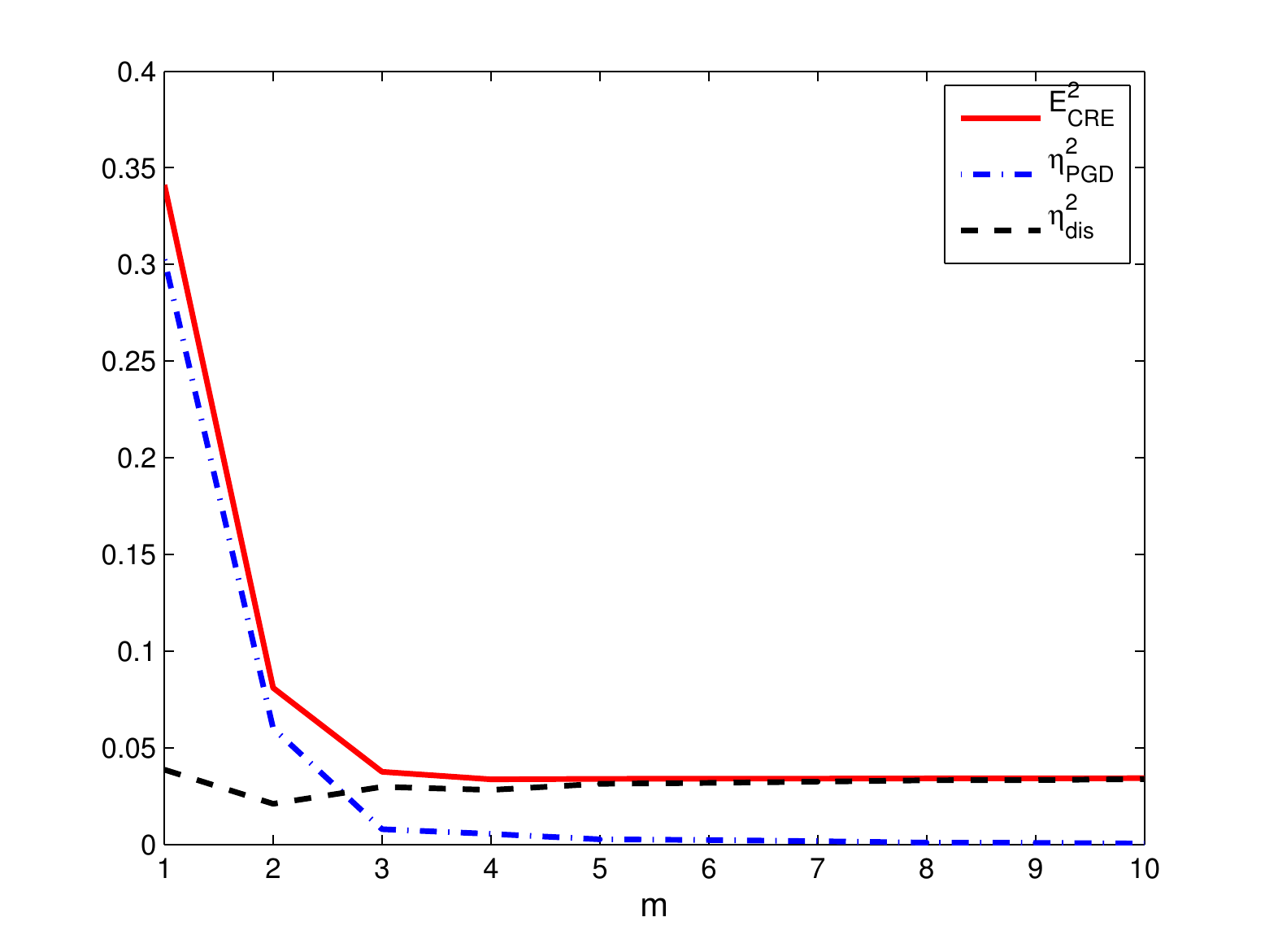}
\caption{Evolutions of the error estimate $E_{\CRE}^2$ and associated error indicators $\eta^2_{\PGD}$ and $\eta_{\dis}^2$ with respect to the number $m$ of PGD modes.}
\label{fig:2D_estimates_global}
\end{figure}

\begin{figure}[h!]
	\centering
	\begin{subfigure}[c]{0.30\textwidth}
		\centering
		\includegraphics[width=6.cm]{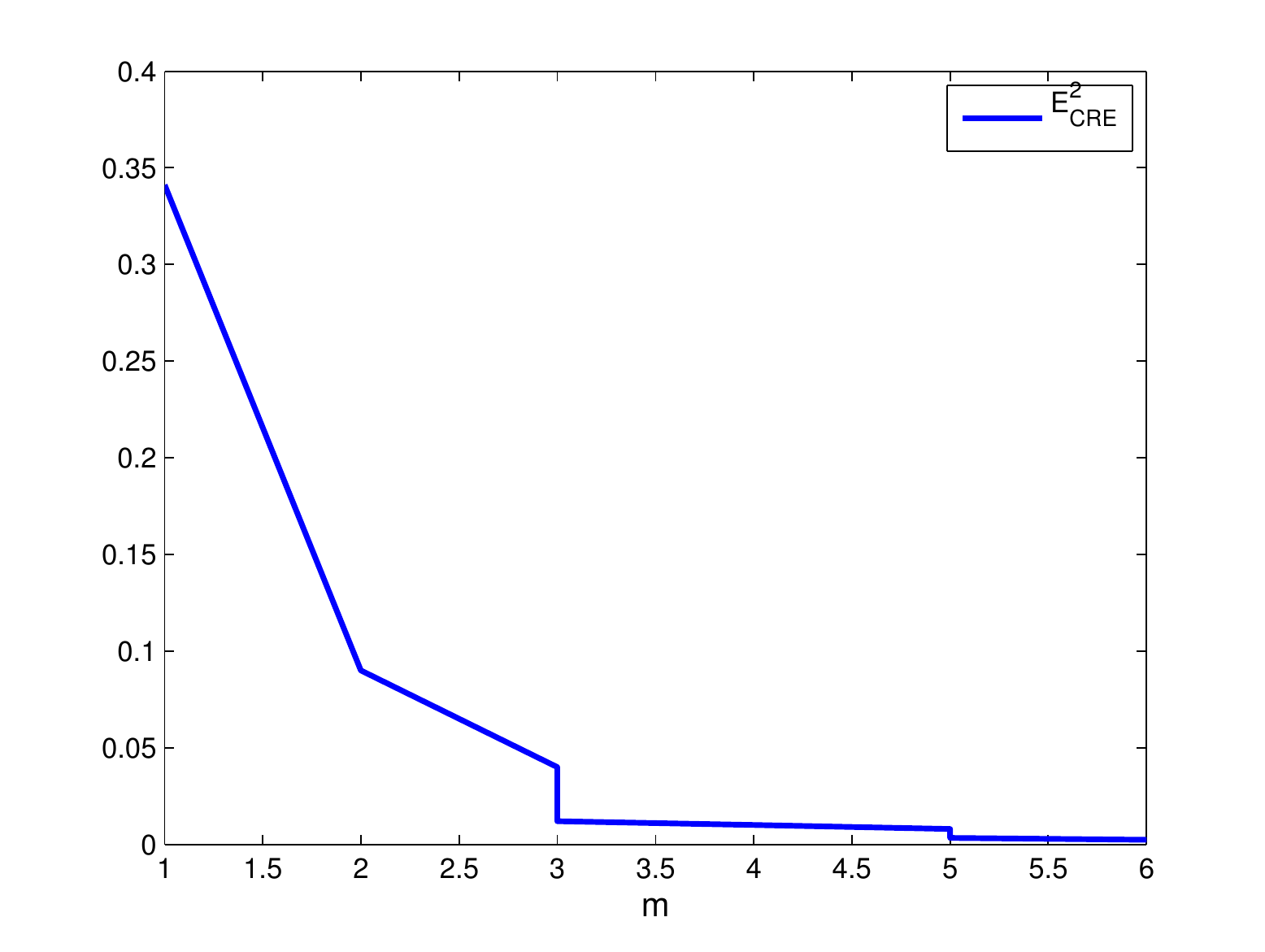}
	\end{subfigure}\hfill
	\begin{subfigure}[c]{0.30\textwidth}
		\centering
		\includegraphics[width=4.cm]{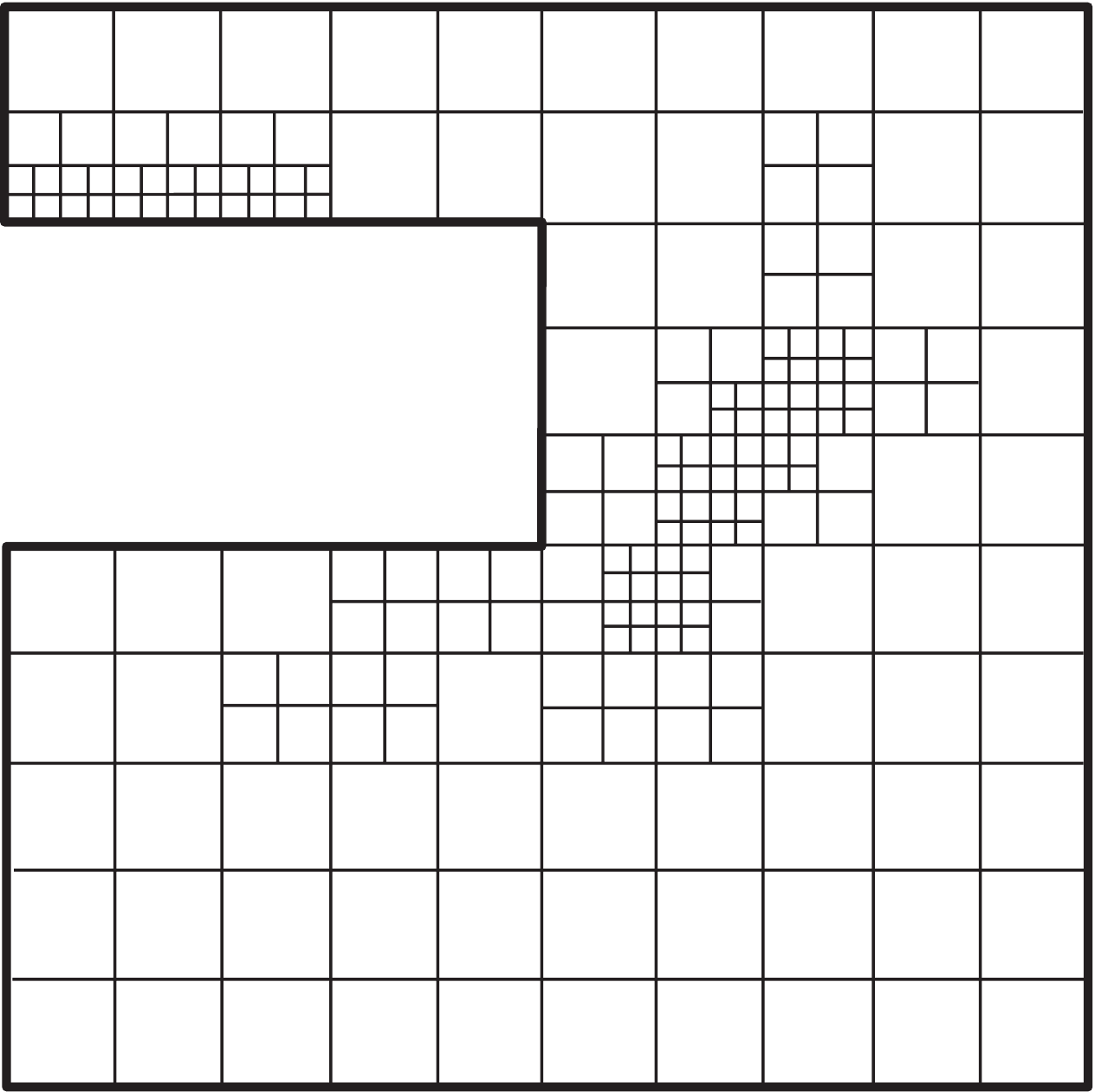}
	\end{subfigure}\hfill
	\begin{subfigure}[c]{0.30\textwidth}
		\centering
		\includegraphics[width=4.cm]{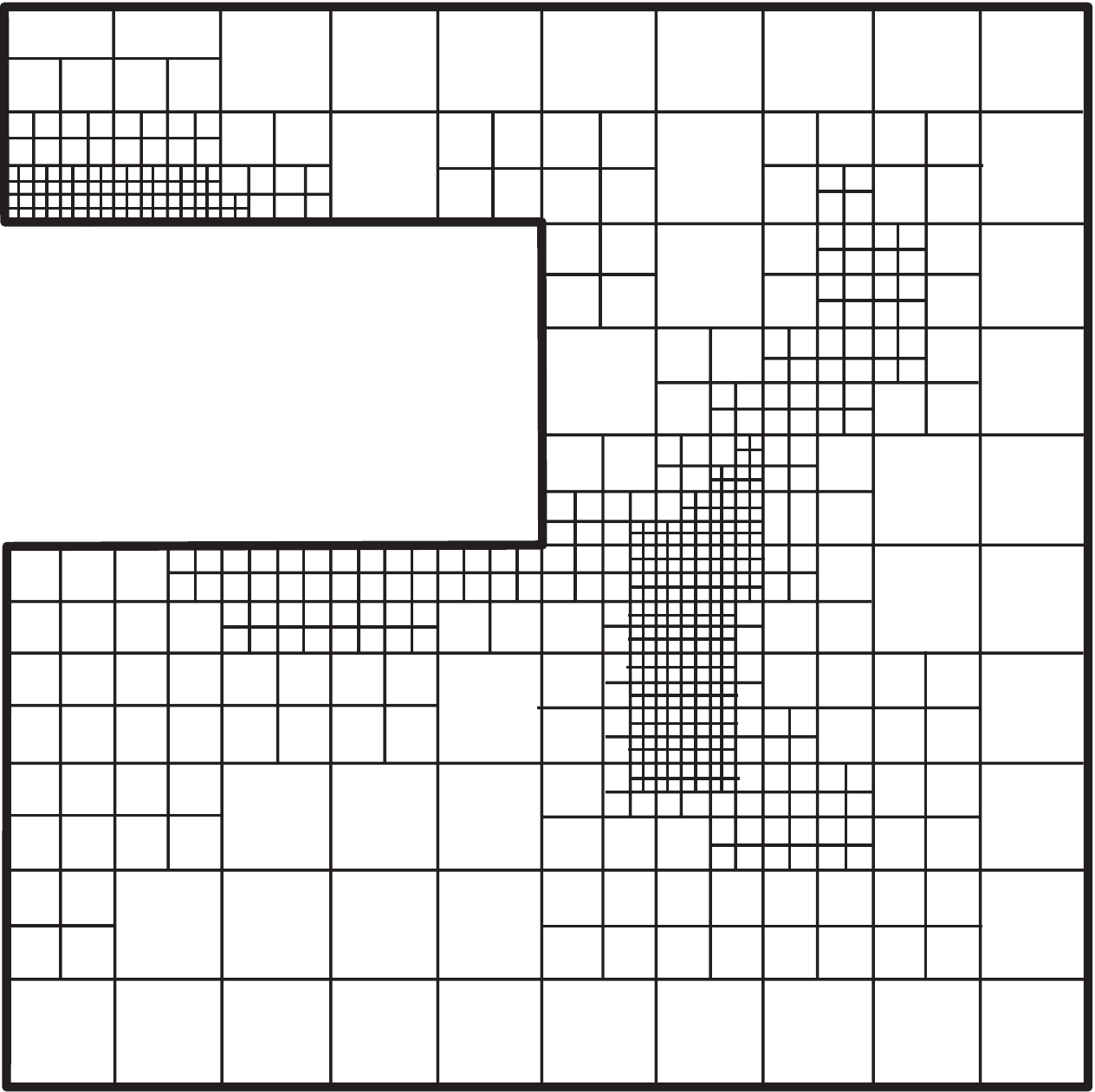}
	\end{subfigure}
\caption{Evolution of the global error estimate $E_{\CRE}^2$ in the adaptive strategy with respect to the number $m$ of PGD modes, and associated refined meshes at order $m=3$ and $m=5$.}
\label{fig:2D_adapt_global}
\end{figure}

We give in Figure~\ref{fig:2D_estimates_global} the convergence of the error estimate $E_{\CRE}$ and associated error indicators $\eta_{\PGD}$ and $\eta_{\dis}$ with respect to the number $m$ of PGD modes, for $m=1,\dots,10$ and for the maximal values obtained with pairs $(k,c)\in P_k\times P_c$. We observe that the error indicator $\eta_{\dis}$ provides a relevant assessment of the discretization error $\Delta_{\dis}$, even for small values of order $m$. Furthermore, we can see that $\eta_{\dis}$ becomes larger than $\eta_{\PGD}$ for $m \geq 3$; this is automatically taken into account in the adaptive process (see Figure~\ref{fig:2D_adapt_global}) in which the number $m$ of computed modes is first increased up to $m=3$, before refining time and space discretizations in order to improve the quality of the approximate PGD solution. The refined meshes obtained after performing mesh adaptation for $m=3$ then $m=5$ are also given in Figure~\ref{fig:2D_adapt_global}. Nested FE meshes based on quad-tree technique are used for practical reasons.

We now consider the control of the local error on a quantity of interest $Q$ defined as the maximal value, for any pair $(k,c)\in P_k\times P_c$, of the average value of the temperature $u$ over a local zone $\om \subset \Om$ and at final time $t=T$:
\begin{equation*}
Q(u)=\max_{(k,c)\in P_k\times P_c}\frac{1}{\abs{\omega}}\into u_{\restrictto{T}}\dom,
\end{equation*}
where subdomain $\omega$ is shown in Figure~\ref{fig:2D_problem}, and $\abs{\omega}$ represents its measure.

\begin{figure}[h!]
	\centering
	\begin{subfigure}[c]{0.30\textwidth}
		\centering
		\includegraphics[height=35mm]{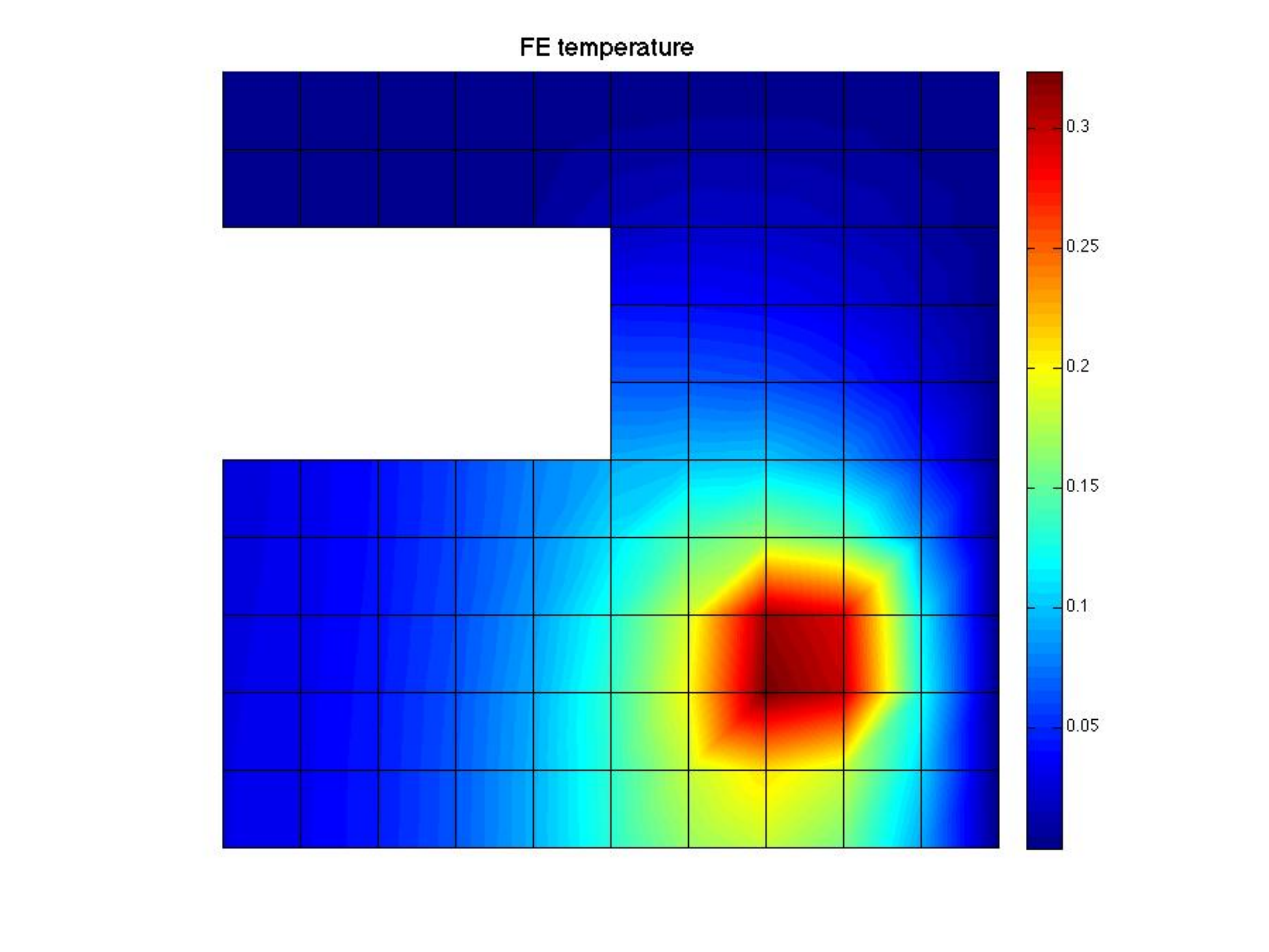}
	\end{subfigure}\hfill
	\begin{subfigure}[c]{0.30\textwidth}
		\centering
		\includegraphics[height=35mm]{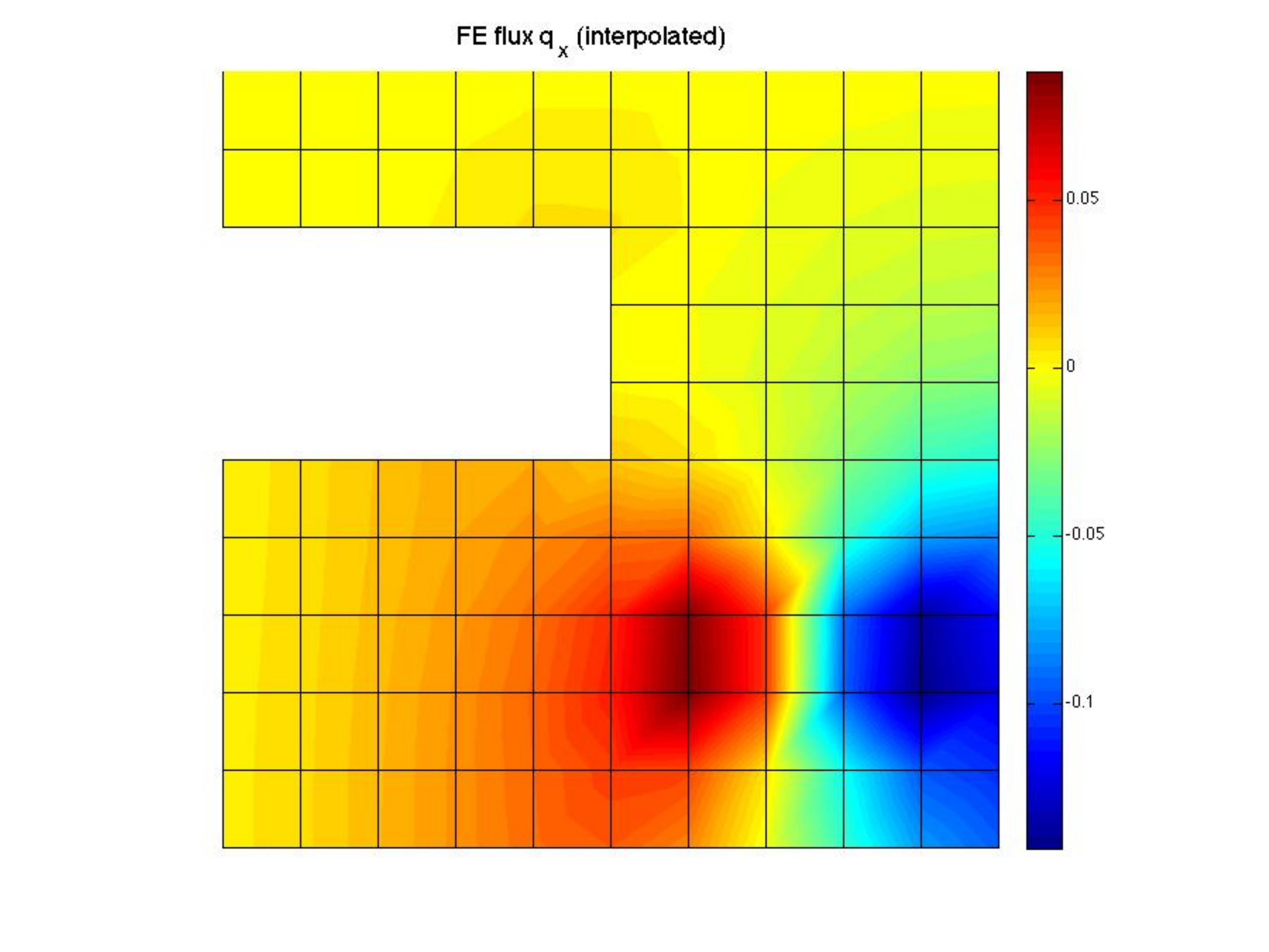}
	\end{subfigure}\hfill
	\begin{subfigure}[c]{0.30\textwidth}
		\centering
		\includegraphics[height=35mm]{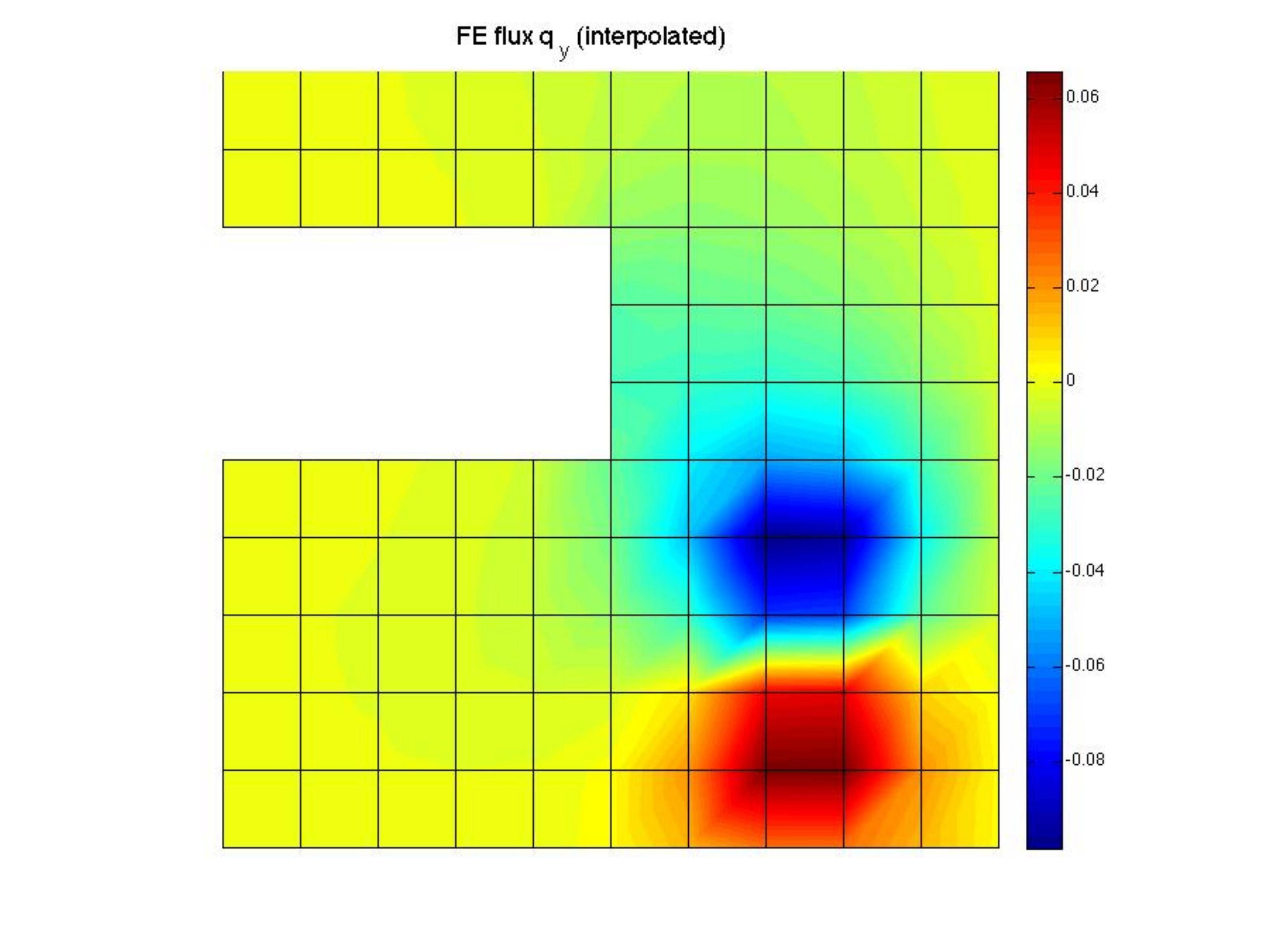}
	\end{subfigure}
	\caption{Adjoint FE solution $\tilde{u}^{h,\Delta t}$ (left) and associated flux components $\tilde{q}_x^{h,\Delta t}$ (center) and $\tilde{q}_y^{h,\Delta t}$ (right).}\label{fig:2D_FEsol_adj}
\end{figure}

We give in Figure~\ref{fig:2D_FEsol_adj} the adjoint FE solution $\tilde{u}^{h,\Delta t}$ and associated flux $\tilde{\qb}^{h,\Delta t}$ for $(k,c)=(1,1)$ and at time $t=T-\Delta t$.

\begin{figure}[h!]
\centering
\includegraphics[width=8.cm]{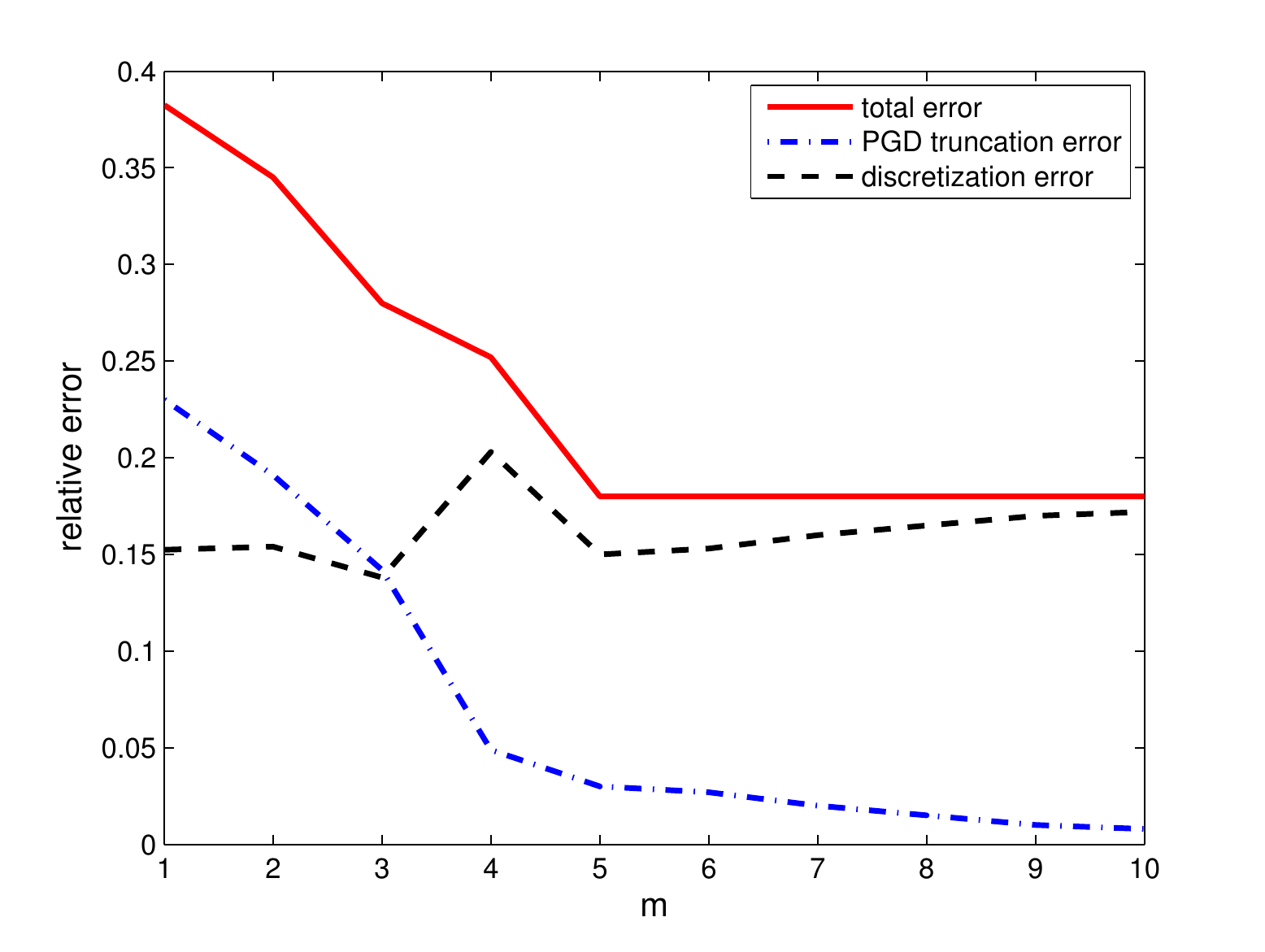}
\caption{Evolutions of the error estimate and indicators with respect to the number $m$ of PGD modes.}\label{fig:2D_estimates_local}
\end{figure}

The evolutions of the normalized upper bound on $\Delta Q - Q_{\corr}$ as well as specific normalized error indicators of $\Delta Q_{\PGD} - Q_{\corr}^{h,\Delta t}$ and $\Delta Q_{\dis} - Q_{\corr} + Q_{\corr}^{h,\Delta t}$ are given in Figure~\ref{fig:2D_estimates_local} with respect to the number $m$ of computed modes.

\begin{figure}[h!]
	\centering
	\begin{subfigure}[c]{0.45\textwidth}
		\centering
		\includegraphics[width=60mm]{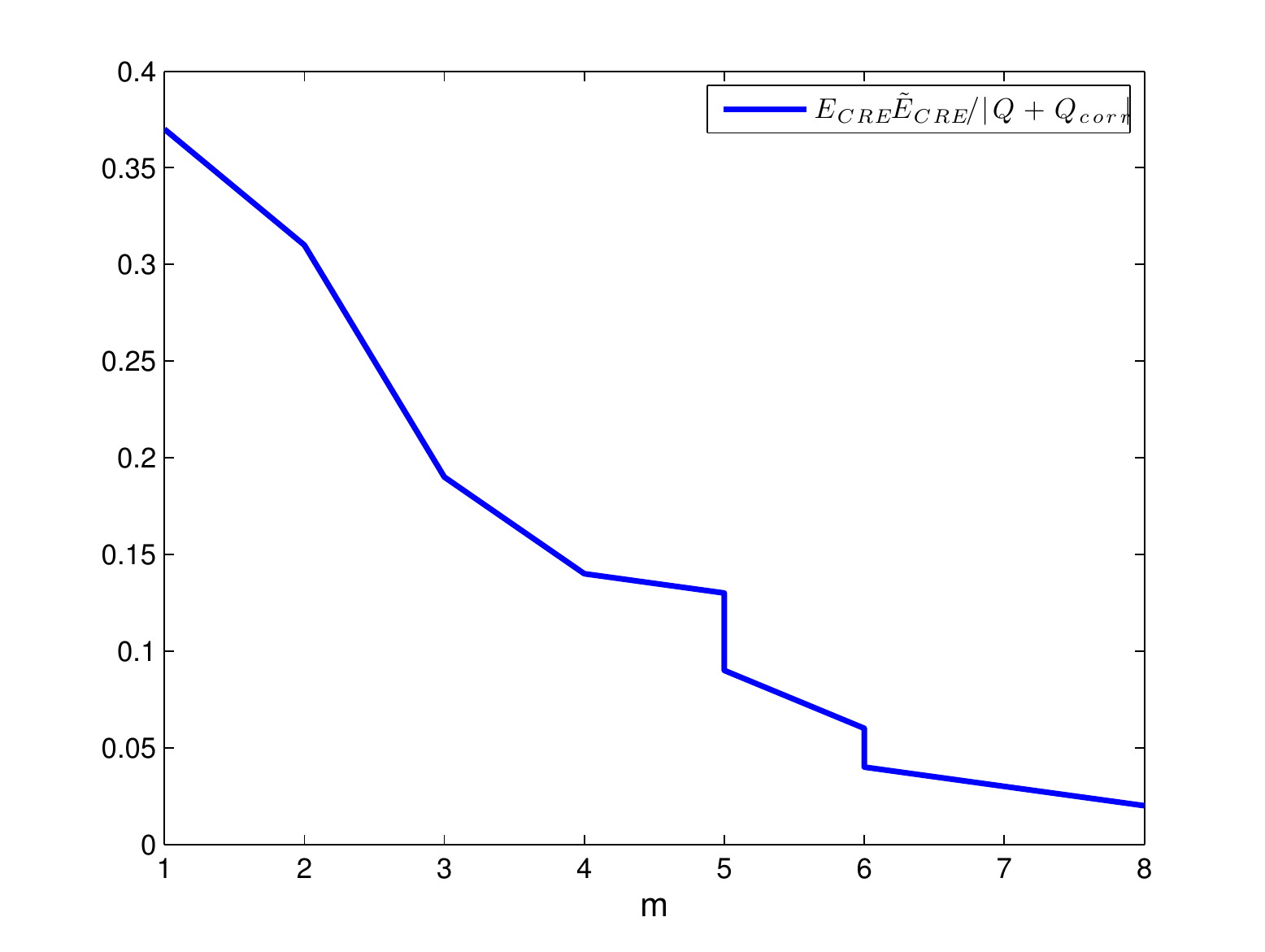}
	\end{subfigure}\hfill
	\begin{subfigure}[c]{0.45\textwidth}
		\centering
		\includegraphics[width=40mm]{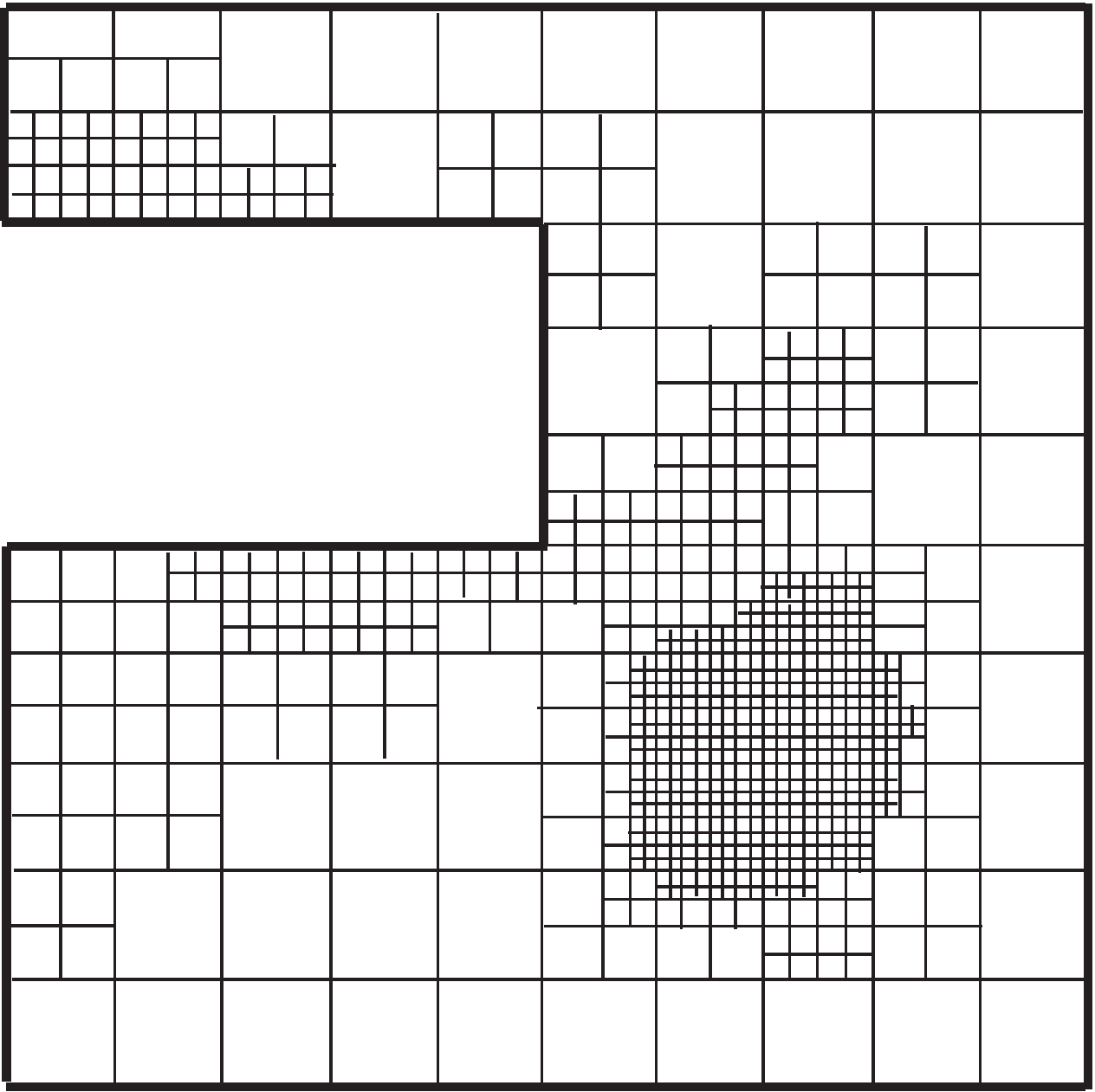}
	\end{subfigure}
	\caption{Evolution of the normalized upper bound on $\Delta Q - Q_{\corr}$ in the adaptive strategy with respect to the number $m$ of PGD modes, and associated refined mesh at order $m=6$.}\label{fig:2D_adapt_local}
\end{figure}

We give in Figure~\ref{fig:2D_adapt_local} the convergence of the normalized upper bound on $\Delta Q - Q_{\corr}$ when performing the adaptive strategy, as well as the refined mesh obtained to compute mode $m=6$.

\begin{figure}[h!]
	\centering
	\begin{subfigure}[c]{0.45\textwidth}
		\centering
		\includegraphics[width=60mm]{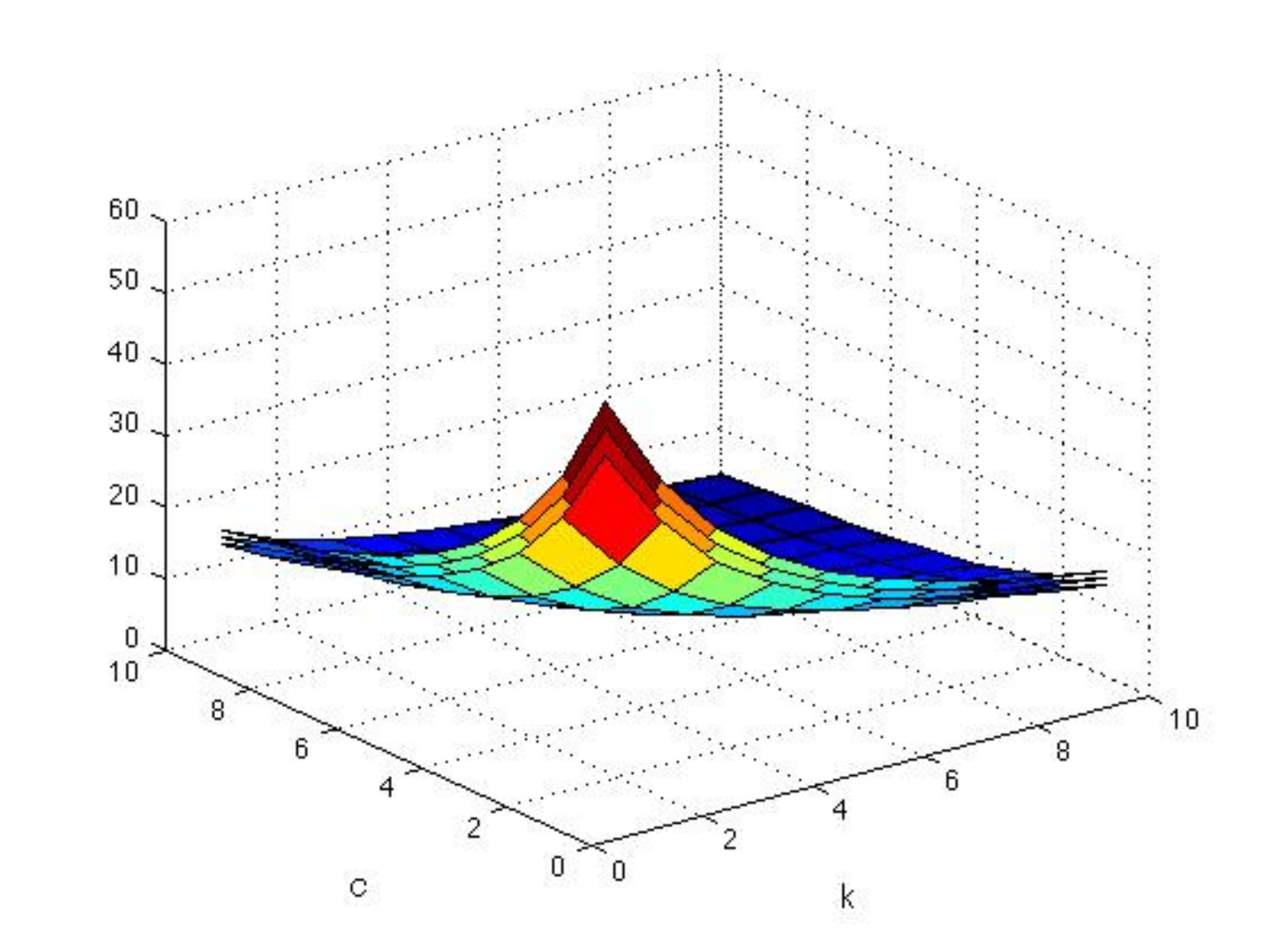}
	\end{subfigure}\hfill
	\begin{subfigure}[c]{0.45\textwidth}
		\centering
		\includegraphics[width=60mm]{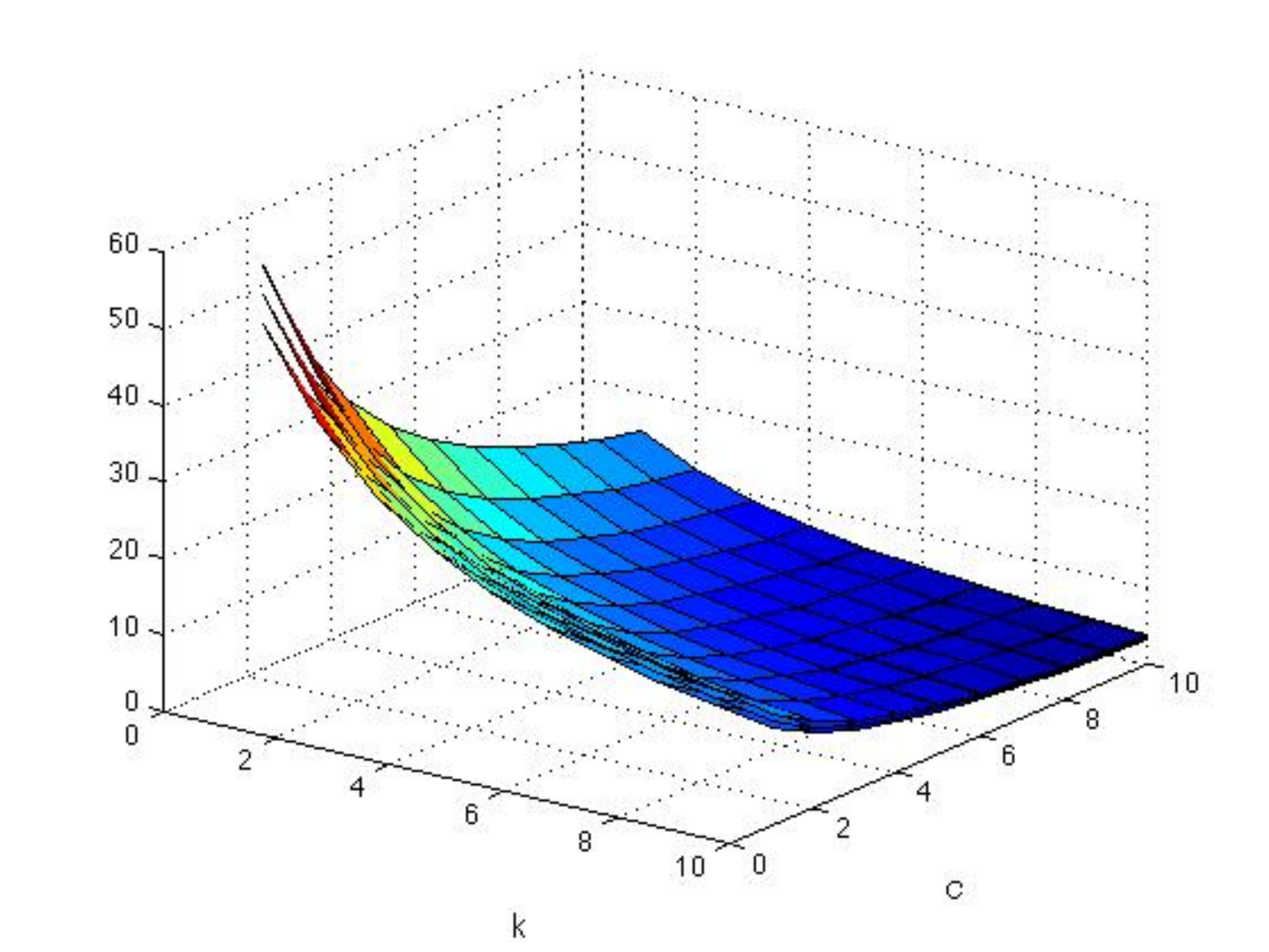}
	\end{subfigure}
	\caption{Evolutions of the predicted value of $Q(u_m^{h,\Delta T})+Q_{\corr}$ and associated error bounds with respect to parameters $k$ and $c$.}\label{fig:2D_evol_Q}
\end{figure}

Eventually, considering $Q(u)
$ as a function of $k$ and $c$, we give in Figure~\ref{fig:2D_evol_Q} the obtained mapping of $Q(u_m^{h,\Delta T})+Q_{\corr}$ and associated guaranteed error bounds over the range of variations of $(k,c) \in P_k\times P_c$, for $m=5$.

\subsection{Three-dimensional elasticity problem}

\begin{figure}[h!]
	\centering
	\begin{subfigure}[c]{0.22\textwidth}
		\centering
		\includegraphics[height=30mm]{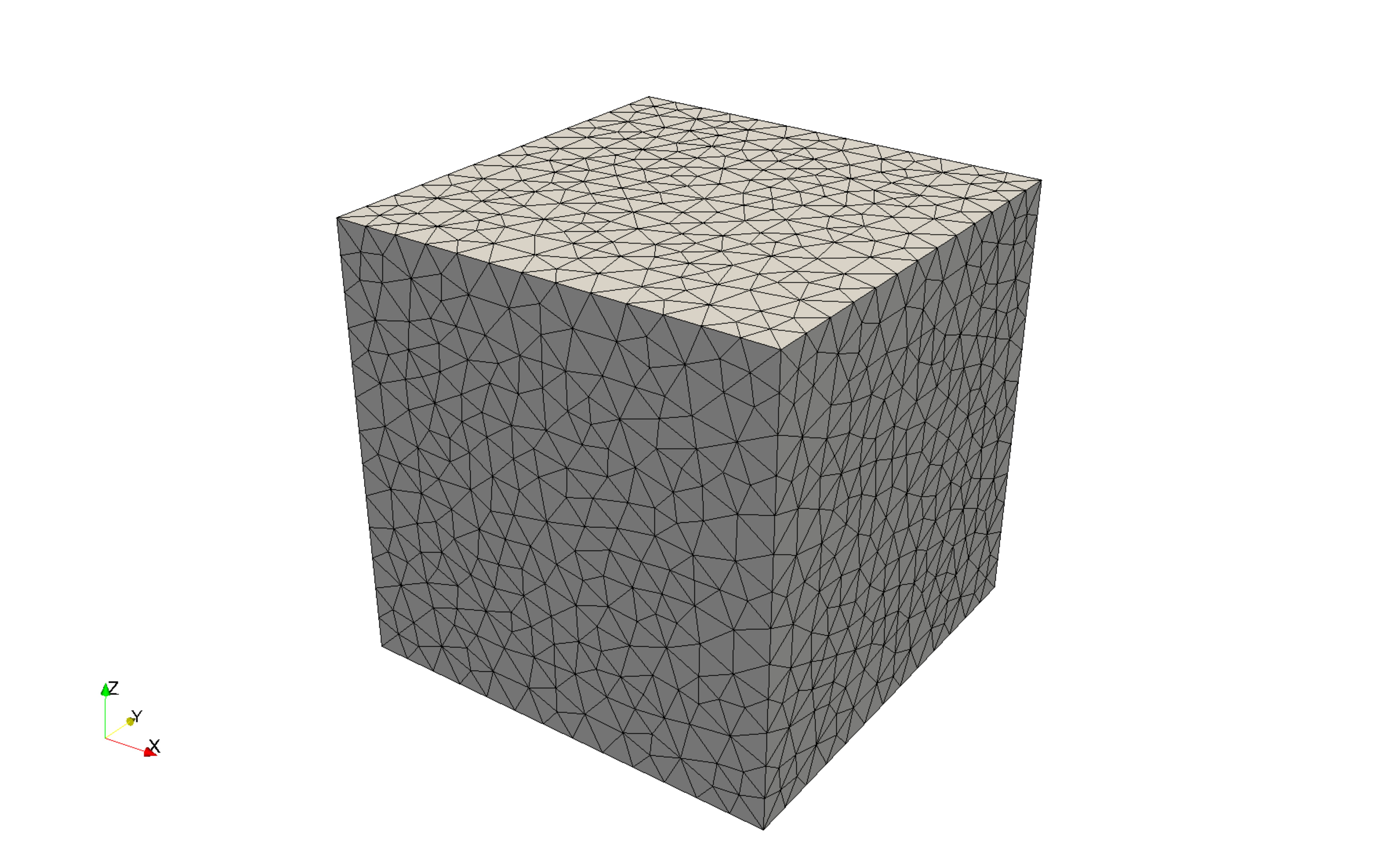}
	\end{subfigure}\hfill
	\begin{subfigure}[c]{0.22\textwidth}
		\centering
		\includegraphics[height=30mm]{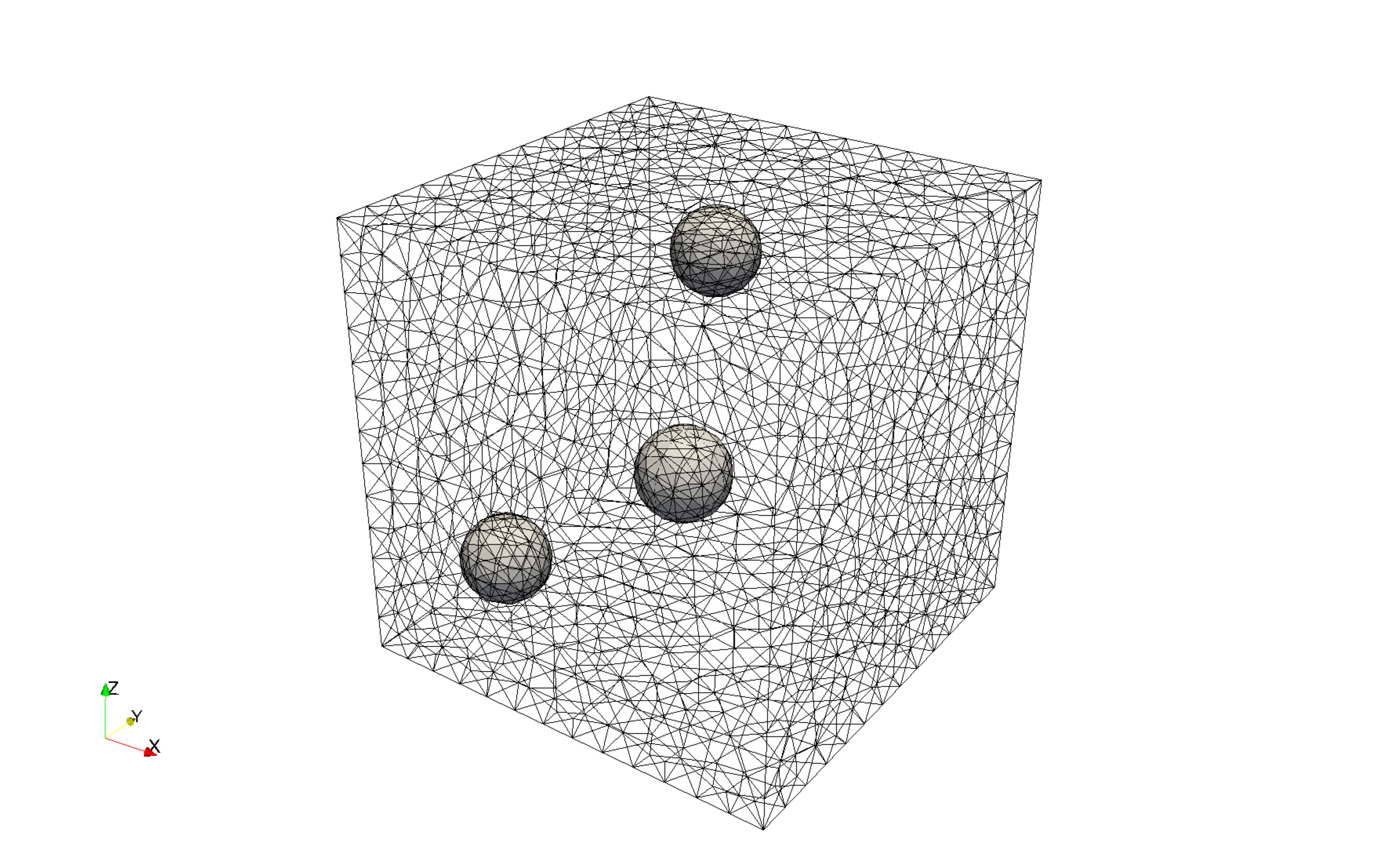}
	\end{subfigure}\hfill
	\begin{subfigure}[c]{0.22\textwidth}
		\centering
		\includegraphics[height=30mm]{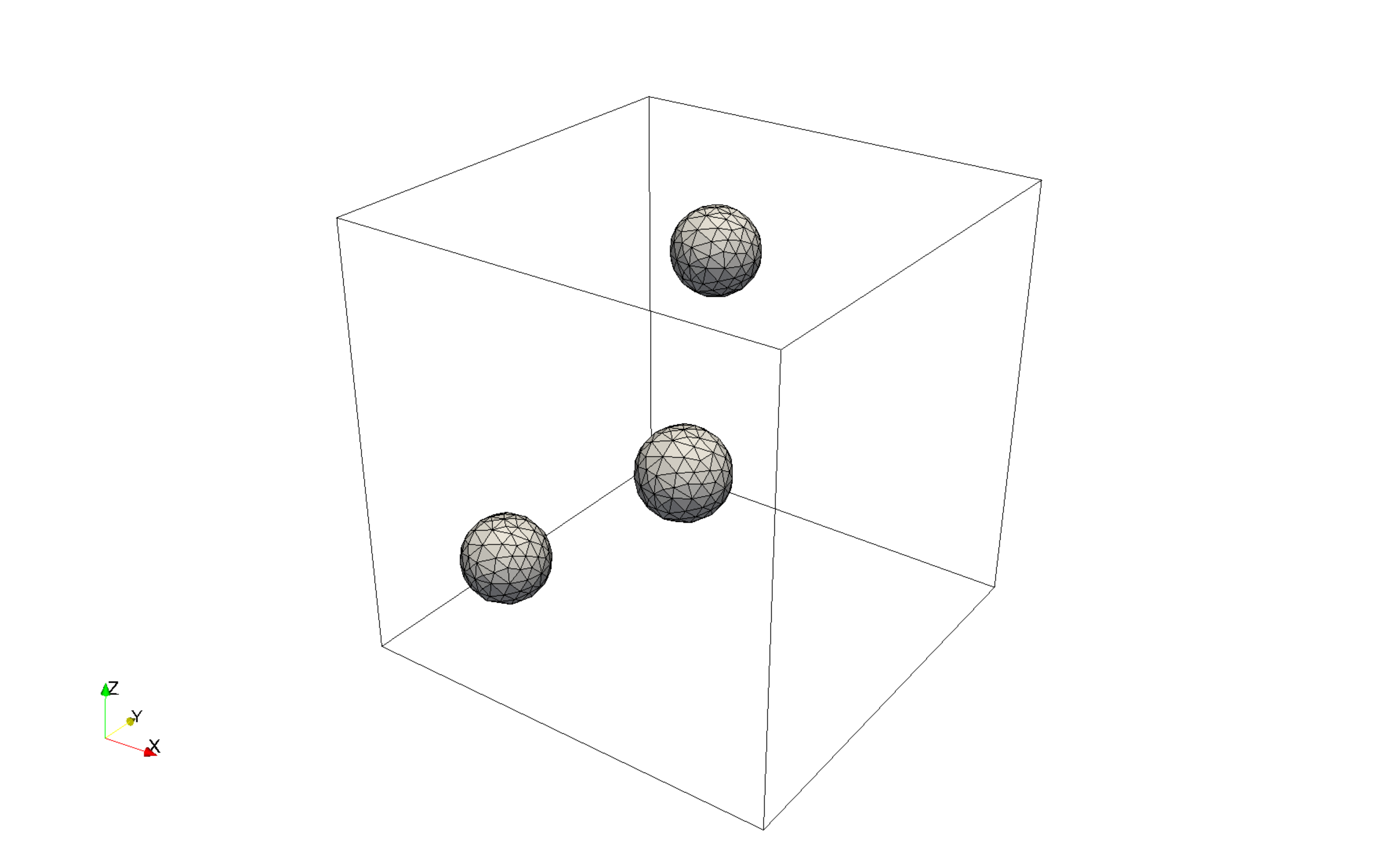}
	\end{subfigure}\hfill
	\begin{subfigure}[c]{0.22\textwidth}
		\centering
		\includegraphics[height=30mm]{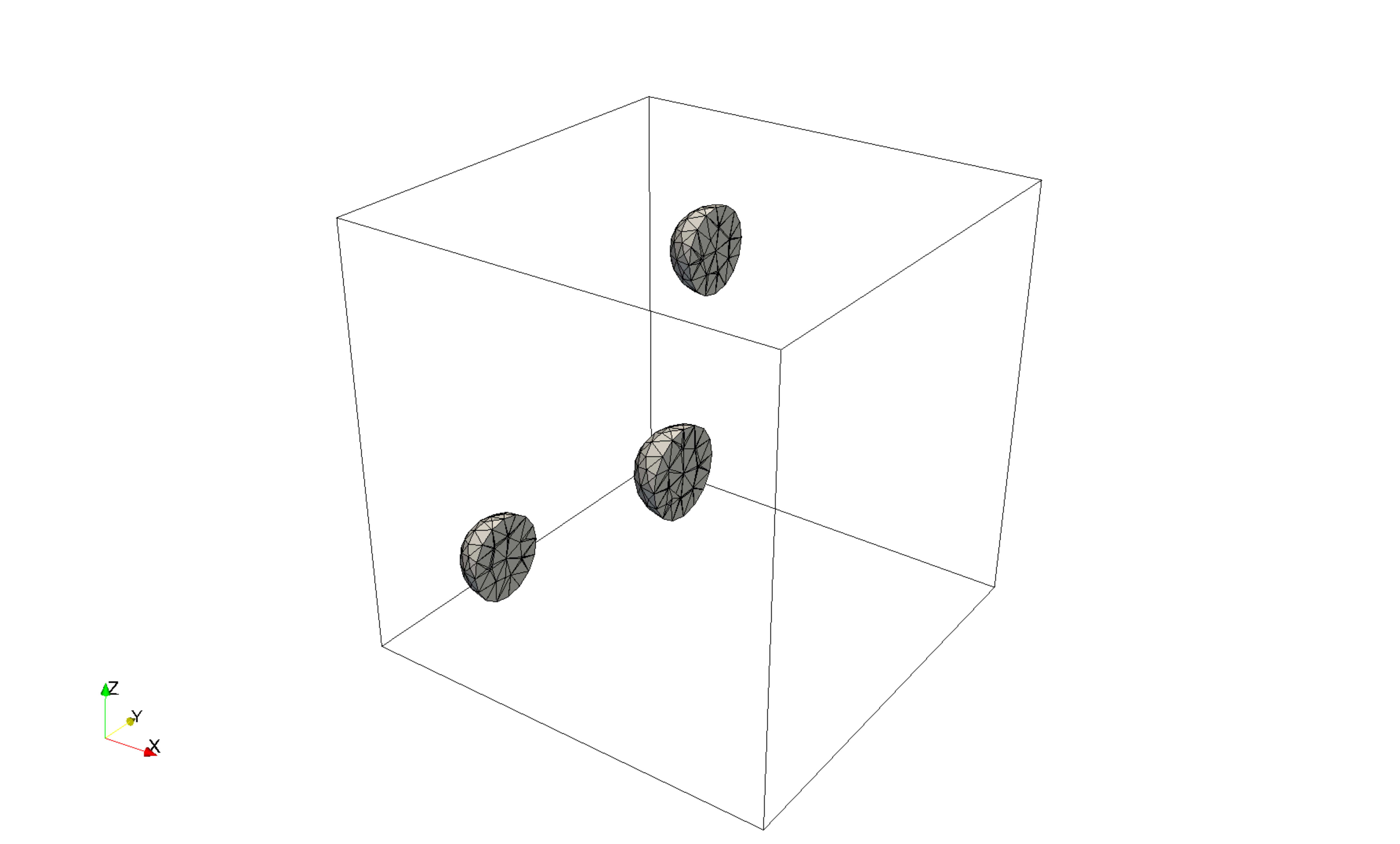}
	\end{subfigure}
	\caption{3D elasticity problem: space domain with three inclusions and associated FE mesh.}\label{fig:3D_problem}
\end{figure}

We eventually consider an elastic cube, of size $1~\text{m}\times 1~\text{m} \times 1~\text{m}$ with three spherical inclusions for which Young's moduli $E_i \in P_E = \intervalcc{1}{10}$ ($1\leq i\leq 3$) are parameters, so that the order~$m$ PGD representation reads $\ub_m(\xb,E_1,E_2,E_3)$. The three inclusions have the same radius $r=0.1$~m, and their centers are respectively located at points $c_1=(0.2,0.2,0.2)$, $c_2=(0.6,0.3,0.5)$ and $c_3=(0.4,0.7,0.8)$ (see Figure~\ref{fig:3D_problem}). The cube is clamped along the plane located at $x=0$ and subjected to a unit traction force $\Fb_d = +\xb$ applied on the plane located at $x=1$.

\begin{figure}[h!]
	\centering
	\rotatebox[origin=c]{90}{$1\textsuperscript{st}$ mode}\hfill
	\begin{subfigure}[c]{0.22\textwidth}
		\centering
		\includegraphics[height=25mm]{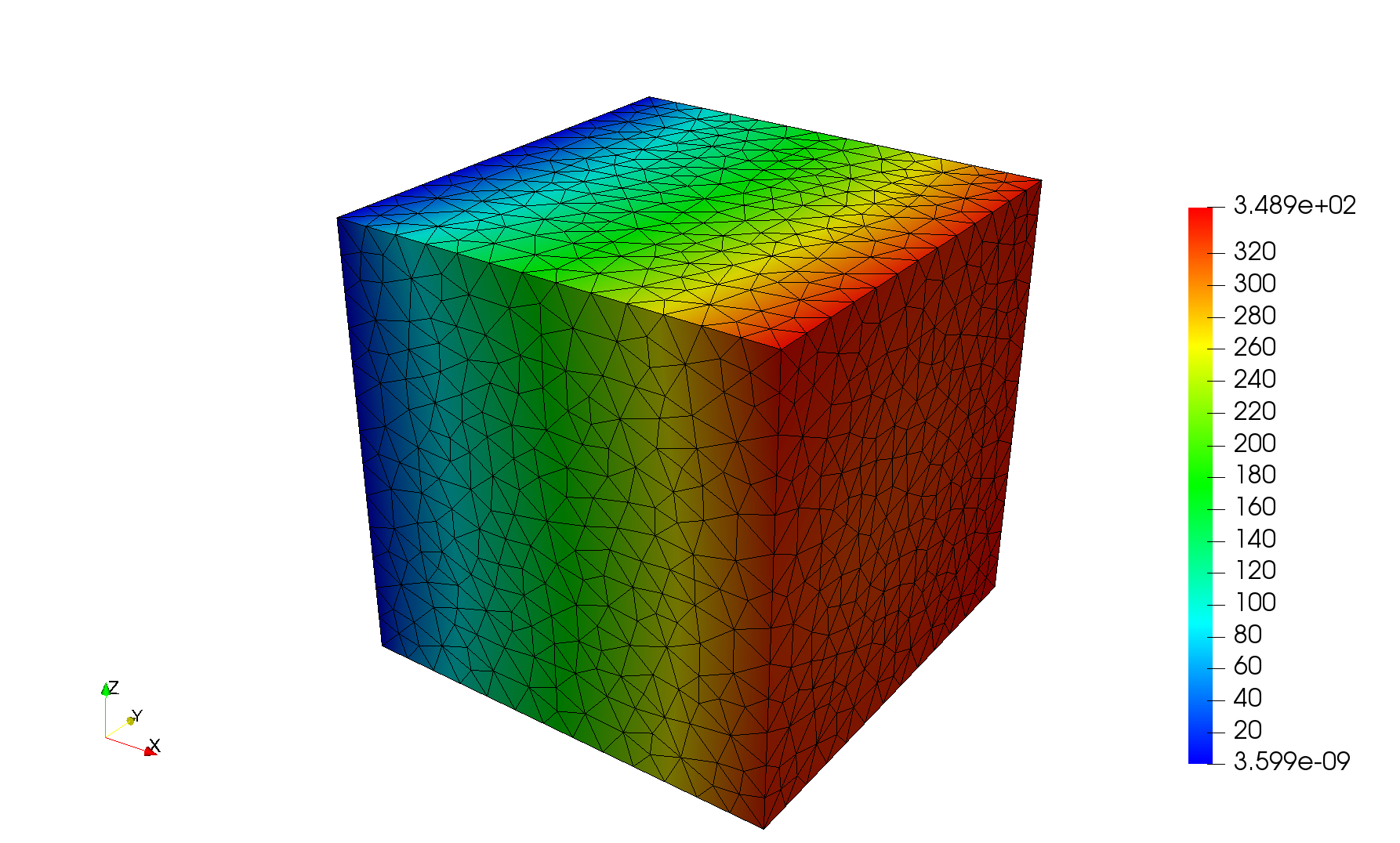}
	\end{subfigure}\hfill
	\begin{subfigure}[c]{0.22\textwidth}
		\centering
		\includegraphics[height=25mm]{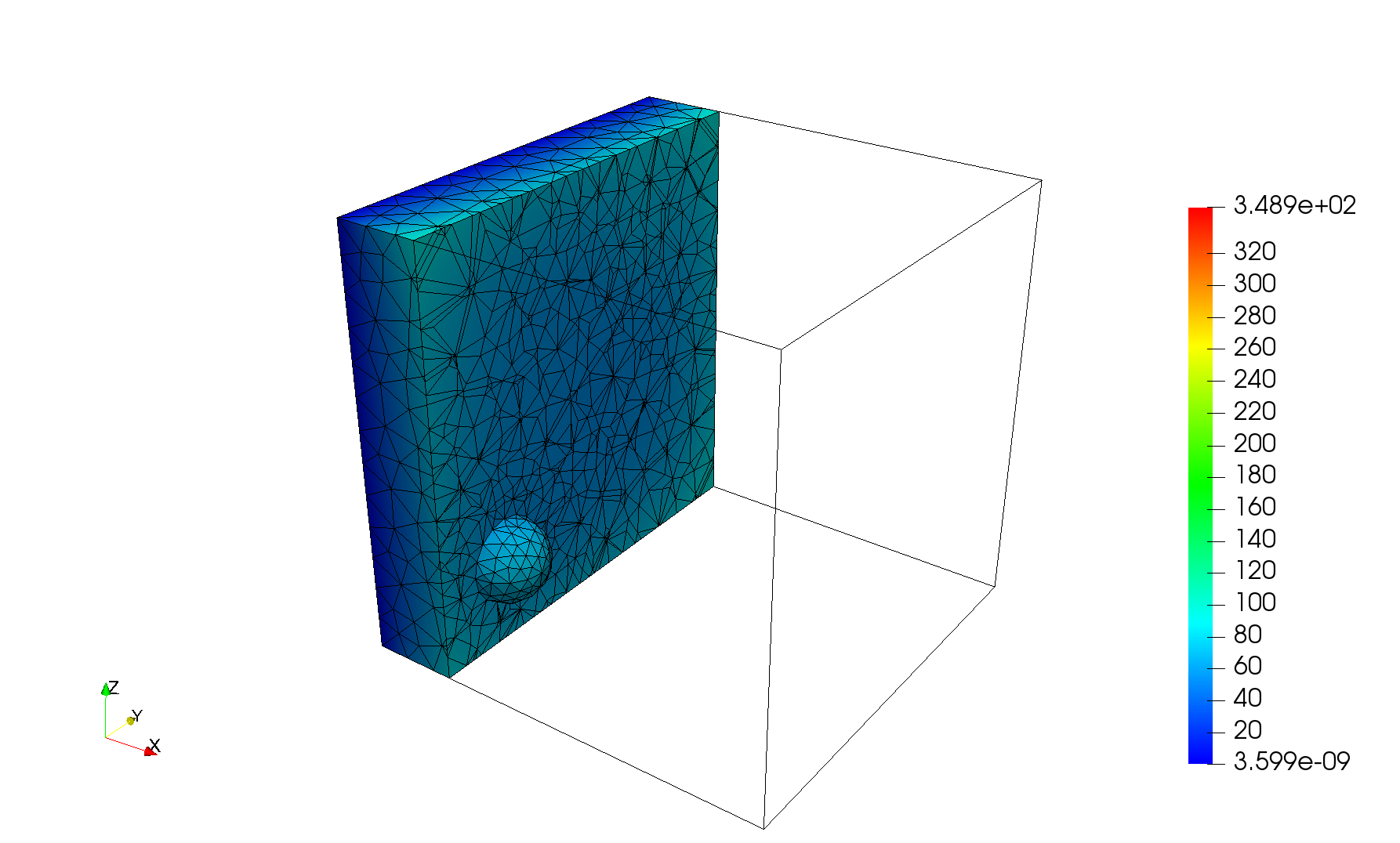}
	\end{subfigure}\hfill
	\begin{subfigure}[c]{0.22\textwidth}
		\centering
		\includegraphics[height=25mm]{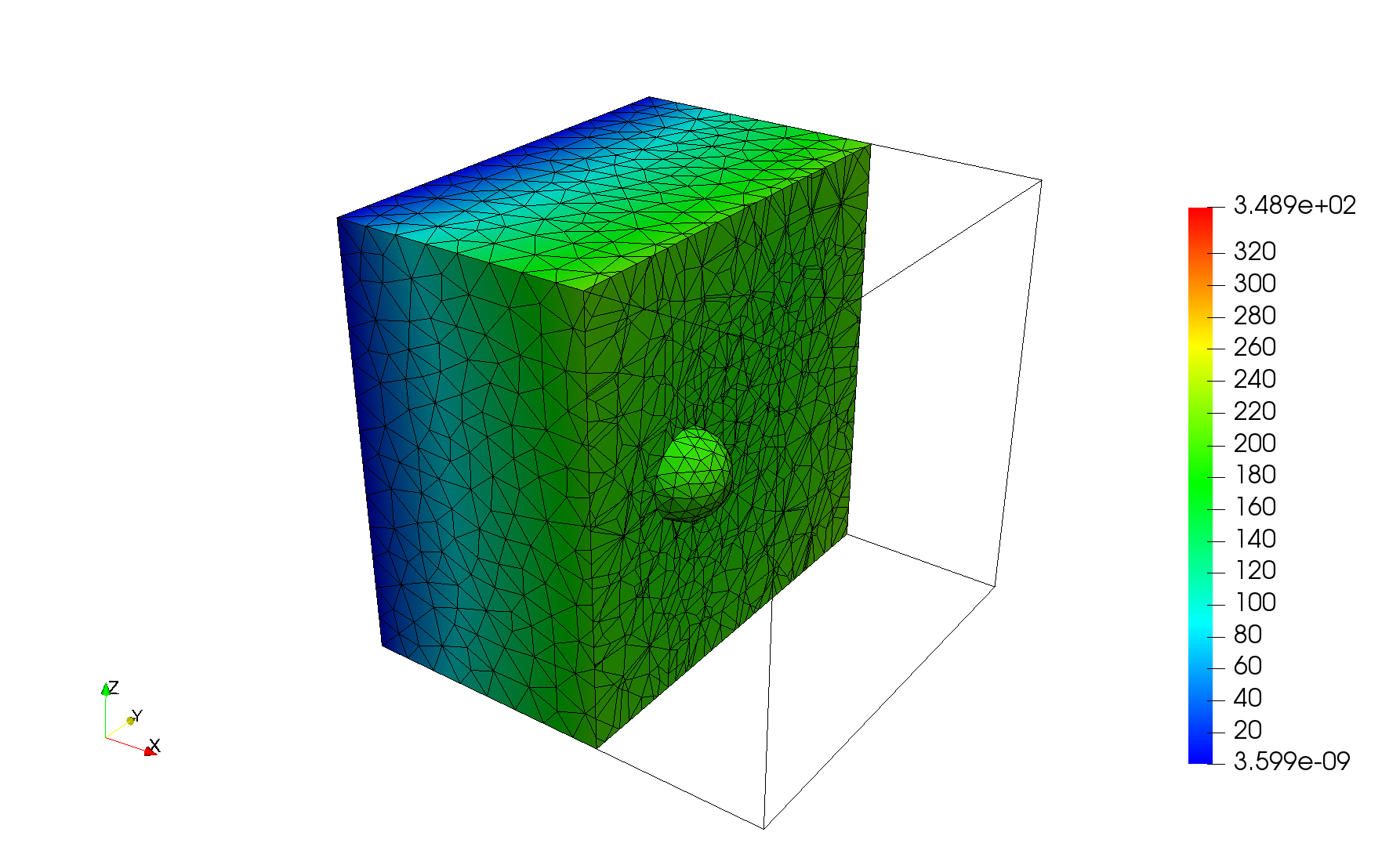}
	\end{subfigure}\hfill
	\begin{subfigure}[c]{0.22\textwidth}
		\centering
		\includegraphics[height=25mm]{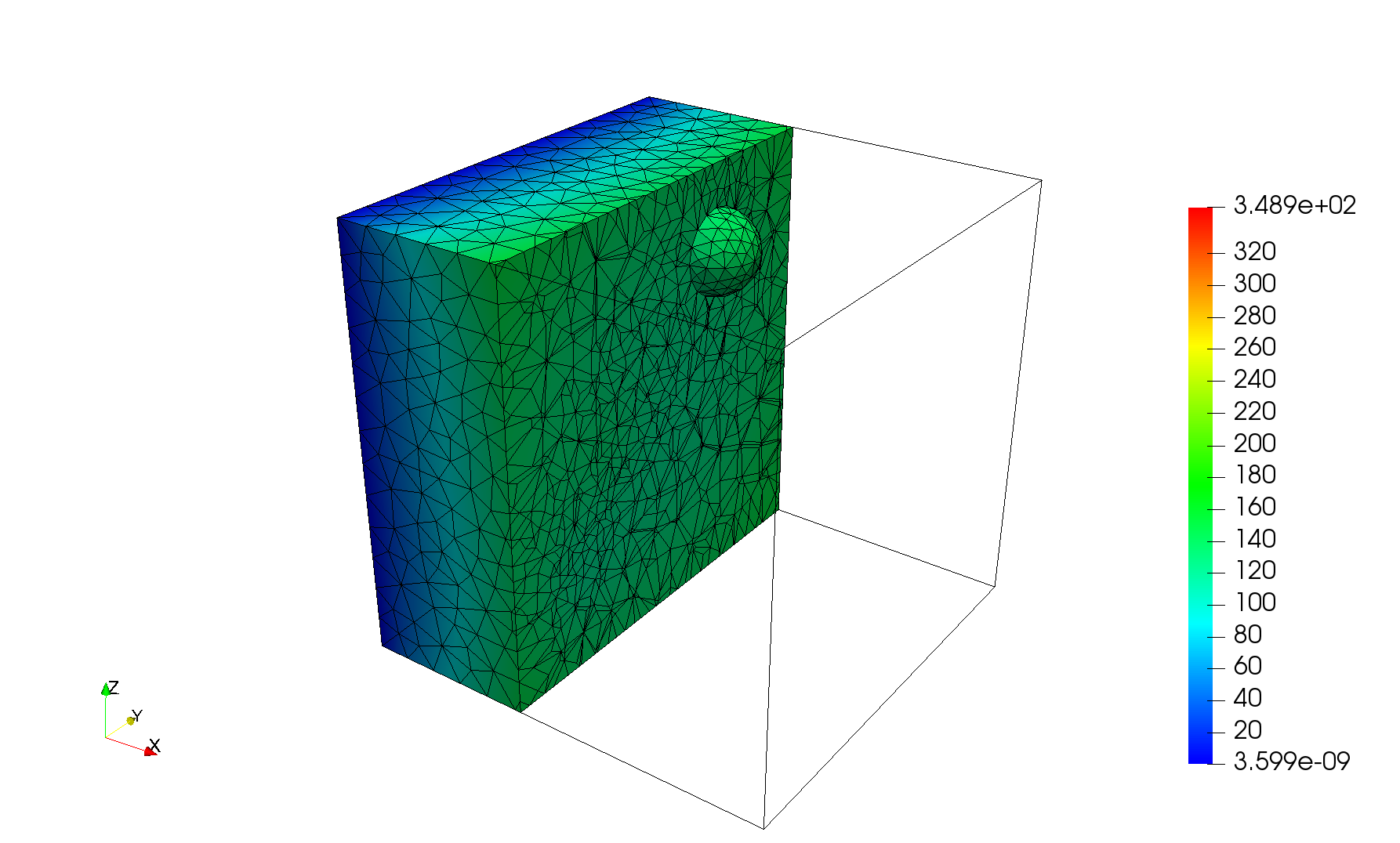}
	\end{subfigure}
	\\
	\rotatebox[origin=c]{90}{$2\textsuperscript{nd}$ mode}\hfill
	\begin{subfigure}[c]{0.22\textwidth}
		\centering
		\includegraphics[height=25mm]{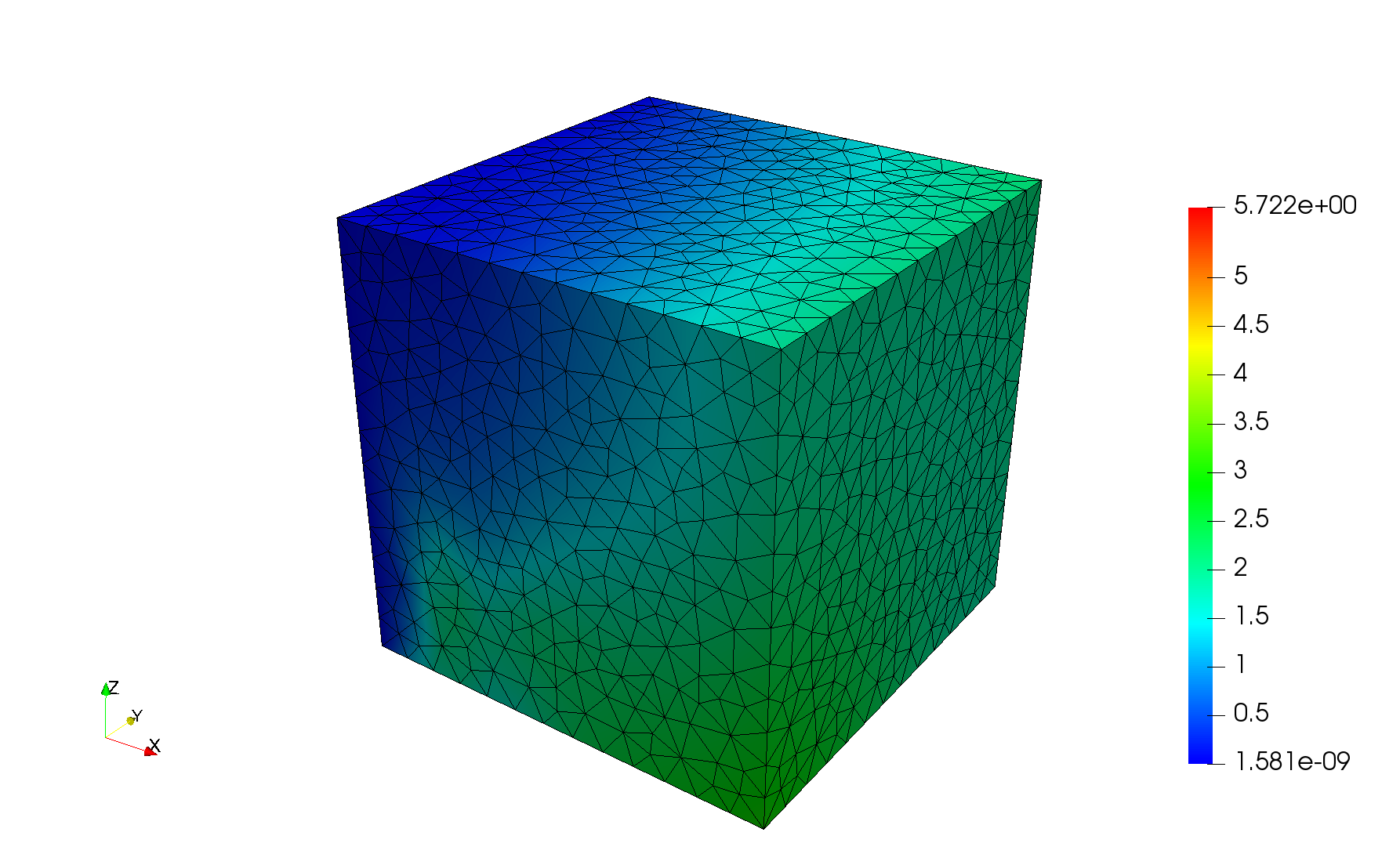}
	\end{subfigure}\hfill
	\begin{subfigure}[c]{0.22\textwidth}
		\centering
		\includegraphics[height=25mm]{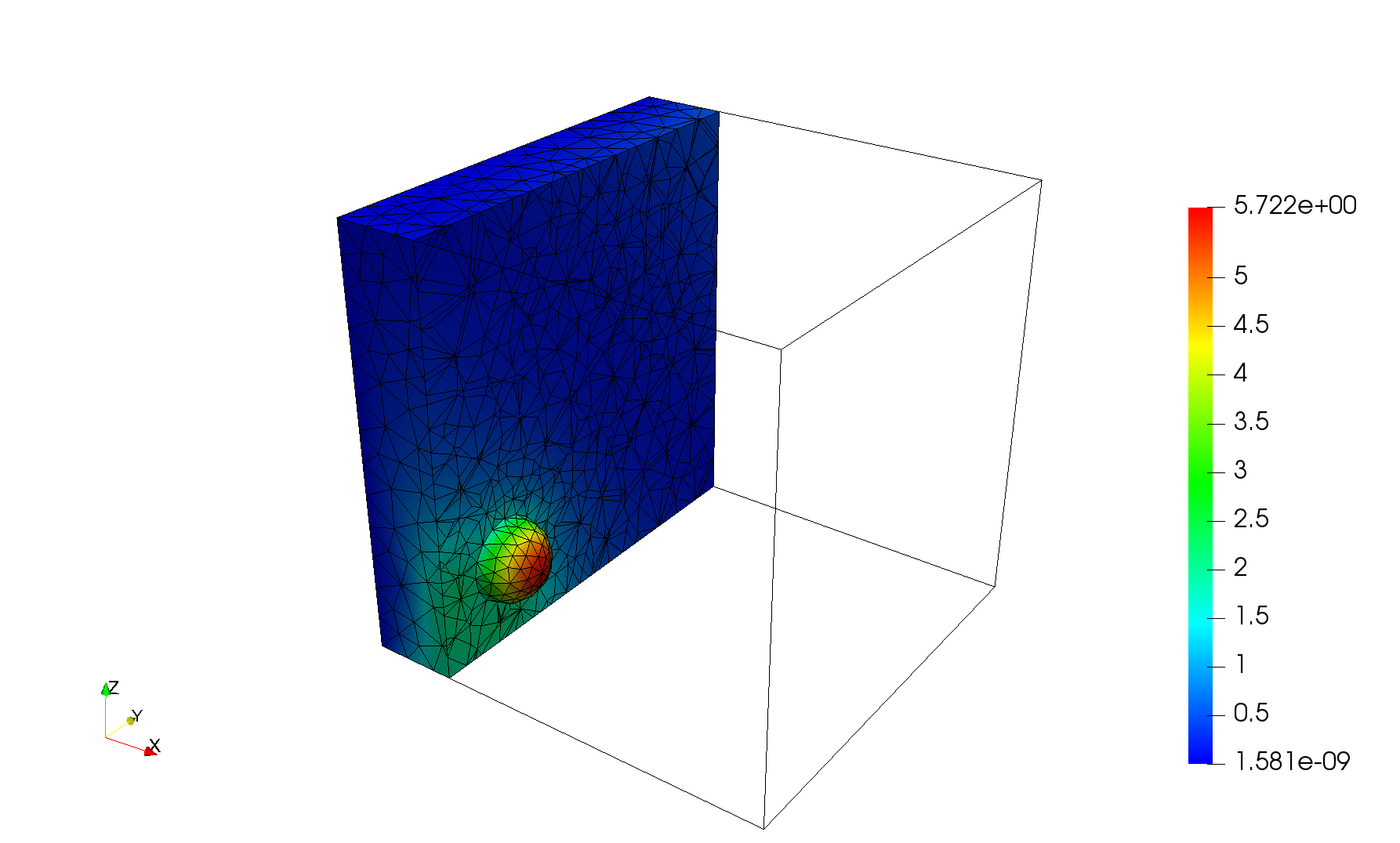}
	\end{subfigure}\hfill
	\begin{subfigure}[c]{0.22\textwidth}
		\centering
		\includegraphics[height=25mm]{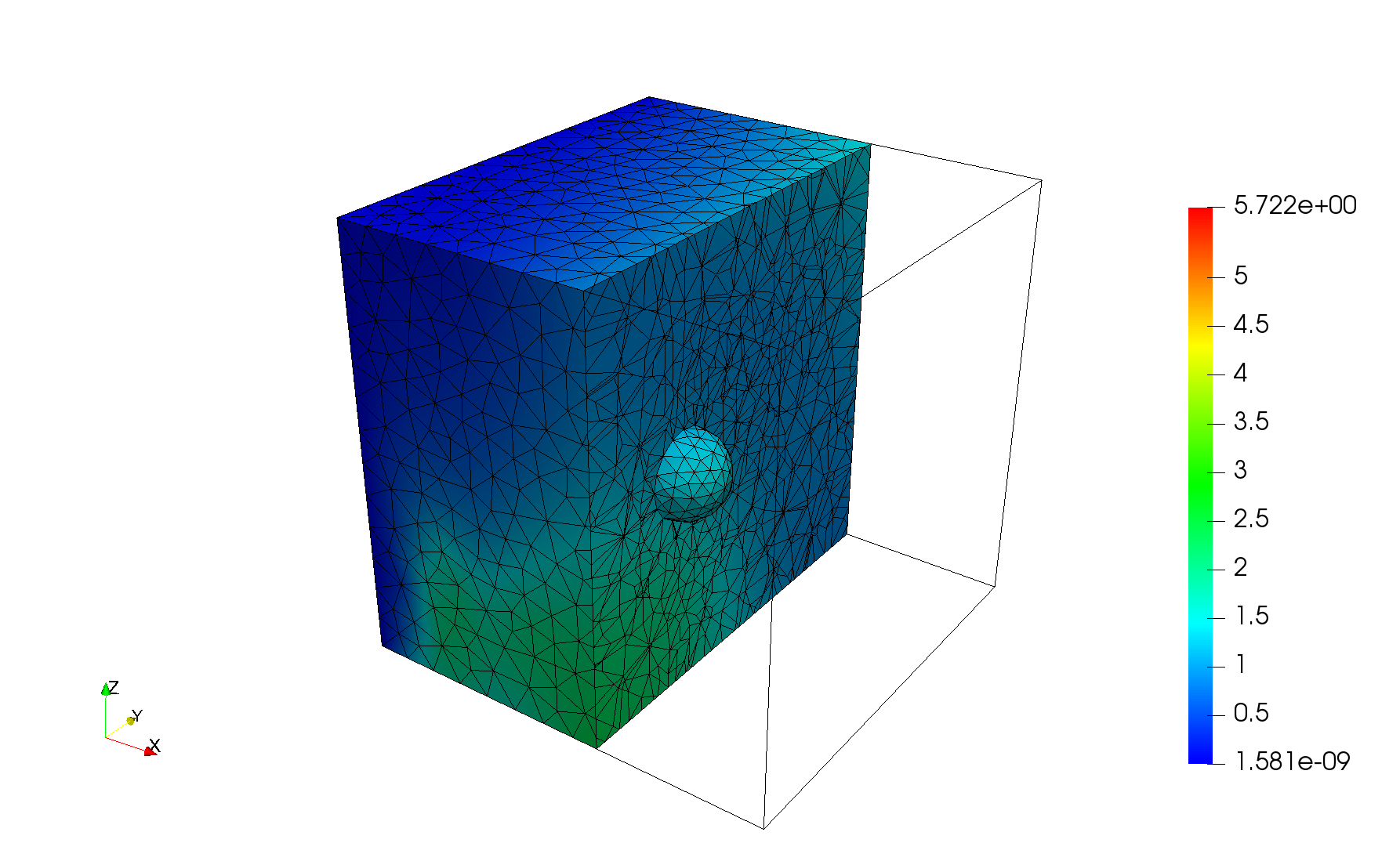}
	\end{subfigure}\hfill
	\begin{subfigure}[c]{0.22\textwidth}
		\centering
		\includegraphics[height=25mm]{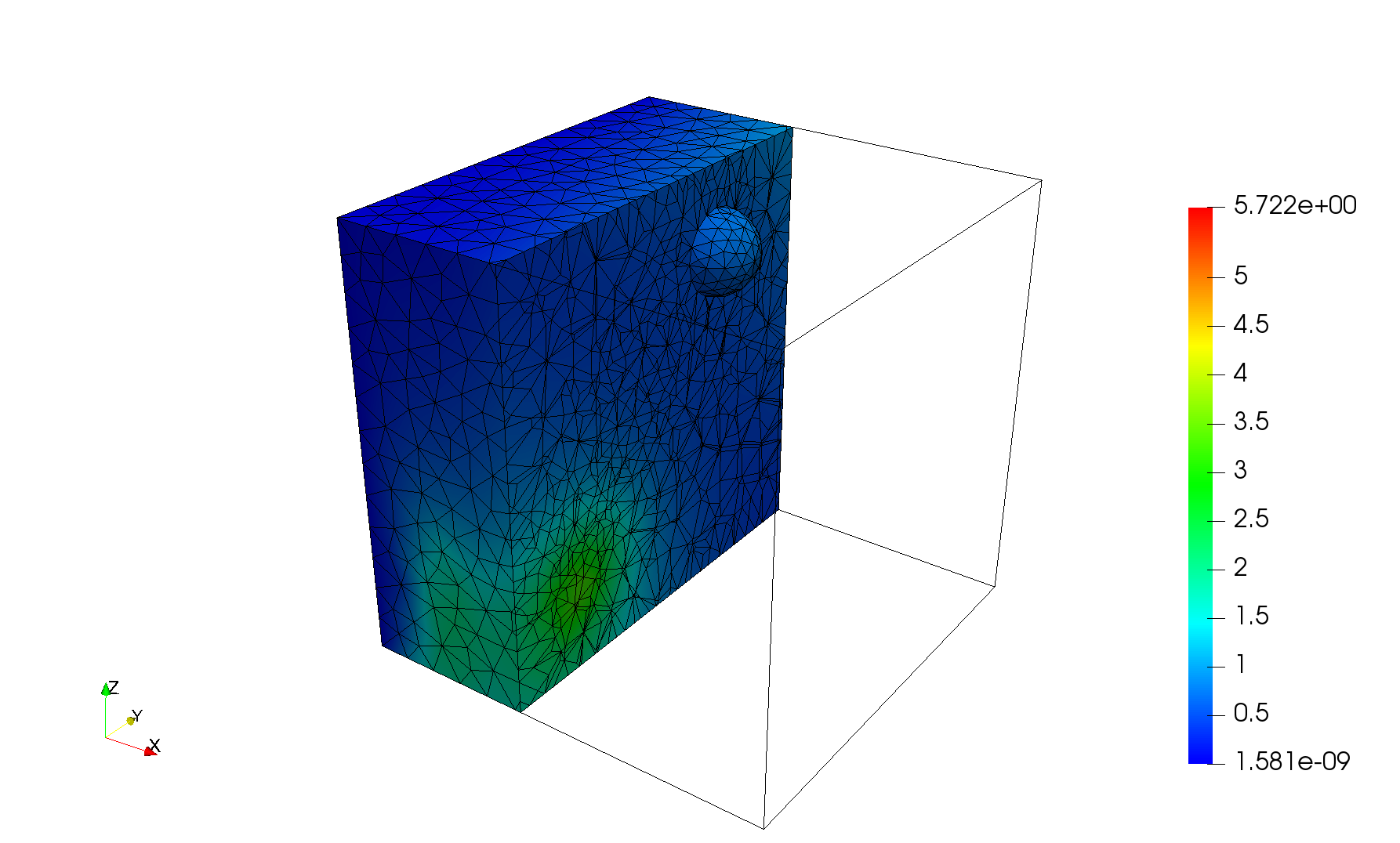}
	\end{subfigure}
	\\
	\rotatebox[origin=c]{90}{$3\textsuperscript{rd}$ mode}\hfill
	\begin{subfigure}[c]{0.22\textwidth}
		\centering
		\includegraphics[height=25mm]{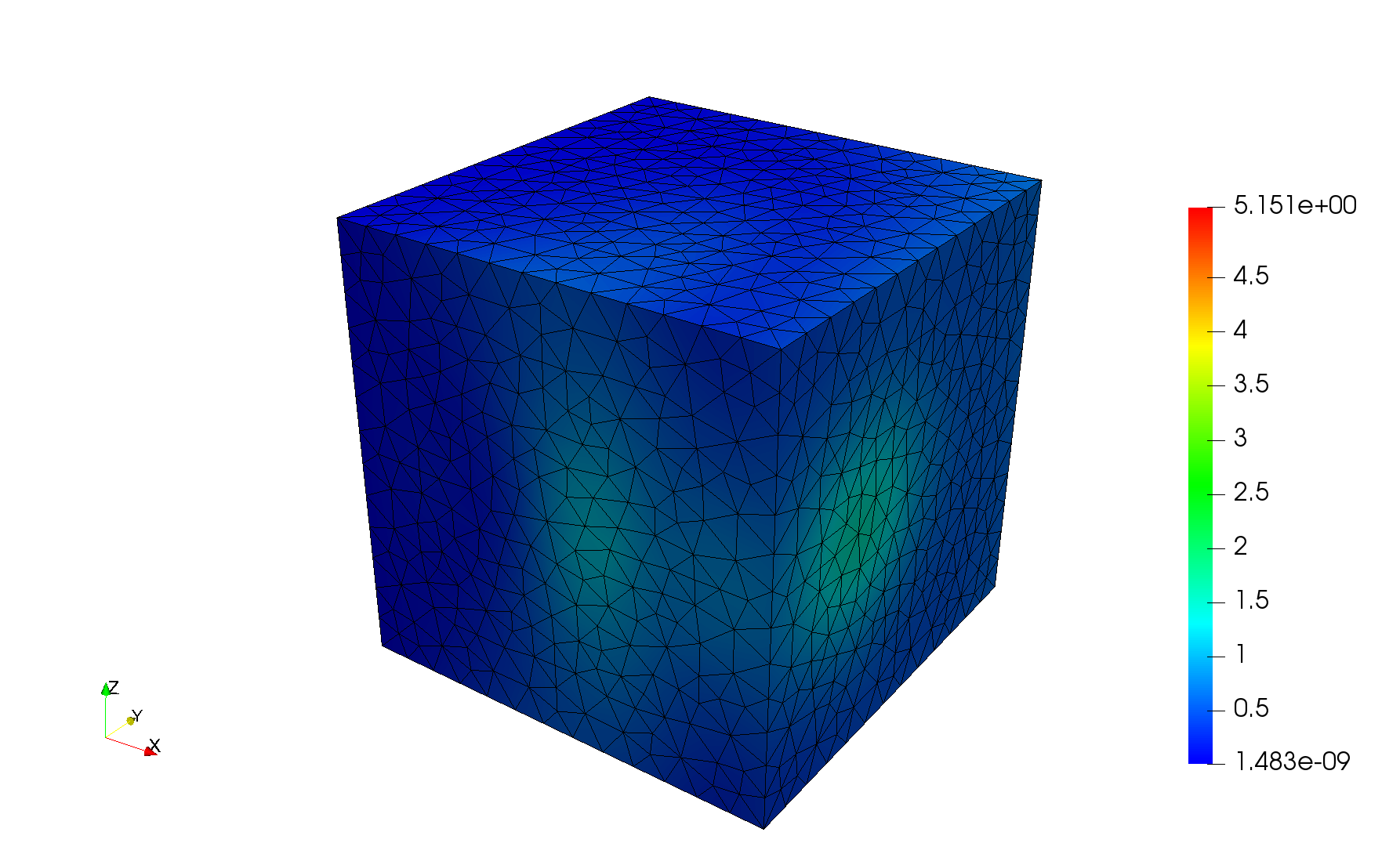}
	\end{subfigure}\hfill
	\begin{subfigure}[c]{0.22\textwidth}
		\centering
		\includegraphics[height=25mm]{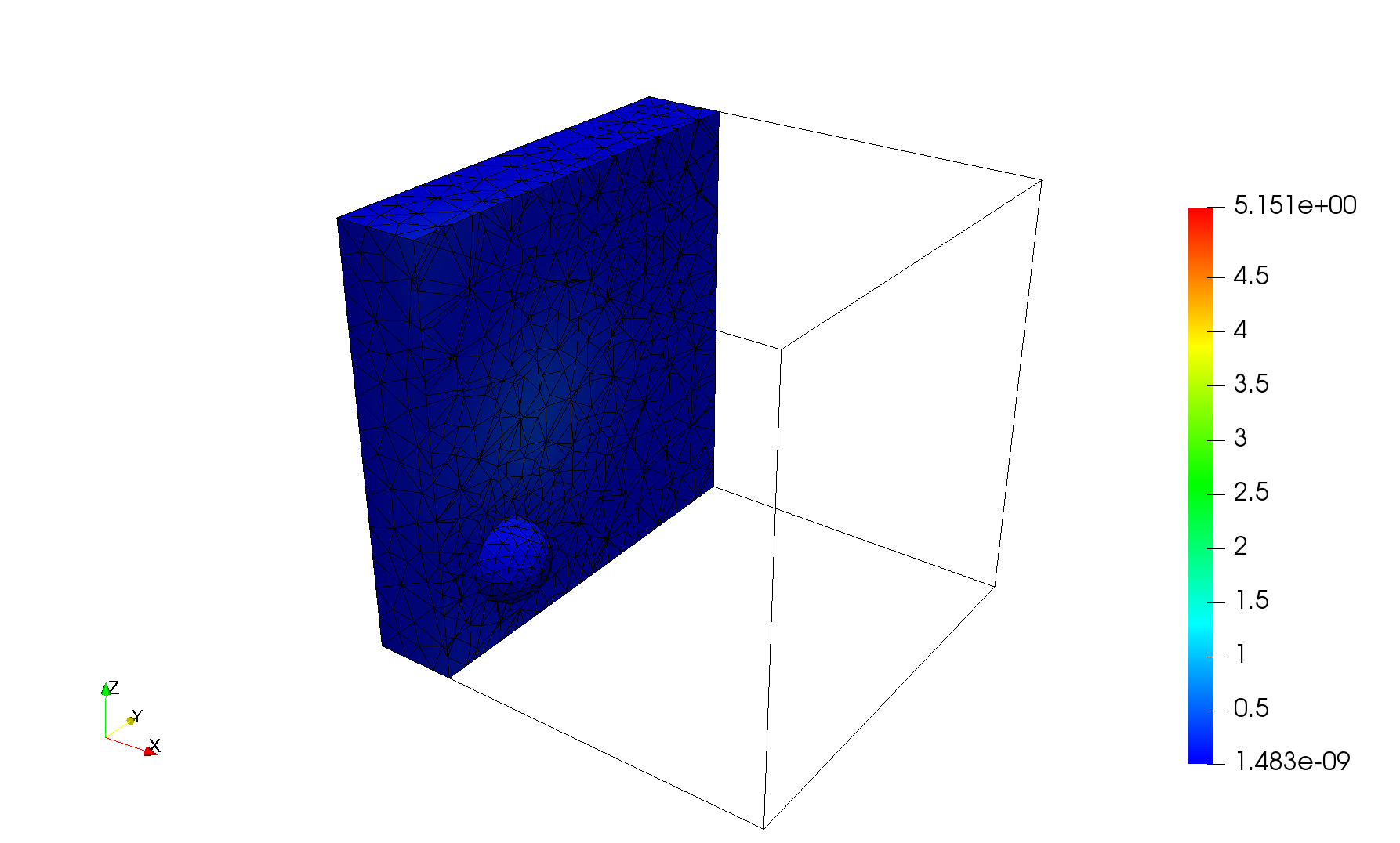}
	\end{subfigure}\hfill
	\begin{subfigure}[c]{0.22\textwidth}
		\centering
		\includegraphics[height=25mm]{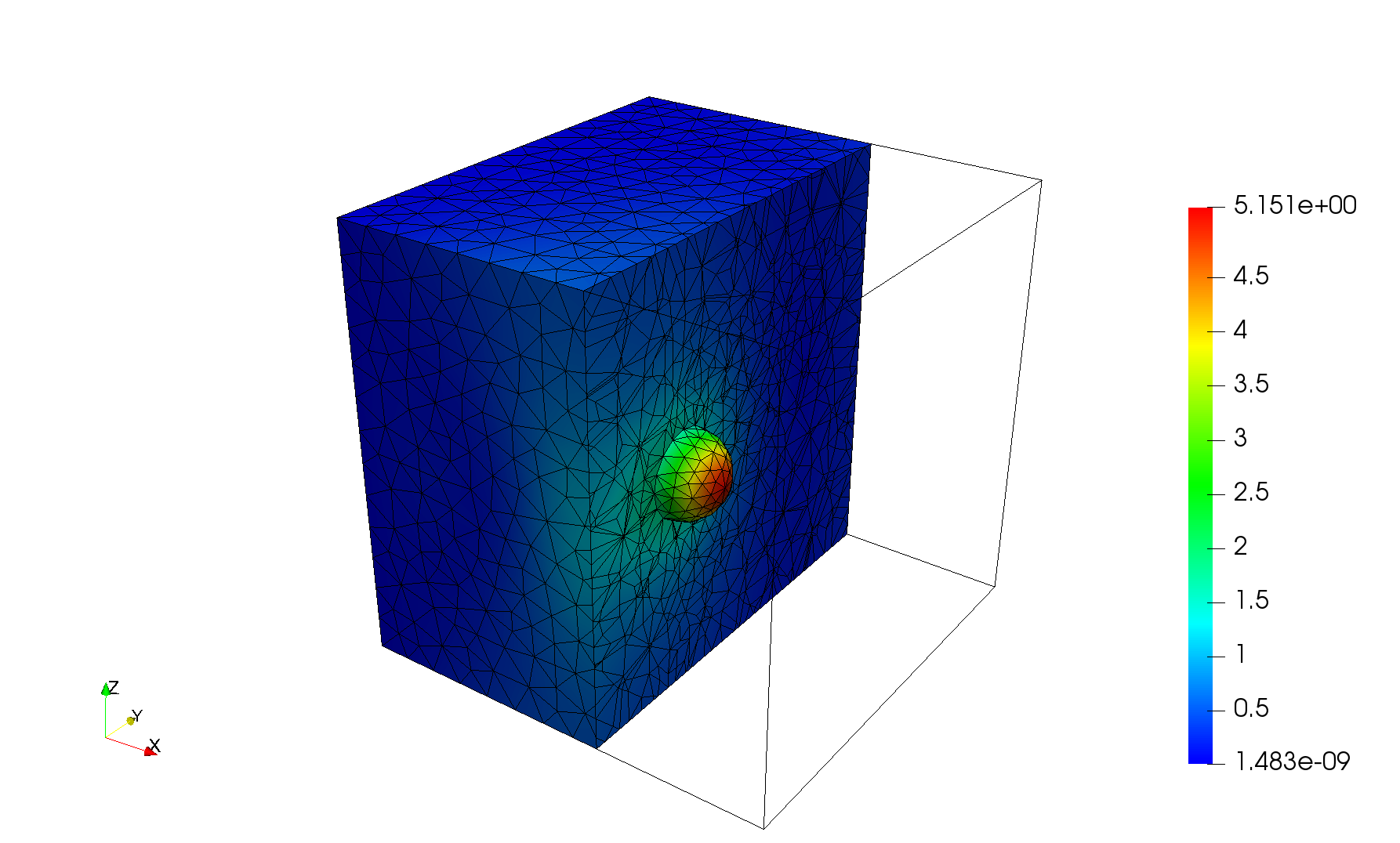}
	\end{subfigure}\hfill
	\begin{subfigure}[c]{0.22\textwidth}
		\centering
		\includegraphics[height=25mm]{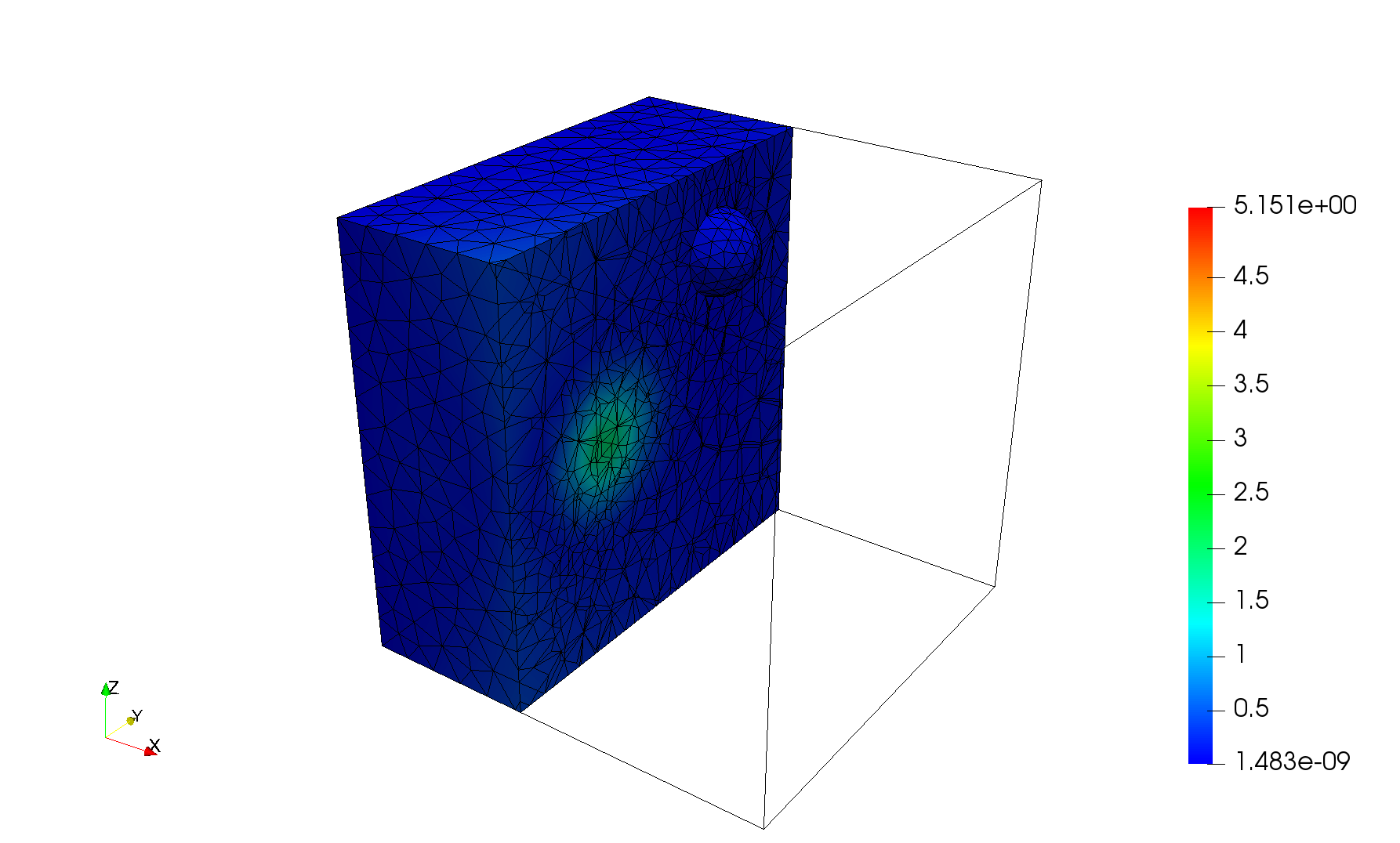}
	\end{subfigure}
	\\
	\rotatebox[origin=c]{90}{$4\textsuperscript{th}$ mode}\hfill
	\begin{subfigure}[c]{0.22\textwidth}
		\centering
		\includegraphics[height=25mm]{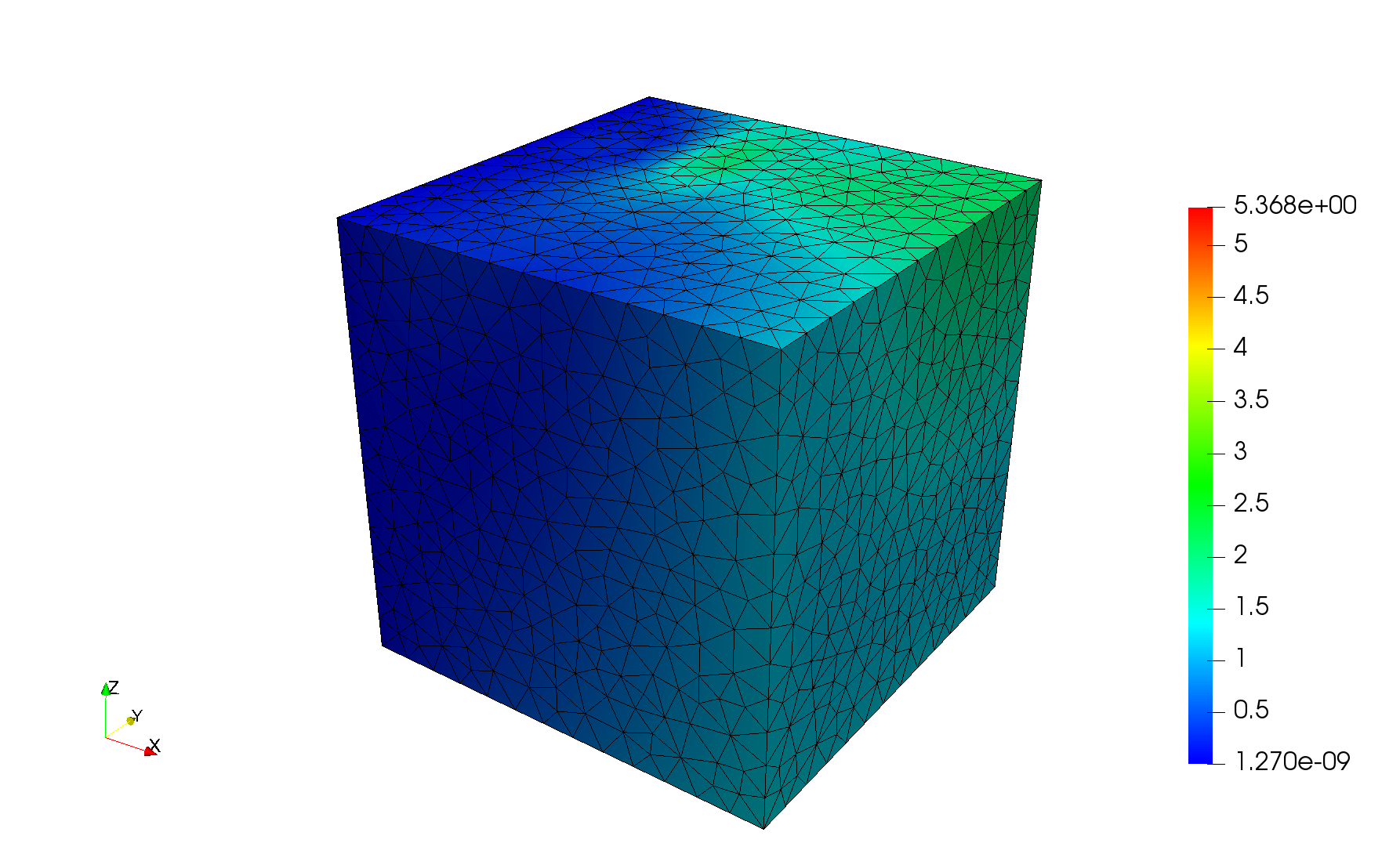}
	\end{subfigure}\hfill
	\begin{subfigure}[c]{0.22\textwidth}
		\centering
		\includegraphics[height=25mm]{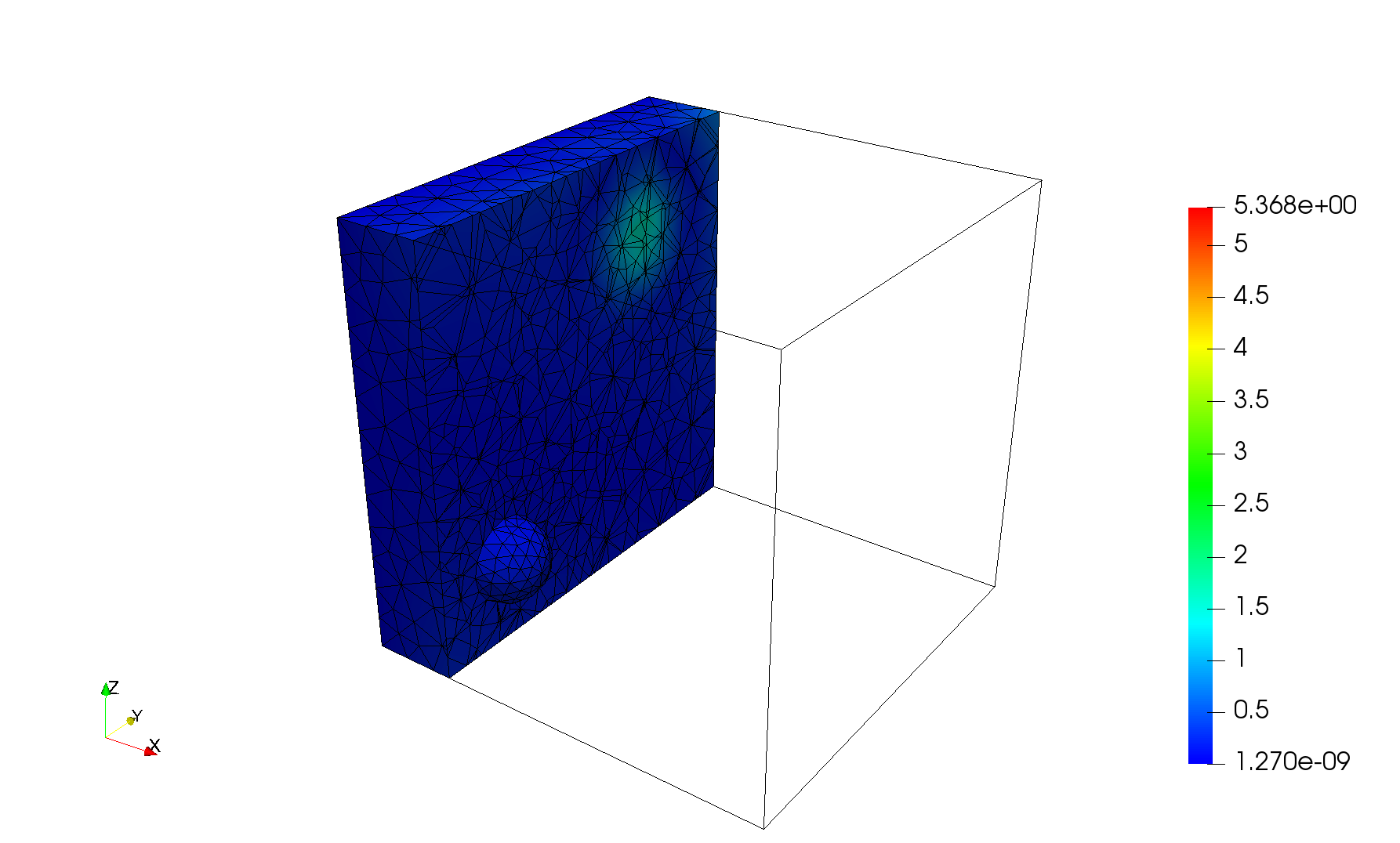}
	\end{subfigure}\hfill
	\begin{subfigure}[c]{0.22\textwidth}
		\centering
		\includegraphics[height=25mm]{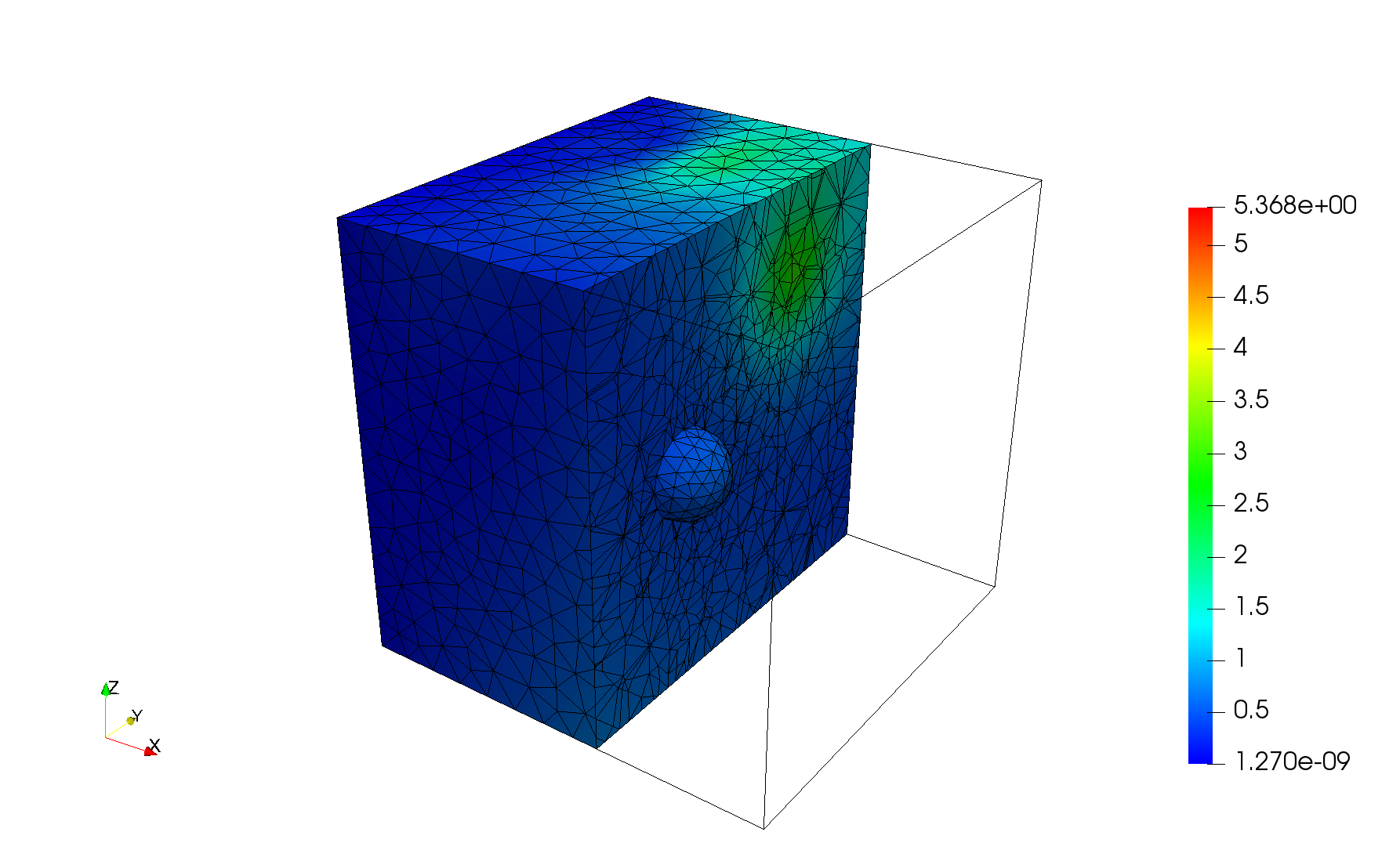}
	\end{subfigure}\hfill
	\begin{subfigure}[c]{0.22\textwidth}
		\centering
		\includegraphics[height=25mm]{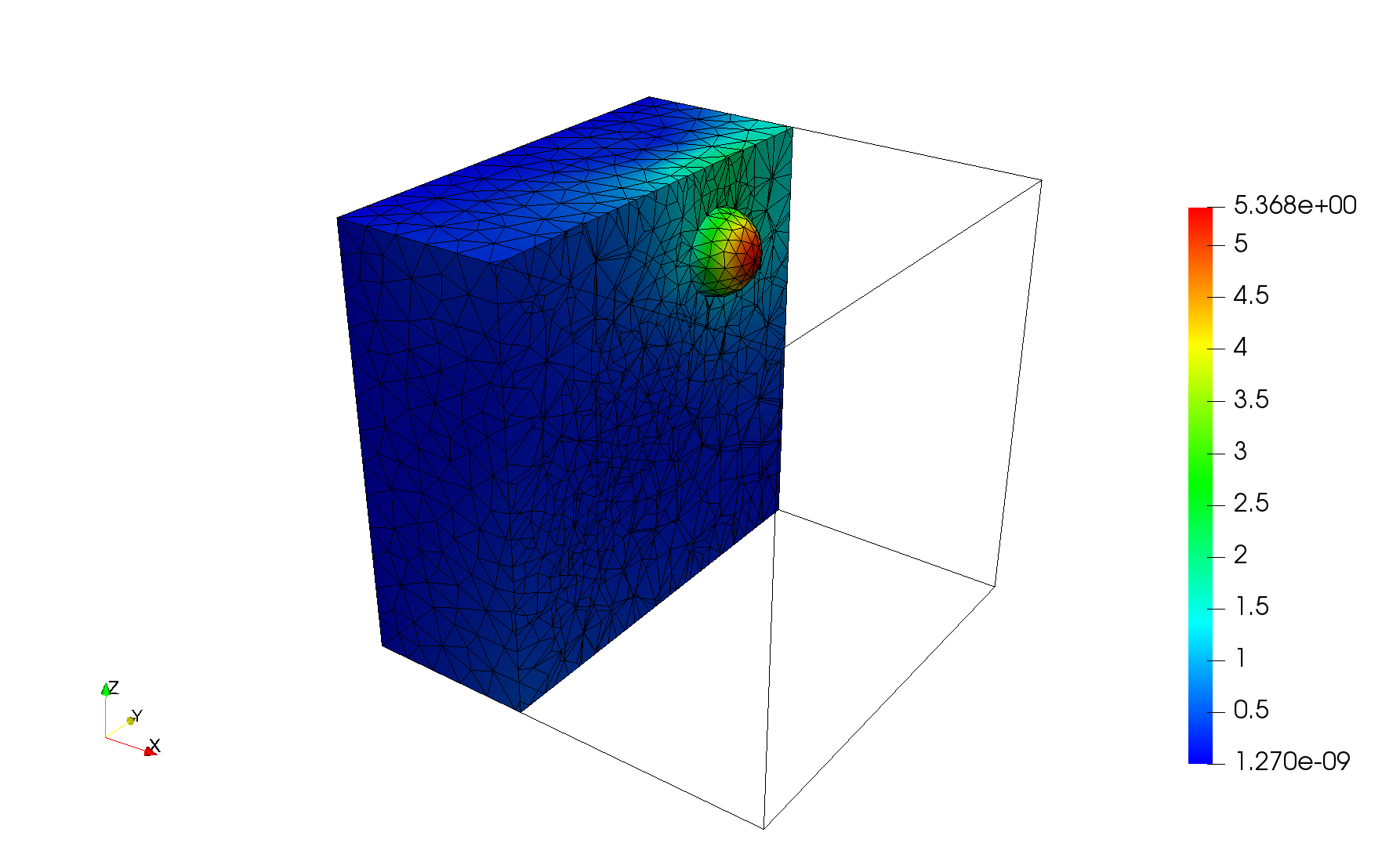}
	\end{subfigure}
	\\
	\rotatebox[origin=c]{90}{$5\textsuperscript{th}$ mode}\hfill
	\begin{subfigure}[c]{0.22\textwidth}
		\centering
		\includegraphics[height=25mm]{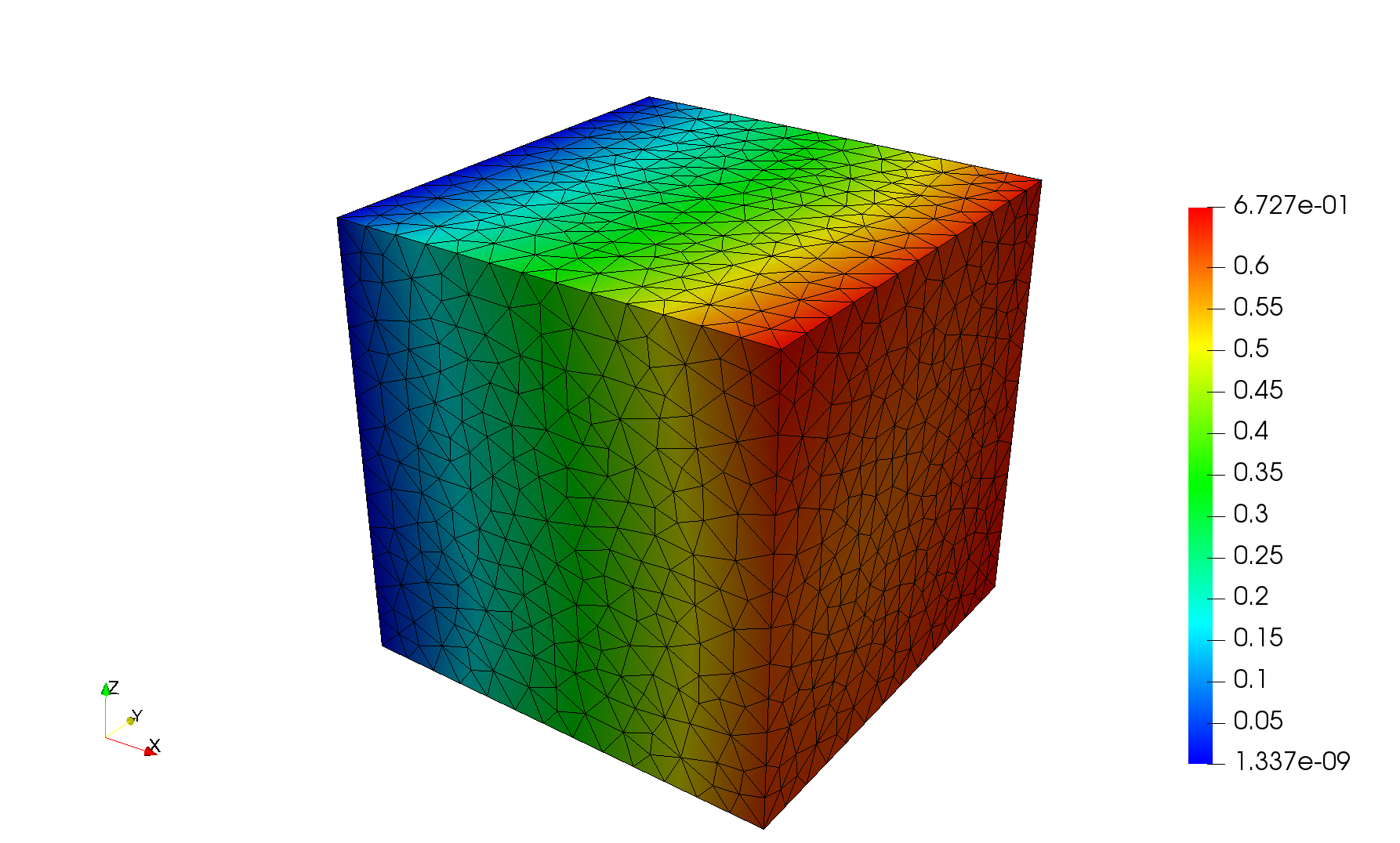}
	\end{subfigure}\hfill
	\begin{subfigure}[c]{0.22\textwidth}
		\centering
		\includegraphics[height=25mm]{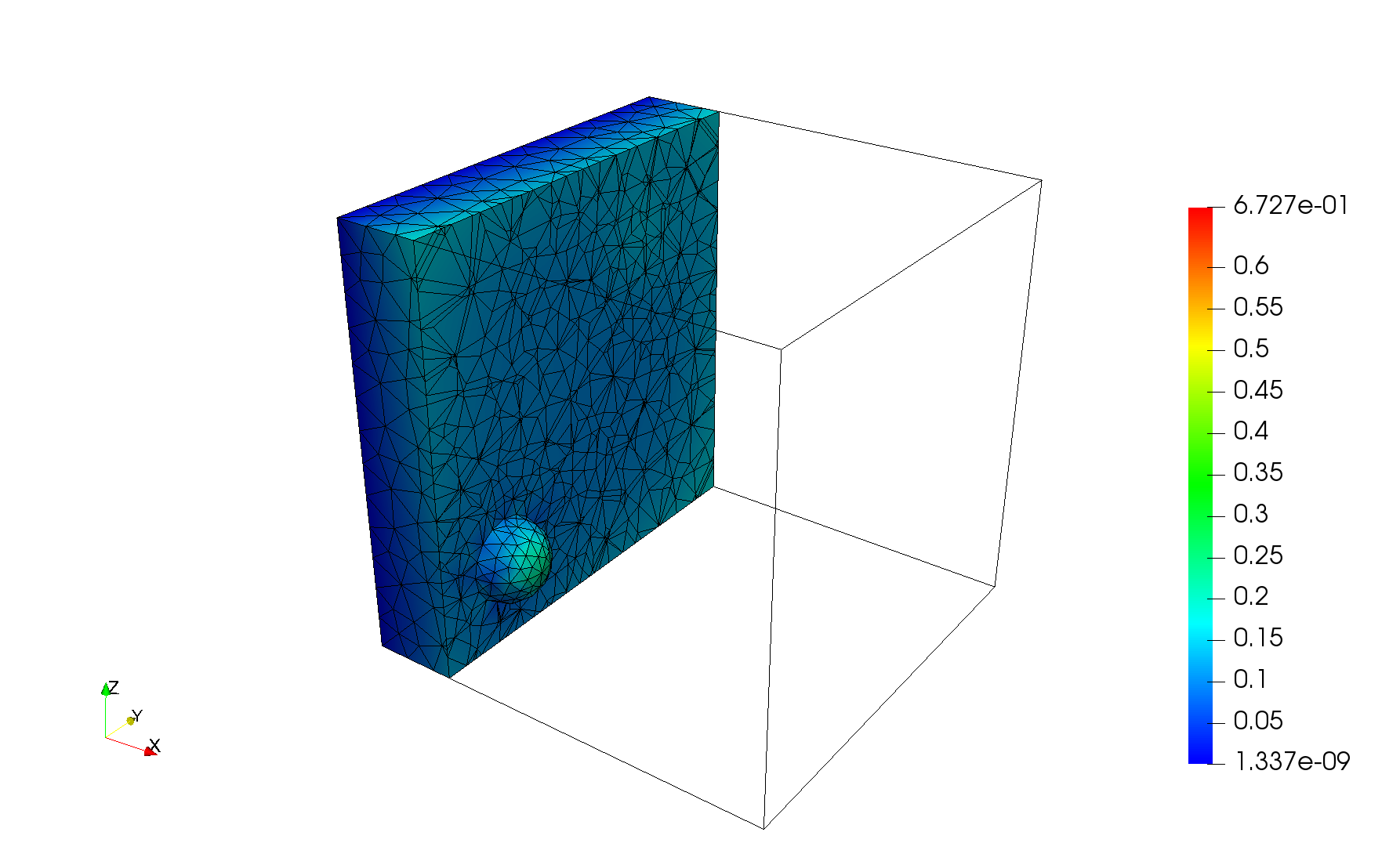}
	\end{subfigure}\hfill
	\begin{subfigure}[c]{0.22\textwidth}
		\centering
		\includegraphics[height=25mm]{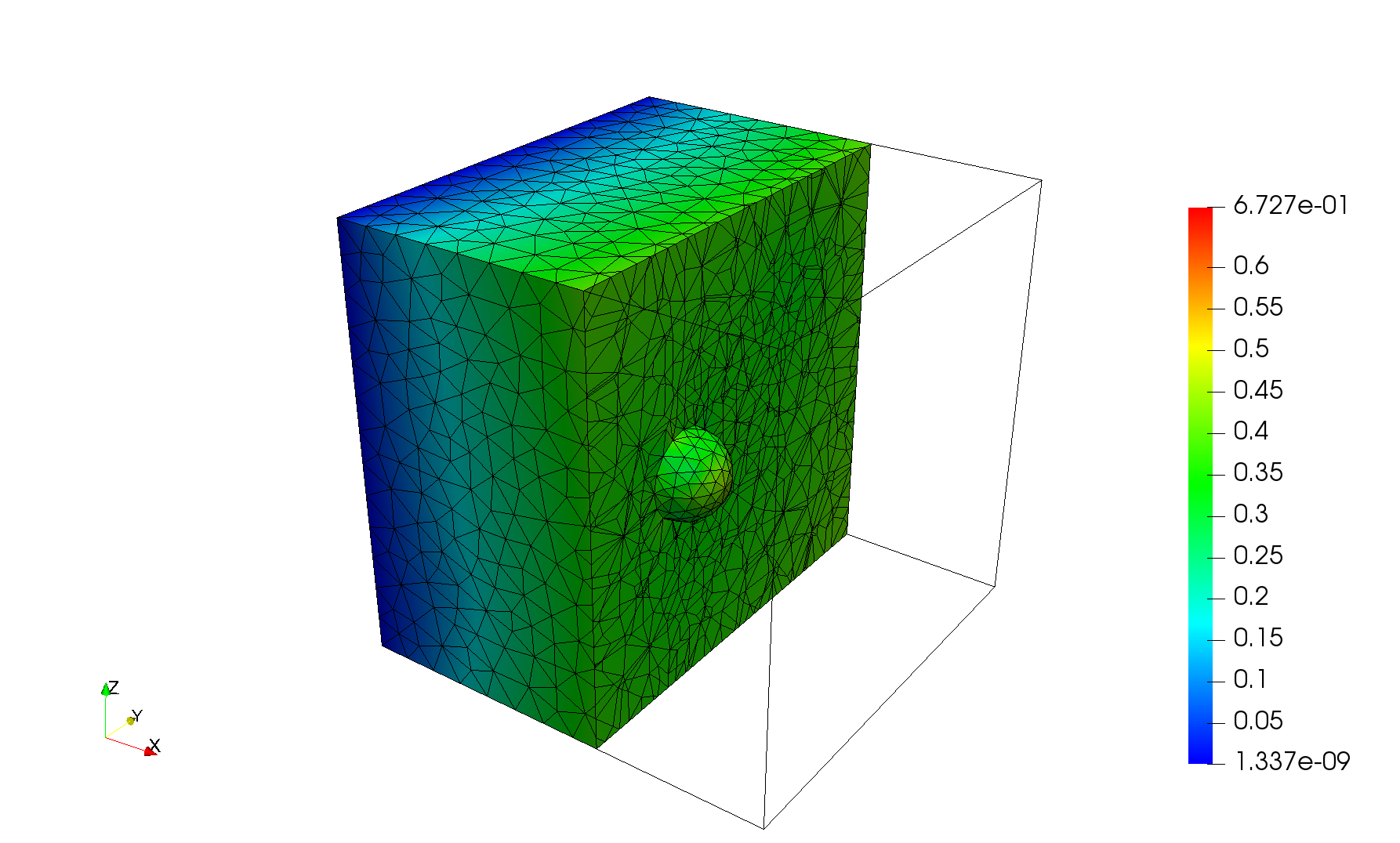}
	\end{subfigure}\hfill
	\begin{subfigure}[c]{0.22\textwidth}
		\centering
		\includegraphics[height=25mm]{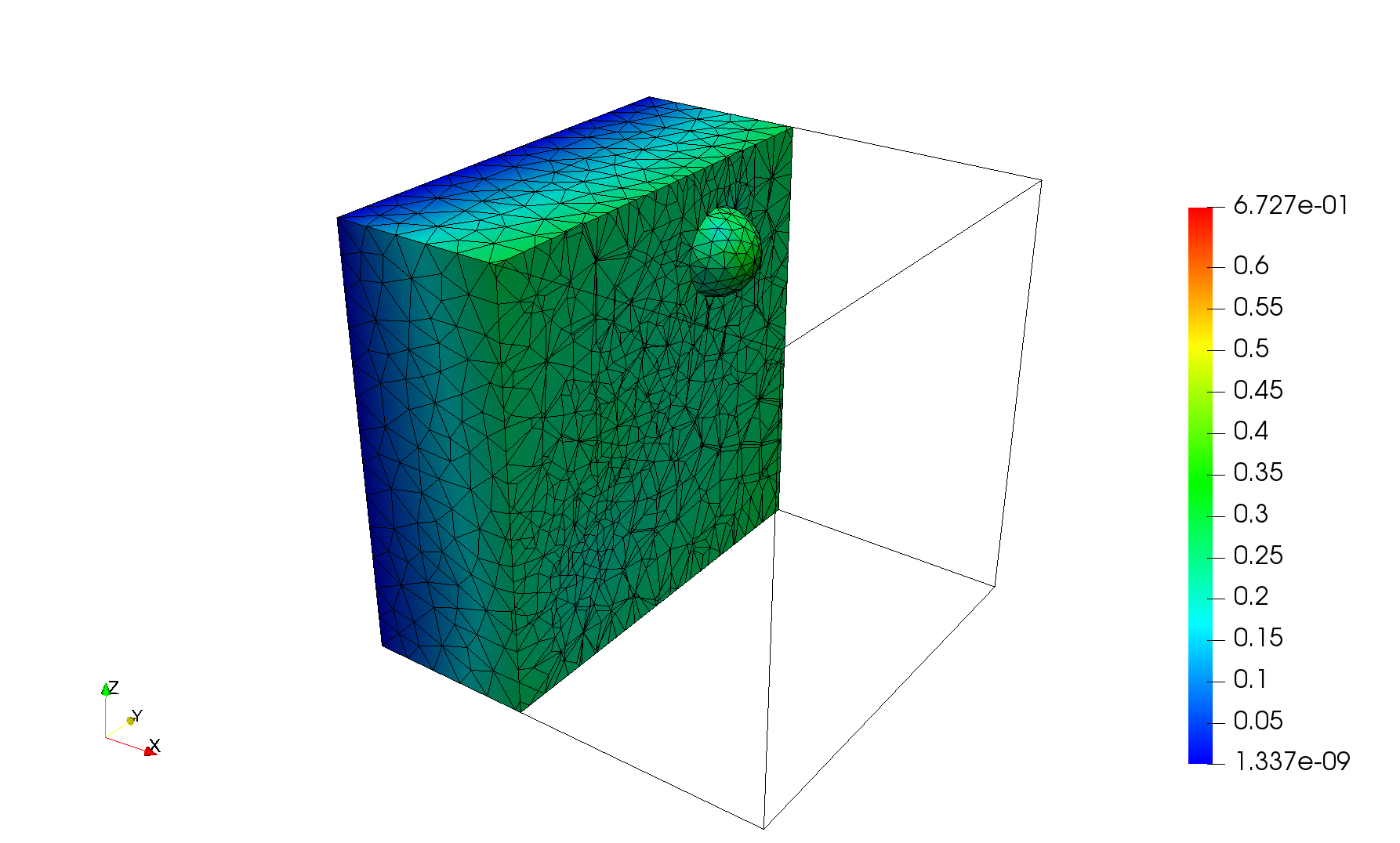}
	\end{subfigure}
	\caption{Magnitude of space functions $\psib_m(\xb)$ obtained for order $m=1,\dots,5$ (from top to bottom).}\label{fig:3D_modes_space}
\end{figure}

\begin{figure}
	\centering
	\rotatebox[origin=c]{90}{$1\textsuperscript{st}$ mode}\hfill
	\begin{subfigure}[c]{0.30\textwidth}
		\centering
		\input{3D_mode1_param1.tex}
	\end{subfigure}\hfill
	\begin{subfigure}[c]{0.30\textwidth}
		\centering
		\input{3D_mode1_param2.tex}
	\end{subfigure}\hfill
	\begin{subfigure}[c]{0.30\textwidth}
		\centering
		\input{3D_mode1_param3.tex}
	\end{subfigure}
	\\
	\rotatebox[origin=c]{90}{$2\textsuperscript{nd}$ mode}\hfill
	\begin{subfigure}[c]{0.30\textwidth}
		\centering
		\input{3D_mode2_param1.tex}
	\end{subfigure}\hfill
	\begin{subfigure}[c]{0.30\textwidth}
		\centering
		\input{3D_mode2_param2.tex}
	\end{subfigure}\hfill
	\begin{subfigure}[c]{0.30\textwidth}
		\centering
		\input{3D_mode2_param3.tex}
	\end{subfigure}
	\\
	\rotatebox[origin=c]{90}{$3\textsuperscript{rd}$ mode}\hfill
	\begin{subfigure}[c]{0.30\textwidth}
		\centering
		\input{3D_mode3_param1.tex}
	\end{subfigure}\hfill
	\begin{subfigure}[c]{0.30\textwidth}
		\centering
		\input{3D_mode3_param2.tex}
	\end{subfigure}\hfill
	\begin{subfigure}[c]{0.30\textwidth}
		\centering
		\input{3D_mode3_param3.tex}
	\end{subfigure}
	\\
	\rotatebox[origin=c]{90}{$4\textsuperscript{th}$ mode}\hfill
	\begin{subfigure}[c]{0.30\textwidth}
		\centering
		\input{3D_mode4_param1.tex}
	\end{subfigure}\hfill
	\begin{subfigure}[c]{0.30\textwidth}
		\centering
		\input{3D_mode4_param2.tex}
	\end{subfigure}\hfill
	\begin{subfigure}[c]{0.30\textwidth}
		\centering
		\input{3D_mode4_param3.tex}
	\end{subfigure}
	\\
	\rotatebox[origin=c]{90}{$5\textsuperscript{th}$ mode}\hfill
	\begin{subfigure}[c]{0.30\textwidth}
		\centering
		\input{3D_mode5_param1.tex}
	\end{subfigure}\hfill
	\begin{subfigure}[c]{0.30\textwidth}
		\centering
		\input{3D_mode5_param2.tex}
	\end{subfigure}\hfill
	\begin{subfigure}[c]{0.30\textwidth}
		\centering
		\input{3D_mode5_param3.tex}
	\end{subfigure}
	\caption{Parameter functions $\gamma_{1,m}(E_1)$, $\gamma_{2,m}(E_2)$ and $\gamma_{3,m}(E_3)$ (from left to right) obtained for order $m=1,\dots,5$ (from top to bottom).}\label{fig:3D_modes_parameters}
\end{figure}
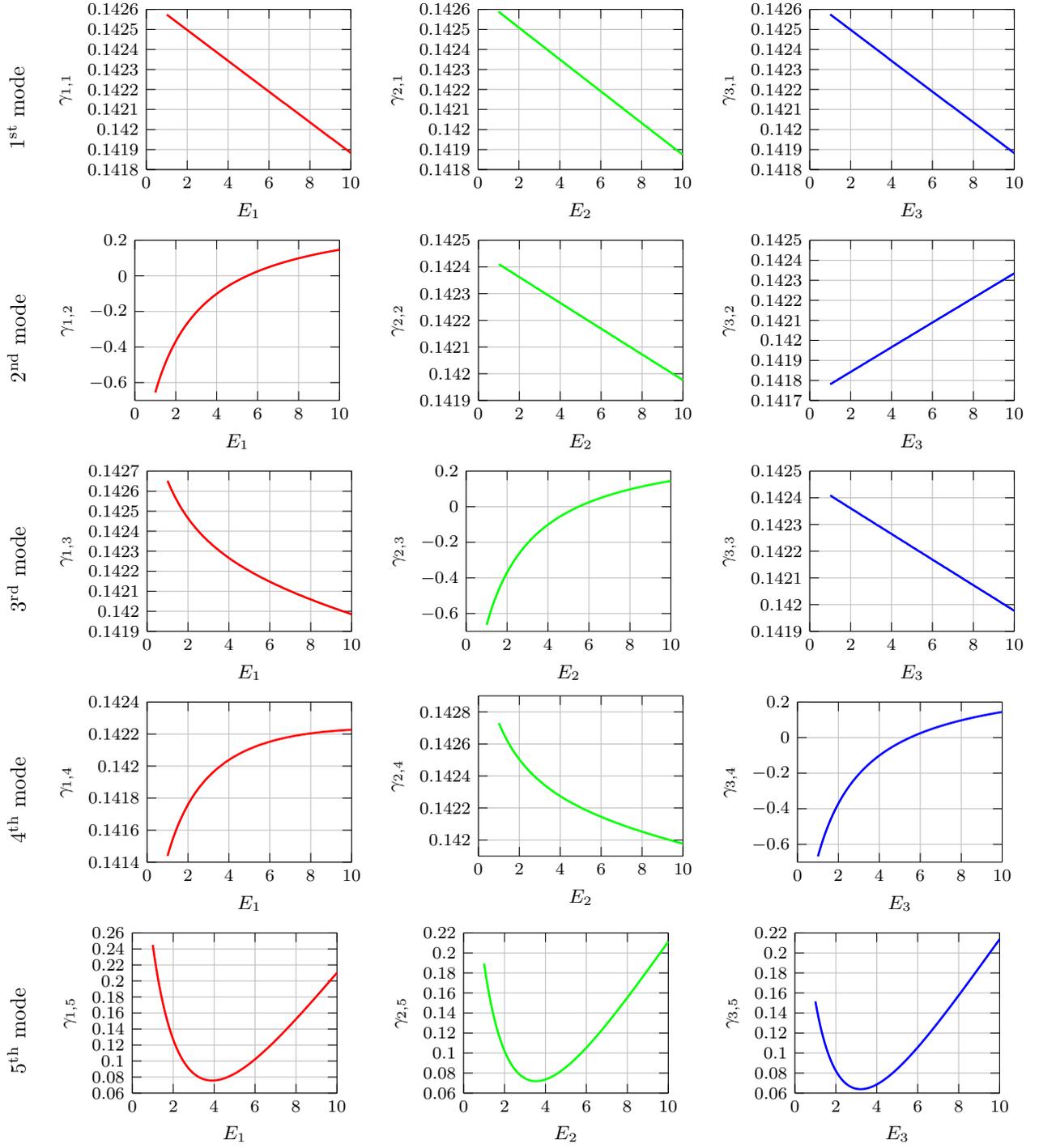

The initial FE mesh contains $17\,731$ $4$-nodes tetrahedral elements and $3\,622$ nodes ($10\,866$ dofs). The first five PGD modes of the PGD approximate solution $\ub^h_m(\xb,E_1,E_2,E_3)$ are given in Figure~\ref{fig:3D_modes_space} for space functions $\boldsymbol{\psi}_m(\xb)$ and in Figure~\ref{fig:3D_modes_parameters} for parameter functions $\gamma_{1,m}(E_1)$, $\gamma_{2,m}(E_2)$ and $\gamma_{3,m}(E_3)$. Note that the first and fifth space functions $\boldsymbol{\psi}_1(\xb)$ and $\boldsymbol{\psi}_5(\xb)$ are global modes, whereas the second, third and fourth space functions $\boldsymbol{\psi}_2(\xb)$, $\boldsymbol{\psi}_3(\xb)$ and $\boldsymbol{\psi}_4(\xb)$ are local modes mostly concentrated around the first, second and third inclusions, respectively.

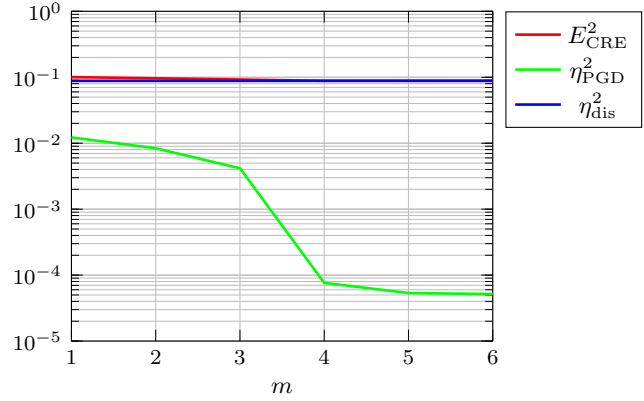
\begin{figure}[h!]
\centering
\input{3D_estimates_global.tex}
\caption{Evolutions of the error estimate $E_{\CRE}^2$ and associated error indicators $\eta^2_{\PGD}$ and $\eta_{\dis}^2$ with respect to the number $m$ of PGD modes.}\label{fig:3D_estimates_global}
\end{figure}

The evolutions of the error estimate $E_{\CRE}$ and associated error indicators $\eta_{\PGD}$ and $\eta_{\dis}$ with respect to the number $m$ of PGD modes are shown in Figure~\ref{fig:3D_estimates_global}, for $m=1,\dots,6$ and for the maximal values obtained with triplets $E_1,E_2,E_3 \in P_E$. The error estimate converges quite fast toward the indicator associated to the discretization error, while the one related to the PGD truncation error decreases toward zero.

\setlength\figureheight{0.20\textheight}

\begin{figure}[h!]
	\centering
	\rotatebox[origin=c]{90}{$m=1$}\hfill
	\begin{subfigure}[c]{0.22\textwidth}
		\centering
		\includegraphics[height=25mm]{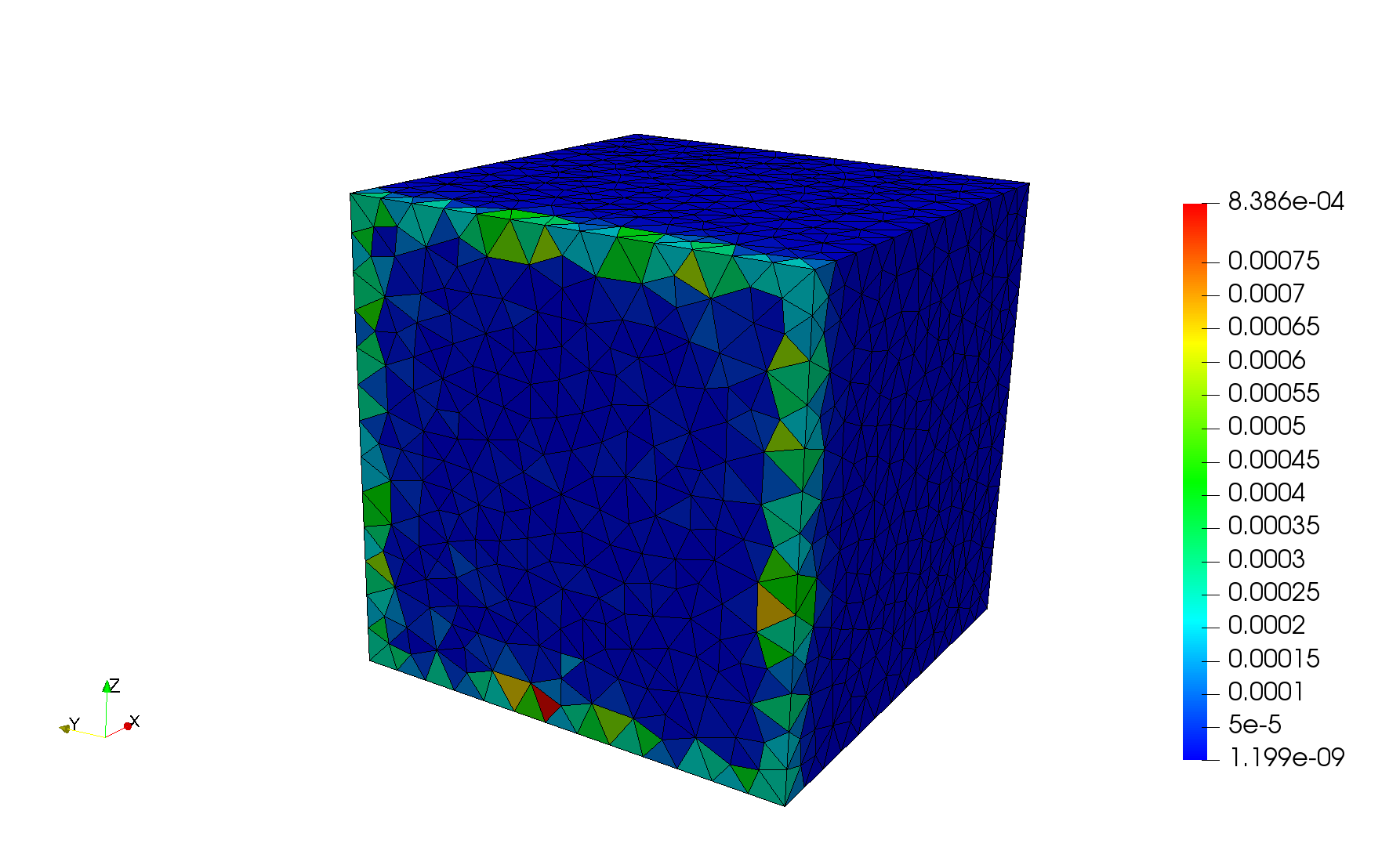}
	\end{subfigure}\hfill
	\begin{subfigure}[c]{0.22\textwidth}
		\centering
		\includegraphics[height=25mm]{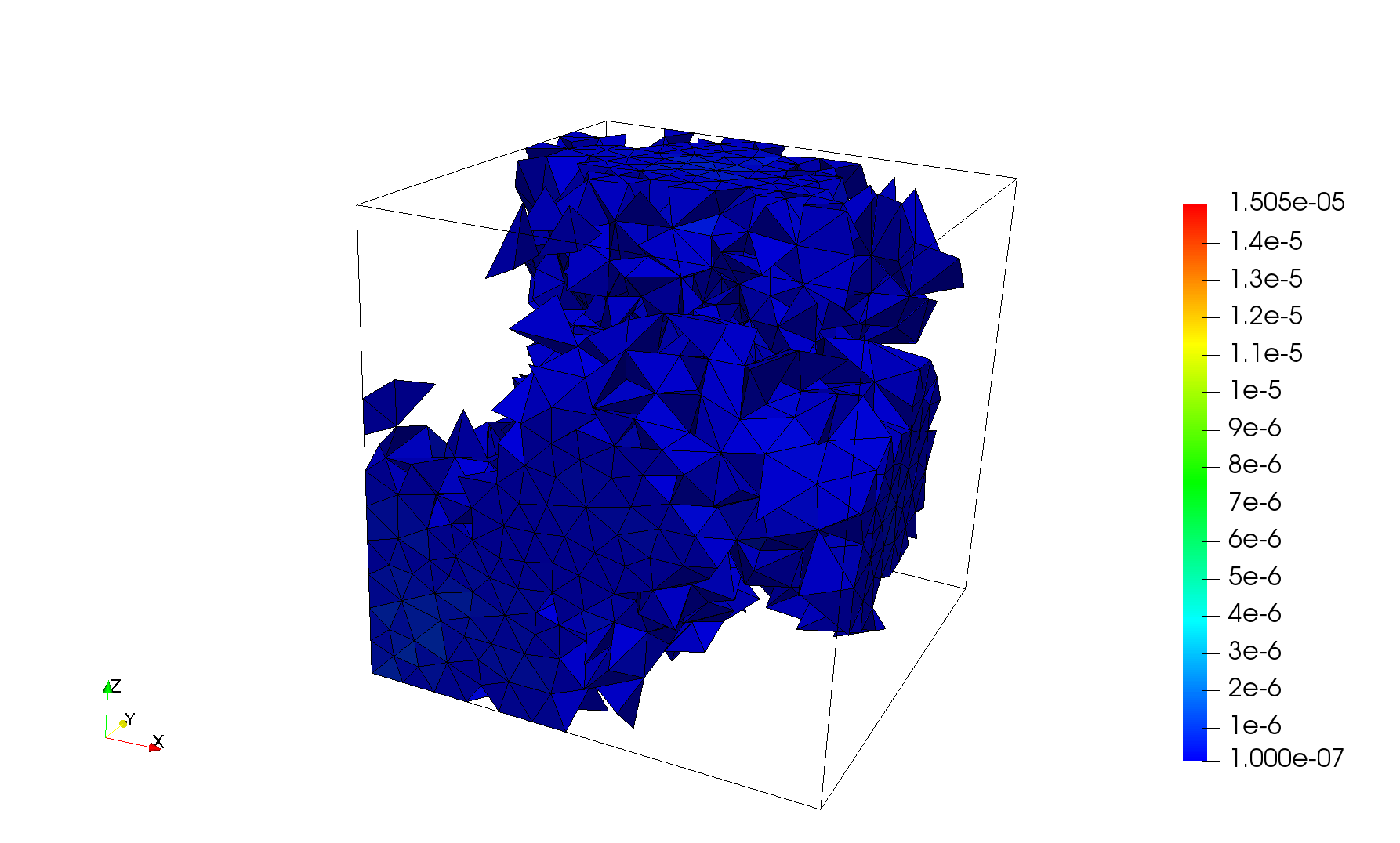}
	\end{subfigure}\hfill
	\begin{subfigure}[c]{0.22\textwidth}
		\centering
		\includegraphics[height=25mm]{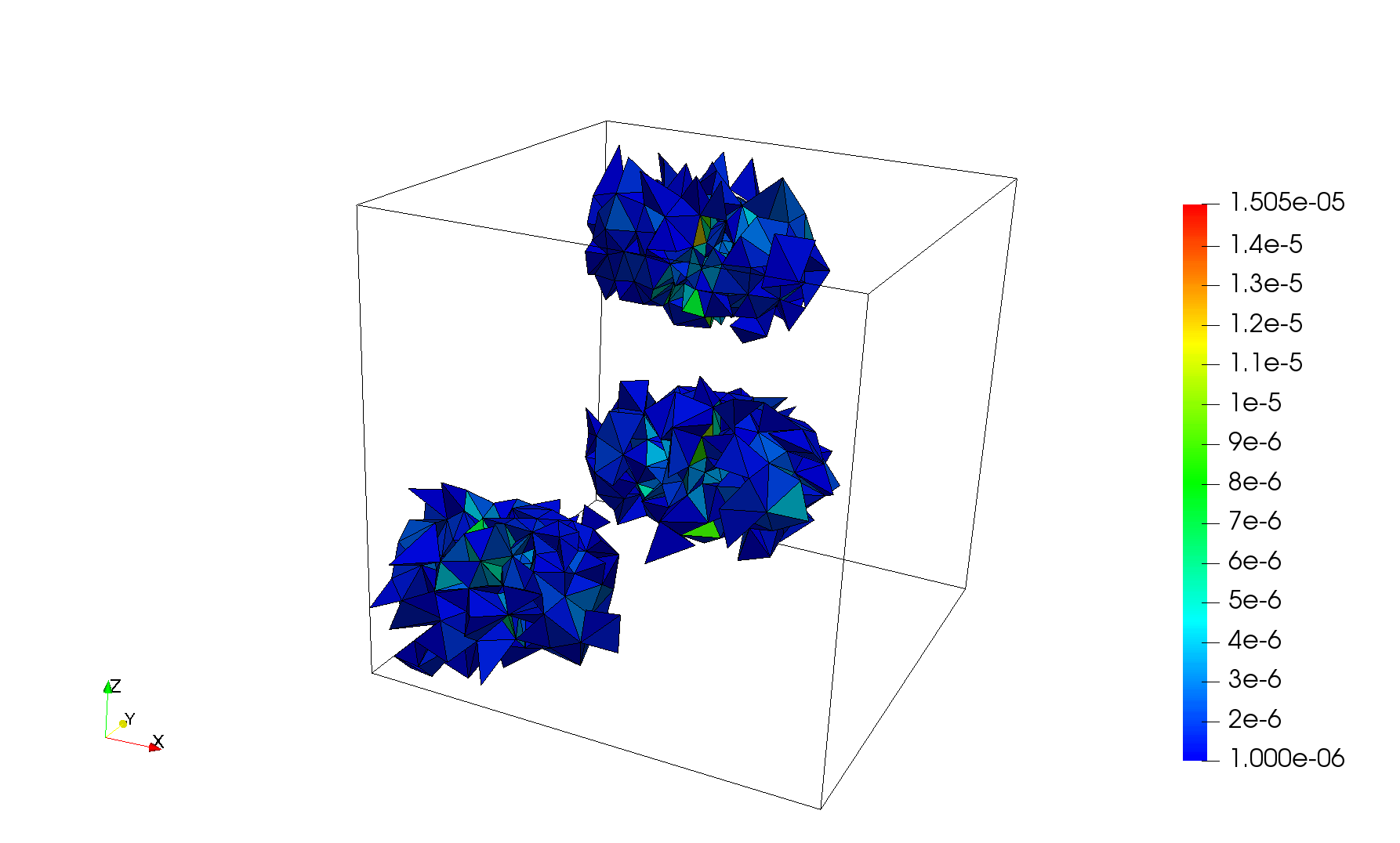}
	\end{subfigure}\hfill
	\begin{subfigure}[c]{0.22\textwidth}
		\centering
		\includegraphics[height=25mm]{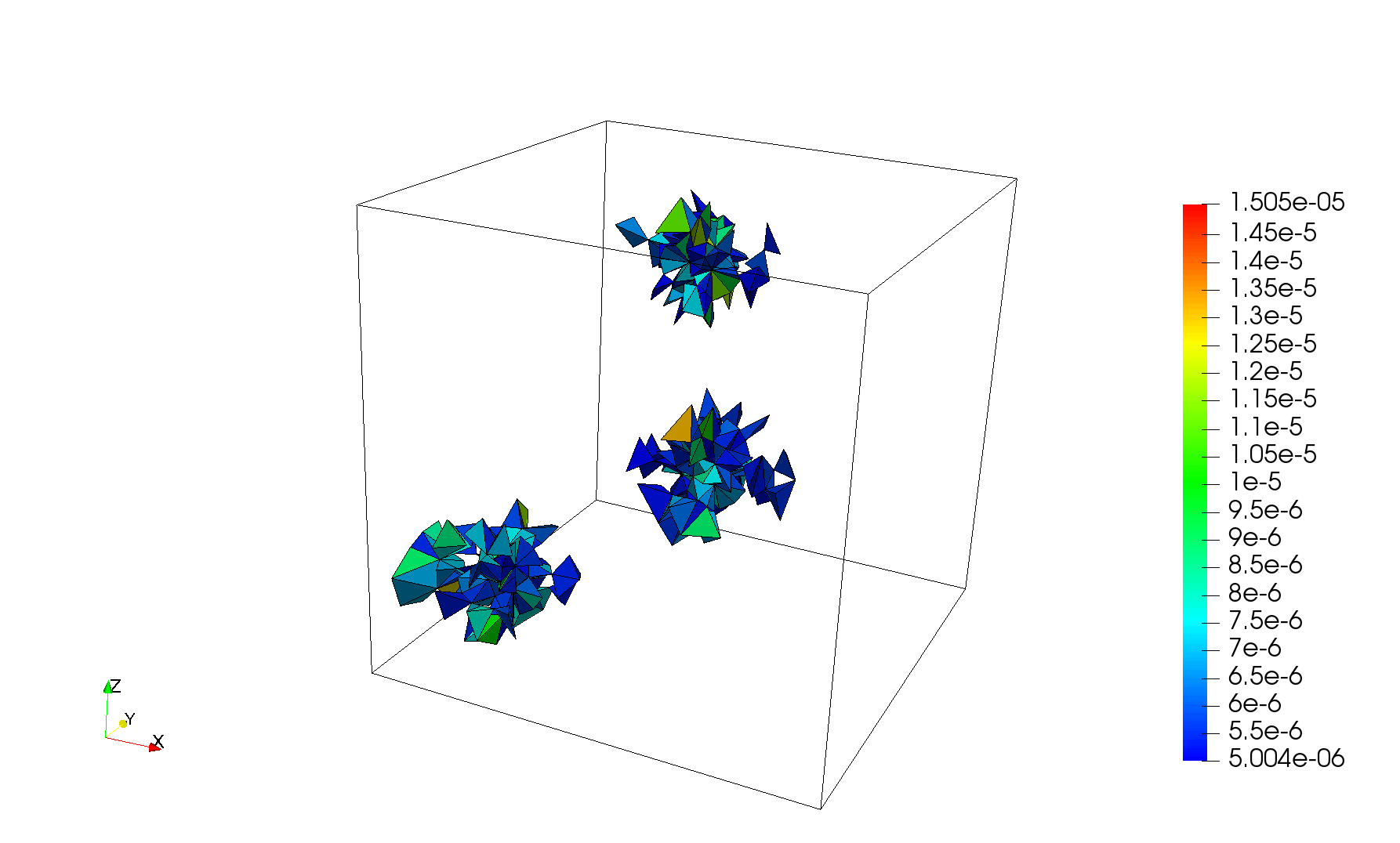}
	\end{subfigure}
	\\
	\rotatebox[origin=c]{90}{$m=2$}\hfill
	\begin{subfigure}[c]{0.22\textwidth}
		\centering
		\includegraphics[height=25mm]{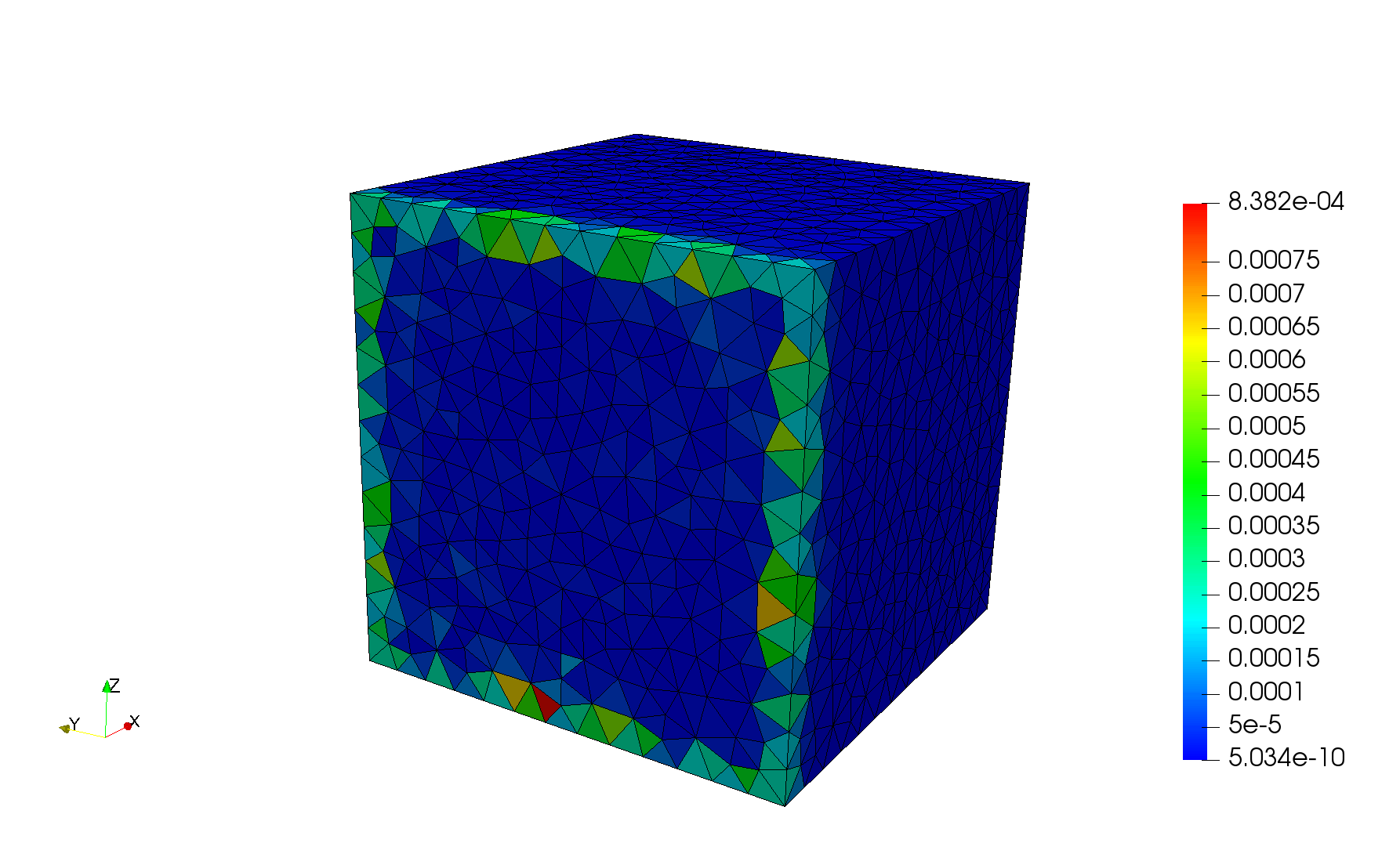}
	\end{subfigure}\hfill
	\begin{subfigure}[c]{0.22\textwidth}
		\centering
		\includegraphics[height=25mm]{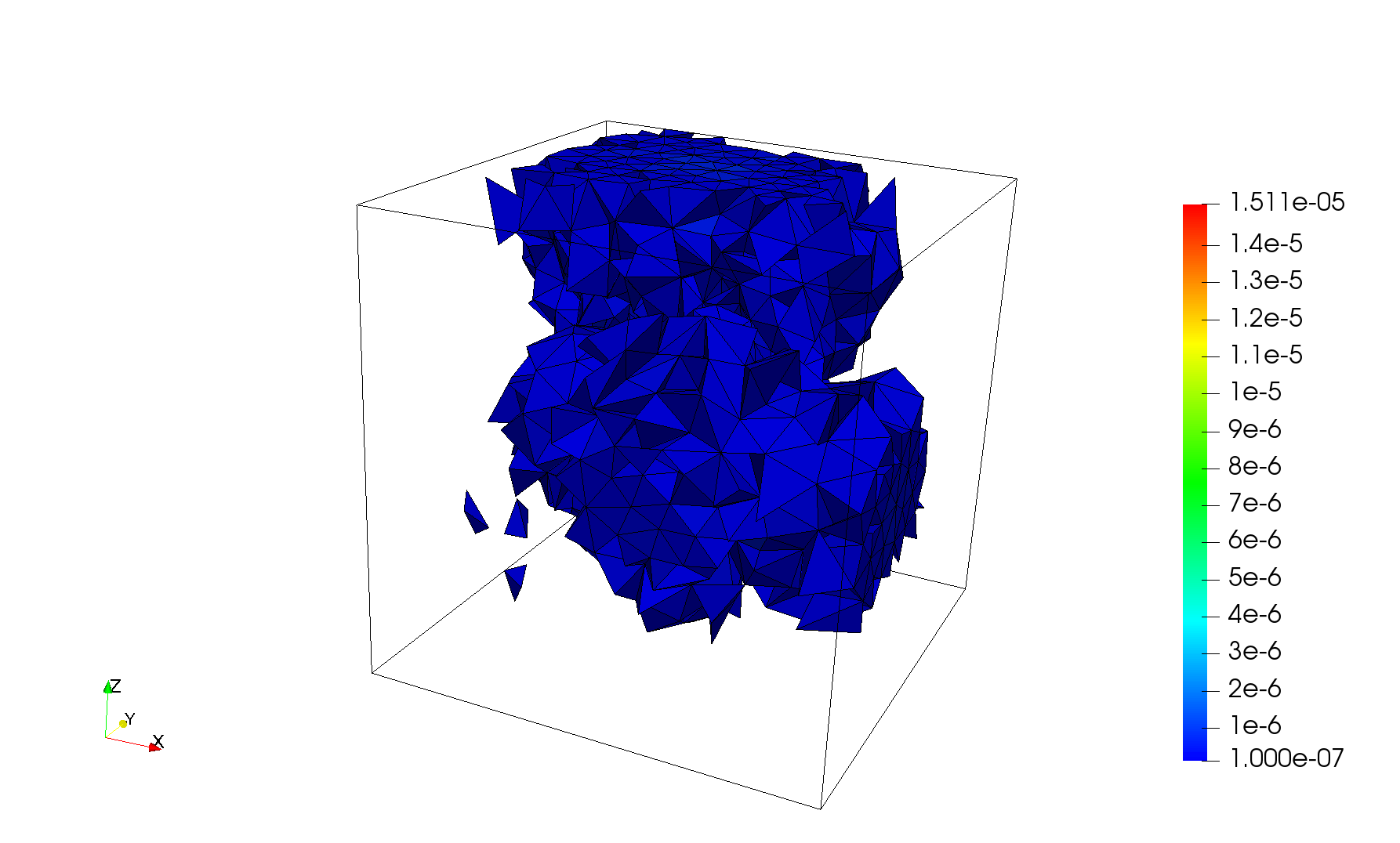}
	\end{subfigure}\hfill
	\begin{subfigure}[c]{0.22\textwidth}
		\centering
		\includegraphics[height=25mm]{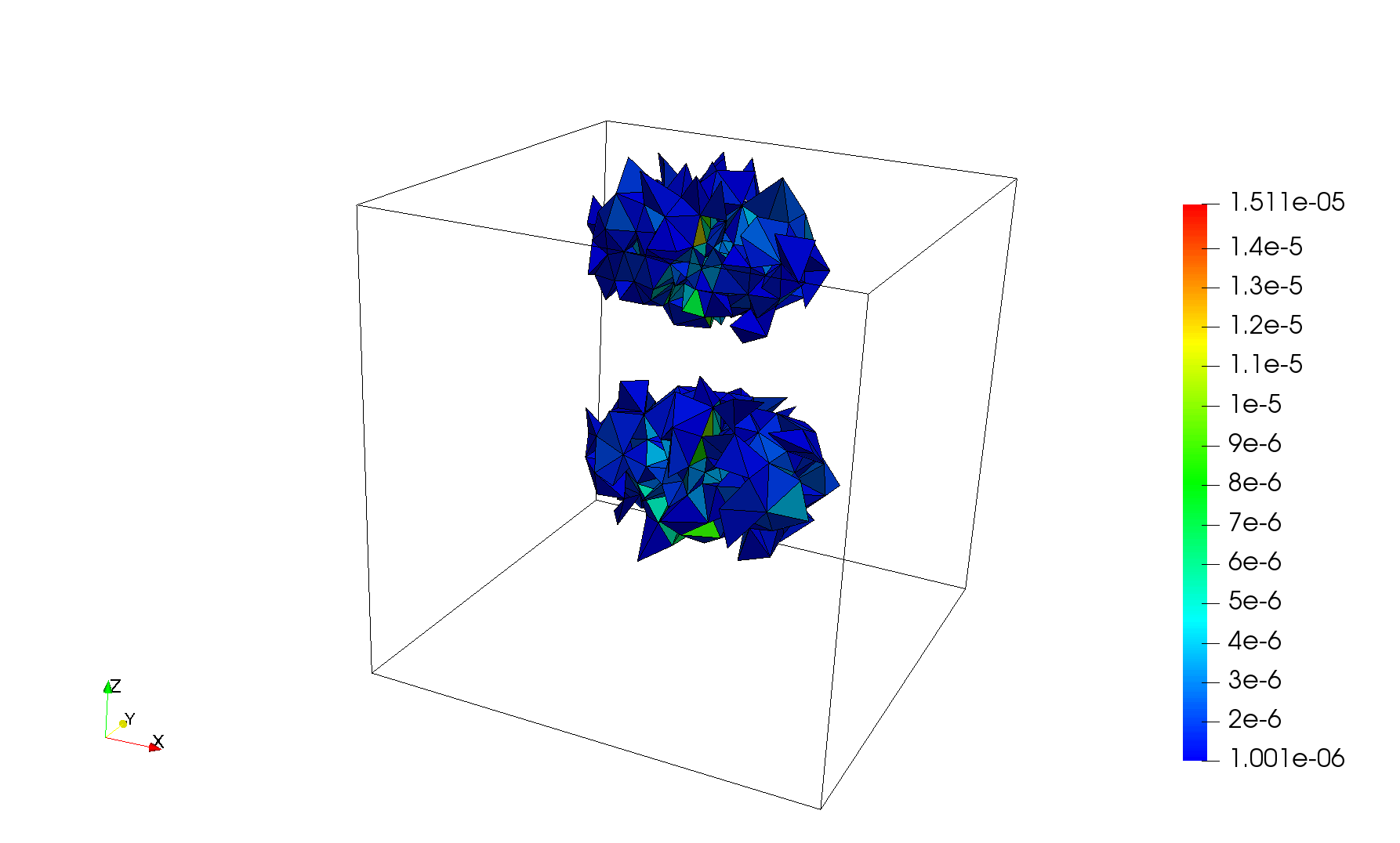}
	\end{subfigure}\hfill
	\begin{subfigure}[c]{0.22\textwidth}
		\centering
		\includegraphics[height=25mm]{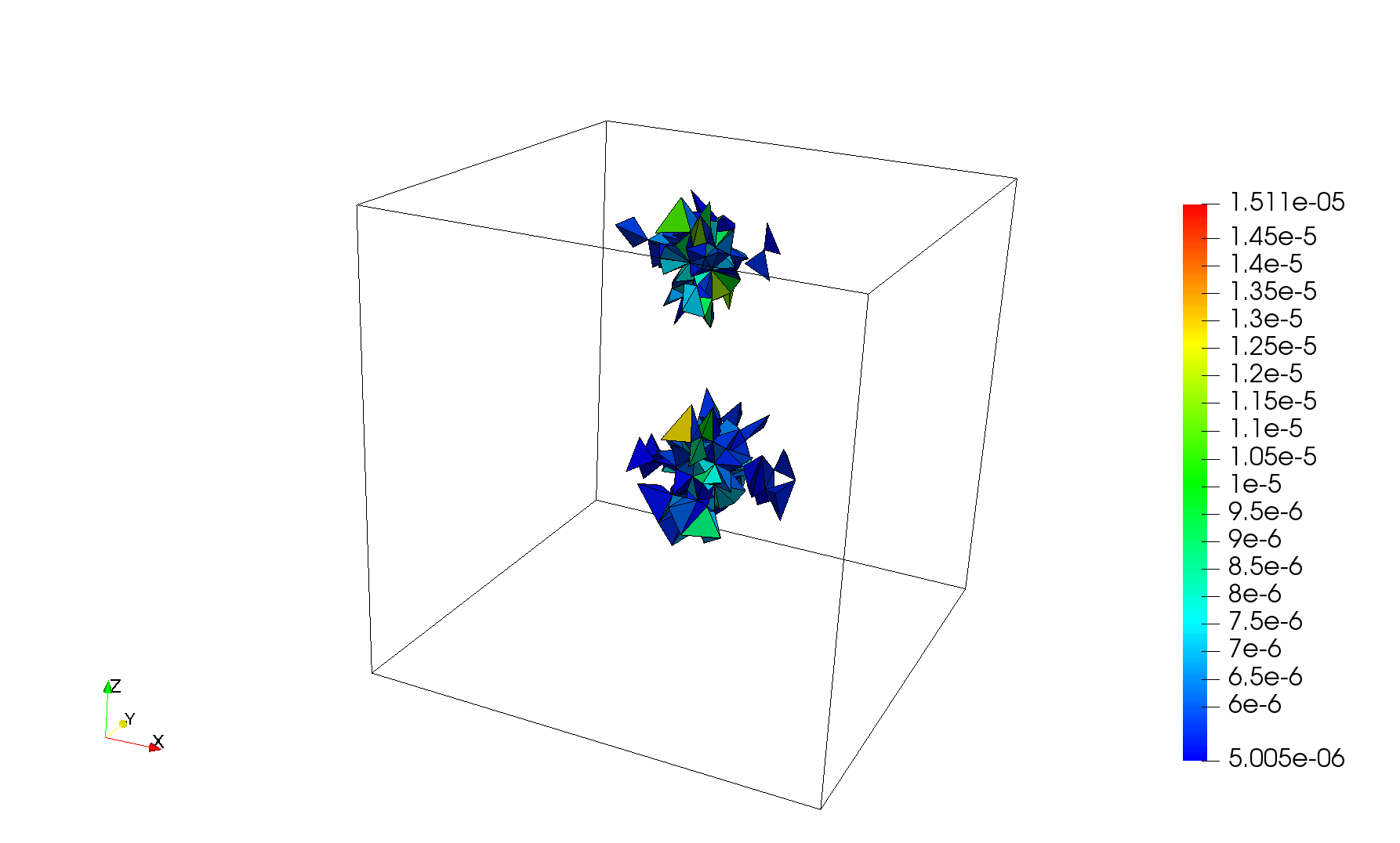}
	\end{subfigure}
	\\
	\rotatebox[origin=c]{90}{$m=3$}\hfill
	\begin{subfigure}[c]{0.22\textwidth}
		\centering
		\includegraphics[height=25mm]{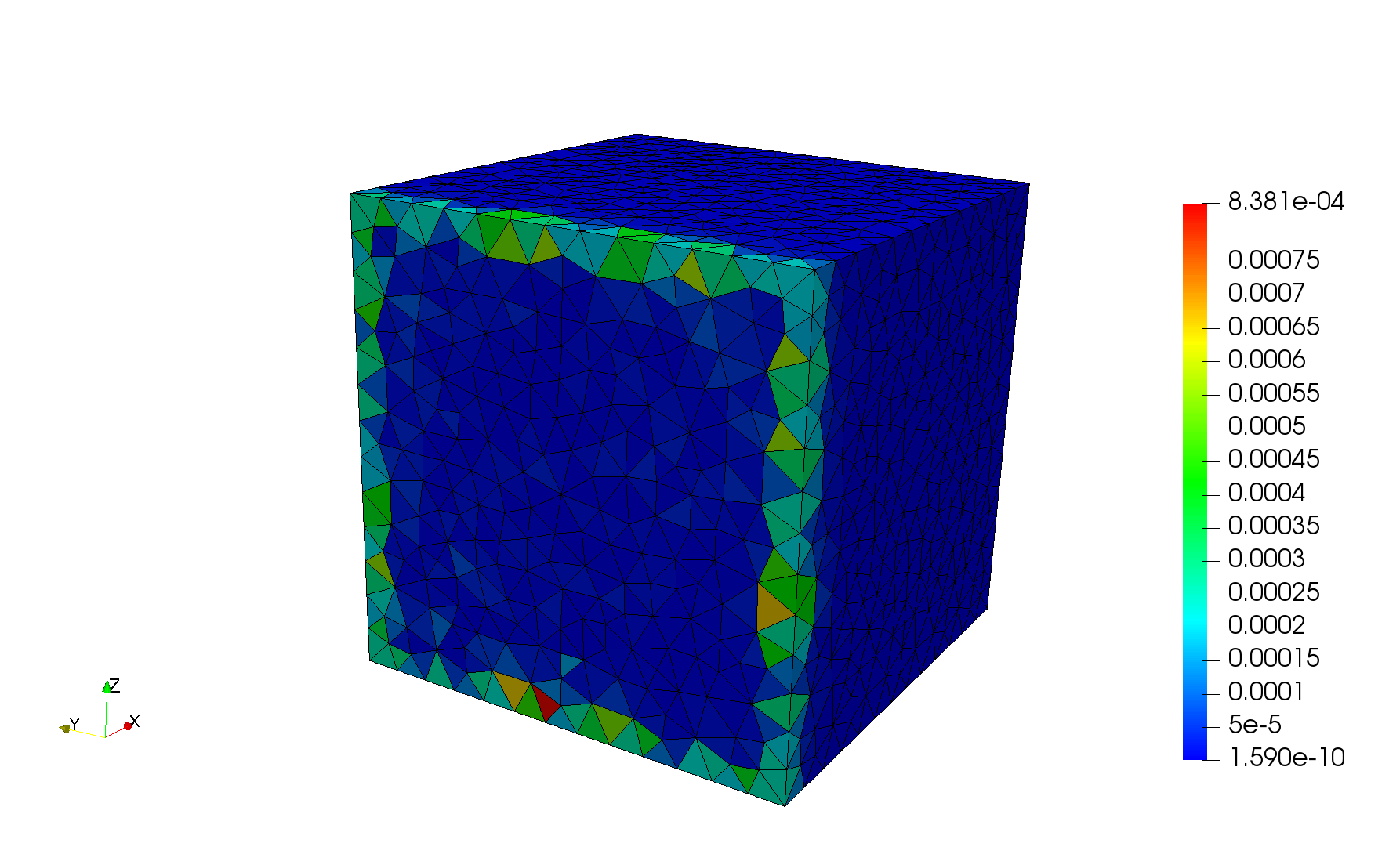}
	\end{subfigure}\hfill
	\begin{subfigure}[c]{0.22\textwidth}
		\centering
		\includegraphics[height=25mm]{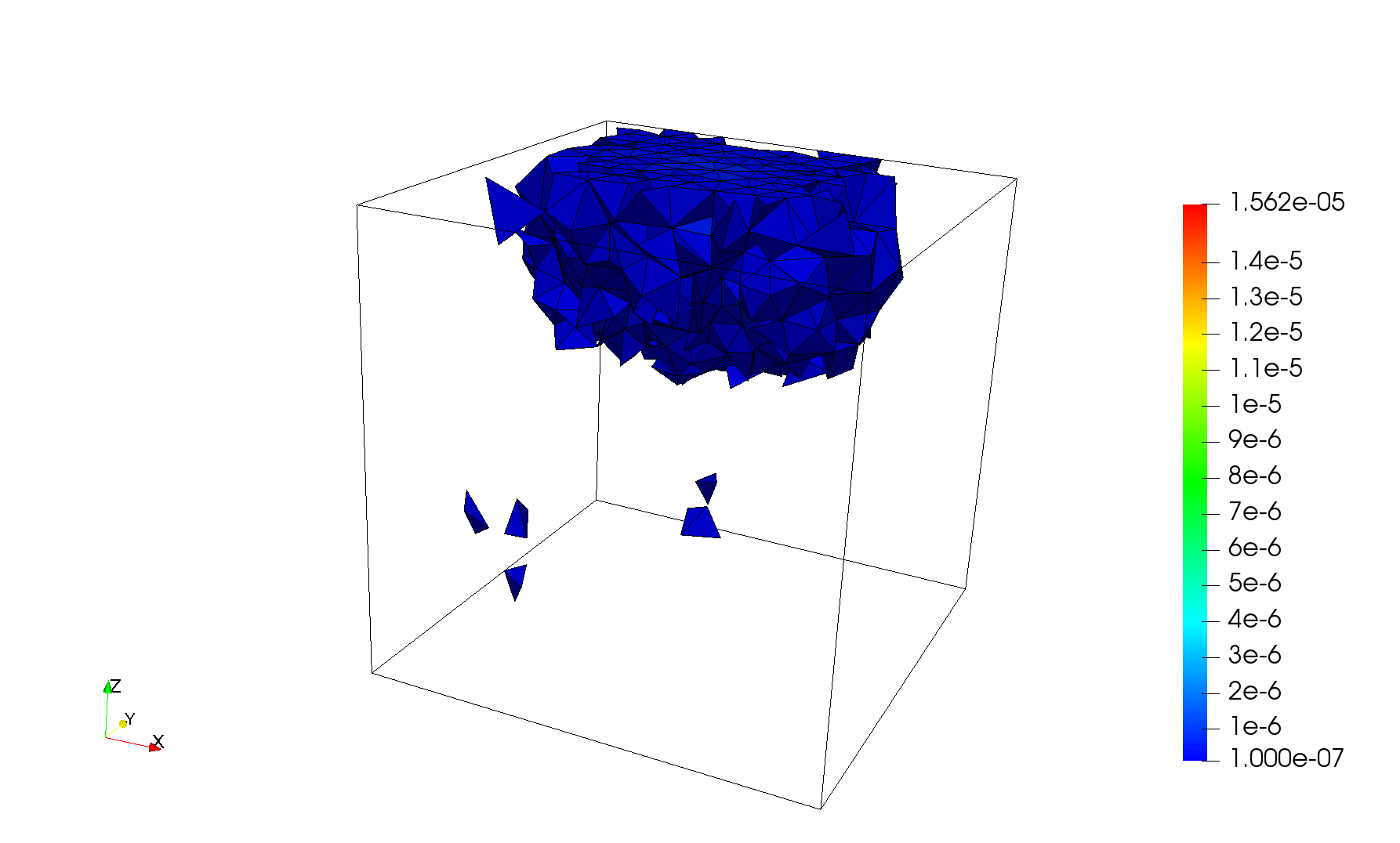}
	\end{subfigure}\hfill
	\begin{subfigure}[c]{0.22\textwidth}
		\centering
		\includegraphics[height=25mm]{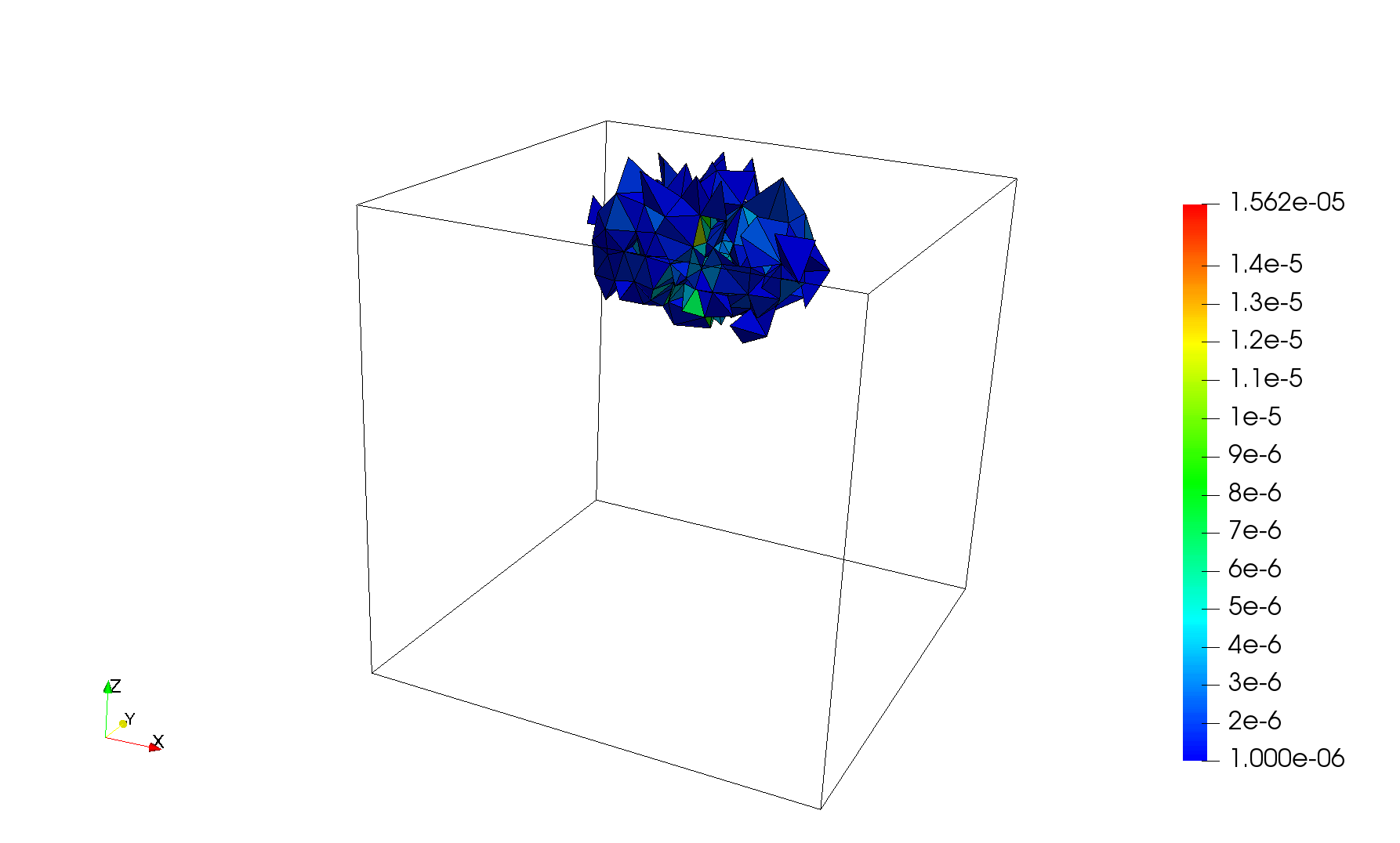}
	\end{subfigure}\hfill
	\begin{subfigure}[c]{0.22\textwidth}
		\centering
		\includegraphics[height=25mm]{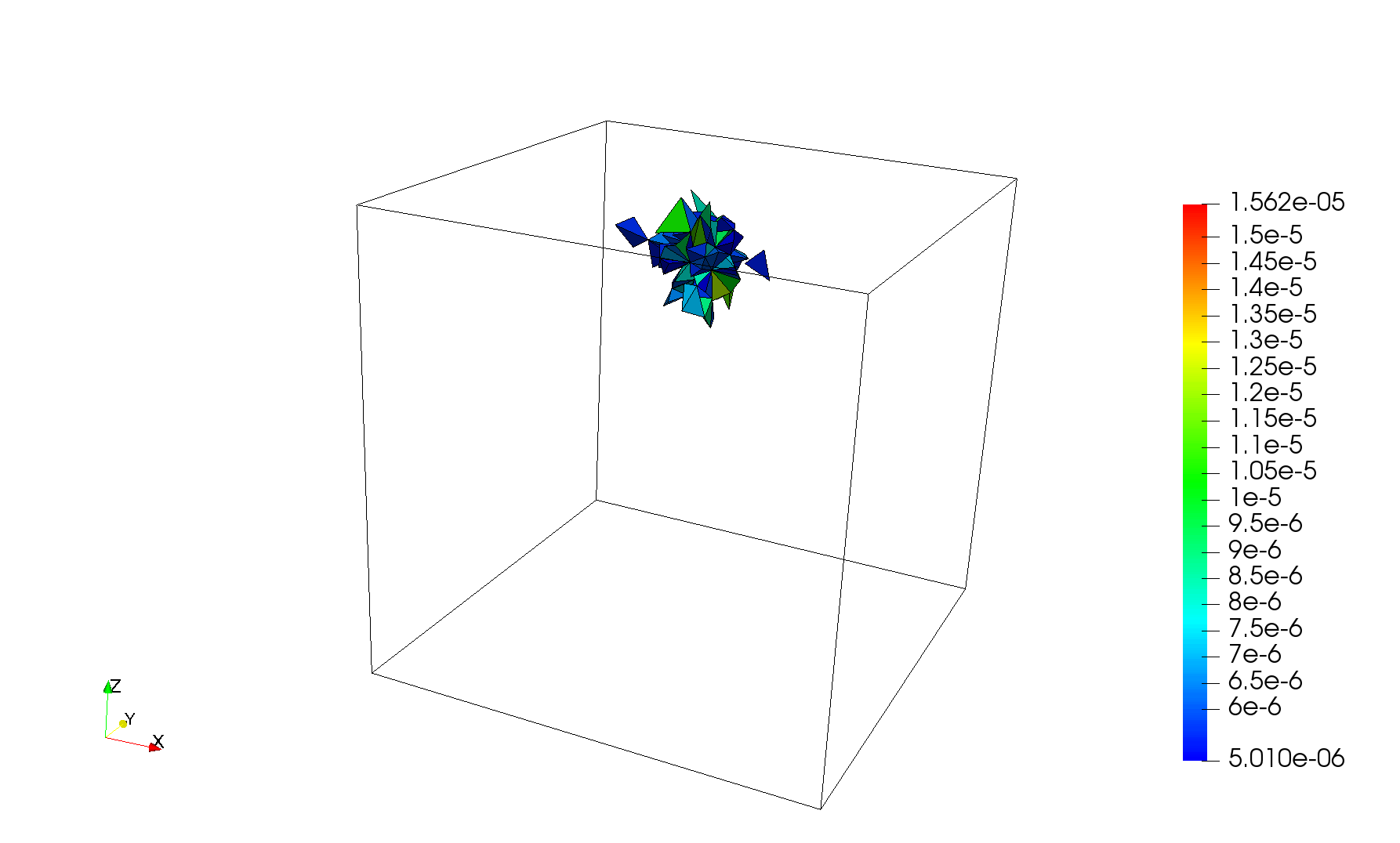}
	\end{subfigure}
	\\
	\rotatebox[origin=c]{90}{$m=4$}\hfill
	\begin{subfigure}[c]{0.22\textwidth}
		\centering
		\includegraphics[height=25mm]{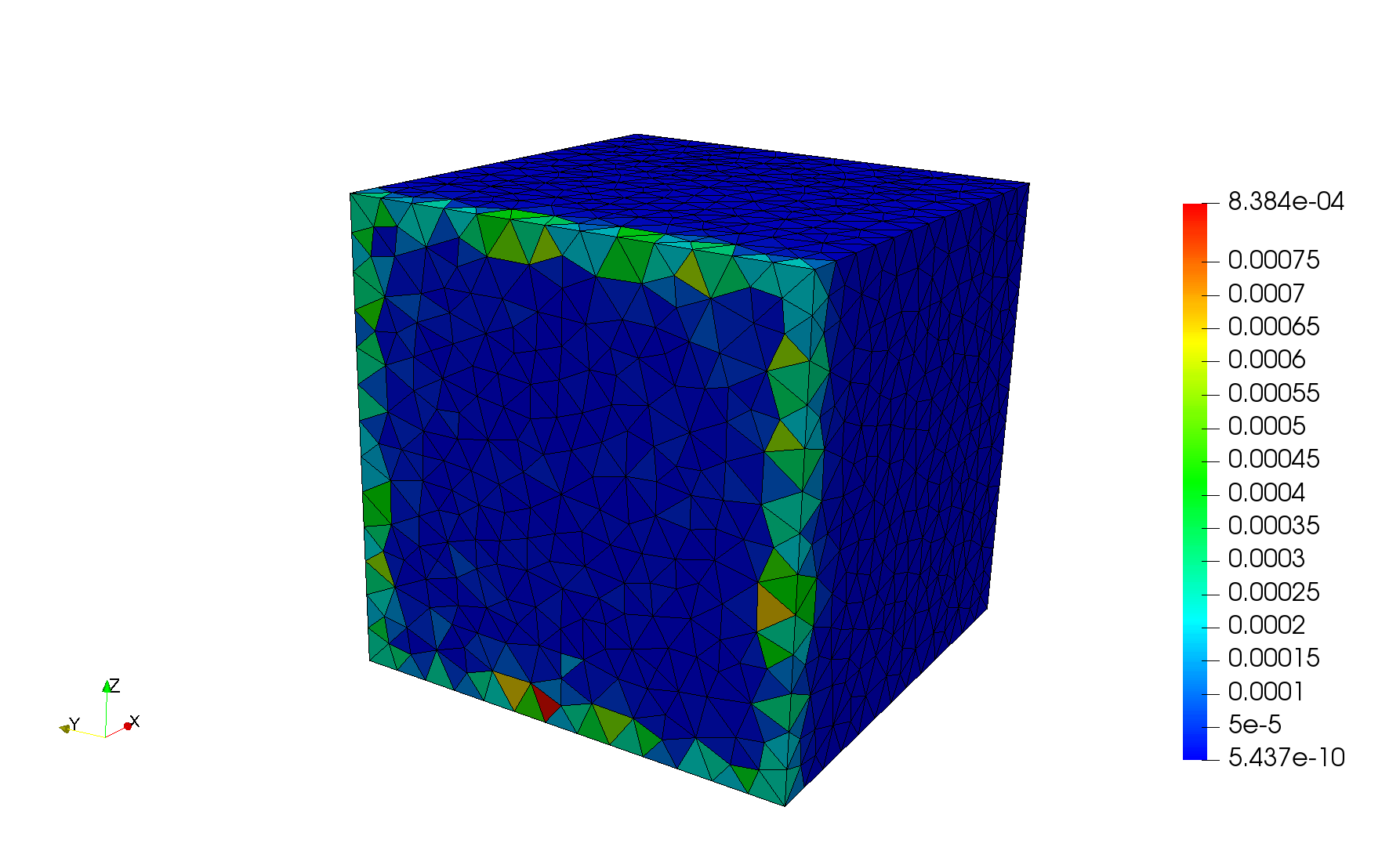}
	\end{subfigure}\hfill
	\begin{subfigure}[c]{0.22\textwidth}
		\centering
		\includegraphics[height=25mm]{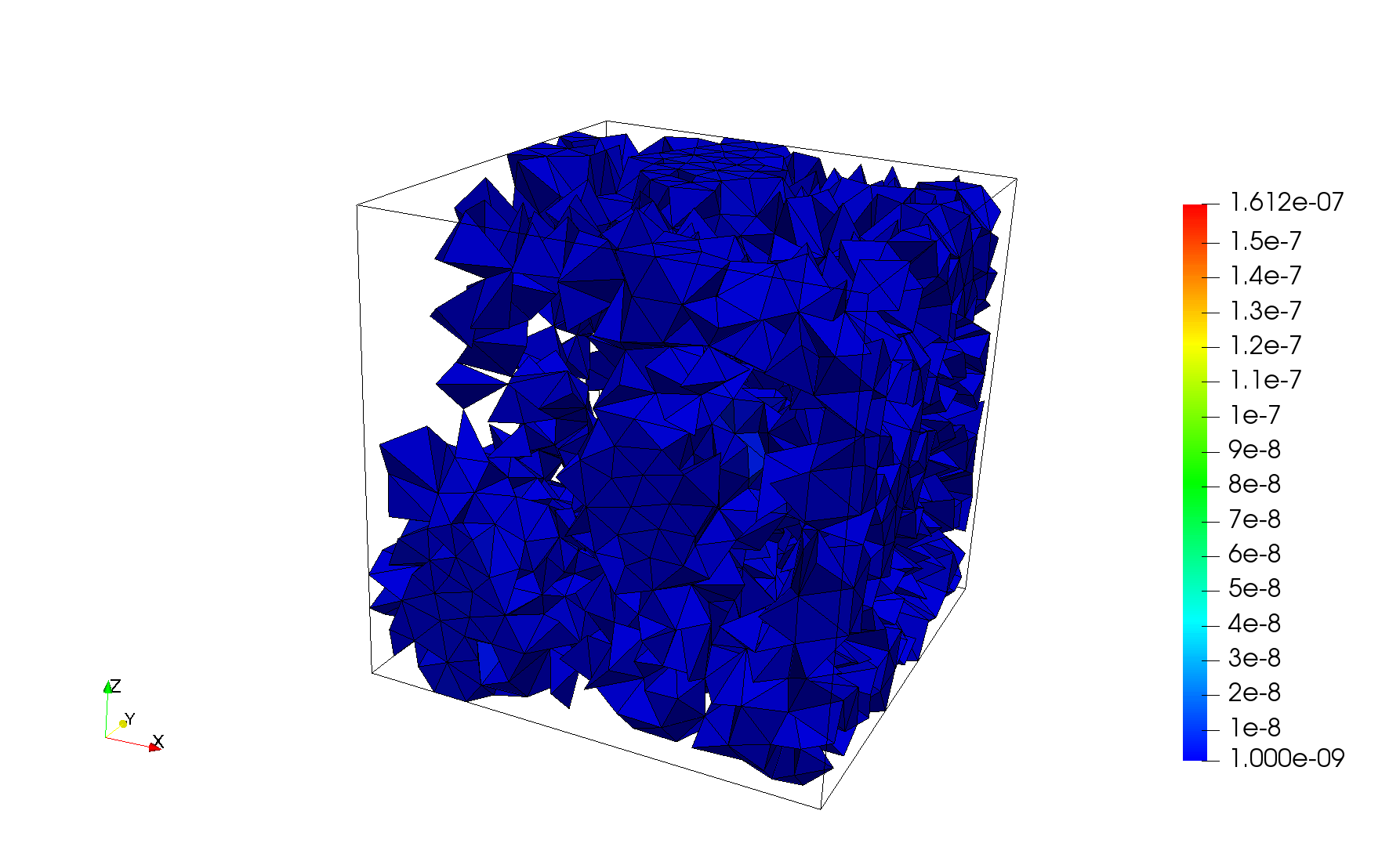}
	\end{subfigure}\hfill
	\begin{subfigure}[c]{0.22\textwidth}
		\centering
		\includegraphics[height=25mm]{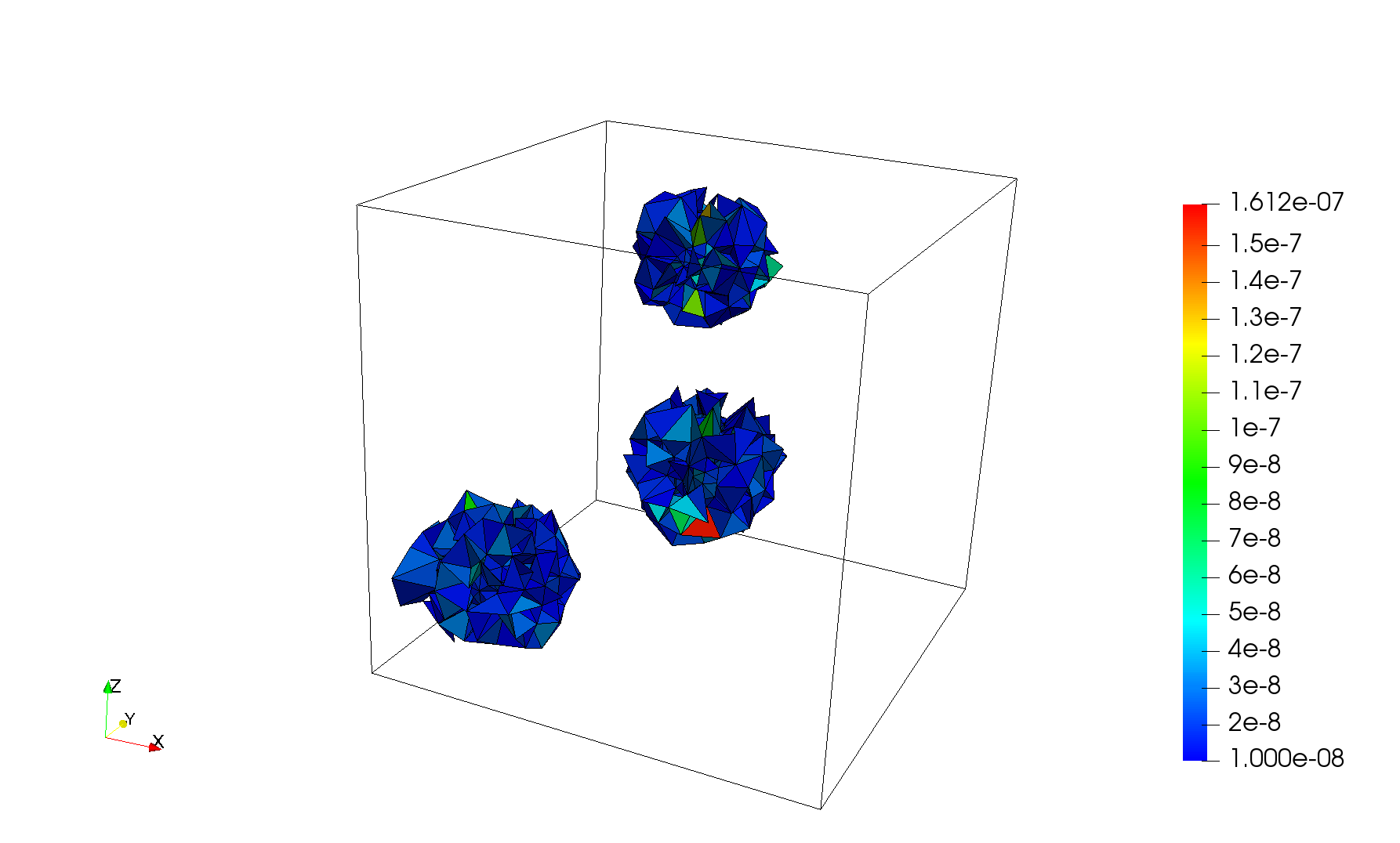}
	\end{subfigure}\hfill
	\begin{subfigure}[c]{0.22\textwidth}
		\centering
		\includegraphics[height=25mm]{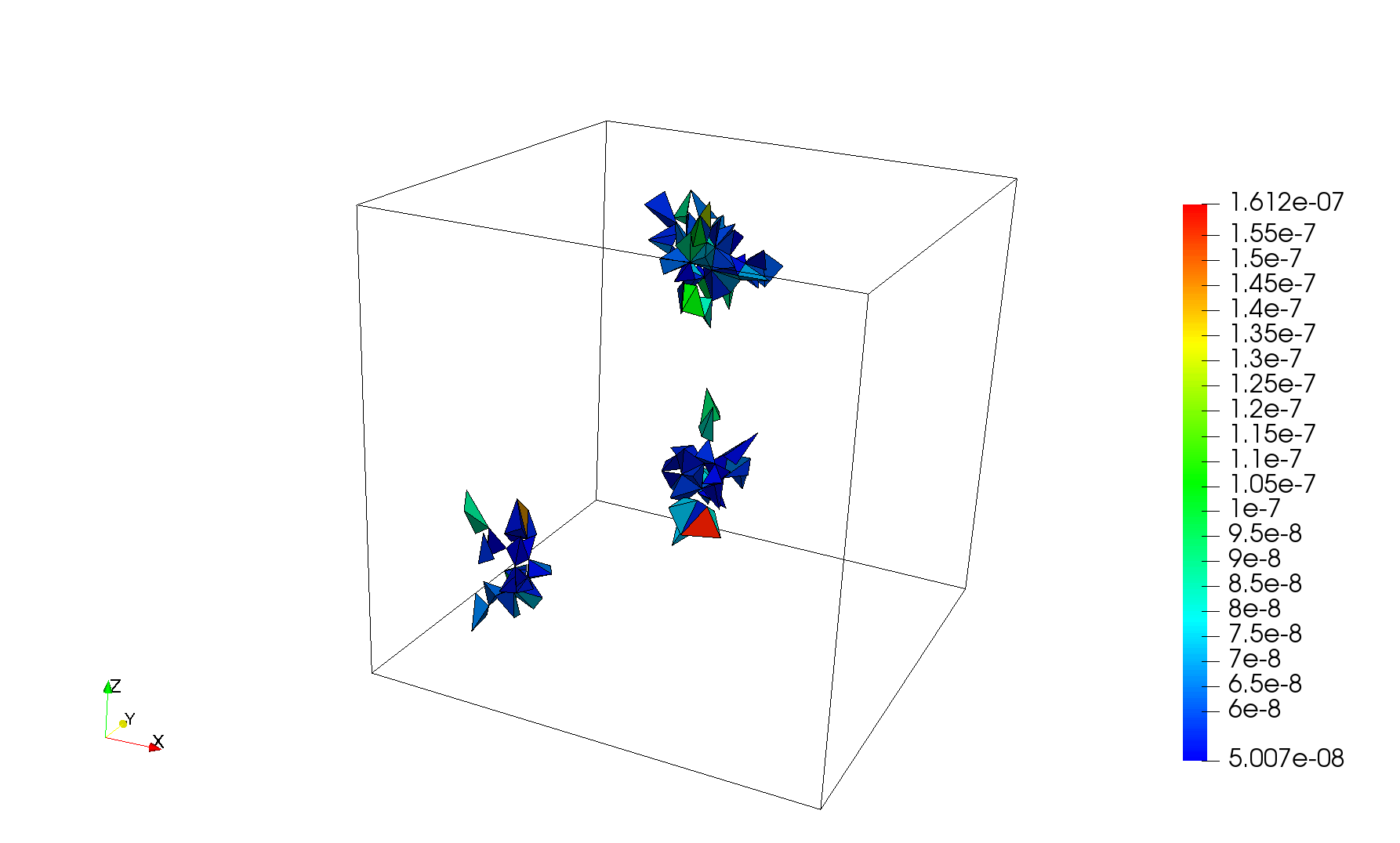}
	\end{subfigure}
	\\
	\rotatebox[origin=c]{90}{$m=5$}\hfill
	\begin{subfigure}[c]{0.22\textwidth}
		\centering
		\includegraphics[height=25mm]{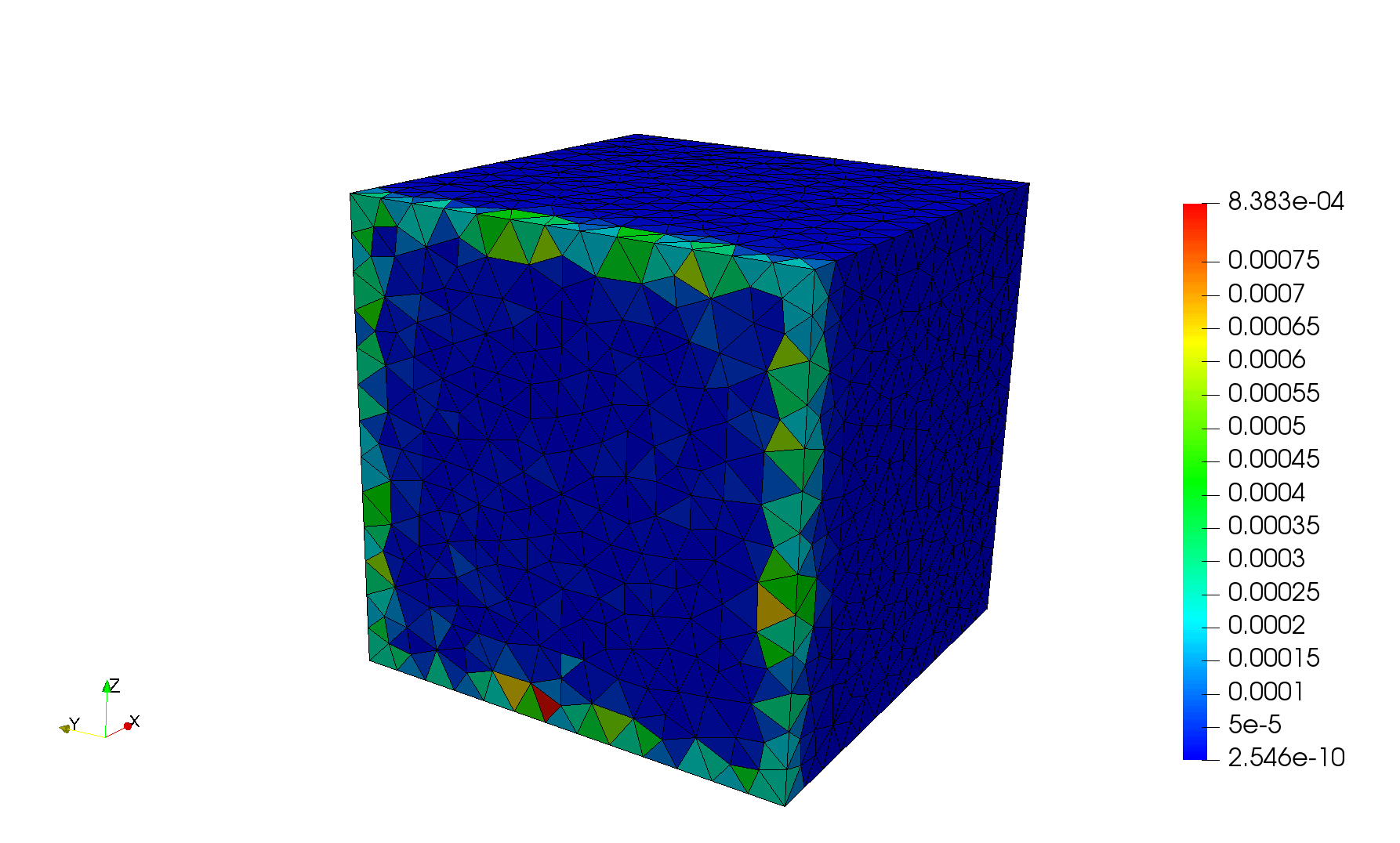}
	\end{subfigure}\hfill
	\begin{subfigure}[c]{0.22\textwidth}
		\centering
		\includegraphics[height=25mm]{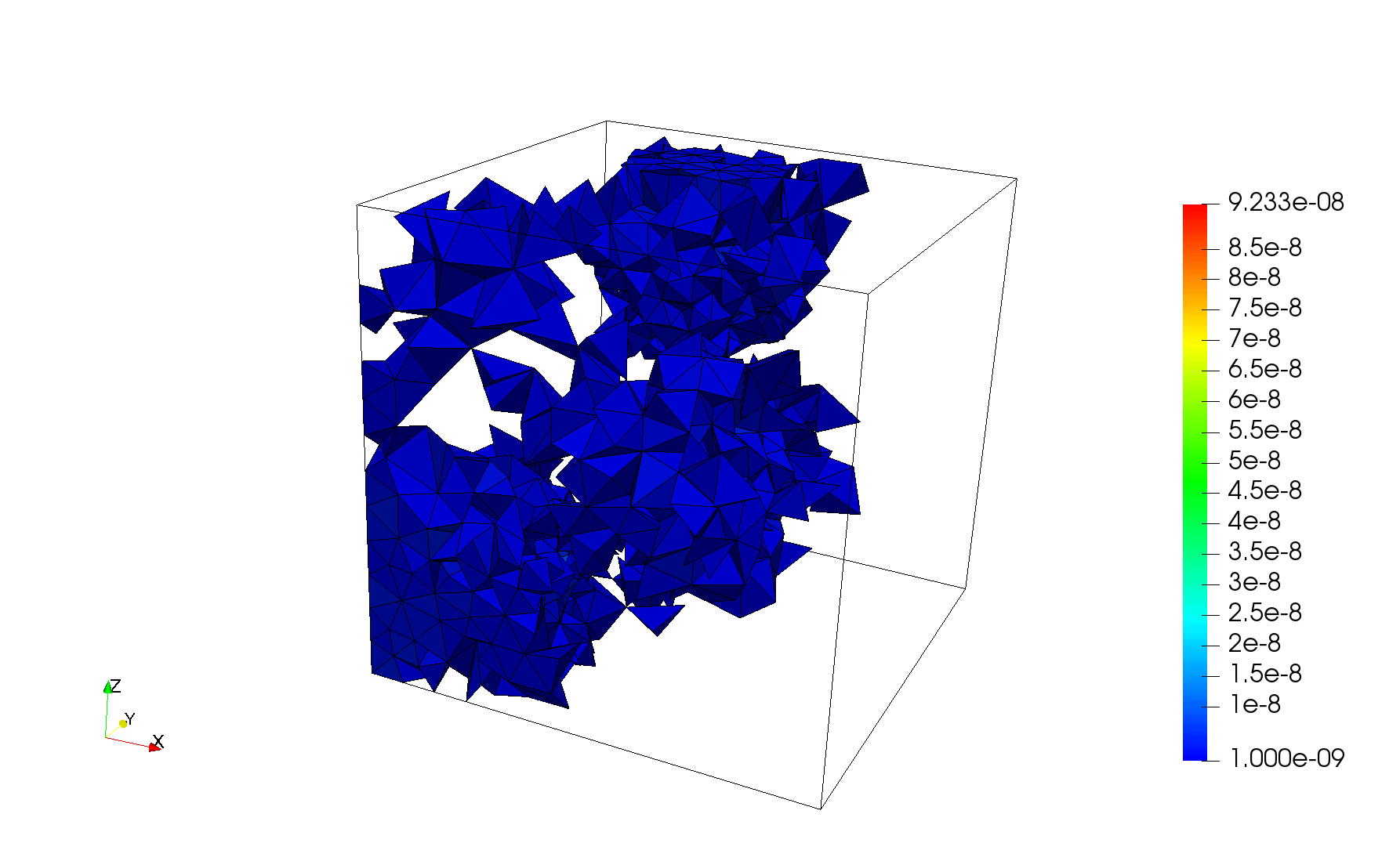}
	\end{subfigure}\hfill
	\begin{subfigure}[c]{0.22\textwidth}
		\centering
		\includegraphics[height=25mm]{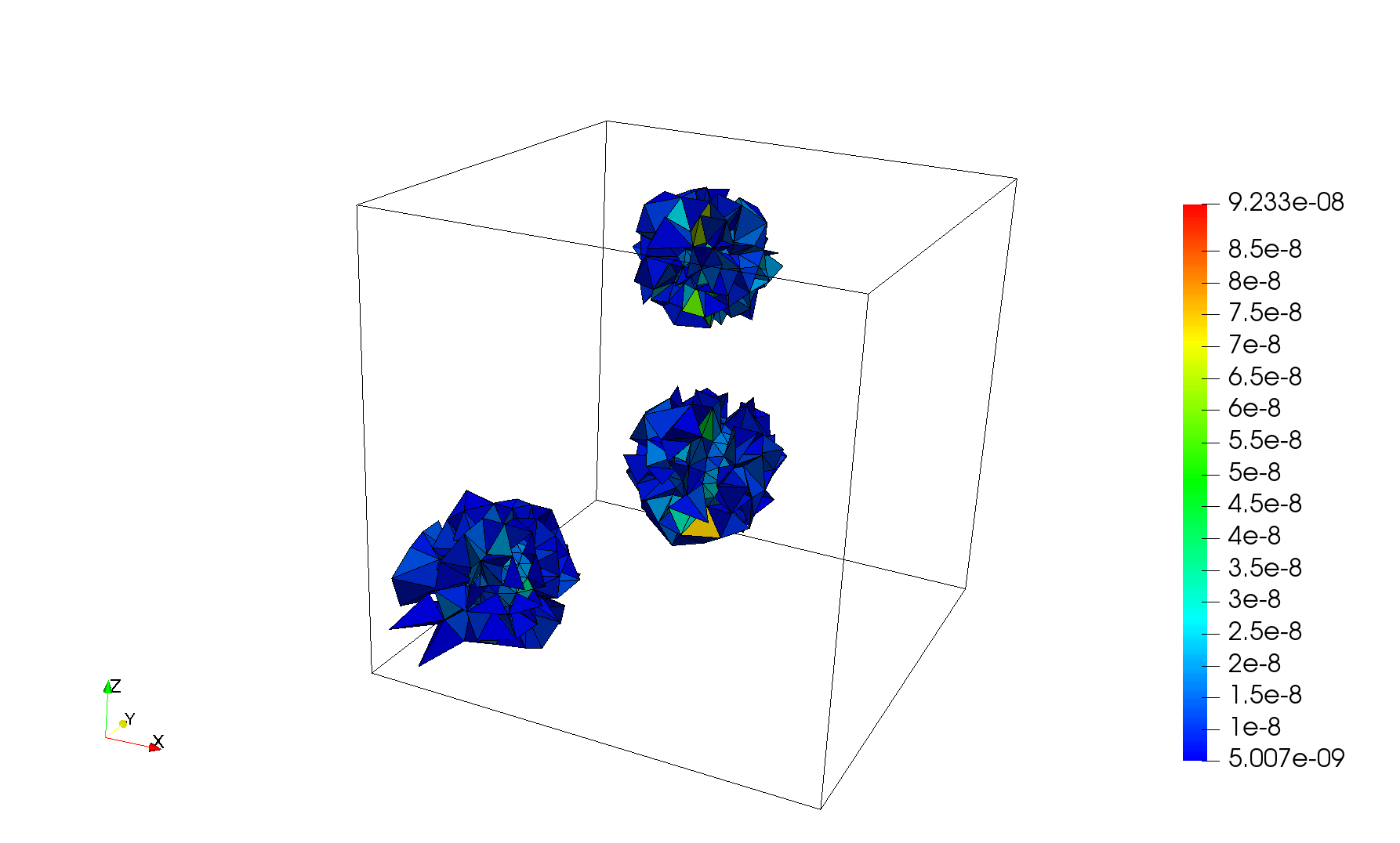}
	\end{subfigure}\hfill
	\begin{subfigure}[c]{0.22\textwidth}
		\centering
		\includegraphics[height=25mm]{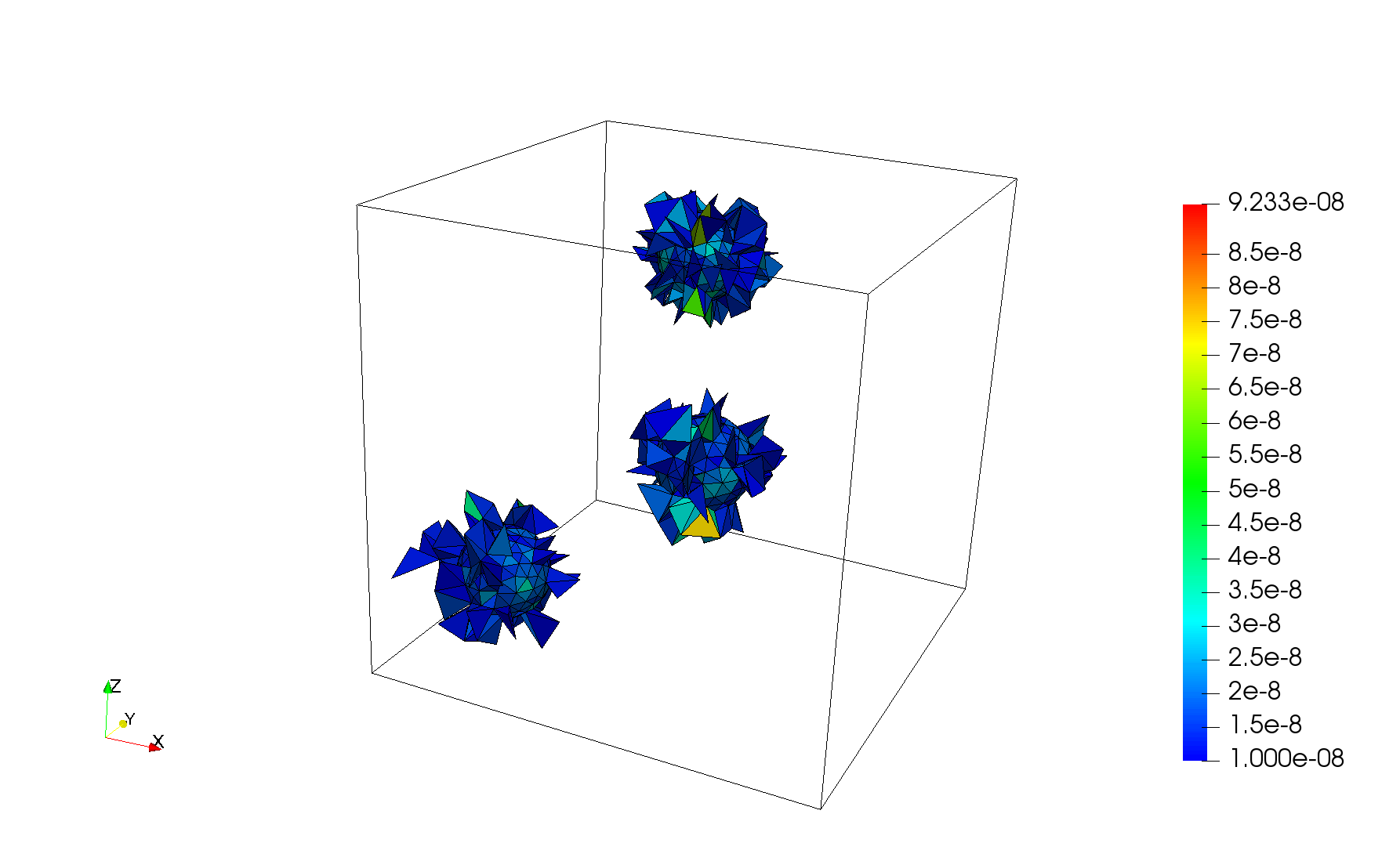}
	\end{subfigure}
	\caption{Spatial distribution of local contributions to the global error estimate $E_{\text{CRE}}^2$ (first column) and to the PGD truncation error indicator $\eta^2_{\PGD}$ (from second to last column) obtained for order $m=1,\dots,5$ (from top to bottom).}\label{fig:3D_distribution_estimates_global}
\end{figure}

Spatial distributions of local contributions to the error estimate $E_{\CRE}^2$ and to the PGD truncation error indicator $\eta^2_{\PGD}$ are shown in Figure~\ref{fig:3D_distribution_estimates_global} for different PGD decompositions ranging from order $1$ to $5$. It should be noted that the highest contributions of the error estimate are concentrated around the clamped boundary, while the ones of the PGD truncation error indicator are mainly located around the three inclusions for the orders~$m=1$, $4$ and $5$, the second and third inclusions for the order~$m=2$ and only the third inclusion for the order~$m=3$. Note that the contributions to the PGD truncation error indicator become negligible compared to the ones of the discretization error indicator for order $m\geq 4$.

The computation of space and parameter modes takes around 3 (resp. 4, 6, 8, 10 and 13)~min for mode 1 (resp. 2, 3, 4, 5 and 6), whereas the computational cost incurred by the error estimation procedure is about 1~min.


\section{Conclusions and prospects}
\label{section:conclusions}

We presented a strategy, based on the Constitutive Relation Error concept, that enables to obtain strict and accurate error estimates and drive adaptive strategies when dealing with the verification of PGD reduced-order models. It takes into account discretization and PGD truncation errors, and is applicable for controlling global error or local error in quantities of interest. Therefore, virtual charts associated with quantities of interest which may be computed from an approximate solution of PGD reduced-order models can now fully benefit from robust verification tools to satisfy a prescribed accuracy.

Future works will deal with the derivation of PGD verification tools for evolution (time-dependent) nonlinear problems in Computational Mechanics. The CRE approach (and associated admissible fields) seems to be a promising way for that purpose, as it has a direct extension to such complex mechanical problems using the concept of dissipation error (associated to the non-verification of the material evolution law) \cite{Lad98c,Lad06,Lad08}. 
This enables to define a residual with strong mechanical foundations.



\bibliographystyle{elsarticle-num}
\bibliography{./Biblio}


%
%
%
\end{document}

%% file: 3D_mode1_param1.tex
%
\begin{tikzpicture}

\begin{axis}[%
width=1.267\figureheight,
height=\figureheight,
at={(0\figureheight,0\figureheight)},
scale only axis,
separate axis lines,
every outer x axis line/.append style={black},
every x tick/.append style={black},
xmin=0,
xmax=10,
xlabel={$E_1$},
every outer y axis line/.append style={black},
every y tick/.append style={black},
ymin=0.1418,
ymax=0.1426,
ylabel={$\gamma_{1,1}$},
axis background/.style={fill=white},
xmajorgrids,
ymajorgrids
]
\addplot [color=red, line width=1.0pt, forget plot]
  table[row sep=crcr]{%
1	0.142574\\
1.09091	0.142567\\
1.18182	0.14256\\
1.27273	0.142553\\
1.36364	0.142546\\
1.45455	0.142539\\
1.54545	0.142532\\
1.63636	0.142525\\
1.72727	0.142518\\
1.81818	0.142511\\
1.90909	0.142504\\
2	0.142497\\
2.09091	0.14249\\
2.18182	0.142483\\
2.27273	0.142476\\
2.36364	0.142469\\
2.45455	0.142462\\
2.54545	0.142455\\
2.63636	0.142448\\
2.72727	0.142441\\
2.81818	0.142434\\
2.90909	0.142427\\
3	0.14242\\
3.09091	0.142413\\
3.18182	0.142406\\
3.27273	0.142399\\
3.36364	0.142392\\
3.45455	0.142385\\
3.54545	0.142378\\
3.63636	0.142371\\
3.72727	0.142364\\
3.81818	0.142357\\
3.90909	0.14235\\
4	0.142343\\
4.09091	0.142336\\
4.18182	0.142329\\
4.27273	0.142322\\
4.36364	0.142315\\
4.45455	0.142308\\
4.54545	0.142301\\
4.63636	0.142294\\
4.72727	0.142287\\
4.81818	0.14228\\
4.90909	0.142273\\
5	0.142266\\
5.09091	0.142259\\
5.18182	0.142252\\
5.27273	0.142245\\
5.36364	0.142238\\
5.45455	0.142231\\
5.54545	0.142224\\
5.63636	0.142217\\
5.72727	0.14221\\
5.81818	0.142203\\
5.90909	0.142196\\
6	0.142189\\
6.09091	0.142182\\
6.18182	0.142175\\
6.27273	0.142168\\
6.36364	0.142161\\
6.45455	0.142154\\
6.54545	0.142147\\
6.63636	0.14214\\
6.72727	0.142133\\
6.81818	0.142127\\
6.90909	0.14212\\
7	0.142113\\
7.09091	0.142106\\
7.18182	0.142099\\
7.27273	0.142092\\
7.36364	0.142085\\
7.45455	0.142078\\
7.54545	0.142071\\
7.63636	0.142064\\
7.72727	0.142057\\
7.81818	0.14205\\
7.90909	0.142043\\
8	0.142036\\
8.09091	0.142029\\
8.18182	0.142022\\
8.27273	0.142015\\
8.36364	0.142008\\
8.45455	0.142001\\
8.54545	0.141994\\
8.63636	0.141987\\
8.72727	0.14198\\
8.81818	0.141973\\
8.90909	0.141966\\
9	0.141959\\
9.09091	0.141953\\
9.18182	0.141946\\
9.27273	0.141939\\
9.36364	0.141932\\
9.45455	0.141925\\
9.54545	0.141918\\
9.63636	0.141911\\
9.72727	0.141904\\
9.81818	0.141897\\
9.90909	0.14189\\
10	0.141883\\
};
\end{axis}
\end{tikzpicture}%

%% file: 3D_mode1_param2.tex
%
\begin{tikzpicture}

\begin{axis}[%
width=1.267\figureheight,
height=\figureheight,
at={(0\figureheight,0\figureheight)},
scale only axis,
separate axis lines,
every outer x axis line/.append style={black},
every x tick/.append style={black},
xmin=0,
xmax=10,
xlabel={$E_2$},
every outer y axis line/.append style={black},
every y tick/.append style={black},
ymin=0.1418,
ymax=0.1426,
ylabel={$\gamma_{2,1}$},
axis background/.style={fill=white},
xmajorgrids,
ymajorgrids
]
\addplot [color=green, line width=1.0pt, forget plot]
  table[row sep=crcr]{%
1	0.142589\\
1.09091	0.142582\\
1.18182	0.142574\\
1.27273	0.142567\\
1.36364	0.14256\\
1.45455	0.142553\\
1.54545	0.142545\\
1.63636	0.142538\\
1.72727	0.142531\\
1.81818	0.142524\\
1.90909	0.142516\\
2	0.142509\\
2.09091	0.142502\\
2.18182	0.142495\\
2.27273	0.142487\\
2.36364	0.14248\\
2.45455	0.142473\\
2.54545	0.142466\\
2.63636	0.142458\\
2.72727	0.142451\\
2.81818	0.142444\\
2.90909	0.142437\\
3	0.14243\\
3.09091	0.142422\\
3.18182	0.142415\\
3.27273	0.142408\\
3.36364	0.142401\\
3.45455	0.142393\\
3.54545	0.142386\\
3.63636	0.142379\\
3.72727	0.142372\\
3.81818	0.142364\\
3.90909	0.142357\\
4	0.14235\\
4.09091	0.142343\\
4.18182	0.142336\\
4.27273	0.142328\\
4.36364	0.142321\\
4.45455	0.142314\\
4.54545	0.142307\\
4.63636	0.142299\\
4.72727	0.142292\\
4.81818	0.142285\\
4.90909	0.142278\\
5	0.142271\\
5.09091	0.142263\\
5.18182	0.142256\\
5.27273	0.142249\\
5.36364	0.142242\\
5.45455	0.142234\\
5.54545	0.142227\\
5.63636	0.14222\\
5.72727	0.142213\\
5.81818	0.142206\\
5.90909	0.142198\\
6	0.142191\\
6.09091	0.142184\\
6.18182	0.142177\\
6.27273	0.14217\\
6.36364	0.142162\\
6.45455	0.142155\\
6.54545	0.142148\\
6.63636	0.142141\\
6.72727	0.142133\\
6.81818	0.142126\\
6.90909	0.142119\\
7	0.142112\\
7.09091	0.142105\\
7.18182	0.142097\\
7.27273	0.14209\\
7.36364	0.142083\\
7.45455	0.142076\\
7.54545	0.142069\\
7.63636	0.142061\\
7.72727	0.142054\\
7.81818	0.142047\\
7.90909	0.14204\\
8	0.142033\\
8.09091	0.142025\\
8.18182	0.142018\\
8.27273	0.142011\\
8.36364	0.142004\\
8.45455	0.141997\\
8.54545	0.141989\\
8.63636	0.141982\\
8.72727	0.141975\\
8.81818	0.141968\\
8.90909	0.141961\\
9	0.141954\\
9.09091	0.141946\\
9.18182	0.141939\\
9.27273	0.141932\\
9.36364	0.141925\\
9.45455	0.141918\\
9.54545	0.14191\\
9.63636	0.141903\\
9.72727	0.141896\\
9.81818	0.141889\\
9.90909	0.141882\\
10	0.141875\\
};
\end{axis}
\end{tikzpicture}%

%% file: 3D_mode1_param3.tex
%
\begin{tikzpicture}

\begin{axis}[%
width=1.267\figureheight,
height=\figureheight,
at={(0\figureheight,0\figureheight)},
scale only axis,
separate axis lines,
every outer x axis line/.append style={black},
every x tick/.append style={black},
xmin=0,
xmax=10,
xlabel={$E_3$},
every outer y axis line/.append style={black},
every y tick/.append style={black},
ymin=0.1418,
ymax=0.1426,
ylabel={$\gamma_{3,1}$},
axis background/.style={fill=white},
xmajorgrids,
ymajorgrids
]
\addplot [color=blue, line width=1.0pt, forget plot]
  table[row sep=crcr]{%
1	0.142575\\
1.09091	0.142568\\
1.18182	0.142561\\
1.27273	0.142554\\
1.36364	0.142547\\
1.45455	0.14254\\
1.54545	0.142533\\
1.63636	0.142526\\
1.72727	0.142519\\
1.81818	0.142512\\
1.90909	0.142505\\
2	0.142498\\
2.09091	0.142491\\
2.18182	0.142484\\
2.27273	0.142477\\
2.36364	0.14247\\
2.45455	0.142463\\
2.54545	0.142456\\
2.63636	0.142449\\
2.72727	0.142442\\
2.81818	0.142435\\
2.90909	0.142428\\
3	0.142421\\
3.09091	0.142414\\
3.18182	0.142407\\
3.27273	0.1424\\
3.36364	0.142393\\
3.45455	0.142386\\
3.54545	0.142379\\
3.63636	0.142372\\
3.72727	0.142365\\
3.81818	0.142357\\
3.90909	0.14235\\
4	0.142343\\
4.09091	0.142336\\
4.18182	0.142329\\
4.27273	0.142322\\
4.36364	0.142315\\
4.45455	0.142308\\
4.54545	0.142301\\
4.63636	0.142294\\
4.72727	0.142287\\
4.81818	0.14228\\
4.90909	0.142273\\
5	0.142266\\
5.09091	0.142259\\
5.18182	0.142252\\
5.27273	0.142245\\
5.36364	0.142238\\
5.45455	0.142231\\
5.54545	0.142224\\
5.63636	0.142217\\
5.72727	0.14221\\
5.81818	0.142203\\
5.90909	0.142196\\
6	0.142189\\
6.09091	0.142182\\
6.18182	0.142175\\
6.27273	0.142168\\
6.36364	0.142161\\
6.45455	0.142154\\
6.54545	0.142147\\
6.63636	0.14214\\
6.72727	0.142133\\
6.81818	0.142126\\
6.90909	0.14212\\
7	0.142113\\
7.09091	0.142106\\
7.18182	0.142099\\
7.27273	0.142092\\
7.36364	0.142085\\
7.45455	0.142078\\
7.54545	0.142071\\
7.63636	0.142064\\
7.72727	0.142057\\
7.81818	0.14205\\
7.90909	0.142043\\
8	0.142036\\
8.09091	0.142029\\
8.18182	0.142022\\
8.27273	0.142015\\
8.36364	0.142008\\
8.45455	0.142001\\
8.54545	0.141994\\
8.63636	0.141987\\
8.72727	0.14198\\
8.81818	0.141973\\
8.90909	0.141966\\
9	0.141959\\
9.09091	0.141952\\
9.18182	0.141945\\
9.27273	0.141938\\
9.36364	0.141931\\
9.45455	0.141924\\
9.54545	0.141917\\
9.63636	0.14191\\
9.72727	0.141903\\
9.81818	0.141896\\
9.90909	0.141889\\
10	0.141882\\
};
\end{axis}
\end{tikzpicture}%

%% file: 3D_mode2_param1.tex
%
\begin{tikzpicture}

\begin{axis}[%
width=1.268\figureheight,
height=\figureheight,
at={(0\figureheight,0\figureheight)},
scale only axis,
separate axis lines,
every outer x axis line/.append style={black},
every x tick/.append style={black},
xmin=0,
xmax=10,
xlabel={$E_1$},
every outer y axis line/.append style={black},
every y tick/.append style={black},
ymin=-0.7,
ymax=0.2,
ylabel={$\gamma_{1,2}$},
axis background/.style={fill=white},
xmajorgrids,
ymajorgrids
]
\addplot [color=red, line width=1.0pt, forget plot]
  table[row sep=crcr]{%
1	-0.655163\\
1.09091	-0.620159\\
1.18182	-0.587452\\
1.27273	-0.556817\\
1.36364	-0.528065\\
1.45455	-0.501027\\
1.54545	-0.475555\\
1.63636	-0.451516\\
1.72727	-0.428793\\
1.81818	-0.40728\\
1.90909	-0.386884\\
2	-0.36752\\
2.09091	-0.349112\\
2.18182	-0.331591\\
2.27273	-0.314893\\
2.36364	-0.298963\\
2.45455	-0.283749\\
2.54545	-0.269203\\
2.63636	-0.255283\\
2.72727	-0.24195\\
2.81818	-0.229166\\
2.90909	-0.216898\\
3	-0.205117\\
3.09091	-0.193793\\
3.18182	-0.182901\\
3.27273	-0.172416\\
3.36364	-0.162317\\
3.45455	-0.152581\\
3.54545	-0.143191\\
3.63636	-0.134127\\
3.72727	-0.125374\\
3.81818	-0.116915\\
3.90909	-0.108736\\
4	-0.100823\\
4.09091	-0.0931642\\
4.18182	-0.0857467\\
4.27273	-0.0785596\\
4.36364	-0.0715923\\
4.45455	-0.064835\\
4.54545	-0.0582782\\
4.63636	-0.0519132\\
4.72727	-0.0457318\\
4.81818	-0.0397261\\
4.90909	-0.0338888\\
5	-0.0282129\\
5.09091	-0.0226918\\
5.18182	-0.0173194\\
5.27273	-0.0120896\\
5.36364	-0.00699691\\
5.45455	-0.00203603\\
5.54545	0.00279807\\
5.63636	0.00751019\\
5.72727	0.0121049\\
5.81818	0.0165864\\
5.90909	0.0209589\\
6	0.0252264\\
6.09091	0.0293924\\
6.18182	0.0334606\\
6.27273	0.0374344\\
6.36364	0.0413171\\
6.45455	0.0451116\\
6.54545	0.0488209\\
6.63636	0.052448\\
6.72727	0.0559954\\
6.81818	0.0594657\\
6.90909	0.0628615\\
7	0.066185\\
7.09091	0.0694386\\
7.18182	0.0726244\\
7.27273	0.0757445\\
7.36364	0.078801\\
7.45455	0.0817957\\
7.54545	0.0847304\\
7.63636	0.087607\\
7.72727	0.0904271\\
7.81818	0.0931924\\
7.90909	0.0959045\\
8	0.0985648\\
8.09091	0.101175\\
8.18182	0.103736\\
8.27273	0.10625\\
8.36364	0.108717\\
8.45455	0.111139\\
8.54545	0.113518\\
8.63636	0.115854\\
8.72727	0.118149\\
8.81818	0.120403\\
8.90909	0.122618\\
9	0.124795\\
9.09091	0.126934\\
9.18182	0.129037\\
9.27273	0.131104\\
9.36364	0.133137\\
9.45455	0.135136\\
9.54545	0.137102\\
9.63636	0.139036\\
9.72727	0.140938\\
9.81818	0.14281\\
9.90909	0.144653\\
10	0.146466\\
};
\end{axis}
\end{tikzpicture}%

%% file: 3D_mode2_param2.tex
%
\begin{tikzpicture}

\begin{axis}[%
width=1.267\figureheight,
height=\figureheight,
at={(0\figureheight,0\figureheight)},
scale only axis,
separate axis lines,
every outer x axis line/.append style={black},
every x tick/.append style={black},
xmin=0,
xmax=10,
xlabel={$E_2$},
every outer y axis line/.append style={black},
every y tick/.append style={black},
ymin=0.1419,
ymax=0.1425,
ylabel={$\gamma_{2,2}$},
axis background/.style={fill=white},
xmajorgrids,
ymajorgrids
]
\addplot [color=green, line width=1.0pt, forget plot]
  table[row sep=crcr]{%
1	0.14241\\
1.09091	0.142406\\
1.18182	0.142402\\
1.27273	0.142397\\
1.36364	0.142393\\
1.45455	0.142388\\
1.54545	0.142384\\
1.63636	0.14238\\
1.72727	0.142375\\
1.81818	0.142371\\
1.90909	0.142366\\
2	0.142362\\
2.09091	0.142358\\
2.18182	0.142353\\
2.27273	0.142349\\
2.36364	0.142344\\
2.45455	0.14234\\
2.54545	0.142336\\
2.63636	0.142331\\
2.72727	0.142327\\
2.81818	0.142322\\
2.90909	0.142318\\
3	0.142314\\
3.09091	0.142309\\
3.18182	0.142305\\
3.27273	0.1423\\
3.36364	0.142296\\
3.45455	0.142292\\
3.54545	0.142287\\
3.63636	0.142283\\
3.72727	0.142278\\
3.81818	0.142274\\
3.90909	0.14227\\
4	0.142265\\
4.09091	0.142261\\
4.18182	0.142256\\
4.27273	0.142252\\
4.36364	0.142248\\
4.45455	0.142243\\
4.54545	0.142239\\
4.63636	0.142234\\
4.72727	0.14223\\
4.81818	0.142226\\
4.90909	0.142221\\
5	0.142217\\
5.09091	0.142213\\
5.18182	0.142208\\
5.27273	0.142204\\
5.36364	0.142199\\
5.45455	0.142195\\
5.54545	0.142191\\
5.63636	0.142186\\
5.72727	0.142182\\
5.81818	0.142177\\
5.90909	0.142173\\
6	0.142169\\
6.09091	0.142164\\
6.18182	0.14216\\
6.27273	0.142156\\
6.36364	0.142151\\
6.45455	0.142147\\
6.54545	0.142142\\
6.63636	0.142138\\
6.72727	0.142134\\
6.81818	0.142129\\
6.90909	0.142125\\
7	0.14212\\
7.09091	0.142116\\
7.18182	0.142112\\
7.27273	0.142107\\
7.36364	0.142103\\
7.45455	0.142099\\
7.54545	0.142094\\
7.63636	0.14209\\
7.72727	0.142085\\
7.81818	0.142081\\
7.90909	0.142077\\
8	0.142072\\
8.09091	0.142068\\
8.18182	0.142064\\
8.27273	0.142059\\
8.36364	0.142055\\
8.45455	0.14205\\
8.54545	0.142046\\
8.63636	0.142042\\
8.72727	0.142037\\
8.81818	0.142033\\
8.90909	0.142029\\
9	0.142024\\
9.09091	0.14202\\
9.18182	0.142016\\
9.27273	0.142011\\
9.36364	0.142007\\
9.45455	0.142002\\
9.54545	0.141998\\
9.63636	0.141994\\
9.72727	0.141989\\
9.81818	0.141985\\
9.90909	0.141981\\
10	0.141976\\
};
\end{axis}
\end{tikzpicture}%

%% file: 3D_mode2_param3.tex
%
\begin{tikzpicture}

\begin{axis}[%
width=1.267\figureheight,
height=\figureheight,
at={(0\figureheight,0\figureheight)},
scale only axis,
separate axis lines,
every outer x axis line/.append style={black},
every x tick/.append style={black},
xmin=0,
xmax=10,
xlabel={$E_3$},
every outer y axis line/.append style={black},
every y tick/.append style={black},
ymin=0.1417,
ymax=0.1425,
ylabel={$\gamma_{3,2}$},
axis background/.style={fill=white},
xmajorgrids,
ymajorgrids
]
\addplot [color=blue, line width=1.0pt, forget plot]
  table[row sep=crcr]{%
1	0.141781\\
1.09091	0.141786\\
1.18182	0.141792\\
1.27273	0.141797\\
1.36364	0.141803\\
1.45455	0.141809\\
1.54545	0.141814\\
1.63636	0.14182\\
1.72727	0.141826\\
1.81818	0.141831\\
1.90909	0.141837\\
2	0.141842\\
2.09091	0.141848\\
2.18182	0.141854\\
2.27273	0.141859\\
2.36364	0.141865\\
2.45455	0.14187\\
2.54545	0.141876\\
2.63636	0.141882\\
2.72727	0.141887\\
2.81818	0.141893\\
2.90909	0.141899\\
3	0.141904\\
3.09091	0.14191\\
3.18182	0.141915\\
3.27273	0.141921\\
3.36364	0.141927\\
3.45455	0.141932\\
3.54545	0.141938\\
3.63636	0.141943\\
3.72727	0.141949\\
3.81818	0.141955\\
3.90909	0.14196\\
4	0.141966\\
4.09091	0.141972\\
4.18182	0.141977\\
4.27273	0.141983\\
4.36364	0.141988\\
4.45455	0.141994\\
4.54545	0.142\\
4.63636	0.142005\\
4.72727	0.142011\\
4.81818	0.142016\\
4.90909	0.142022\\
5	0.142028\\
5.09091	0.142033\\
5.18182	0.142039\\
5.27273	0.142044\\
5.36364	0.14205\\
5.45455	0.142056\\
5.54545	0.142061\\
5.63636	0.142067\\
5.72727	0.142072\\
5.81818	0.142078\\
5.90909	0.142084\\
6	0.142089\\
6.09091	0.142095\\
6.18182	0.1421\\
6.27273	0.142106\\
6.36364	0.142112\\
6.45455	0.142117\\
6.54545	0.142123\\
6.63636	0.142128\\
6.72727	0.142134\\
6.81818	0.14214\\
6.90909	0.142145\\
7	0.142151\\
7.09091	0.142156\\
7.18182	0.142162\\
7.27273	0.142167\\
7.36364	0.142173\\
7.45455	0.142179\\
7.54545	0.142184\\
7.63636	0.14219\\
7.72727	0.142195\\
7.81818	0.142201\\
7.90909	0.142207\\
8	0.142212\\
8.09091	0.142218\\
8.18182	0.142223\\
8.27273	0.142229\\
8.36364	0.142234\\
8.45455	0.14224\\
8.54545	0.142246\\
8.63636	0.142251\\
8.72727	0.142257\\
8.81818	0.142262\\
8.90909	0.142268\\
9	0.142274\\
9.09091	0.142279\\
9.18182	0.142285\\
9.27273	0.14229\\
9.36364	0.142296\\
9.45455	0.142301\\
9.54545	0.142307\\
9.63636	0.142313\\
9.72727	0.142318\\
9.81818	0.142324\\
9.90909	0.142329\\
10	0.142335\\
};
\end{axis}
\end{tikzpicture}%

%% file: 3D_mode3_param1.tex
%
\begin{tikzpicture}

\begin{axis}[%
width=1.267\figureheight,
height=\figureheight,
at={(0\figureheight,0\figureheight)},
scale only axis,
separate axis lines,
every outer x axis line/.append style={black},
every x tick/.append style={black},
xmin=0,
xmax=10,
xlabel={$E_1$},
every outer y axis line/.append style={black},
every y tick/.append style={black},
ymin=0.1419,
ymax=0.1427,
ylabel={$\gamma_{1,3}$},
axis background/.style={fill=white},
xmajorgrids,
ymajorgrids
]
\addplot [color=red, line width=1.0pt, forget plot]
  table[row sep=crcr]{%
1	0.142652\\
1.09091	0.14263\\
1.18182	0.142609\\
1.27273	0.14259\\
1.36364	0.142571\\
1.45455	0.142554\\
1.54545	0.142537\\
1.63636	0.142521\\
1.72727	0.142506\\
1.81818	0.142492\\
1.90909	0.142478\\
2	0.142465\\
2.09091	0.142452\\
2.18182	0.14244\\
2.27273	0.142429\\
2.36364	0.142417\\
2.45455	0.142407\\
2.54545	0.142396\\
2.63636	0.142386\\
2.72727	0.142376\\
2.81818	0.142367\\
2.90909	0.142358\\
3	0.142349\\
3.09091	0.142341\\
3.18182	0.142332\\
3.27273	0.142324\\
3.36364	0.142316\\
3.45455	0.142309\\
3.54545	0.142301\\
3.63636	0.142294\\
3.72727	0.142287\\
3.81818	0.14228\\
3.90909	0.142273\\
4	0.142266\\
4.09091	0.14226\\
4.18182	0.142253\\
4.27273	0.142247\\
4.36364	0.142241\\
4.45455	0.142235\\
4.54545	0.142229\\
4.63636	0.142224\\
4.72727	0.142218\\
4.81818	0.142212\\
4.90909	0.142207\\
5	0.142202\\
5.09091	0.142196\\
5.18182	0.142191\\
5.27273	0.142186\\
5.36364	0.142181\\
5.45455	0.142176\\
5.54545	0.142171\\
5.63636	0.142166\\
5.72727	0.142162\\
5.81818	0.142157\\
5.90909	0.142152\\
6	0.142148\\
6.09091	0.142143\\
6.18182	0.142139\\
6.27273	0.142134\\
6.36364	0.14213\\
6.45455	0.142126\\
6.54545	0.142122\\
6.63636	0.142117\\
6.72727	0.142113\\
6.81818	0.142109\\
6.90909	0.142105\\
7	0.142101\\
7.09091	0.142097\\
7.18182	0.142093\\
7.27273	0.142089\\
7.36364	0.142085\\
7.45455	0.142081\\
7.54545	0.142078\\
7.63636	0.142074\\
7.72727	0.14207\\
7.81818	0.142066\\
7.90909	0.142063\\
8	0.142059\\
8.09091	0.142056\\
8.18182	0.142052\\
8.27273	0.142048\\
8.36364	0.142045\\
8.45455	0.142041\\
8.54545	0.142038\\
8.63636	0.142034\\
8.72727	0.142031\\
8.81818	0.142028\\
8.90909	0.142024\\
9	0.142021\\
9.09091	0.142018\\
9.18182	0.142014\\
9.27273	0.142011\\
9.36364	0.142008\\
9.45455	0.142004\\
9.54545	0.142001\\
9.63636	0.141998\\
9.72727	0.141995\\
9.81818	0.141992\\
9.90909	0.141988\\
10	0.141985\\
};
\end{axis}
\end{tikzpicture}%

%% file: 3D_mode3_param2.tex
%
\begin{tikzpicture}

\begin{axis}[%
width=1.268\figureheight,
height=\figureheight,
at={(0\figureheight,0\figureheight)},
scale only axis,
separate axis lines,
every outer x axis line/.append style={black},
every x tick/.append style={black},
xmin=0,
xmax=10,
xlabel={$E_2$},
every outer y axis line/.append style={black},
every y tick/.append style={black},
ymin=-0.7,
ymax=0.2,
ylabel={$\gamma_{2,3}$},
axis background/.style={fill=white},
xmajorgrids,
ymajorgrids
]
\addplot [color=green, line width=1.0pt, forget plot]
  table[row sep=crcr]{%
1	-0.664065\\
1.09091	-0.627932\\
1.18182	-0.59424\\
1.27273	-0.562741\\
1.36364	-0.533231\\
1.45455	-0.505526\\
1.54545	-0.479466\\
1.63636	-0.454908\\
1.72727	-0.431727\\
1.81818	-0.40981\\
1.90909	-0.389057\\
2	-0.369376\\
2.09091	-0.350688\\
2.18182	-0.332919\\
2.27273	-0.316003\\
2.36364	-0.29988\\
2.45455	-0.284496\\
2.54545	-0.269801\\
2.63636	-0.25575\\
2.72727	-0.242301\\
2.81818	-0.229417\\
2.90909	-0.217062\\
3	-0.205206\\
3.09091	-0.193818\\
3.18182	-0.182871\\
3.27273	-0.17234\\
3.36364	-0.162203\\
3.45455	-0.152436\\
3.54545	-0.143021\\
3.63636	-0.133939\\
3.72727	-0.125172\\
3.81818	-0.116704\\
3.90909	-0.108521\\
4	-0.100608\\
4.09091	-0.0929517\\
4.18182	-0.0855405\\
4.27273	-0.0783626\\
4.36364	-0.0714071\\
4.45455	-0.0646639\\
4.54545	-0.0581235\\
4.63636	-0.0517769\\
4.72727	-0.0456156\\
4.81818	-0.0396316\\
4.90909	-0.0338174\\
5	-0.0281659\\
5.09091	-0.0226704\\
5.18182	-0.0173245\\
5.27273	-0.0121223\\
5.36364	-0.007058\\
5.45455	-0.00212625\\
5.54545	0.00267808\\
5.63636	0.00735983\\
5.72727	0.0119236\\
5.81818	0.0163739\\
5.90909	0.0207147\\
6	0.0249501\\
6.09091	0.0290838\\
6.18182	0.0331195\\
6.27273	0.0370606\\
6.36364	0.0409103\\
6.45455	0.0446718\\
6.54545	0.0483481\\
6.63636	0.051942\\
6.72727	0.0554562\\
6.81818	0.0588933\\
6.90909	0.0622559\\
7	0.0655464\\
7.09091	0.0687669\\
7.18182	0.0719198\\
7.27273	0.075007\\
7.36364	0.0780307\\
7.45455	0.0809927\\
7.54545	0.0838949\\
7.63636	0.0867391\\
7.72727	0.089527\\
7.81818	0.0922602\\
7.90909	0.0949404\\
8	0.097569\\
8.09091	0.100148\\
8.18182	0.102677\\
8.27273	0.10516\\
8.36364	0.107596\\
8.45455	0.109988\\
8.54545	0.112336\\
8.63636	0.114642\\
8.72727	0.116907\\
8.81818	0.119131\\
8.90909	0.121316\\
9	0.123463\\
9.09091	0.125573\\
9.18182	0.127647\\
9.27273	0.129686\\
9.36364	0.13169\\
9.45455	0.133661\\
9.54545	0.135599\\
9.63636	0.137505\\
9.72727	0.13938\\
9.81818	0.141225\\
9.90909	0.14304\\
10	0.144826\\
};
\end{axis}
\end{tikzpicture}%

%% file: 3D_mode3_param3.tex
%
\begin{tikzpicture}

\begin{axis}[%
width=1.267\figureheight,
height=\figureheight,
at={(0\figureheight,0\figureheight)},
scale only axis,
separate axis lines,
every outer x axis line/.append style={black},
every x tick/.append style={black},
xmin=0,
xmax=10,
xlabel={$E_3$},
every outer y axis line/.append style={black},
every y tick/.append style={black},
ymin=0.1419,
ymax=0.1425,
ylabel={$\gamma_{3,3}$},
axis background/.style={fill=white},
xmajorgrids,
ymajorgrids
]
\addplot [color=blue, line width=1.0pt, forget plot]
  table[row sep=crcr]{%
1	0.142409\\
1.09091	0.142405\\
1.18182	0.1424\\
1.27273	0.142396\\
1.36364	0.142391\\
1.45455	0.142387\\
1.54545	0.142383\\
1.63636	0.142378\\
1.72727	0.142374\\
1.81818	0.142369\\
1.90909	0.142365\\
2	0.142361\\
2.09091	0.142356\\
2.18182	0.142352\\
2.27273	0.142348\\
2.36364	0.142343\\
2.45455	0.142339\\
2.54545	0.142334\\
2.63636	0.14233\\
2.72727	0.142326\\
2.81818	0.142321\\
2.90909	0.142317\\
3	0.142313\\
3.09091	0.142308\\
3.18182	0.142304\\
3.27273	0.142299\\
3.36364	0.142295\\
3.45455	0.142291\\
3.54545	0.142286\\
3.63636	0.142282\\
3.72727	0.142278\\
3.81818	0.142273\\
3.90909	0.142269\\
4	0.142265\\
4.09091	0.14226\\
4.18182	0.142256\\
4.27273	0.142251\\
4.36364	0.142247\\
4.45455	0.142243\\
4.54545	0.142238\\
4.63636	0.142234\\
4.72727	0.14223\\
4.81818	0.142225\\
4.90909	0.142221\\
5	0.142216\\
5.09091	0.142212\\
5.18182	0.142208\\
5.27273	0.142203\\
5.36364	0.142199\\
5.45455	0.142195\\
5.54545	0.14219\\
5.63636	0.142186\\
5.72727	0.142182\\
5.81818	0.142177\\
5.90909	0.142173\\
6	0.142168\\
6.09091	0.142164\\
6.18182	0.14216\\
6.27273	0.142155\\
6.36364	0.142151\\
6.45455	0.142147\\
6.54545	0.142142\\
6.63636	0.142138\\
6.72727	0.142134\\
6.81818	0.142129\\
6.90909	0.142125\\
7	0.142121\\
7.09091	0.142116\\
7.18182	0.142112\\
7.27273	0.142108\\
7.36364	0.142103\\
7.45455	0.142099\\
7.54545	0.142094\\
7.63636	0.14209\\
7.72727	0.142086\\
7.81818	0.142081\\
7.90909	0.142077\\
8	0.142073\\
8.09091	0.142068\\
8.18182	0.142064\\
8.27273	0.14206\\
8.36364	0.142055\\
8.45455	0.142051\\
8.54545	0.142047\\
8.63636	0.142042\\
8.72727	0.142038\\
8.81818	0.142034\\
8.90909	0.142029\\
9	0.142025\\
9.09091	0.142021\\
9.18182	0.142016\\
9.27273	0.142012\\
9.36364	0.142007\\
9.45455	0.142003\\
9.54545	0.141999\\
9.63636	0.141994\\
9.72727	0.14199\\
9.81818	0.141986\\
9.90909	0.141981\\
10	0.141977\\
};
\end{axis}
\end{tikzpicture}%

%% file: 3D_mode4_param1.tex
%
\begin{tikzpicture}

\begin{axis}[%
width=1.267\figureheight,
height=\figureheight,
at={(0\figureheight,0\figureheight)},
scale only axis,
separate axis lines,
every outer x axis line/.append style={black},
every x tick/.append style={black},
xmin=0,
xmax=10,
xlabel={$E_1$},
every outer y axis line/.append style={black},
every y tick/.append style={black},
ymin=0.1414,
ymax=0.1424,
ylabel={$\gamma_{1,4}$},
axis background/.style={fill=white},
xmajorgrids,
ymajorgrids
]
\addplot [color=red, line width=1.0pt, forget plot]
  table[row sep=crcr]{%
1	0.141439\\
1.09091	0.141479\\
1.18182	0.141515\\
1.27273	0.14155\\
1.36364	0.141582\\
1.45455	0.141612\\
1.54545	0.141641\\
1.63636	0.141667\\
1.72727	0.141692\\
1.81818	0.141716\\
1.90909	0.141738\\
2	0.14176\\
2.09091	0.14178\\
2.18182	0.141799\\
2.27273	0.141817\\
2.36364	0.141834\\
2.45455	0.14185\\
2.54545	0.141866\\
2.63636	0.14188\\
2.72727	0.141894\\
2.81818	0.141908\\
2.90909	0.141921\\
3	0.141933\\
3.09091	0.141945\\
3.18182	0.141956\\
3.27273	0.141967\\
3.36364	0.141977\\
3.45455	0.141987\\
3.54545	0.141996\\
3.63636	0.142005\\
3.72727	0.142014\\
3.81818	0.142023\\
3.90909	0.142031\\
4	0.142038\\
4.09091	0.142046\\
4.18182	0.142053\\
4.27273	0.14206\\
4.36364	0.142066\\
4.45455	0.142073\\
4.54545	0.142079\\
4.63636	0.142085\\
4.72727	0.142091\\
4.81818	0.142096\\
4.90909	0.142101\\
5	0.142107\\
5.09091	0.142111\\
5.18182	0.142116\\
5.27273	0.142121\\
5.36364	0.142125\\
5.45455	0.14213\\
5.54545	0.142134\\
5.63636	0.142138\\
5.72727	0.142141\\
5.81818	0.142145\\
5.90909	0.142149\\
6	0.142152\\
6.09091	0.142156\\
6.18182	0.142159\\
6.27273	0.142162\\
6.36364	0.142165\\
6.45455	0.142168\\
6.54545	0.142171\\
6.63636	0.142173\\
6.72727	0.142176\\
6.81818	0.142179\\
6.90909	0.142181\\
7	0.142183\\
7.09091	0.142186\\
7.18182	0.142188\\
7.27273	0.14219\\
7.36364	0.142192\\
7.45455	0.142194\\
7.54545	0.142196\\
7.63636	0.142198\\
7.72727	0.142199\\
7.81818	0.142201\\
7.90909	0.142203\\
8	0.142204\\
8.09091	0.142206\\
8.18182	0.142207\\
8.27273	0.142209\\
8.36364	0.14221\\
8.45455	0.142211\\
8.54545	0.142213\\
8.63636	0.142214\\
8.72727	0.142215\\
8.81818	0.142216\\
8.90909	0.142217\\
9	0.142218\\
9.09091	0.142219\\
9.18182	0.14222\\
9.27273	0.142221\\
9.36364	0.142222\\
9.45455	0.142223\\
9.54545	0.142223\\
9.63636	0.142224\\
9.72727	0.142225\\
9.81818	0.142226\\
9.90909	0.142226\\
10	0.142227\\
};
\end{axis}
\end{tikzpicture}%

%% file: 3D_mode4_param2.tex
%
\begin{tikzpicture}

\begin{axis}[%
width=1.267\figureheight,
height=\figureheight,
at={(0\figureheight,0\figureheight)},
scale only axis,
separate axis lines,
every outer x axis line/.append style={black},
every x tick/.append style={black},
xmin=0,
xmax=10,
xlabel={$E_2$},
every outer y axis line/.append style={black},
every y tick/.append style={black},
ymin=0.1419,
ymax=0.1429,
ylabel={$\gamma_{2,4}$},
axis background/.style={fill=white},
xmajorgrids,
ymajorgrids
]
\addplot [color=green, line width=1.0pt, forget plot]
  table[row sep=crcr]{%
1	0.14273\\
1.09091	0.142703\\
1.18182	0.142678\\
1.27273	0.142654\\
1.36364	0.142632\\
1.45455	0.14261\\
1.54545	0.14259\\
1.63636	0.142571\\
1.72727	0.142553\\
1.81818	0.142536\\
1.90909	0.14252\\
2	0.142504\\
2.09091	0.142489\\
2.18182	0.142475\\
2.27273	0.142461\\
2.36364	0.142448\\
2.45455	0.142435\\
2.54545	0.142423\\
2.63636	0.142411\\
2.72727	0.1424\\
2.81818	0.142389\\
2.90909	0.142379\\
3	0.142368\\
3.09091	0.142359\\
3.18182	0.142349\\
3.27273	0.14234\\
3.36364	0.142331\\
3.45455	0.142322\\
3.54545	0.142314\\
3.63636	0.142305\\
3.72727	0.142297\\
3.81818	0.142289\\
3.90909	0.142282\\
4	0.142274\\
4.09091	0.142267\\
4.18182	0.14226\\
4.27273	0.142253\\
4.36364	0.142246\\
4.45455	0.14224\\
4.54545	0.142233\\
4.63636	0.142227\\
4.72727	0.142221\\
4.81818	0.142215\\
4.90909	0.142209\\
5	0.142203\\
5.09091	0.142197\\
5.18182	0.142191\\
5.27273	0.142186\\
5.36364	0.142181\\
5.45455	0.142175\\
5.54545	0.14217\\
5.63636	0.142165\\
5.72727	0.14216\\
5.81818	0.142155\\
5.90909	0.14215\\
6	0.142145\\
6.09091	0.14214\\
6.18182	0.142135\\
6.27273	0.142131\\
6.36364	0.142126\\
6.45455	0.142122\\
6.54545	0.142117\\
6.63636	0.142113\\
6.72727	0.142108\\
6.81818	0.142104\\
6.90909	0.1421\\
7	0.142096\\
7.09091	0.142091\\
7.18182	0.142087\\
7.27273	0.142083\\
7.36364	0.142079\\
7.45455	0.142075\\
7.54545	0.142071\\
7.63636	0.142067\\
7.72727	0.142064\\
7.81818	0.14206\\
7.90909	0.142056\\
8	0.142052\\
8.09091	0.142049\\
8.18182	0.142045\\
8.27273	0.142041\\
8.36364	0.142038\\
8.45455	0.142034\\
8.54545	0.142031\\
8.63636	0.142027\\
8.72727	0.142024\\
8.81818	0.14202\\
8.90909	0.142017\\
9	0.142013\\
9.09091	0.14201\\
9.18182	0.142007\\
9.27273	0.142003\\
9.36364	0.142\\
9.45455	0.141997\\
9.54545	0.141993\\
9.63636	0.14199\\
9.72727	0.141987\\
9.81818	0.141984\\
9.90909	0.141981\\
10	0.141978\\
};
\end{axis}
\end{tikzpicture}%

%% file: 3D_mode4_param3.tex
%
\begin{tikzpicture}

\begin{axis}[%
width=1.268\figureheight,
height=\figureheight,
at={(0\figureheight,0\figureheight)},
scale only axis,
separate axis lines,
every outer x axis line/.append style={black},
every x tick/.append style={black},
xmin=0,
xmax=10,
xlabel={$E_3$},
every outer y axis line/.append style={black},
every y tick/.append style={black},
ymin=-0.7,
ymax=0.2,
ylabel={$\gamma_{3,4}$},
axis background/.style={fill=white},
xmajorgrids,
ymajorgrids
]
\addplot [color=blue, line width=1.0pt, forget plot]
  table[row sep=crcr]{%
1	-0.667228\\
1.09091	-0.630707\\
1.18182	-0.596676\\
1.27273	-0.564882\\
1.36364	-0.535112\\
1.45455	-0.507179\\
1.54545	-0.480918\\
1.63636	-0.456185\\
1.72727	-0.432848\\
1.81818	-0.410794\\
1.90909	-0.389919\\
2	-0.370132\\
2.09091	-0.351349\\
2.18182	-0.333496\\
2.27273	-0.316507\\
2.36364	-0.300319\\
2.45455	-0.284877\\
2.54545	-0.270132\\
2.63636	-0.256036\\
2.72727	-0.242549\\
2.81818	-0.229631\\
2.90909	-0.217247\\
3	-0.205366\\
3.09091	-0.193956\\
3.18182	-0.182991\\
3.27273	-0.172445\\
3.36364	-0.162294\\
3.45455	-0.152517\\
3.54545	-0.143093\\
3.63636	-0.134004\\
3.72727	-0.125233\\
3.81818	-0.116762\\
3.90909	-0.108577\\
4	-0.100664\\
4.09091	-0.0930084\\
4.18182	-0.0855991\\
4.27273	-0.078424\\
4.36364	-0.0714722\\
4.45455	-0.0647335\\
4.54545	-0.0581983\\
4.63636	-0.0518574\\
4.72727	-0.0457025\\
4.81818	-0.0397254\\
4.90909	-0.0339185\\
5	-0.0282748\\
5.09091	-0.0227875\\
5.18182	-0.0174501\\
5.27273	-0.0122566\\
5.36364	-0.00720137\\
5.45455	-0.00227891\\
5.54545	0.00251593\\
5.63636	0.00718801\\
5.72727	0.011742\\
5.81818	0.0161822\\
5.90909	0.020513\\
6	0.0247382\\
6.09091	0.0288617\\
6.18182	0.032887\\
6.27273	0.0368177\\
6.36364	0.040657\\
6.45455	0.0444081\\
6.54545	0.0480738\\
6.63636	0.0516572\\
6.72727	0.0551609\\
6.81818	0.0585876\\
6.90909	0.0619397\\
7	0.0652197\\
7.09091	0.0684298\\
7.18182	0.0715723\\
7.27273	0.0746491\\
7.36364	0.0776625\\
7.45455	0.0806142\\
7.54545	0.0835062\\
7.63636	0.0863402\\
7.72727	0.089118\\
7.81818	0.0918412\\
7.90909	0.0945114\\
8	0.0971301\\
8.09091	0.0996987\\
8.18182	0.102219\\
8.27273	0.104692\\
8.36364	0.107118\\
8.45455	0.109501\\
8.54545	0.111839\\
8.63636	0.114136\\
8.72727	0.116391\\
8.81818	0.118606\\
8.90909	0.120782\\
9	0.12292\\
9.09091	0.125021\\
9.18182	0.127086\\
9.27273	0.129116\\
9.36364	0.131111\\
9.45455	0.133073\\
9.54545	0.135002\\
9.63636	0.1369\\
9.72727	0.138767\\
9.81818	0.140603\\
9.90909	0.142409\\
10	0.144187\\
};
\end{axis}
\end{tikzpicture}%

%% file: 3D_mode5_param1.tex
%
\begin{tikzpicture}

\begin{axis}[%
width=1.268\figureheight,
height=\figureheight,
at={(0\figureheight,0\figureheight)},
scale only axis,
separate axis lines,
every outer x axis line/.append style={black},
every x tick/.append style={black},
xmin=0,
xmax=10,
xlabel={$E_1$},
every outer y axis line/.append style={black},
every y tick/.append style={black},
ymin=0.06,
ymax=0.26,
ylabel={$\gamma_{1,5}$},
axis background/.style={fill=white},
xmajorgrids,
ymajorgrids
]
\addplot [color=red, line width=1.0pt, forget plot]
  table[row sep=crcr]{%
1	0.245087\\
1.09091	0.228889\\
1.18182	0.214112\\
1.27273	0.200615\\
1.36364	0.188279\\
1.45455	0.176996\\
1.54545	0.166674\\
1.63636	0.157229\\
1.72727	0.148587\\
1.81818	0.140682\\
1.90909	0.133455\\
2	0.126853\\
2.09091	0.120827\\
2.18182	0.115335\\
2.27273	0.110336\\
2.36364	0.105796\\
2.45455	0.101682\\
2.54545	0.0979637\\
2.63636	0.094615\\
2.72727	0.0916109\\
2.81818	0.0889288\\
2.90909	0.0865476\\
3	0.0844484\\
3.09091	0.0826134\\
3.18182	0.0810263\\
3.27273	0.079672\\
3.36364	0.0785365\\
3.45455	0.0776068\\
3.54545	0.076871\\
3.63636	0.0763179\\
3.72727	0.075937\\
3.81818	0.0757186\\
3.90909	0.0756537\\
4	0.0757338\\
4.09091	0.075951\\
4.18182	0.0762979\\
4.27273	0.0767676\\
4.36364	0.0773536\\
4.45455	0.0780497\\
4.54545	0.0788503\\
4.63636	0.0797498\\
4.72727	0.0807434\\
4.81818	0.0818261\\
4.90909	0.0829934\\
5	0.0842412\\
5.09091	0.0855653\\
5.18182	0.0869621\\
5.27273	0.0884278\\
5.36364	0.0899592\\
5.45455	0.0915529\\
5.54545	0.0932061\\
5.63636	0.0949157\\
5.72727	0.096679\\
5.81818	0.0984934\\
5.90909	0.100357\\
6	0.102266\\
6.09091	0.10422\\
6.18182	0.106215\\
6.27273	0.108251\\
6.36364	0.110324\\
6.45455	0.112434\\
6.54545	0.114578\\
6.63636	0.116755\\
6.72727	0.118964\\
6.81818	0.121202\\
6.90909	0.123468\\
7	0.125761\\
7.09091	0.12808\\
7.18182	0.130423\\
7.27273	0.13279\\
7.36364	0.135178\\
7.45455	0.137587\\
7.54545	0.140016\\
7.63636	0.142464\\
7.72727	0.144929\\
7.81818	0.147412\\
7.90909	0.149911\\
8	0.152424\\
8.09091	0.154953\\
8.18182	0.157494\\
8.27273	0.160049\\
8.36364	0.162616\\
8.45455	0.165194\\
8.54545	0.167783\\
8.63636	0.170382\\
8.72727	0.172991\\
8.81818	0.175609\\
8.90909	0.178235\\
9	0.180869\\
9.09091	0.18351\\
9.18182	0.186158\\
9.27273	0.188813\\
9.36364	0.191473\\
9.45455	0.194139\\
9.54545	0.196809\\
9.63636	0.199485\\
9.72727	0.202164\\
9.81818	0.204848\\
9.90909	0.207535\\
10	0.210224\\
};
\end{axis}
\end{tikzpicture}%

%% file: 3D_mode5_param2.tex
%
\begin{tikzpicture}

\begin{axis}[%
width=1.268\figureheight,
height=\figureheight,
at={(0\figureheight,0\figureheight)},
scale only axis,
separate axis lines,
every outer x axis line/.append style={black},
every x tick/.append style={black},
xmin=0,
xmax=10,
xlabel={$E_2$},
every outer y axis line/.append style={black},
every y tick/.append style={black},
ymin=0.06,
ymax=0.22,
ylabel={$\gamma_{2,5}$},
axis background/.style={fill=white},
xmajorgrids,
ymajorgrids
]
\addplot [color=green, line width=1.0pt, forget plot]
  table[row sep=crcr]{%
1	0.189387\\
1.09091	0.177098\\
1.18182	0.165952\\
1.27273	0.15583\\
1.36364	0.146636\\
1.45455	0.138282\\
1.54545	0.130693\\
1.63636	0.123799\\
1.72727	0.117543\\
1.81818	0.111869\\
1.90909	0.106731\\
2	0.102086\\
2.09091	0.0978939\\
2.18182	0.0941216\\
2.27273	0.0907372\\
2.36364	0.0877124\\
2.45455	0.0850212\\
2.54545	0.0826403\\
2.63636	0.0805482\\
2.72727	0.0787253\\
2.81818	0.0771536\\
2.90909	0.0758167\\
3	0.0746996\\
3.09091	0.0737881\\
3.18182	0.0730696\\
3.27273	0.0725321\\
3.36364	0.0721647\\
3.45455	0.0719572\\
3.54545	0.0719002\\
3.63636	0.0719849\\
3.72727	0.0722032\\
3.81818	0.0725474\\
3.90909	0.0730107\\
4	0.0735862\\
4.09091	0.0742679\\
4.18182	0.07505\\
4.27273	0.0759271\\
4.36364	0.0768942\\
4.45455	0.0779464\\
4.54545	0.0790794\\
4.63636	0.0802889\\
4.72727	0.0815709\\
4.81818	0.0829219\\
4.90909	0.0843382\\
5	0.0858166\\
5.09091	0.087354\\
5.18182	0.0889474\\
5.27273	0.0905941\\
5.36364	0.0922913\\
5.45455	0.0940367\\
5.54545	0.0958279\\
5.63636	0.0976626\\
5.72727	0.0995387\\
5.81818	0.101454\\
5.90909	0.103407\\
6	0.105396\\
6.09091	0.107419\\
6.18182	0.109474\\
6.27273	0.111559\\
6.36364	0.113674\\
6.45455	0.115817\\
6.54545	0.117987\\
6.63636	0.120182\\
6.72727	0.1224\\
6.81818	0.124642\\
6.90909	0.126905\\
7	0.129189\\
7.09091	0.131493\\
7.18182	0.133815\\
7.27273	0.136155\\
7.36364	0.138512\\
7.45455	0.140884\\
7.54545	0.143272\\
7.63636	0.145675\\
7.72727	0.148091\\
7.81818	0.15052\\
7.90909	0.152961\\
8	0.155414\\
8.09091	0.157878\\
8.18182	0.160352\\
8.27273	0.162836\\
8.36364	0.16533\\
8.45455	0.167831\\
8.54545	0.170342\\
8.63636	0.172859\\
8.72727	0.175384\\
8.81818	0.177916\\
8.90909	0.180454\\
9	0.182998\\
9.09091	0.185547\\
9.18182	0.188102\\
9.27273	0.190661\\
9.36364	0.193224\\
9.45455	0.195792\\
9.54545	0.198363\\
9.63636	0.200938\\
9.72727	0.203515\\
9.81818	0.206096\\
9.90909	0.208679\\
10	0.211264\\
};
\end{axis}
\end{tikzpicture}%

%% file: 3D_mode5_param3.tex
%
\begin{tikzpicture}

\begin{axis}[%
width=1.268\figureheight,
height=\figureheight,
at={(0\figureheight,0\figureheight)},
scale only axis,
separate axis lines,
every outer x axis line/.append style={black},
every x tick/.append style={black},
xmin=0,
xmax=10,
xlabel={$E_3$},
every outer y axis line/.append style={black},
every y tick/.append style={black},
ymin=0.06,
ymax=0.22,
ylabel={$\gamma_{3,5}$},
axis background/.style={fill=white},
xmajorgrids,
ymajorgrids
]
\addplot [color=blue, line width=1.0pt, forget plot]
  table[row sep=crcr]{%
1	0.151593\\
1.09091	0.141535\\
1.18182	0.132465\\
1.27273	0.124281\\
1.36364	0.116897\\
1.45455	0.110237\\
1.54545	0.104235\\
1.63636	0.098833\\
1.72727	0.0939777\\
1.81818	0.0896229\\
1.90909	0.0857273\\
2	0.0822537\\
2.09091	0.0791689\\
2.18182	0.0764427\\
2.27273	0.0740479\\
2.36364	0.0719602\\
2.45455	0.070157\\
2.54545	0.0686183\\
2.63636	0.0673254\\
2.72727	0.0662615\\
2.81818	0.0654112\\
2.90909	0.0647602\\
3	0.0642955\\
3.09091	0.0640053\\
3.18182	0.0638783\\
3.27273	0.0639045\\
3.36364	0.0640744\\
3.45455	0.0643792\\
3.54545	0.064811\\
3.63636	0.0653621\\
3.72727	0.0660257\\
3.81818	0.0667951\\
3.90909	0.0676644\\
4	0.0686279\\
4.09091	0.0696803\\
4.18182	0.0708166\\
4.27273	0.0720324\\
4.36364	0.0733232\\
4.45455	0.074685\\
4.54545	0.0761139\\
4.63636	0.0776064\\
4.72727	0.0791592\\
4.81818	0.080769\\
4.90909	0.0824329\\
5	0.0841481\\
5.09091	0.0859118\\
5.18182	0.0877217\\
5.27273	0.0895753\\
5.36364	0.0914704\\
5.45455	0.0934049\\
5.54545	0.0953768\\
5.63636	0.0973841\\
5.72727	0.0994251\\
5.81818	0.101498\\
5.90909	0.103601\\
6	0.105734\\
6.09091	0.107893\\
6.18182	0.110079\\
6.27273	0.112289\\
6.36364	0.114522\\
6.45455	0.116778\\
6.54545	0.119055\\
6.63636	0.121351\\
6.72727	0.123667\\
6.81818	0.126\\
6.90909	0.128351\\
7	0.130717\\
7.09091	0.133099\\
7.18182	0.135495\\
7.27273	0.137905\\
7.36364	0.140328\\
7.45455	0.142763\\
7.54545	0.145209\\
7.63636	0.147666\\
7.72727	0.150134\\
7.81818	0.152611\\
7.90909	0.155096\\
8	0.157591\\
8.09091	0.160093\\
8.18182	0.162603\\
8.27273	0.165119\\
8.36364	0.167642\\
8.45455	0.170171\\
8.54545	0.172706\\
8.63636	0.175246\\
8.72727	0.17779\\
8.81818	0.180339\\
8.90909	0.182892\\
9	0.185449\\
9.09091	0.188009\\
9.18182	0.190572\\
9.27273	0.193137\\
9.36364	0.195705\\
9.45455	0.198276\\
9.54545	0.200848\\
9.63636	0.203421\\
9.72727	0.205996\\
9.81818	0.208572\\
9.90909	0.211148\\
10	0.213726\\
};
\end{axis}
\end{tikzpicture}%

%% file: 3D_estimates_global.tex
%
\begin{tikzpicture}

\begin{axis}[%
width=1.268\figureheight,
height=\figureheight,
at={(0\figureheight,0\figureheight)},
scale only axis,
separate axis lines,
every outer x axis line/.append style={black},
every x tick/.append style={black},
xmin=1,
xmax=6,
xlabel={$m$},
every outer y axis line/.append style={black},
every y tick/.append style={black},
ymode=log,
ymin=1e-05,
ymax=1,
yminorticks=true,
axis background/.style={fill=white},
xmajorgrids,
ymajorgrids,
yminorgrids,
]
\addplot [color=red, line width=1.0pt]
  table[row sep=crcr]{%
1	0.100485\\
2	0.0966414\\
3	0.0924295\\
4	0.0883598\\
5	0.0883955\\
6	0.0883642\\
};
\addlegendentry{$E^2_{\mathrm{CRE}}$}

\addplot [color=green, line width=1.0pt]
  table[row sep=crcr]{%
1	0.0122554\\
2	0.00838311\\
3	0.00416868\\
4	7.63145e-05\\
5	5.35344e-05\\
6	5.15998e-05\\
};
\addlegendentry{$\eta^2_{\mathrm{PGD}}$}

\addplot [color=blue, line width=1.0pt]
  table[row sep=crcr]{%
1	0.0882293\\
2	0.0882583\\
3	0.0882608\\
4	0.0882835\\
5	0.088342\\
6	0.0883126\\
};
\addlegendentry{$\eta^2_{\mathrm{dis}}$}

\end{axis}
\end{tikzpicture}%